\documentclass[hypertex,11pt]{article}
\usepackage{amsfonts,amsmath,amscd,amssymb,amsthm}
\usepackage{fullpage}
\usepackage[alphabetic,initials]{amsrefs}
\usepackage[all]{xy}
\usepackage{hyperref} 
\usepackage{mathrsfs}
\input xy%
\usepackage{fullpage}
\usepackage{amssymb}

\def\convf{\hbox{\space \raise-2mm\hbox{$\textstyle      \bigotimes \atop \scriptstyle \omega$} \space}}
\def\0{{\bar 0}}
\def\1{{\bar 1}}
\def\C{{\mathbb C}}
\def\Z{{\mathbb Z}}
\def\Q{{\mathbb Q}}

\def\vv{{\mathbb V}}
\def\vw{{\mathbb W}}

\def\N{{\mathbb N}}

\def\ds{\stackrel{\cdot}{\cup}}

\def\ad{{\operatorname{ad \;}}}

\def\trp{{\operatorname{transpose}}}

\def\hgt{{\operatorname{ht}}}

\def\ST{{\operatorname{Str}\;}}

\def\Ind{{\operatorname{Ind}}}

\def\ch{{\operatorname{ch}\:}}

\def\pwg1{{\operatorname{PWG}}}
\def\pwg{{\operatorname{pwg}}}

\def\wt{{\operatorname{wt}}}
\def\rank{{\operatorname{rank}}}
\def\span{{\operatorname{span}}}

\def\CS{{\operatorname{C-span}}\;}

\def\Hom {{\operatorname{Hom}}}
\def\Ker {{\operatorname{Ker}\;}}

\def\dist{{\operatorname{dist}}}
\def\ann {{\operatorname{ann}}}

\def\Img{{\operatorname{Im}\;}}

\newcommand{\ttk}{\mathtt{k}}
\newcommand{\tte}{\mathtt{e}}
\newcommand{\ttd}{\mathtt{d}}
\newcommand{\ttS}{\mathtt{S}}
\newcommand{\ttT}{\mathtt{T}}

\newcommand{\ovd}{{\Delta}}

\newcommand{\gL}{\Lambda}

\newcommand{\itemi}{\item[{{\rm(i)}}]}
\newcommand{\itemii}{\item[{{\rm(ii)}}]}
\newcommand{\itemiii}{\item[{{\rm(iii)}}]}
\newcommand{\itemiv}{\item[{{\rm(iv)}}]}

\newcommand{\itema}{\item[{{\rm(a)}}]}
\newcommand{\itemb}{\item[{{\rm(b)}}]}
\newcommand{\itemc}{\item[{{\rm(c)}}]}
\newcommand{\itemd}{\item[{{\rm(d)}}]}
\newcommand{\iteme}{\item[{{\rm(e)}}]}
\newcommand{\itemf}{\item[{{\rm(f)}}]}

\newcommand{\da}{\downarrow}
\newcommand{\itemo}{\item[{}]}
\newcommand{\noi}{\noindent}
\newcommand{\ga}{\alpha}
\newcommand{\gb}{\beta}
\newcommand{\gc}{\gamma}

\newcommand{\Gc}{\Gamma}
\newcommand{\Gd}{\Delta}
\newcommand{\Gt}{\Theta}
\newcommand{\gd}{\delta}

\newcommand{\gs}{\sigma}

\newcommand{\go}{\omega}
\newcommand{\gO}{\Omega}
\newcommand{\cg}{{\rm Cdeg}}
\newcommand{\gt}{\tau}
\newcommand{\gz}{\zeta}
\newcommand{\gi}{\iota}
\newcommand{\gl}{\lambda}
\newcommand{\la}{\lambda}

\newcommand{\gr}{\rho}
\newcommand{\gk}{\kappa}
\newcommand{\gep}{\epsilon}
\newcommand{\gth}{\theta}
\newcommand{\op}{\oplus}
\newcommand{\gve}{\varepsilon}

\def\LT{\operatorname{LT}}
\def\Supp{{\operatorname{Supp}}}

\newcommand{\ot}{\otimes}

\newcommand{\fg}{\mathfrak{g}}\newcommand{\fgl}{\mathfrak{gl}}
\newcommand{\fsl}{\mathfrak{sl}}\newcommand{\fpsl}{\mathfrak{psl}}\newcommand{\osp}{\mathfrak{osp}}\newcommand{\spo}{\mathfrak{spo}}

\newcommand{\fr}{\mathfrak{r}}
\newcommand{\fu}{\mathfrak{u}}
\newcommand{\fv}{\mathfrak{v}}
\newcommand{\fw}{\mathfrak{w}}
\newcommand{\fs}{\mathfrak{s}}

\newcommand{\fh}{\mathfrak{h}}

\newcommand{\fb}{\mathfrak{b}}
\newcommand{\fm}{\mathfrak{m}}
\newcommand{\fn}{\mathfrak{n}}

\newcommand{\fc}{\mathfrak{c}}
\newcommand{\fo}{\mathfrak{o}}
\newcommand{\fsp}{\mathfrak{sp}}
\newcommand{\fk}{\mathfrak{k}}
\newcommand{\fp}{\mathfrak{p}}
\newcommand{\fq}{\mathfrak{q}}
\newcommand{\fl}{\mathfrak{l}}

\newcommand{\sfq}{\small{\mathfrak{q}}}
\newcommand{\sfb}{\small{\mathfrak{b}}}

\newcommand{\sfg}{\small{\mathfrak{g}}}

\newcommand{\ff}{\footnote}
\newfont{\eufm}{eufm10 scaled\magstep1}

 \newcommand{\ti}{\times}

\newcommand{\cO}{\mathcal{O}}
\newcommand{\cC}{\mathcal{C}}

\newcommand{\cI}{\mathcal{I}}
\newcommand{\cN}{\mathcal{N}}

\newcommand{\cU}{\mathcal{U}}
\newcommand{\cB}{\mathcal{B}}

\newcommand{\cA}{\mathcal{A}}

\newcommand{\cE}{\mathcal{E}}

\newcommand{\cS}{\mathcal{S}}

\newcommand{\cQ}{\mathcal{Q}}

\newcommand{\cV}{\mathcal{V}}
\newcommand{\cH}{\mathcal{H}}

\newcommand{\bbZ}{\mathbb{Z}}

\newcommand{\ey}{\end{eqnarray}}
\newcommand{\by}{\begin{eqnarray}}
\newcommand{\nn}{\nonumber}

\newcommand{\bco}{\begin{conjecture}}
\newcommand{\ba}{\begin{alg}}
\newcommand{\ea}{\end{alg}}
\newcommand{\eco}{\end{conjecture}}
\newcommand{\bpf}{\begin{proof}}
\newcommand{\epf}{\end{proof}}
\newcommand{\bt}{\begin{theorem}}

\newcommand{\et}{\end{theorem}}

\newcommand{\br}{\begin{rem}}
\newcommand{\er}{\end{rem}}
\newcommand{\brs}{\begin{rems}}
\newcommand{\ers}{\end{rems}}
\newcommand{\bi}{\begin{itemize}}
\newcommand{\ei}{\end{itemize}}
\newcommand{\bl}{\begin{lemma}}
\newcommand{\bsul}{\begin{sublemma}}
\newcommand{\esul}{\end{sublemma}}
\newcommand{\bp}{\begin{proposition}}
\newcommand{\be}{\begin{equation}}
\newcommand{\bc}{\begin{corollary}}
\newcommand{\bexs}{\begin{examples}}
\newcommand{\eexs}{\end{examples}}
\newcommand{\bexa}{\begin{example}}
\newcommand{\eexa}{\end{example}}
\newcommand{\bex}{\begin{exercise}}
\newcommand{\eex}{\end{exercise}}
\newcommand{\btab}{\begin{tab}}
\newcommand{\etab}{\end{tab}}
\newcommand{\bg}{\begin{fig}}
\newcommand{\eg}{\end{fig}}
\newcommand{\el}{\end{lemma}}
\newcommand{\ep}{\end{proposition}}
\newcommand{\ee}{\end{equation}}
\newcommand{\ec}{\end{corollary}}
\newcommand{\Bc}{\begin{center}}
\newcommand{\Ec}{\end{center}}
\newcommand{\bh}{\begin{hyp}}
\newcommand{\eh}{\end{hyp}}
\newcommand{\bhs}{\begin{hyps}}
\newcommand{\ehs}{\end{hyps}}
\newcommand{\bd}{\begin{dfn}}
\newcommand{\ed}{\end{dfn}}

\begin{document}
\title{Table of Contents}

\newtheorem{thm}{Theorem}[section]
\newtheorem{hyp}[thm]{Hypothesis}
 \newtheorem{hyps}[thm]{Hypotheses}
  \newtheorem{rems}[thm]{Remarks}

\newtheorem{conjecture}[thm]{Conjecture}
\newtheorem{theorem}[thm]{Theorem}
\newtheorem{theorem a}[thm]{Theorem A}
\newtheorem{example}[thm]{Example}
\newtheorem{examples}[thm]{Examples}
\newtheorem{corollary}[thm]{Corollary}
\newtheorem{rem}[thm]{Remark}
\newtheorem{lemma}[thm]{Lemma}
\newtheorem{sublemma}[thm]{Sublemma}
\newtheorem{cor}[thm]{Corollary}
\newtheorem{proposition}[thm]{Proposition}
\newtheorem{exs}[thm]{Examples}
\newtheorem{ex}[thm]{Example}
\newtheorem{exercise}[thm]{Exercise}
\numberwithin{equation}{section}%
\setcounter{part}{0}
\newcommand{\drar}{\rightarrow}
\newcommand{\lra}{\longrightarrow}
\newcommand{\rra}{\longleftarrow}
\newcommand{\dra}{\Rightarrow}
\newcommand{\dla}{\Leftarrow}

\newtheorem{Thm}{Main Theorem}


\newtheorem*{thm*}{Theorem}
\newtheorem{lem}[thm]{Lemma}
\newtheorem{fig}[thm]{Figure}
\newtheorem*{lem*}{Lemma}
\newtheorem*{prop*}{Proposition}
\newtheorem*{cor*}{Corollary}
\newtheorem{dfn}[thm]{Definition}
\newtheorem*{defn*}{Definition}
\newtheorem{notadefn}[thm]{Notation and Definition}
\newtheorem*{notadefn*}{Notation and Definition}
\newtheorem{nota}[thm]{Notation}
\newtheorem*{nota*}{Notation}
\newtheorem{note}[thm]{Remark}
\newtheorem*{note*}{Remark}
\newtheorem*{notes*}{Remarks}
\newtheorem{hypo}[thm]{Hypothesis}
\newtheorem*{ex*}{Example}
\newtheorem{prob}[thm]{Problems}
\newtheorem{conj}[thm]{Conjecture}

\title{\v Sapovalov elements and the Jantzen sum formula for contragredient Lie superalgebras.}
\author{Ian M. Musson\ff{Research partly supported by  NSA Grant H98230-12-1-0249, and Simons Foundation grant 318264.} \\Department of Mathematical Sciences\\
University of Wisconsin-Milwaukee\\ email: {\tt
musson@uwm.edu}}
\maketitle
\begin{abstract}
If $\fg$ is a contragredient Lie superalgebra
and $\gc$ is a root of $\fg,$ we prove the existence and uniqueness of \v Sapovalov elements for $\gc$ and give  upper bounds on the degrees of their coefficients.
Then we use \v Sapovalov elements
 to define some new highest weight modules.

If $X$ is a set of orthogonal isotropic roots and $\gl \in \fh^*$ is such that $\gl +\gr$ is orthogonal to all
roots in $X$, 	we construct a highest weight module $M^X(\gl)$ with character $\tte^\gl{p}_X$. Here $p_X$ is a partition function that counts partitions not involving roots in $X$. Examples of such modules can be constructed via
parabolic induction provided $X$ is contained in the set of simple roots of some Borel subalgebra.  However our construction works without this condition and provides a highest weight module for the distinguished Borel subalgebra.
The main results are analogs of the \v Sapovalov determinant and the Jantzen sum formula  for $M^X(\gl)$ when $\fg$ has type A.

 We also explore the behavior of \v Sapovalov elements when the Borel subalgebra is changed,  relations between \v Sapovalov elements for different roots,
and the survival of \v{S}apovalov elements in factor modules of Verma modules. In type A we give a closed formula
for \v Sapovalov elements and give a new 
approach to results of Carter and Lusztig \cite{CL}.

For the proof of the main results it is enough to study the behavior for  certain ``relatively general" highest weights.  Using an  equivalence of categories due to Cheng,  Mazorchuk and  Wang \cite{CMW}, the information we require is
deduced from the behavior of the modules $M^X(\gl)$ when $\fg=\fgl(2,1)$ or $\fgl(2,2)$. These low dimensional
cases are studied in detail in an appendix.
\end{abstract}
\noi \ref{uv.1}. {Introduction.} \\
\ref{1s.1}. {Uniqueness of \v Sapovalov elements.} \\
\ref{1s.5}. {Proof of Theorems \ref{1Shap} and \ref{1aShap}.}\\
\ref{1sscbs}. {Changing the Borel subalgebra.} \\
\ref{RS}. {Relations between \v Sapovalov elements.} \\
\ref{jaf}. {Highest weight modules with prescribed characters.}\\
\ref{SV}. {The submodule structure of Verma modules.} \\
\ref{1cosp}. {An (ortho) symplectic example.} \\
\ref{1s.8}. {The Type A Case.} \\
\ref{1surv}. {{Survival of \v{S}}apovalov elements in factor modules.}\\
\noi \ref{sf}. {The Jantzen sum formula. }\\
\noi  {Appendix \ref{pip}: Anti-distinguished Borel subalgebras.}\\
\noi {Appendix \ref{ldc}: Low Dimensional Cases}.

\section{Introduction.} \label{uv.1} 
Throughout this paper we work over an  algebraically closed field $\ttk $
  of characteristic zero. If  $\fg$ is a simple Lie algebra
necessary and sufficient conditions for the existence of a non-zero homomorphism from $M(\mu)$ to $M(\gl)$ can be obtained by combining work of Verma
\cite{Ve} with work of Bernstein, Gelfand and  Gelfand \cite{BGG1}, \cite{BGG2}.  Such maps can be described explicitly in terms of certain elements introduced by N.N.\v Sapovalov in \cite{Sh}.
Verma modules are fundamental objects in the study of category $\cO$, a study that has blossomed into
an extremely rich theory in the years since these early papers appeared.
We refer to the
 book by Humphreys \cite{H2} for a survey.
\\ \\
Significant advances have been made in the study of finite dimensional modules, and more generally modules in the category $\cO$ for classical simple Lie superalgebras using a variety of techniques.  After the early work of Kac  \cite{K}, \cite{Kac3} the first major advance was made by Serganova who used geometric techniques to obtain a character formula for finite dimensional simple modules over $\fgl(m,n)$, \cite{S2}.  The next development was Brundan's approach to the same problem using a combination of algebraic and  combinatorial  techniques \cite{Br}.  
For more recent developments concerning the category $\cO$ for Lie superalgebras, 
see the survey article \cite{B} by Brundan and the book \cite{CW} by Cheng and Wang.  
 \\ \\The \v Sapovalov determinant, also introduced in \cite{Sh} has been developed in a variety of contexts, such as Kac-Moody algebras \cite{KK}, quantum groups \cite{Jo1}, 
generalized Verma modules \cite{KM},
and Lie superalgebras \cite{G4}, \cite{G}, \cite{G2}.  
The factorization of the \v Sapovalov determinant is the key ingredient in the proof of the original Jantzen sum formula, \cite{J1} Satz 5.3, see also 
\cite{KK}, \cite{MP} Corollary to Theorem 6.6.1 for the Kac-Moody case. 
\\ \\
\v Sapovalov elements  have also appeared in a number of situations in representation theory.  
Though not given this name, they appear in the work  of Carter and Lusztig \cite{CL}.  
Indeed determinants similar to those in our  Theorem \ref{1shgl2} were introduced in \cite{CL} Equation (5), and   our Corollary  \ref{11.5} may be viewed as a version of \cite{CL} Theorem 2.7. Carter and Lusztig use their result to study tensor powers of the defining representation of $GL(V)$, and homomorphisms between Weyl modules in positive characteristic, see also \cite{CP} and \cite{F}. Later Carter \cite{Car} used \v Sapovalov  
elements to construct raising and lowering operators for $\fsl(n,\C)$, see also \cite{Br3}, \cite{Carlin}.  In \cite{Car}, these operators are  used to construct orthogonal bases for non-integral Verma modules, and all finite dimensional modules for $\fsl(n,\C)$. 
\\ \\
More recently Kumar and Letzter gave degrees on the coefficients of \v Sapovalov  
elements. In the Lie algebra case, our Theorem \ref{1Shap} is roughly equivalent to \cite{KL}  Propositions 5.2 and 5.6.
Kumar and Letzter use their result to obtain a
new proof of the irreducibility of the Steinberg module for restricted enveloping algebras and their quantized cousins. 
They also apply their results in these cases to derive versions of the Jantzen sum formula originally obtained by Andersen, Jantzen and Soergel, \cite{AJS} Proposition 6.6. The Strong Linkage Principle may be deduced from this sum formula. 
\\ \\
The purpose of this paper is to initiate the study of
these closely related topics in the super case. New phenomena arise due to the presence of isotropic roots.   We note that in order to pass to positive characteristic, the results in \cite{CL} and \cite{KL} 
are formulated using the Kostant $\Z$-form of $U(\fg).$ 
It is easily seen that our main results can be so formulated.  However we do not go in that direction. 
\\ \\
I would like to thank Jon Brundan for suggesting the use of noncommutative determinants to write \v Sapovalov elements for $\fgl(m,n)$ in Section    \ref{1s.8}, and raising the possibility of using Theorem \ref{1shgl} to prove Theorem \ref{1shapel}.
I also thank Kevin Coulembier, Volodymyr Mazorchuk and Vera Serganova for some helpful correspondence.
\subsection{\bf Preliminaries.}  \label{sss7.1}
\noi
Before we state the main results, we introduce some notation. Throughout the paper $[n]$ denotes the set of integers
$\{1,2,\ldots,n\}$.  Let $\fg= \fg(A,\gt)$ be a finite dimensional contragredient Lie superalgebra with Cartan subalgebra $\fh$, and set of simple roots $\Pi$.
 The superalgebras $\fg(A,\gt)$ coincide with the basic classical simple Lie superalgebras, except that instead of $\fpsl(n,n)$ we obtain $\fgl(n,n).$
 Let $\Delta^{+}$ and  %
\be \label{gtri}\mathfrak{g} = \mathfrak{n}^- \oplus \mathfrak{h}
\oplus \mathfrak{n}^+\ee
be the set of positive roots  containing $\Pi$, and  the corresponding triangular decomposition  of $\fg$ respectively. We use the Borel subalgebras $\mathfrak{b} =  \mathfrak{h}
\oplus \mathfrak{n}^+$ and
$\fb^- =
\mathfrak{n}^- \oplus \mathfrak{h}$.
 The Verma module
$M(\gl)$ with highest weight $\gl \in \fh^*$, and highest weight vector $v_\lambda$ is induced from $\mathfrak{b}$.
Denote the unique simple factor of $M(\gl)$ by $L(\gl)$.
\\ \\
We use the definition of partitions from \cite{M} Remark 8.4.3. Set $Q^+=\sum_{\ga\in \Pi} \N\ga$.
If $\eta \in Q^+$, a {\it
partition} of $\eta$ is a map
$\pi: \Delta^+ \longrightarrow
\mathbb{N} $ such that
$\pi(\alpha) = 0$ or $1$ for all isotropic roots $\alpha$,
 $\pi(\alpha) = 0$ for all even roots $\alpha$ such that $\ga/2$ is a root, and
\be  \label{suma}\sum_{\alpha \in {{\Delta^+}}} \pi(\alpha)\alpha = \eta.\ee
For $\eta \in Q^+$, we denote by $\bf{\overline{P}(\eta)}$ the set of partitions of $\eta$. (Unlike the Lie algebra case, there can be $\eta \in Q^+$ for which $\bf{\overline{P}(\eta)}$ is empty).
If $\pi \in \bf{\overline{P}(\eta)}$  the {\it degree} of $\pi$ is defined to be $|\pi| = \sum_{\alpha \in \Delta^+} \pi(\alpha).$ If $\gs, \pi$ are partitions we say that $\gs + \pi$ is a partition if
$\gs(\ga)+\pi(\alpha) \le 1$ for all isotropic roots  $\ga$. In this case   $\gs + \pi$ is defined by
$(\gs + \pi)(\gc) = \gs(\gc)+\pi(\gc) $  for all positive roots  $\gc$.
\\ \\
Next we  introduce generating functions for certain kinds of partitions.
If $X$ is any  set of
positive roots and $\eta\in \Gd^+$,
set
\[{\bf \overline{P}}_{X}(\eta) = \{\pi \in {\bf \overline{P}}(\eta) | \pi(\alpha) = 0 \mbox{ for all } \ga \in X\}.\]
and ${\bf p}_{X}(\eta) =|{\bf \overline{P}}_{X}(\eta)|$. Usually $X$ will be a set of pairwise orthogonal isotropic roots, but in
the Jantzen sum formula $X$ will sometimes include even roots, see \eqref{yew5}. Let $X_0$ (resp. $X_1$) be the set of even (resp. odd)
roots contained in $X$, and
set $p_X = \sum {\bf p}_X(\eta)\tte^{-\eta}$.
Then
\be \label{pfun}  p_X = \prod_{\alpha \in \Delta^+_{1}\backslash X_1} (1 +
\tte^{- \alpha})/ \prod_{\alpha \in \Delta^+_{0}\backslash X_0} (1 -
\tte^{- \alpha}).\ee
If $X$ is empty, set $p = p_X$, and if $X =\{\gc\}$ is a singleton write

\be \label{sing}{\bf \overline{P}}_{\gc}(\eta),\quad{\bf p}_{\gc}(\eta),\; \mbox{ and } \;p_\gc
\ee
instead of
${\bf \overline{P}}_{X}(\eta),
{\bf p}_{X}(\eta),\mbox{ and } p_X.$
In Section \ref{jaf} we also use the usual definition of a partition.  Thus
for $\eta \in Q^+$, let ${\bf P}_X(\eta)$ denote the set of all functions $\pi:\Gd^+\lra \N$ such that   \eqref{suma} holds, $\pi(\alpha) = 0$ or 1 for $\ga$ an odd root and  $\pi({\alpha}) =0$ for $\ga\in X$. The main difference between the two notations arises in the definition of the elements $e_{-\pi}$ in \eqref{negpar} below, see \cite{M} Example 8.4.5 for the case where $\fg=\osp(1,2)$.
\noi 
We denote the usual BGG category $\cO$ of $\fg_0$ modules here by $\cO_0$ and reserve $\cO$ for the category of 
$\Z_2$-graded $\fg$-modules which are objects of $\cO_0$ when regarded as $\fg_0$-modules. Morphisms in $\cO$ preserve the grading. For a module $M$ in the category $\cO_0$,   the character of $M$ is defined by $\ch M= \sum_{\eta \in \fh^*}\dim_\ttk M^{\eta}\tte^\eta.$  Recall that the Verma module $M(\gl)$ has character $\tte^\gl p.$
\\ \\
\noi Fix a non-degenerate invariant symmetric bilinear form $(\;,\;)$ on $\fh^*$, and for all $\ga \in \fh^*$, let $h_\ga \in \fh$ be the unique element such that $(\ga,\gb) = \gb(h_\ga)$ for all $\gb \in \fh^*$. Then for all $\alpha
\in \Delta^+$, choose elements $e_{\pm \alpha} \in
\mathfrak{g}^{\pm \alpha}$
 such that

\[ [e_{\alpha}, e_{-\alpha}] = h_{\alpha}.\]
It follows that if
$v_\mu$ is a highest weight vector of weight $\mu$ then
\be \label{ebl} e_\ga e_{-\ga} v_\mu = h_\ga v_\mu =(\mu, \ga)v_\mu . \ee
\noi  \noi Fix an order on the set $\Delta^+$, and for $\pi$ a partition,
set

\be \label{negpar} e_{-\pi} = \prod_{\alpha \in \Delta^+} e^{\pi (\alpha)}_{-\alpha},\ee
the product being taken with respect to this order. In
addition set \be \label{pospar}e_\pi = 
\prod_{\alpha \in \Delta^+} e^{\pi(\alpha)}_\alpha,\ee 
where the product is taken in the opposite order. 
The elements $e_{- \pi},$ with
$\pi \in \bf{\overline{P}}(\eta)$ form a basis of
$U(\mathfrak{n}^-)^{- \eta},$
\cite{M} Lemma 8.4.1.
\\ \\
  \noi For a non-isotropic root $\ga,$ we set $\alpha^\vee = 2\alpha /
(\alpha, \alpha)$, and denote the reflection corresponding to $\alpha$ by $
s_\ga.$
As usual the  Weyl group  $W$ is the subgroup of $GL({\fh}^*)$
generated by all reflections.  For $u \in W$ set
 \be \label{nu} N(u) = \{ \alpha \in \Delta_0^+ | u \alpha < 0 \},\qquad \ell(u) = |N(u)|.\ee
We use the following well-known fact several times.
\bl
If $w  = s_\ga u$ with  $\ell(w)>\ell(u)$ and $\ga$ is a isotropic root which is simple for $\fg_0$, then we have a disjoint union
\begin{eqnarray} \label{Nw}
N(w^{-1}) = s_\ga N(u^{-1}) \ds \{\ga\}.
\end{eqnarray}
\el
\bpf See for example \cite{H3} Chapter 1.\epf
\noi Set \be\label{rde}\gr_0=\frac{1}{2}\sum_{\ga \in \Gd_0^+}\ga , \quad
\gr_1=\frac{1}{2}\sum_{\ga \in \Gd_1^+}\ga ,\quad  \gr=
\gr_0 -\gr_1.\ee
Usually  we work with a fixed Borel subalgebra 
$\fb$ which we take to be either the distinguished Borel from \cite{K} Table VI, or an anti-distinguished Borel,   defined in Appendix \ref{pip}. If necessary to emphasize the role played by the Borel subalgebra, we write $\gr_1 =
\gr_1({\fb})$  and $\gr=\gr({\fb}).$ Note that $\fb_0$ is fixed throughout. 
A finite dimensional contragredient Lie superalgebra $\fg$ has, in general several conjugacy classes of Borel subalgebras, and this both complicates and enriches the representation theory of $\fg$.  The complications are partially resolved by at first fixing a Borel subalgebra (or equivalently a basis of simple roots for $\fg$) with special properties.  Then in Section  \ref{1sscbs}  we study the effect of changing the Borel subalgebra.

\subsection{Main Themes.} \label{s.2}
{\bf A. Coefficients of \v Sapovalov  elements.}\\
\\
Fix a positive root $\gc$  and a  positive   integer $m$. If $\gc$ is isotropic, assume $m=1$, and if $\gc$ is odd non-isotropic, assume that $m$ is odd.
There are two special partitions of $m\gc$. Let
$\pi^0 \in {\overline{\bf P}}(m\gc)$ be the unique partition of $m\gc$ such
that $\pi^0(\ga) = 0$ if $\ga \in \Delta^+ \backslash \Pi.$ The partition $m\pi^{\gamma}$ of $m\gc$ is given by $m\pi^{\gamma}(\gc)=m,$ and $m\pi^{\gamma}(\ga)= 0$ for all positive roots $\ga$ different from $\gc.$
We say that $\gth = \theta_{\gamma,m}\in U({\mathfrak b}^{-} )^{- m\gc}$ is a
{\it \v Sapovalov element for the pair} $(\gc,m)$ if it has the form
has the form \be \label{rat}
\theta = \sum_{\pi \in {\overline{\bf P}}(m\gamma)} e_{-\pi}
H_{\pi},\ee where $H_{\pi} \in U({\mathfrak h})$, $H_{\pi^0} = 1,$  and
\be \label{boo} e_{\ga} \theta \in U({\mathfrak g})(h_{\gc} + \rho(h_{\gc})-m(\gc,\gc)/2)+U({\mathfrak g}){\mathfrak n}^+ , \; \rm{ for \; all }\;\ga \in \Delta^+. \ee
We call the $H_\pi$ in \eqref{rat} the {\it coefficients} of $\gth$. For a semisimple Lie algebra, the existence of such elements was shown by \v Sapovalov, \cite{Sh} Lemma 1.
The \v Sapovalov element
$\theta_{\gamma,m}$
has the important property that if $\gl$ lies on a certain hyperplane then
$\theta_{\gc, m}v_{\gl}$ is a highest weight vector  in $M(\gl).$  Indeed set  \be \label{vat}{\mathcal H}_{\gc, m} = \{ \lambda \in  {\mathfrak h}^*|(\lambda + \rho, \gc) = m(\gc, \gc)/2  \}.\ee
Then
$\theta_{\gc, m} v_\lambda$ is a highest weight vector of weight $\lambda -m\gc$ in
$M(\gl)$ for all
$\gl \in \mathcal H_{\gc, m}$.
The normalization condition $H_{\pi^0}=1$
guarantees that $\theta_{\gc, m}v_\gl$ is never zero.
If $\gc$ is isotropic, to simplify notation we set $\mathcal{H}_{\gamma} = \mathcal{H}_{\gamma,1},$ and
denote a  \v Sapovalov element for the pair $(\gc, 1)$ by $\gth_{\gamma}$.  When $X$ is an orthogonal set of isotropic roots, we set $\cH_X= \bigcap_{\gc\in X} \cH_\gc$.
\\ \\
 We give bounds on the degrees of the coefficients $H_{\pi}$ in (\ref{rat}).  There is always a unique coefficient  of highest degree, and we determine the leading term of this coefficient up to a scalar multiple. These results appear to be new even for simple Lie algebras.
The exact form of the coefficients depends on the way the positive roots are ordered.
Suppose $\Pi_{\rm nonisotropic},$ (resp. $\Pi_{\rm even}$) is the set of nonisotropic (resp. even) simple roots, and let $W_{\rm nonisotropic}$ (resp. $W_{\rm even}$) be
the subgroup of $W$ generated by the reflections $s_\ga,$ where $\ga
\in \Pi_{\rm nonisotropic}$ (resp. $\ga \in \Pi_{\rm even}$).
\ff{\rm
In general we have
\be \label{gag1} W_{\rm even} \subseteq W_{\rm nonisotropic}\subseteq W.
\ee
Suppose we  use the distinguished set of simple roots.  Then 
$ W_{\rm even} =W_{\rm nonisotropic}$ if and only if $\fg\neq\osp(1,2n)$, and $W_{\rm nonisotropic} =W$ if and  only if 
$\fg$ has type $A, C$ or $B(0,n).$
}
Consider the following hypotheses.
\be \label{i} \mbox{ The set of simple roots of } \Pi \mbox{ is  either distinguished or anti-distinguished.} \ee
\be \label{iii} \gc = w\gb \mbox{ for a simple root } \gb \mbox{ and } w \in W_{\rm even}.\ee
\be \label{ii}
\gc = w\gb \mbox{ for a simple root } \gb \mbox{ and } w \in W_{\rm nonisotropic}.\ee
When (\ref{i}) and either of (\ref{iii}) or (\ref{ii}) holds we always assume that $\ell(w)$ is minimal, and for $\alpha \in N(w^{-1}),$ we define $q(w,\ga) = (w\gb, \ga^\vee).$ If $\Pi = \{\ga_i|i   = 1, \ldots, t \}$ is the set of simple roots, and $\gc = \sum_{i=1}^t a_i\ga_i,$ then the {\it height} $\hgt \gc$ of $\gc$ is defined to be
$\hgt \gc = \sum_{i=1}^t a_i$.
Let ${\mathcal I}(\mathcal H_{\gc, m})$ be the ideal of $S(\fh)$ consisting of  functions vanishing on $\mathcal{H}_{\gamma,m}$, and $\cO (\mathcal H_{\gc, m})= S(\fh)/{\mathcal I}(\mathcal H_{\gc, m}).$

\bt \nn\label{1Shap}
Suppose $\fg$ is semisimple or a contragredient Lie superalgebra and 
\bi \itema $\gc$ is a positive root such that $(\ref{i})$ and $(\ref{iii})$ hold.   
\itemb If $\gc$ is isotropic assume that $m=1.$ 
\itemc If $\fg= G_2$ or $G(3)$ assume that $\gc$ is  not the highest short root of $\fg_0$.\ei  Then
there exists a \v Sapovalov  element
$\theta_{\gamma,m} \in U({\mathfrak b}^{-} )^{- m\gamma}$, which is unique modulo the left ideal
$U({\mathfrak b}^{-} ){\mathcal I}(\mathcal H_{\gamma,m})$.
The coefficients of $\theta_{\gc, m}$ satisfy
\be \label{x1} |\pi|+\deg  H_{\pi} \le m\hgt \gc, \ee and
\be\label{hig} H_{m\pi^{\gamma}} \mbox{ has leading term } \prod_{\ga \in N(w^{-1})}h_\ga^{mq(w,\ga)}.\ee
\et

\noi 
If $\fg$ is a semisimple Lie algebra this result is due to Kumar and Letzter, see \cite{KL} where the case  $\fg=G_2$ is also covered. 
\\ \\
We show that the exponents in \eqref{hig} are always positive.
\bl If $\gc =w \gb$ with $\ell(w)$ minimal, then $q(w,\ga)$ is a positive integer. 
\el 
\bpf By considering sub-expressions of $w$, it is enough to show this  
when  $w  = s_\ga u$ with  
$\ell(w)>\ell(u)$ and $\ga$ is a simple non-isotropic root. 
Here $\ga_1= -w^{-1}\ga$ is a positive root, so  $(w\gb,\ga^\vee) = -(\gb,\ga_1)$ is a non-negative integer since $\gb$ is a simple root.
If $(w\gb,\ga^\vee)=0,$ then $w\gb = u\gb$ so $w$ does not have minimal length.
\epf
\noi If we assume hypothesis $(\ref{ii})$ instead of $(\ref{iii})$,
it seems difficult to obtain the same estimates on \v Sapovalov elements as in Theorem \ref{1Shap}.
However it is still possible to obtain a reasonable estimate using a different definition of the  degree of a partition, at least if $m=1$. 
Information on $\gth_{\gc, m}$ in general can then be deduced using Theorem \ref{1calu}.
In the Theorem below we assume that $\Pi$ contains an odd non-isotropic root, since otherwise (\ref{iii}) holds and  the
situation is covered by Theorem \ref{1Shap}.  This assumption is essential for Lemma \ref{1gag}. Likewise if
$\gc$ is odd and non-isotropic, then again (\ref{iii}) holds, so we assume that $\gc=w\gb$  with $w \in W_{\rm nonisotropic}$ and
$\gb \in \overline{\Delta}^+_{0} \cup {\overline{\Delta}}^+_{1}$, where
\be \label{Kacroot}  {\overline{\Delta}}^+_{0} = \{ \alpha \in \Delta^+_{0} | \alpha
/ 2 \not\in \Delta^+_{1}\}, \quad {\overline{\Delta}}^+_{1}= \{ \alpha
\in \Delta^+_{1}|2\alpha \not\in \Delta^+_{0}\}.\ee
\noi For $\ga$ a positive root, and then for $\pi$ a partition, we define the {\it Clifford degree} of $\ga, \pi$ by $$\cg (\ga) = 2-i, \mbox{ for } \ga \in \Delta^+_{i},\quad   \cg (\pi) = \sum_{\alpha \in \Delta^+} \pi(\alpha)\cg (\ga).$$  The reason for this terminology is that if we set $U_n = \mbox{ span } \{e_{-\pi}|\cg(\pi)\le n\},$ then $\{U_n\}_{n \ge 0}$ is the {\it Clifford filtration} on $U(\fn^-)$ as in \cite{M} Section 6.5. The associated graded algebra  is a Clifford algebra over $S(\fn^-_0)$.

\begin{theorem} \label{1aShap} Suppose that $\fg$ is  finite dimensional contragredient, and that $\Pi$ contains an odd non-isotropic root. Assume $\gc$ is a positive root such that $\eqref{i}$ and $\eqref{ii}$ hold.  
If $\gc$ is isotropic assume that $m=1,$ and if   $\fg= G_2$ or $G(3)$ assume that $\gc$ is  not the highest short root of $\fg_0$. 
Then
\bi \itema there exists a \v Sapovalov  element
$\theta_{\gamma,m} \in U({\mathfrak b}^{-} )^{- m\gamma}$, which is unique modulo the left ideal
$U({\mathfrak b}^{-} ){\mathcal I}(\mathcal H_{\gc, m})$.
\itemb If $m=1$,
the coefficients of $\theta_{\gc}$ satisfy
\be \label{x2}
 \cg (\pi)+2\deg  H_{\pi} \le 2\ell(w)+\cg(\gc)  ,\ee
 and \[H_{\pi^{\gamma}} \mbox{
has leading term } \prod_{\ga \in N(w^{-1})}h_\ga.\]
\ei
\end{theorem}

\brs \label{rmd}{\rm  (a) \noi By  induction using (\ref{Nw}), we have
\be \label{xxx}\hgt \gc = 1 +\sum_{\ga \in N(w^{-1})} q(w,\ga).\ee
It is interesting to compare the inequalities (\ref{x1}) and (\ref{x2}).  If $\Pi$ does not contain an odd non-isotropic root, then $\cg(\pi) = 2|\pi|$ for all ${\pi \in {\overline{\bf P}}(\gamma)}$, compare Lemma \ref{1pig1} (a).   
Moreover from (\ref{x1}) with $m=1$ and \eqref{xxx}, we have
$$\cg(\pi) + 2\deg  H_{\pi} \le 2+ 2\sum_{\ga \in N(w^{-1})}{q(w,\ga)}.$$
Thus if $q(w,\ga) = 1$ for all  $\ga \in N(w^{-1})$ we have
$$\cg(\pi) + 2\deg  H_{\pi}  \le 2+2\ell(w) ,$$ and if $\gc$ is even this equals the right side of \eqref{x2}.  For example all
the above assumptions hold if $\fg$ is a simply laced semisimple Lie algebra. However if $\fg =\osp(3,2)$ with (odd) 
simple roots $\gb=\gep-\gd, \ga=\gd$, where $\ga$ is non-isotropic and $\gc = w(\gb)$ with $w=s_\ga$, then $q(w,\ga)=2$ and the upper bound for $\deg H_{\pi^{\gamma}}$ given by  (\ref{x2}) is sharper than \eqref{x1}. If instead we take $\gb=\gd-\gep, \ga=\gep$ as simple roots, and $\gc = w(\gb)$ with $w=s_\ga$, then again $q(w,\ga)=2$, but (\ref{x2}) does not hold.  Note in this case $\Pi$ does not contain an  odd non-isotropic root.
\\ \\
(b) It is often the case that equality holds in Equation \eqref{x1} for all   ${\pi \in {\overline{\bf P}}(m\gc)}$, for example this is the case in
Type A, and in the case $\fg=\osp(2,4)$ when $\gc$ is isotropic, see Sections \ref{1s.8} and \ref{1cosp}  respectively.
Furthermore in these cases, it is possible to order the positive roots so that all coefficients are products of linear factors. Equality also holds for all $\pi$ in  \eqref{x2} when $\fg=\osp(3,2)$. However if $\fg=\osp(1,4)$  with simple roots $\ga, \gb$ where $\gb$ is odd, then using \eqref{1f11}, it follows that the coefficient of $e_{- \ga - \gb}
e_{- \gb }$ in $\gth_{\ga  +2\gb}$ is constant.  This means that the inequality in \eqref{x2} can  be strict.
These matters deserve further investigation.
}\ers

\noi Next we mention some important consequences of Theorems \ref{1Shap}  and  \ref{1aShap}.

\bc \label{lter} In Theorems \ref{1Shap} and  \ref{1aShap} $(b),$ $H_{m\pi^{\gamma}}$ is the unique coefficient of highest degree in $\gth_{\gc, m}.$
\ec
\bpf This follows easily from the given degree estimates, and the statements about the leading terms.\epf

\noi Suppose that
$\{U(\fn^-)_j\}$ is either the standard or the Clifford
filtration on $U(\fn^-)$. Then if $\cU = U(\fn^-)\ot \ttk[T]$ we define a filtration $\{\cU_m\}$ on $\cU$ by setting
\be\label{nuf} \cU_m=\sum_{i+j\le m} T^iU(\fn^-)_j.\ee
For $u\in \cU$, set $\deg_{\cU} u= N$ if $N$ is minimal such that $u\in \cU_N$. Choose $\xi \in \fh^*$ such that $(\xi, \gamma) = 0$ and $(\xi, \ga) \neq 0$ for all roots $\ga\in N(w^{-1})$.
Then for $\gl \in \cH_{m\gc}$,
$\widetilde{\lambda} = \lambda + T \xi$ satisfies
$(\widetilde{\lambda} + \rho, \gamma) = m(\gc, \gc)/2  $. Let $\gth=\gth_{\gc, m}$ be as in Theorem \ref{1Shap} or \ref{1aShap} (b).
We obtain the following result for the evaluation of  $\gth$ at ${{\widetilde{\lambda} }}$.

\bc \label{1.5} If $N = m\hgt \gc,$ we have  $\gth_{\gc, m}({{\widetilde{\lambda} }}) = \sum_{\pi \in {\overline{\bf P}}(m\gc)} c_\pi e_{-\pi}\in \cU_N$ with $c_\pi \in\ttk[T]$, and if $\pi =m\pi^\gc$ or $\pi^0$ we have  $\deg_{\cU} c_\pi e_{-\pi}=N.$
\ec
\bpf   This follows easily from \eqref{xxx}.
\epf
\subsubsection{The case where $\gc$ is isotropic.}
We show that positive isotropic roots satisfy Hypothesis \eqref{ii}.
\bl \label{vi} 
Suppose that 
 $\Pi$ is either distinguished or anti-distinguished and that $\gc$ is a positive isotropic root which is not simple.
If $\fg=G(3)$ suppose that 
$\Pi$ is distinguished. 
Then for some non-isotropic root $\ga$ we have $(\gc,\ga^\vee)> 0$, and $s_\ga \gc$ is a positive odd root.  
\el    
\bpf This can be shown by carrying  out a  case-by-case check. For types A-D, see the proof of Lemma \eqref{sim}. 
The Lie superalgebra $D(2,1;\ga)$ is treated in the same way as $D(2,1)= \osp(4,2).$ 
For the other exceptional Lie superalgebras it is convenient to use the quiver $\cQ$ defined as follows.   
Set  $\cQ_0=\Gd^1_+$ the set  of odd positive roots of $\fg$. There is
an arrow $\gc' \longrightarrow \gc$ whenever 
some non-isotropic root $\ga$ we have $(\gc,\ga^\vee)> 0$, and $\gc'=s_\ga \gc$.\
Then $\cQ$ is a sub-quiver of the 
Borel quiver defined in \cite{M} Section 3.6.   
 When $\fg=G(3)$ this quiver appears in Exercise 4.7.10 of \cite{M}; the first and last diagrams at the bottom of page 92 correspond to the distinguished and anti-distinguished Borel subalgebra respectively. For the case of $\fg = F(4)$ see Exercise 4.7.12.  
Examination of these quivers easily gives the result.
\epf
\bc If $\gc$ is as in the statement of Lemma \ref{vi}, then $\gc$ satisfies Hypothesis \eqref{ii}. \ec

\noi  \v Sapovalov elements corresponding to non-isotropic roots for  a basic classical simple  Lie superalgebra
were constructed  in \cite{M} Chapter 9.  This closely parallels the semisimple
case.
Properties of the coefficients of these elements were announced in \cite{M} Theorem 9.2.10. However the bounds on the degrees of the coefficients claimed in \cite{M} are incorrect
if $\Pi$ contains a non-isotropic odd root.
They are corrected by Theorem \ref{1aShap}.\\

\noi {\bf B. Modules with prescribed characters.}\\ \\
 The existence of a unique coefficient  of highest degree in the \v Sapovalov  element
$\theta_{\gamma},$ when $\gc$ is isotropic, is useful in the construction of  some new highest weight modules $M^X(\gl)$, as in the  next result.  We assume that $\Pi$ is a basis of simple roots satisfying hypothesis (\ref{i}).
Let $X$ be an isotropic set of roots and set ${\cH}_{X}= \bigcap_{\gc\in X} {\cH}_{\gc}$.
 Since ${\cH}_{\gc}={\cH}_{-\gc}$ for $\gc\in X$,
we may assume that all roots in $X$ are positive with respect to the basis $\Pi.$
If $\Pi$ contains (resp. does not contain) an odd non-isotropic root,
we may assume that each $\gc\in X$ satisfies \eqref{ii} (resp. \eqref{iii})
and then apply Theorem \ref{1aShap} (resp. Theorem \ref{1Shap}) to each  \v Sapovalov  element $\theta_{\gamma}$.
These assumptions will remain in place for the remainder of this introduction, and then from Section \ref{jaf} onwards.

\bt \label{newmodgen}
Suppose that  $X$ is an isotropic set of positive roots and $\lambda \in {\cH}_{X}$. Then there
exists a factor module $M^{X}(\lambda)$ of
$M(\lambda)$ such that
\[ \ch M^X(\lambda) = \tte^{\lambda} p_{X}.\]
\et\noi
If $X$ is contained in the set of simple roots of some Borel subalgebra $\fb$, it is possible to a construct module with
character as in the Theorem by parabolic induction.  If $\fp = \fb \op \bigoplus_{\gc\in X} \fg^{-\gc}$, then
provided $(\gl+ \gr_\fb, \gc)=0$ for all $\gc\in X,$ the module $\Ind_{\;\fp}^{\;\fg} \; \ttk_\gl$ induced from the one dimensional
$\fp$-module $\ttk_\gl$  with weight $\gl$ will have character $\tte^{\lambda} p_{X}.$ Our construction has however two features lacking in the approach via parabolic induction.   First we do not need the condition that $X$ is contained in the set of simple roots for a Borel subalgebra.  Secondly we provide a uniform construction of a highest weight module for  the distinguished or anti-distinguished Borel subalgebra.  This allows us to compare $M^{X}(\lambda)$ to $M^{Y}(\lambda)$ when, for example $Y\subset X$, or $Y=s_\ga(X)$ for a reflection
$s_\ga.$
\\
\\
The construction of the modules in Theorem \ref{newmodgen} involves a process of deformation and specialization. First we  extend scalars to $A = \ttk[T]$ and $B = \ttk(T)$.  If $R$ is either of these algebras we set
$U(\fg)_R = U(\fg) \otimes R$. Let $\fh_{\mathbb Q}^*$ be the rational subspace of $\fh^*$  spanned by the roots of $\fg$.  Since ${\mathbb Q}X$  is an isotropic subspace of $\fh_{\mathbb Q}^*,$ it cannot contain any non-isotropic root.  Thus we can choose $\xi \in \fh^*$ such that $(\xi, \gamma) = 0$ for all $\gc \in X$, and $(\xi, \ga) \neq 0$ for all even roots $\ga$.
Then for $\gl \in \cH_X$,
\be \label{wld} \widetilde{\lambda} = \lambda + T \xi\ee satisfies
$(\widetilde{\lambda} + \rho, \gamma) = 0$ for all $\gc \in X$.\ff{Replacing $T$ by $c\in \ttk$ in  ${{\widetilde{\lambda} }}$ gives the  point
${{\widetilde{\lambda} }} =\gl+c\xi$
on the line $\gl+ \ttk \xi$  and we think of ${{\widetilde{\lambda} }} =\gl+T \xi$  as a generic point on this line.}
Consider the $U(\fg)_B$-module $M({{\widetilde{\lambda} }})_{B}$ with highest weight ${{\widetilde{\lambda} }} $.
\ff{
Some general remarks on base change for Lie superalgebras and their modules can be found in
\cite{M} subsection 8.2.6.}
The next step is to consider the factor module $M^X({{\widetilde{\lambda} }})_{B}$ of $M({{\widetilde{\lambda} }})_{B}$ obtained
by setting $\gth_\gc v_{{{\widetilde{\lambda} }}}$ equal to zero for $\gc \in X$.
Then we take the $U(\fg)_{A}$-submodule of this factor module  generated by the  highest weight vector $v_{{{\widetilde{\lambda} }}}$  and reduce mod $T$ to obtain the module $M^X(\lambda)$.
\\ \\
 We explain the idea behind the proof that $M^{X}(\lambda)$ has the character asserted in the Theorem.
 If $X =\{\gc\}$ we write $M^{\gc}(\lambda)$ in place of $M^{X}(\lambda)$ throughout the paper.  Set $\gth=\gth_\gc.$
Let
$L =U(\fg)_Bv$  a module with highest weight $\widetilde{\lambda}$ and highest weight vector $v$, such that 
$\gth v = 0$.
Evaluation of \v Sapovalov elements at elements of $\fh^*\ot B$ is defined at the start of Section \ref{1s.5}.
Based on the fact that
$\gth(\widetilde{\lambda})$
has a unique coefficient in $A$ of highest degree we obtain an upper bound on the dimension over $B$ of the weight space $L^{\widetilde{\lambda} -\eta}$, see 
Corollary \ref{ratss}. Now $\gth^2 v_{\widetilde{\lambda}} = 0$ by Theorem \ref{1zprod},
 so the upper bound holds for the first and third terms in the short exact sequence
\be \label{bim} 0\lra N_B=U(\fg)_B \gth v_{\widetilde{\lambda}}\lra M({{\widetilde{\lambda} }})_{B} \lra M_B=M({{\widetilde{\lambda} }})_{B}/U(\fg)_B \gth v_{\widetilde{\lambda}} \lra 0.\ee
 This easily implies that the $U(\fg)_B$-modules $N_B$ and $M_B$ have characters $\tte^{\gl-\gc}{p}_\gc$ and $\tte^\gl{p}_\gc$ respectively.   Now let $M_A=U(\fg)_Av_{\widetilde{\lambda}}$ be the  $U(\fg)_A$-submodule of $M_B$ generated by
$v_{\widetilde{\lambda}}$.  It follows easily that $M_A/TM_A$
is a highest weight module with weight $\gl$ that has character
$\tte^{\gl}{p}_\gc$.
\\ \\
The construction of $M^X(\gl)$ in the general case is similar, but we need an extra ingredient which we illustrate in the case
$X = \{\gc,\gc'\}$.
Here we write $\gth, \gth', U_B$ in place of $\gth_\gc$, $\gth_{\gc'}$ and $U(\fg)_B$ respectively.
We want to apply Corollary \ref{ratss} to the four factor modules
arising from the series

\[ U_B\gth \gth' v_{{\widetilde{\lambda} }}   \subset  U_B\gth v_{{\widetilde{\lambda} }}  \subset U_B\gth' v_{{\widetilde{\lambda} }}   +U_B\gth v_{{\widetilde{\lambda} }} \subset M({{\widetilde{\lambda} }})_B,\]
but to do this we need to know that $\gth\gth' v_{{\widetilde{\lambda} }} \in U_B\gth v_{{\widetilde{\lambda} }}$, $\gth' \gth v_{{\widetilde{\lambda} }} \in U_B\gth \gth' v_{{\widetilde{\lambda} }}$ and $\gth' \gth \gth' v_{{\widetilde{\lambda} }}= 0.$
All this will follow if we know that if $\mu  \in \cH_\gc\cap \cH_{\gc'} $, and $v_{{\widetilde{\mu} }}$ is the highest weight vector in the $U(\fg)_B$ Verma module $M({{\widetilde{\mu} }})_B,$
then  $\gth' \gth v_{{\widetilde{\mu} }}$  and $\gth \gth' v_{{\widetilde{\mu} }}$ are equal up to a scalar multiple.
\\ \\
We give two proofs of this fact. 
The first, Theorem \ref{fiix2} shows that
\be \label{pee} \gth_{\gc'}(\mu - \gc) \gth_\gc(\mu)=a(\mu)\gth_ \gc(\mu - {\gc'})  \gth_ {\gc'}(\mu),\ee
for a rational function $a(\mu)$, which is a
ratio of linear polynomials that are equal except for their constant terms. The computation relies on the fact that any orthogonal set of isotropic roots is $W$-conjugate to a set of simple roots in some Borel subalgebra, see \cite{DS}.
\\ \\ Before turning to the second proof, we need some notation which will also be important for other results.
For $\lambda \in
\mathfrak{h}^*$ define
\by \label{rue} A(\lambda)_{0}  &=&  \{ \alpha \in \overline{\Delta}^+_{0} | (\lambda + \rho,
\alpha^\vee) \in \mathbb{N} \backslash \{0\} \}\nn \\
A(\lambda)_{1} &=& \{ \alpha \in \Delta^+_{1} \backslash \overline{\Delta}^+_{1} |
(\lambda + \rho, \alpha^\vee ) \in 2\mathbb{N} + 1 \} \\
 A(\lambda) &=& A(\lambda)_{0} \cup  A(\lambda)_{1}\nn \\
B(\gl) &=& \{ \alpha \in \overline{\Delta}^+_{1} | (\lambda + \rho,
\alpha) = 0 \} .\nn \ey
Throughout the paper we say that a property holds {\it for general } $\mu$ in a closed subset $X$ of $\fh^*$, if it holds for all $\mu$ in a Zariski dense subset of $X$. (Such $\mu$ are the relatively general highest weights referred to in the abstract.)
The second proof, which is representation theoretic, is based on the fact that for general $\mu$ such that  $B(\mu) = \{ {\gc}, {\gc'}\}$  (that is for general $\mu\in\cH_\gc\cap \cH_{\gc'} $),   the space of highest weight vectors with weight $\mu-\gc-\gc'$ in the Verma module ${M}(\mu)$ has dimension one, 
see Corollary \ref{kt1}. 
This implies \eqref{pee} without however the explicit computation of $a(\mu).$ 
\\ \\
\noi {\bf C. The \v Sapovalov determinant and the Jantzen Sum Formula.}\\ \\
First we recall the Jantzen sum formula for Verma modules given in \cite{M} Theorem 10.3.1.   The Jantzen filtration $\{{M}_{i}(\lambda)\}_{i\ge1}$ on ${M}(\lambda)$ satisfies
\be \label{1jfn}
 \sum_{i > 0} \ch {M}_{i}(\lambda) = \sum_{\alpha \in A(\lambda)}
\ch {M}(s_\alpha \cdot  \lambda) + \sum_{\alpha \in
B(\lambda)} \tte^{\lambda -\alpha}p_\ga.\ee
We use the modules  $M^{\gc}(\lambda)$  to
obtain  an improvement to \eqref{1jfn}.
\ff{This improvement is mentioned by Brundan at the end of section 2 of \cite{B}.}
At the same time, rather than using characters, it will be useful to rewrite the result using the Grothendieck group $K(\cO)$ of the category $\cO$. 
We define $K(\cO)$ to be the free abelian group generated by the symbols $[L(\gl)]$ for $\gl \in \fh^*$.  If $M\in \cO$, the class of $M$ in $K(\cO)$ is defined as
$[M]=\sum_{\gl \in \fh^*} |{M}:L(\gl)|  [L(\lambda)] ,$
where $|{M}:L(\gl)| $ is the multiplicity of the composition factor $L(\gl)$ in $M$.

\bt \label{Jansum}
For all $\lambda \in \mathfrak{h}^*$

\be \label{lb} \sum_{i > 0} [{M}_{i}(\lambda)] = \sum_{\alpha \in A(\lambda)}[{M}(s_{\alpha}\cdot \lambda)] +
\sum_{\gc \in B(\lambda)} [M^{\gc}(\lambda -\gc)].\ee
\et

\bpf Combine Theorem \ref{newmodgen} with \eqref{1jfn}.\epf
\noi The advantage of using this version of the formula is that $K(\cO)$ has a natural partial order.  For $A, B \in \cO$ we write $A
\ge B$ if $[A]-[B]$ is a linear combination of classes of simple modules with non-negative integer coefficients.  Clearly if $B$ is a subquotient of $A$
we have $A \ge B$. \ff{In order not to lengthen this introduction, we postpone a discussion of the merits of working with characters until the start of Section  \ref{SV}. }
\\ \\
 This raises the question of whether there is a similar sum formula for the modules $M^\gc(\gl)$ and more generally for the modules $ M^X(\gl).$  In Theorem \ref{Jansum101} we obtain such a formula for Lie superalgebras of Type A.
As  in the classical case the sum formula follows from the
calculation of the \v Sapovalov determinant  $  \det F_\eta^X $ for the modules $M^X(\gl)$, see Theorem \ref{shapdet}.
These results depend on most of the other results in the paper. As usual $  \det F_\eta^X $ factors into a product linear polynomials, each of which has an interpretation in terms of representation theory.  However there is a rather subtle point about the computation of the leading term of $  \det F_\eta^X $, see Subsection \ref{com}.
The issue is that there
need not be a natural basis for $M^X(\gl)$, indexed by partitions and independent of $\gl$ as there is in the Verma module case. To remedy this we introduce another determinant $\det G^X_{\eta}$ whose leading term is more readily computed and then using \v Sapovalov elements we compare the two determinants.
\\ \\
There are also three
conditions on the highest weight which result in behavior that has no analog for semisimple Lie algebras.
\\ \\
$\bullet$  \indent
The first should be expected based on Theorem \ref{Jansum}. The set $X$ could be contained in a larger isotropic set of simple roots
$Y = X\cup \{\gc\}$, where $\gc \notin X.$
Here  we show that for general $\gl \in \cH_Y$, there is a non-split exact sequence
\[0 \lra  M^{Y}(\gl-\gc) \lra  M^{X}(\gl) \lra M^{Y}(\gl)\lra 0,\]
where the two end terms are simple, compare \eqref{bim}.\\
$\bullet$  \indent Secondly there could be distinct
non-isotropic roots $\ga_1, \ga_2$ such that $\ga_1^\vee\equiv \ga_2^\vee \mod \Q X.$
Then  for  general $\gl \in \cH_X$ such that  $\gl$ is dominant integral, the module $M^X(\gl)$ has length four, except
in one exceptional case, where it has length five, see Theorems \ref{wane} and \ref{name}.\\
$\bullet$  \indent Thirdly suppose  there is no pair of non-isotropic roots $\ga_1, \ga_2$  as in second case.  There could be a
non-isotropic root $\ga$, such that there is a unique
isotropic root $\gc\in X$ with
$(\gc,\ga^\vee)\neq0.$ Then set $\gc' =s_\ga \gc$.
We show that if $(\gc,\ga^\vee)>0$, $(\gl+\gr,\gc')=0$ and $Y =s_\ga X$, then
for general $\gl\in \cH_X\cap\cH_Y$, there is an onto map $M^X(\gl)\lra M^{Y}(\gl)$ having a simple kernel,
Theorem  \ref{eon}.
\\ \\
The Jantzen sum formula (Theorem \ref{Jansum101}) follows easily once the \v Sapovalov determinant  is known.  We note that,  like Theorem \ref{Jansum}, the sum of characters  $\sum_{i>0} \ch {M}^X_i({\lambda})$ is expressed as a sum of other modules with {\it positive} coefficients. 
In the case of semisimple Lie algebras there are other modules which have Jantzen sum formulas, see \cite{AL}, \cite{AJS} Proposition 6.6 and \cite{J1} Satz 5.14. However these sum formulas involve terms with negative coefficients. 
\\ \\
{\bf D. Low Dimensional Examples and Twisting Functors.}
\\ \\
In an appendix to this paper we study the modules $M^{X}(\gl) $ in detail for the cases of $\fgl(2,1)$ and $\fgl(2,2)$. This serves two purposes.  First it
illustrates many of the phenomena that arise elsewhere in the paper.
For example the behavior described after the second (resp. third) bullet above occurs already in the case of $\fgl(2,2)$ (resp. $\fgl(2,1)$).
Secondly the results of the appendix play an important role
in the evaluation of the \v Sapovalov determinant and the proof of the Jantzen sum formula  for the modules $ M^X(\gl)$ in Type A in general. Namely the representation theory needed for these results  can be reduced to the case of $\fgl(2,1)$ and $\fgl(2,2)$ using an equivalence of categories due to Cheng,  Mazorchuk and  Wang \cite{CMW}.
A similar study of low dimensional cases will presumably be necessary to extend these results to the orthosymplectic case.\\ \\
The equivalence in \cite{CMW} is obtained using twisting functors, parabolic induction and  odd  reflections.
In Subsection \ref{agm}, we adapt this result to our needs, and present some other results on twisting functors which may be of independent interest.
\\ \\Finally  our study of $\fgl(2,1)$ and $\fgl(2,2)$ leads us to make the following conjecture.
\bco \label{conje}
For any orthogonal set of isotropic roots $X$ for $\fgl(m,n)$, and any dominant weight $\gl\in\cH_X$, the maximal finite dimensional quotient of $M^X(\gl)$ is simple.
\eco
\noi {\bf Organization of the paper.}\\ \\
\noi In the next Section we discuss the uniqueness of \v Sapovalov elements. Theorems \ref{1Shap}
and \ref{1aShap} are proved in subections \ref{zzw} and \ref{1s.51} respectively.
The proofs depend on a rather subtle cancellation property which is illustrated in
Section \ref{1cosp} in the cases of $\fsp(6)$ and $\osp(2,4)$. Changing the Borel subalgebra and relations between \v Sapovalov elements form the subject of Sections \ref{1sscbs} and \ref{RS} respectively. The modules $M^X(\gl)$ are constructed and Theorem \ref{newmodgen} is proved in Section \ref{jaf}. Applications to the submodule structure of Verma modules follow in Section \ref{SV}.
In Section \ref{1s.8} we give a closed formula for \v Sapovalov elements in Type A.
By definition \v Sapovalov elements give rise  to highest weight vectors in  Verma modules.
The question of when the images of these highest weight vectors in
various factor modules is  non-zero is studied in Section \ref{1surv}. The \v Sapovalov determinant and the Jantzen sum formula are treated in Section \ref{sf}.
\\ \\
The material on the coefficients of \v Sapovalov  elements has appeared in preliminary form in the unpublished preprint \cite{M1}. There is a survey of some of the results from Sections 3-8 in \cite{M2}.
\section{Uniqueness of \v Sapovalov elements.} \label{1s.1}
If $\fk$ is a simple Lie algebra with Cartan subalgebra $\fh$ and  $\mu, \gl \in \fh^*$ we have, for $\fk$ Verma modules
\be \label{1kit}\dim \Hom_\fk(M(\mu),M(\gl))\le 1,\ee
by \cite{D} Theorem 7.6.6,
and it follows that the \v Sapovalov element $\theta_{\gc, m}$ for the pair $(\gc,m)$ is unique modulo the left ideal
$U({\mathfrak b}^{-} ){\mathcal I}(\mathcal H_{\gc, m})$.
\\ \\
\noi We do not know  whether (\ref{1kit}) holds in general for Lie superalgebras, but we note that the analog of (\ref{1kit}) fails for parabolic Verma modules over simple Lie algebras, \cite{IS}, \cite{IS1}.
Here we use an alternative argument to show the uniqueness of \v{S}apovalov elements.
\bt \label{17612} Suppose $\gth_1, \gth_2$ are \v Sapovalov elements for the pair $(\gc,m)$.  Set $\cH = \mathcal H_{\gc, m}.$ Then
\bi \itema for all $\gl \in \mathcal H$ we have $\gth_1 v_\gl =\gth_2 v_\gl$
\itemb $\gth_1- \gth_2
\in U({\mathfrak g}){\mathfrak n}^+  + U({\mathfrak g}){\mathcal I}(\mathcal H).
$ \ei\et \bpf
Set
\[\gL = \{\gl \in \mathcal H| A(\gl) = \{\gc\}, \;B(\gl) = \emptyset \},\]
if $\gc$ is non-isotropic, and
\[\gL = \{\gl \in \mathcal H| B(\gl) = \{\gc\}, \;A(\gl) = \emptyset \},\]
if $\gc$ is isotropic.  If $\gl \in \gL$
it follows from the sum formula \eqref{1jfn} that
$ M_1(\gl)^{\gl-m\gc}$ is one-dimensional.
Because ${M}_{1}(\lambda)$ is the unique maximal submodule of
${M}(\lambda)$, $\gth_1 v_\gl$ and $\gth_2 v_\gl$ are proportional. Then from the requirement that $e_{-\pi^0}$ occurs with coefficient 1 in
a \v{S}apovalov element we have $\gth_1 v_\gl =\gth_2 v_\gl$.  Since $\gL$ is Zariski dense in $\mathcal H$, (a) holds and (b) follows from (a) because
by \cite{M} Lemma 9.4.1 we have
\be \label{1rut} \bigcap_{\gl \in \Lambda} \ann_{U(\sfg)}v_\gl =
 U({\mathfrak g}){\mathfrak n}^+  + U({\mathfrak g}){\mathcal I}(\mathcal H).
 \ee
\epf

\section{Proof of Theorems \ref{1Shap} and  \ref{1aShap}.}\label{1s.5}
 \subsection{Outline of the Proof and Preliminary Lemmas.}\label{1s.4}
Theorems \ref{1Shap} and \ref{1aShap}
are proved by looking at the proofs given in \cite{H2} or \cite{M} and keeping track of the coefficients. 
Given $\gl \in {\mathfrak h}^* $ we define  {\it evaluation at}  $\gl$ to be the
map
$$\varepsilon^\gl:U({\mathfrak b}^{-}) = U({\mathfrak n}^-)\otimes S({\mathfrak h})  \longrightarrow M (\lambda),\;  \quad
\sum_{i} a_i \otimes b_i \longrightarrow \sum_{i} a_i  b_i(\gl)v
_{\lambda}.$$ If $B$ is a $\ttk$-algebra, then evaluation at elements of $\fh^*\ot B$ is defined similarly.
Let  $(\gc,m)$ be as in the statement of the Theorems and set ${\mathcal H} = {\mathcal H}_{\gc, m}$.
If $\gth =\gth_{\gc, m}$ is as in the conclusion of the Theorem, then for any $\gl \in  {\mathcal H},$ $\theta(\gl)v_\gl$
 is a highest weight vector
in $M (\lambda)^{\gl-m\gamma}$. Conversely suppose that $\Lambda$
is a dense subset of ${\mathcal H} $  and that for all $\gl \in
\Lambda$ we have constructed $\theta^\gl \in U({\mathfrak n})^{-m\gamma}$
such that $\theta^\gl v_\gl$ is  a highest weight vector in
$M(\lambda)^{\gl-m\gamma}$ and that
$$\theta^\gl = \sum_{\pi \in {\overline{\bf P}}(m\gc)}a_{\pi, \gl}e_{-\pi}.$$
where $a_{\pi, \gl} $ is a polynomial function of $\gl \in \Lambda$ satisfying suitable conditions.
For $\pi \in {\overline{\bf P}}(m\gc)$,
the assignment $\gl
 \rightarrow a_{\pi, \gl}$ for $\gl \in \Lambda$ determines a polynomial map from ${\mathcal H}$ to
 $U({\mathfrak n}^-)^{-\gamma},$ so  there exists  an element  $H_\pi \in U({\mathfrak h})$ uniquely determined modulo ${\mathcal I}(\mathcal H)$
such that $H_{\pi}(\gl) = a_{\pi, \gl}$ for all $\gl \in \Lambda$.   We define the element
$\theta \in U(\fb^-)$  by
setting
$$\theta= \sum_{\pi \in {\overline{\bf P}}(m\gc)}e_{- \pi}H_{\pi}.$$
\noi Note that $\gth$ is uniquely determined modulo the left ideal
$U({\mathfrak b^-}){\mathcal I}(\mathcal H),$ and that
$\theta(\gl) = \theta^\gl$.
Also, for $\alpha \in \Delta^+$ and
 $\gl \in \Lambda$ we have $e_\alpha \theta v_\gl = e_\alpha
 \theta^\gl v_\gl = 0$, because $\theta^\gl v_\gl = 0$ is a highest weight vector, so  $e_\alpha \theta \in \bigcap_{\gl \in \Lambda} \ann_{U(\sfg)}v_\gl.$   Thus (\ref{boo}) follows from (\ref{1rut}).\\ \\
\noi We need to examine the polynomial nature of the coefficients of $\gth_{\gc,m}$. 
The following easy observation (see \cite{D} Lemma 7.6.9),  is the key to doing this.
For any associative algebra and  $e, a \in A$
     and all $r \in \mathbb{N} $ we have,
\be \label{1cow} e^{r}a = \sum^r_{i=0} \left( \begin{array}{c}
                r \\
                i \end{array}\right) ((\ad e)^i a)e^{r - i} .\ee
\noi Here we should interpret $(\ad e)a $ as $ea-ae$.  The following consequence is well-known, \cite{BR} Satz 1.4, \cite{Ma} Lemma 4.2.  We give the short proof for completeness.
\bc \label{1oreset} Suppose that $S$ is a multiplicatively closed  subset of $A$, consisting of non-zero divisors,  which is  generated by 
locally $\ad $-nilpotent  elements. Then    
$S$ is an Ore set in $A$. \ec 
\noi We write $A_S$ for the resulting Ore localization, or $A_e$ in the case $S=\{e^n|n \in \mathbb{N}\} $. \bpf It suffices to show that for a locally ad-nilpotent generator $e\in S$, the set $\{e^n|n \in \mathbb{N}\} $  
satisfies the  Ore   condition. Given $a \in A$  and $n \in \mathbb{N}$, suppose that
$(\ad e)^{k+1} a = 0$.  Then $e^{k+n}a = a'e^n$, where \[a' = \sum^k_{i=0} \left( \begin{array}{c}
                n+k \\
                i \end{array}\right)
                ((\ad e)^i a)e^{k-i}.\]
 \epf
\noi Now suppose $\ga \in \Pi_{\rm nonisotropic}$, and set $e = e_{-\ga}.$ Then $e$ is a nonzero divisor in
$U = U(\fn^-)$, and the set $\{e^n|n \in \mathbb{N}\} $ is an Ore set in $U $ (and in $U(\fg)$) by Corollary \ref{1oreset}.  Given $a \in U^{-\eta}$, and $\pi$ a partition, the coefficient of $e_{-\pi}$ in \eqref{cows} depends polynomially on $r\in\ttk$, and hence this coefficient is determined by its values on $r\in\N$. Also,
the adjoint action of $\fh$ on $U$ extends to $U_e$, and in the next result we give a basis for the weight spaces of $U_e$.
Let ${\widehat{\bf P}}(\eta)$ be the set of pairs $(k,\pi)$ such that $k \in \mathbb{Z}, \pi \in {\overline{\bf P}}(\eta -k\ga)$ and $\pi(\ga) = 0.$  Then we have \bl \label{1uebasis}  \begin{itemize} \item[{}]
\itema The set $\{e_{-\pi} e^k|(k,\pi ) \in {\widehat{\bf P}}(\eta)\}$ forms a $\ttk $-basis for the weight space $U_e^{-\eta}.$
 \itemb If $u = \sum_{(k,\pi ) \in {\widehat{\bf P}}(\eta)} c_{(k,\pi )} e_{-\pi} e^k \in U_e^{-\eta}$ with $c_{(k,\pi )} \in \ttk ,$ then $u \in U$ if and only if $c_{(k,\pi )} \neq 0$ implies $k \ge 0$.
\end{itemize}
 \el \bpf Throughout the proof we order the set $\Gd^+$ so that for any partition $\gs$ we have
$e_{-\gs}= e_{-\pi}e^{\ell}$ where $\pi(\ga)=0.$\\
(a) Suppose $u \in U_e^{-\eta}$. We need to show that $u$ is uniquely expressible in the form
 \be \label{1china} u = \sum_{(k,\pi ) \in {\widehat{\bf P}}(\eta)} c_{(k,\pi )} e_{-\pi} e^k \ee
 We have $ue^N \in U^{-(N\ga+\eta)}$ for some $N$.  Hence by the PBW Theorem for $U$ we have a unique expression $$ue^N = \sum_{\gs \in {\bf\overline{P}}(\eta+N\ga)}a_\gs  e_{-\gs}  .$$ Now if $a_\gs \neq 0$, then $e_{-\gs}= e_{-\pi}e^{\ell}$ where $\gs(\ga) = {\ell}$ and
  where $\pi \in {\bf\overline{P}} (\eta+N\ga-{\ell}\ga)$ satisfies $\pi(\ga) = 0$. Then $\pi$ and ${k = \ell-N}$ are uniquely determined by $\gs$, so we set $c_{(k,\pi )} =a_\gs$.  Then clearly (\ref{1china})  holds. Given (a), (b) follows from the PBW Theorem.
  \epf

\noi  We remark that if $\ga$ is a non-isotropic odd root, then we can use $e^2$ in place of $e = e_{-\ga}$ in the above Corollary and Lemma.  \noi   However we will need a version of Equation (\ref{1cow}) when $e$ is replaced by an odd element $x$ of a $\Z_2$-graded algebra $A$.
Suppose that $z$ is homogeneous, and define $z^{[j]} = (\ad x)^j z$.  Set $e = x^2$ and apply (\ref{1cow}) to $a = xz =[x, z] + (-1)^{\overline z}zx,$  to obtain
\be \label{1f11}  x^{2\ell+1}z = \sum^\ell _{i=0} \left( \begin{array}{c}
                \ell \\
                i \end{array}\right) z^{[2i+1]}x^{2\ell - 2i}
 + (-1)^{\overline z}\sum^\ell _{i=0} \left( \begin{array}{c}
                \ell \\
i \end{array}\right) z^{[2i]} x^{2\ell - 2i +1} .\ee
\noi The \v Sapovalov elements in Theorems \ref{1Shap} and  \ref{1aShap} are constructed inductively using the next Lemma.

\begin{lemma}\label{11768}
Suppose that  $\mu \in \cH_{\gc',m}, \;\alpha \in\Pi\cap A(\mu)$ and set
\be \label{121c}\gl = s_\alpha\cdot \mu,\;\;\gc = s_\alpha\gc',\;\;p = (\mu + \rho, \alpha^\vee),\;\;q = (\gc, \alpha^\vee).\ee
Assume that $ q, m  \in \mathbb{N}\backslash \{0\}$, and
\begin{itemize}
\itema  $\theta' \in
U({\mathfrak n}^-)^{-m\gc'}$ is such that $v = \theta'v_\mu \in
M(\mu)$ is a highest weight vector.
\itemb  If  ${\ga\in \Delta^+_{1} \backslash \overline{\Delta}^+_{1} }$, then $q = 2$.
\end{itemize}
\noi Then there is a unique $\theta \in U({\mathfrak n}^-)^{-m\gc}$ such that
\begin{equation} \label{121nd}
e^{p + mq}_{- \alpha}\theta' = \theta e^p_{- \alpha}.
\end{equation}
\end{lemma}
\bpf This is well-known, see for example  \cite{H2} Section 4.13 or \cite{M} Theorem 9.4.3.\epf
\noi Equation \eqref{121nd} is the basis for the proof of many properties of \v Sapovalov elements. We note the following variations. First under the hypothesis of the Lemma the \v Sapovalov elements
$\theta_{\gamma',m}$ and
$\theta_{\gamma,m}$ are related by
\begin{equation} \label{121a}
e^{p + mq}_{- \alpha}\theta_{\gc',m}(\mu) = \theta_{\gc,m}(\gl) e^p_{- \alpha}.
\end{equation}
Now suppose that  
instead of the hypothesis $\alpha \in \Pi\cap A(\mu)$, we have $\alpha \in \Pi\cap A(\gl -m\gc)$, and set $r  = (\gl + \rho-m\gc, \alpha^\vee)$. Then
\begin{equation} \label{121b}
\theta_{\gc',m}(\mu)
e^{r + mq}_{- \alpha} = e^r_{- \alpha} \theta_{\gc,m}(\gl).
\end{equation} Note that \eqref{121b} becomes  formally equivalent to \eqref{121a} when we set
$r=-{(p + mq)}$.\\ \\
\noi In the proofs of Theorems \ref{1Shap} and \ref{1aShap} we write $\gc = w\gb$ for $\gb \in \Pi$ and $w\in  W$.
We use the Zariski dense  subset $\gL$ of  ${\mathcal H}_{\beta,m}$ defined by (see \eqref{rue} for the notation)
\be \label{1tar}\Lambda =
\{\nu \in {\mathcal H}_{\beta,m}| \Pi_{\rm nonisotropic}\subseteq A(\nu)
\}
.\ee
\br {\rm Before giving the proofs we remark that it  
is also  possible to construct \v Sapovalov elements  directly in $U(\fb^-)$ instead of their  evaluations in $U(\fn^-)$. To do this note that if $e=e_{-\ga},$ $h\in \fh$ and $f(h)$ is a polynomial in $\fh$, then 
$ef(h) = f(h+\ga(h))e$.  It follows that $\{e^n|n\in \N\}$ is an Ore set in $U(\fb^-)$, and we have 
$U(\fb^-)_e= U(\fn^-)_e\ot S(\fh).$ }\er

\subsection{Proof of Theorem \ref{1Shap}.}\label{zzw}
In this section we assume $\fg$ is contragredient and hypotheses $(\ref{i})$ and $(\ref{iii})$ hold.
In particular $\gc=w\gb$ where $\gb$ is a simple root, and $w\in W.$
Since the statement of Theorem \ref{1Shap} involves precise but somewhat lengthy
conditions on the coefficients, we introduce the following definition as a shorthand.
\bd \label{1sag} {\rm We say that a family of elements $\theta^\gl_{\gamma,m} \in
U({\mathfrak n}^-)^{- m\gamma}$ is {\it in good position for $w$} if  for all
$\gl \in w\cdot  \Lambda$ we have} \ed
\begin{equation} \label{12}
\theta^\gl_{\gamma,m} = \sum_{\pi \in {\overline{\bf P}}(m\gamma)} a_{\pi,\gl}
e_{-\pi},
\end{equation}
where
the coefficients $a_{\pi,\gl} \in \ttk $ depend polynomially on $\gl \in
w\cdot \Lambda$, and \bi%
\itema $\deg  \; a_{\pi,\gl} \leq m\hgt \gc - |\pi|$
\itemb $a_{m\pi^{\gamma},\gl}$ is a polynomial function of
$\gl$ of degree $m(\hgt \gc - 1)$ with highest term equal to $c\prod_{\ga
\in N(w^{-1})}(\gl,\ga)^{mq(w,\ga)}$ for a nonzero constant $c$.
\end{itemize}
We show that the conditions on the coefficients in this definition
are independent of the order on the positive roots $\Delta^+$ used
to define the $e_{-\pi}.$
  Consider  two orders on $\Delta^+$, and for $\pi \in
\overline{\bf P}(m\gc)$, set $e_{-\pi} = \prod_{\alpha \in \Delta^+} e^{\pi
(\alpha)}_{- \alpha}$, and $\overline{e}_{-\pi} = \prod_{\alpha \in
\Delta^+} e^{\pi (\alpha)}_{- \alpha}$, the product being taken with
respect to the given orders. \bl \label{1delorder1}%
Fix a total order on the set ${\overline{\bf P}}(m\gc)$ such that if $\pi, \gs \in {\overline{\bf P}}(m\gc)$ and $|\pi| > |\gs|$ then $\pi$ precedes $\gs,$ and use this order on partitions to induce orders on the bases ${\bf B_1} = \{
e_{-\pi}|\pi \in {\overline{\bf P}}(m\gc)\}$  and ${\bf B_2} = \{\overline{e}_{-\pi}|\pi \in {\overline{\bf P}}(m\gc)\}$  for $U(\fn^-)^{-m\gc}$. Then the change of basis matrix from the basis ${\bf B_1}$ to ${\bf B_2}$ is upper triangular with all diagonal entries equal to $\pm1.$
\el %
\noindent \bpf Let $\{ U_n =U_n(\fn^-) \}$ be the standard filtration on
$U=U({\mathfrak n}^-).$  Note that if $\pi \in {\overline{\bf P}}(m\gc),$ then $e_{-\pi}, \overline{e}_{-\pi} \in
U_{|\pi|}({\mathfrak n}^-)^{-m\gc}.$  Also the factors of $e_{-\pi}$ supercommute
modulo lower degree terms, so for all $\pi \in {\overline{\bf P}}(m\gc),$ $e_{-\pi} \pm\overline{e}_{-\pi} \in
U_{|\pi| - 1}({\mathfrak n}^-)^{-m\gc}$.  The result follows easily.
 \epf
\bl \label{1delorder}%
For $x \in U({\mathfrak n}^-)^{-m\gc}
\otimes S({\mathfrak h}),$ write
\be \label{1owl} x \; = \; \sum_{\pi \in {\overline{\bf P}}(m\gc)}
e_{-\pi}f_\pi \; = \; \sum_{\pi \in {\overline{\bf P}}(m\gc)}
\overline{e}_{-\pi}g_\pi.\ee Suppose that $f_{m\pi^{\gamma}}$ has
degree $m(\hgt \gc - 1),$ and that for all $\pi \in {\overline{\bf P}}(m\gc),$ we have $\deg  f_\pi \le m\hgt \gc - |\pi|.$  Then
$g_{m\pi^{\gamma}}$ has the same degree and leading term as
$f_{m\pi^{\gamma}}$
 and
 for all $\pi \in {\overline{\bf P}}(m\gc),$ we
have $\deg  \;g_\pi \le m\hgt \gc - |\pi|.$
\el %
\noindent \bpf By Lemma \ref{1delorder1} we can write
\[e_{-\pi} = \sum_{\zeta
 \in {\overline{\bf P}}(m\gc)}
c_{\pi,\zeta}\overline{e}_{-\zeta},\] where $c_{\pi,\zeta} \in \ttk ,
c_{\pi,\pi} = \pm1$ and if $c_{\pi,\zeta} \neq 0$ with $\zeta \neq
\pi,$ then $|\gz| < |\pi|$. Thus (\ref{1owl}) holds with
$$g_\zeta = \sum_{\pi \in {\overline{\bf P}}(m\gc)}
c_{\pi,\zeta}f_{\pi}.$$ It follows that $g_\zeta$ is a linear combination of polynomials of
degree less than or equal to $m\hgt \gc - |\gz|.$ Also $|m\pi^{\gamma}| = m,$ and for $\zeta \in {\overline{\bf P}}(m\gc), \zeta \neq
m\pi^{\gamma},$ we have $|\gz|>m.$  Therefore $$g_{m\pi^{\gamma}}
= f_{m\pi^{\gamma}} + \mbox{ a linear combination of polynomials of smaller
degree.}$$ The result follows easily from this.\epf
\noi If $\gc$ is a simple root, then $\theta_{\gamma,m} = e_{-\gc}^m$ satisfies the conditions of Theorem \ref{1Shap}.
Otherwise we have $\gc = w\gb$ for some  $w \in W_{\rm even}, w
\neq 1$.
Write \be \label{1wuga} w =
s_{\ga}u, \quad   \gamma' = u\gb, \quad \gamma = w\gb = s_{\alpha} \gamma',
\ee
 with $\alpha \in \Pi_{\rm even}$ and $\ell(w) =\ell(u) + 1.$
Suppose $\nu \in \Lambda$ and set  \be \label{1wuga1}
\mu =
u\cdot \nu, \quad \lambda = w\cdot \nu
= s_{\alpha}\cdot \mu.
 \ee
\noi For the remainder of this subsection we assume that the positive roots are ordered so that for any partition $\pi$ we have $e_{-\pi} = e_{-\gs}e_{-\ga}^k$ for a non-negative integer $k$ and
 a partition  $\gs$ with $\gs(\ga)=0$. The next Lemma is the key step in establishing the degree estimates in the proof of Theorem \ref{1Shap}.
The idea is to use Equation (\ref{121nd}) and the fact that
$\theta \in U({\mathfrak n}^-)$, rather than a localization of $U({\mathfrak n}^-)$, to show that certain coefficients cancel.
Then using induction and (\ref{121nd}) we obtain the required degree estimates. Since the proof of the Lemma is rather long we break it into
a number of steps.
\bl \label{1wpfg} Suppose that $p, m, q$ are as in Lemma \ref{11768}, $\ga \in \Pi_{\rm even}$ and
\begin{equation} \label{12nd} e^{p+mq}_{-
\alpha}\theta^{\mu}_{\gamma',m} = \theta^\lambda_{\gamma,m}e^p_{-
\alpha}.
\end{equation} 
Then the family $\theta^\gl_{\gamma,m} $ is in good position for $w$ if
the family $\theta^\mu_{\gamma',m}$ is in good position for $u$. \el \bpf
{\it Step 1. Setting the stage.}\\  \\
 Suppose that \begin{equation} \label{14cha7} \theta^{\mu}_{\gamma',m} =
\sum_{\pi' \in {\overline{\bf P}}(m\gamma')} a'_{\pi',\mu} e_{-\pi'},
\end{equation}
and let
\be \label{1qwer} e_{- \pi'}^{(j)} = (\ad e_{- \alpha})^j e_{- \pi'} \in U_{|\pi'|}(\fn^-)^{-(m\gc' +
j \ga )},\ee
for all $j \geq 0,$ and $\pi' \in {\overline{\bf P}}(m\gamma')$. Then by Equation (\ref{1cow})
\begin{eqnarray} \label{149}
e^{p + mq}_{- \alpha} e_{- \pi'} = \sum_{i \geq 0} \left(
\begin{array}{c}
p + mq \\
j
\end{array}
\right) e^{(j)}_{- \pi'} \; e^{p + mq-j}_{- \alpha}.
\end{eqnarray}
Choose $N$ so that
$e^{(N+1)}_{- \pi'} = 0,$ for all $\pi' \in {\overline{\bf P}}(m\gc').$ Then for
all such $\pi' $ and $j = 0, \ldots, N$ we can write
\begin{eqnarray} \label{199}
e_{- \pi'}^{(j)}e^{N-j}_{- \alpha} = \sum_{\zeta \in {\overline{\bf P}}(m\gc' +
N \ga )} b_{j, \zeta}^{\pi'}e_{- \zeta} ,
\end{eqnarray}
with $ b_{j, \zeta}^{\pi'} \in \ttk .$ Furthermore if  $ b_{j,
\zeta}^{\pi'} \neq 0,$ then
since
$e_{- \pi'}^{(j)}e^{N-j}_{- \alpha} \in  U_{|\pi'| + N - j}$,
 (\ref{1qwer}) gives
\begin{eqnarray} \label{198}
|\gz| \leq |\pi'| + N - j.
\end{eqnarray}
{\it Step 2. The cancellation step.}\\ \\
By Equations (\ref{14cha7}) and (\ref{149})
\begin{eqnarray} \label{1199} e^{p + mq}_{- \alpha}
\theta^{\mu}_{\gamma',m}  & = & \sum _{\pi' \in
{\overline{\bf P}}(m\gamma')} a'_{\pi',\mu} 
e^{p + mq}_{- \alpha}e_{-\pi'}\\
& = & \sum_{ \pi' \in {\overline{\bf P}}(m\gamma'),\;j \geq 0} \left(
\begin{array}{c}
p + mq \nonumber\\
j
\end{array}
\right) a'_{\pi',\mu} e_{- \pi'}^{(j)}e^{p + mq-j}_{- \alpha}.
\end{eqnarray}
Now collecting coefficients, set \be\label{159} c_{ \zeta,\lambda}
 = \sum_{\pi' \in {\overline{\bf P}}(m\gamma'),\;j \geq 0} \left(
\begin{array}{c}
p + mq \\
j
\end{array}
\right) a'_{\pi',\mu} 
b_{j, \zeta}^{\pi'}. \ee Then using Equations (\ref{199}) and (\ref{1199}), we have in $U_e,$ where $e = e_{-\ga}$,
\begin{eqnarray} \label{1lab} e^{p + mq}_{- \alpha} \theta^{\mu}_{\gamma',m} &
= & \sum_{\zeta \in {\overline{\bf P}}(m\gc' + N\ga )}%
c_{ \zeta, \lambda}e_{-\zeta}e^{p + mq-N}_{- \alpha}.
\end{eqnarray}
By (\ref{12nd}) and Lemma \ref{1uebasis}, $c_{ \zeta, \lambda} = 0$ unless $\zeta(\ga) \geq
N-mq.$\\ \\
{\it Step 3. The coefficients $a_{\pi,\lambda}$.}\\ \\
It remains to deal with the nonzero terms $c_{ \zeta, \lambda}$. There is a bijection
\[f:{\overline{\bf P}}(m\gamma) \longrightarrow \{\zeta \in {\overline{\bf P}}(m\gamma' +N\ga)\;|\;\zeta(\ga) \geq N-mq\},\]
defined by \be \label{1yam}
 (f\pi)(\gs)=\left\{ \begin{array}
  {ccl}\pi(\gs)&\mbox{if} \;\; \gs \neq \ga,
\\\pi(\ga) + N - mq &\mbox{if} \;\; \gs = \ga.
\end{array} \right. \ee
Moreover if $f\pi = \zeta,$ then
\begin{eqnarray} \label{1119}|\gz| = |\pi| + N -mq
\end{eqnarray} and $e_{-\pi} = e_{-\zeta}e^{mq-N}_{-
\alpha}$.
Thus in Equation (\ref{12nd}) the coefficients $a_{\pi,\lambda}$
of $\theta^\lambda_{\gamma,m} $ (see (\ref{12})) are given by \be
\label{1cwp2} a_{\pi,\lambda} = c_{f(\pi),\lambda}.\ee
{\it Step 4. Completion of the proof.}\\ \\
We now show that the family $\theta^\gl_{\gamma,m} $ is in good position for $w.$ For this we
use Equations (\ref{159}) and  (\ref{1cwp2}), noting that $p =
-(\gl + \rho,\ga)$ depends linearly on $\gl.$   It is clear
that the coefficients $a_{\pi,\lambda}$ are polynomials in
$\lambda$.
 By induction $\deg  a'_{\pi',\mu} \le m\hgt \gc' - |\pi'|.$ Thus using
 (\ref{159}), \be \label{1fox} \; \mbox{deg} \;\; a_{\pi,\gl} = \;
\mbox{deg} \; c_{\zeta,\gl}  \leq  \max\{j + \; \mbox{deg} \;
a'_{\pi',\mu}\;|\; b_{j, \zeta}^{\pi'}
\neq 0\}. \ee Now if $b_{j, \zeta}^{\pi'}
\neq 0$ then Equation (\ref{198}) holds. Therefore by Equation (\ref{1119}) 
\begin{eqnarray} 
\mbox{deg} \;\; a_{\pi,\gl} \nonumber &\leq&
\mbox{deg} \;\; a'_{\pi',\mu}+ |\pi'| + N - |\gz|
\\ &=&  \mbox{deg} \;\; a'_{\pi',\mu}+ |\pi'| - |\pi| +mq\nonumber 
\end{eqnarray}
Finally since $\gamma = \gc' +
q\ga$, induction gives
(a) in Definition \ref{1sag}.\\ \\
\noi  Also, modulo terms of lower
degree
\be \label{1Nw1}  a_{m\pi^{\gamma},\gl} \equiv \left( \begin{array}{c}
p + mq\\
mq
\end{array}   \right)
a'_{m\pi^{\gamma'},\mu} .\ee
Note that the above binomial coefficient is a polynomial of degree $mq$ in $p$.
 By induction 
 $a'_{m\pi^{\gamma'},\mu}$ has
highest term $c'\prod_{\gt \in N(u^{-1})}(\mu,\gt)^{mq(u,\tau)}$ as a polynomial in $\mu$, for
a nonzero constant $c'$. Now $(\mu + \gr,\gt) = (\gl +\gr, s_\ga \gt)$, and $(\mu + \gr,\gt)- (\mu,\gt)$, $(\gl +\gr, s_\ga \gt) - (\gl, s_\ga \gt)$ are independent of $\mu, \gl$.  Therefore  as a polynomial in $\gl$, $a'_{m\pi^{\gamma'},\mu}$ has
highest term $$c'\prod_{\gt \in N(u^{-1})}(\gl,s_\ga\gt)^{mq(u,\tau)} = c'\prod_{\gt \in N(u^{-1})}(\gl,s_\ga\gt)^{mq(w,s_\ga\tau)}.$$
 Since we assume that $\gc$ is  not the highest short root of $\fg_0$ if $\fg= G_2$ or $G(3)$, it follows from 
\cite{M} Table 3.4.1, that the $\ga$-string through $\gc'$ has length $q+1$, and $\gc'-\ga$ is not a root. Thus from the representation theory of $\fsl(2)$, it follows that $(\ad e_{- \alpha})^q(\fg^{-\gc'}) = \fg^{-\gc}$. Hence, as $e_{-\gc}$ is not used in the construction of $\theta_{\gamma',m},$ we can choose the notation so that
$e^{(mq)}_{- m\gamma'} = e_{- m\gamma}$. Then $e^{(mq+1)}_{- m\gamma'} = 0.$ Now $q(w,\ga) = (\gc, \alpha^\vee) = q,$ and the degree of the binomial coefficient in (\ref{1Nw1}) as a polynomial in $p$ is $mq$, so the
claim about the leading term of $a_{m\pi^{\gamma},\gl} $ in Definition \ref{1sag} (b) follows from
Equations (\ref{1Nw1}) and (\ref{Nw}).\epf
\begin{theorem}  \label{1localshap}
Suppose $\gamma = w\beta$ with $w \in W_{\rm even}$ and $\gb$ simple. 
Define $\gL$ as in \eqref{1tar}. Then there exists a family of elements
$\theta^\gl_{\gamma,m} \in U({\mathfrak n}^-)^{- m\gamma}$ for all $\gl \in
w\cdot  \Lambda$  which is in good position for $w$, such that
\[\theta^\gl_{\gamma,m}v_{\gl} \; \mbox{is a highest weight vector in}\;
M(\gl)^{\gl - m\gamma}.\]
\end{theorem}
\noindent \bpf We use induction on the length of $w.$ If $w = 1,$ we
take $\theta^\gl_{\gamma,m} = e_{-\gb}^m$ for all $\gl.$
Now assume that $w
\neq 1,$ and use the notation of Equations (\ref{1wuga}) and (\ref{1wuga1}).
Suppose $\lambda =s_{\alpha}\cdot \mu =w\cdot \nu\in w\cdot\Lambda,$ and set \[p = (\mu + \rho,
\alpha^{\vee}) = (\nu + \rho, u^{-1}\alpha^{\vee}), \quad(\gamma, \alpha^{\vee}) = q.\] Then $p$ and $q$ are positive integers, and 
$\gl = \mu - p\alpha$ and $\gamma = \gamma' +q\alpha.$
Now $U({\mathfrak n}^-)e^p_{- \alpha} v_{\mu}$ 
is a submodule of  $M(\mu)$ which is isomorphic to 
$M(\gl)$.  Also 
$M(\gl)$ is uniquely embedded in $M(\mu)$, by \cite{M} Theorem 9.3.2, so we set $M(\gl)=U({\mathfrak n}^-)e^p_{- \alpha} v_{\mu}$.
Induction gives  elements $\theta^{\mu}_{\gamma',m} \in U({\mathfrak n}^-)^{-
\gamma'}$ which are in good position for $u,$ such that
\[v = \theta^{\mu}_{\gamma',m}v_{\mu} \in M(\mu)^{\mu -
m\gamma'} \mbox{is a highest weight vector}.
\]
By Lemma \ref{11768} there exists a unique element
$\theta^{\lambda}_{\gamma,m} \in U({\mathfrak n}^-)^{- m\gamma}$ such that
(\ref{12nd})  holds
and therefore
\[e^{p+mq}_{- \alpha} v
=
\theta^{\lambda}_{\gamma,m}
e^{p}_{- \alpha} v_{\mu}
\in U({\mathfrak n}^-)e^p_{- \alpha} v_{\mu} =
M(\lambda).\]
It follows from Lemma \ref{1wpfg} that the family
$\theta^\gl_{\gamma,m}$ is in good position for $w$. \epf \noindent {\it Proof of Theorem \ref{1Shap}.} Let $\theta^\gl_{\gamma,m}$ be the family of elements from Theorem \ref{1localshap}. The existence of the elements
$$\theta_{\gc, m} = \sum_{\pi \in {\overline{\bf P}}(m\gc)} e_{-\pi} H_{\pi} \in U(\fb^-),$$ with $\theta_{\gc, m}(\gl) = \theta_{\gc, m}^\gl$ for all $ \lambda \in w\cdot\Lambda$
follows since $ w \cdot\Lambda$ is
Zariski dense in ${\mathcal H}_{\gc,m}.$
The claims about the coefficients
$H_{\pi}$ hold since  the family $\theta^\gl_{\gamma,m}$ is in good position for $w$.
\hfill  $\Box$
\subsection{Proof of Theorem \ref{1aShap}.} \label{1s.51}
Now suppose that $\fg$ is contragredient, with $\Pi$ as in $(\ref{i})$. We assume that
$\Pi$ contains an odd non-isotropic root and $\gc=w\gb$  with $w \in W_{\rm nonisotropic}$,
$\;\gb \in \overline{\Delta}^+_{0}\cup {\overline{\Delta}}^+_{1}$. Let $\{ U_n =U_n(\fn^-) \}$ be the Clifford filtration on
$U=U({\mathfrak n}^-).$
Our assumptions have the following consequence.
\bl \label{1gag} Let   $\gc$ be as above and $\ga \in \Pi_{\rm nonisotropic}$ is such that $(\gc,\ga)\neq 0.$ Then
\bi \itema If $\ga$ is even then $(\gc,\ga^\vee) = \pm 1$
\itemb  If $\ga$ is odd then $(\gc,\ga^\vee) = \pm 2$. \ei\el \bpf Left to the reader. The assumption that
$\Pi$ contains an odd non-isotropic root is crucial to (a).  Without this $\osp(3,2)$ would be a counterexample.\epf
\bl \label{1it} The Clifford filtration on $U(\fn^-)$  is stable under the adjoint action of $\fn_0^-,$ and satisfies $\ad \fn_1^- (U_n) \subseteq U_{n+1}$. \el
\bpf Left to the reader.\epf\noi
\noi Fix a total order on the set ${\overline{\bf P}}(\gc)$ such that if $\pi, \gs \in {\overline{\bf P}}(\gc)$ and $|\pi| > |\gs|,$ or if
$|\pi| = |\gs|$ and $\cg(\pi) > \cg(\gs)$ then $\pi$ precedes $\gs$.
\\ \\
\noi Next we prove a Lemma relating $\cg(\pi)$ to $|\pi|$ and the order defined above.
\bl
Set $\Xi = {\Delta}^+_{0} \backslash \overline{\Delta}^+_{0},$ and suppose $\pi, \gs \in {\overline{\bf P}}(\gc)$.
Let $a(\pi)  = \sum_{2\gd \in \Xi}\pi(\gd).$
\bi \label{1pig1}\itema We have   $2|\pi| -\cg(\pi) = a(\pi)$.\ei
\bi \itemb $a(\pi) \le 2.$
\itemc If $\gs$ precedes $\pi$, then $\cg(\pi) \le \cg(\gs).$
\ei \el
\bpf (a) follows since 
\[ |\pi| =a(\pi)+ \sum_{\ga \in \overline{\Delta}_{0}^+}\pi(\ga)  + \sum_{\ga \; \rm non-isotropic}\pi(\ga) ,\]
and \[ \cg(\pi) =a(\pi)+ 2\sum_{\ga \in \overline{\Delta}_{0}^+}\pi(\ga)+ \sum_{\ga \; \rm non-isotropic}\pi(\ga) .\]
For  (b) we note that the Lie superalgebras that have 
 an odd non-isotropic root $\gd$ are $G(3)$ and the family $\osp(2m+1,2n)$.
     Define a group homomorphism $f:\bigoplus_{\ga \in \Pi}\Z\ga \lra \mathbb{Z}$ by setting  $f( \gd) = 1$ and $f(\ga) = 0$ for any $\ga \in \Pi$, $\ga \neq \gd.$
It can be checked on a case-by-case basis that if $\gd\in \Pi$, then $\gd$ occurs with coefficient at most two when a positive root $\gc$ is written as a linear combination of simple roots. Since $a(\pi)  = f(\gc)$  for $\pi\in  {\overline{\bf P}}(\gc)$,
(b) follows. If (c) is false, then  by definition of the order, we must have $|\pi| < |\gs|$ and
$\cg(\pi) > \cg(\gs).$ But then by (a) this implies that $$a(\gs) =2|\gs| -\cg(\gs)  \ge 2|\pi| -\cg(\pi)+3 = a(\pi) + 3\ge 3$$ which contradicts (b). \epf
\br {\rm
We remark that if $\fg=G(3)$ and we work with the anti-distinguished system of simple roots, then
$W'=W_{\rm nonisotropic}$ has order 4
and is generated by the commuting reflections $s_{-\gd}$ and $s_{-3\ga_1-2\ga_2}$ (notation as in \cite{M} section 4.4).
Thus in particular if $\gb$ is the simple isotropic root, then $|W'\gb|=4$.
We leave it to the reader to check the assertions in Theorem  \ref{1aShap}  in this case,  using Equation  \eqref{1f11} and Lemma \ref{11768}. \ff{If $\gb$ is simple non-isotropic, then $W'\gb\cap \Gd^+ =\{\gb\}$.  If $\fg=G(3)$ it is more interesting to use the distinguished system of simple roots.} For the rest of this section we assume that
$\fg\neq G(3)$. Thus since we assume that
$\Pi$ contains an odd non-isotropic root, $\fg=\osp(2m+1,2n)$ for some $m, n$.  In the proof of Lemma \ref{1ant} we use the explicit description
of the roots of $\fg$ given in, for example \cite{M} Section 2.3.}
\er

\bd \label{1sog} We say that a family of elements $\theta^\gl_{\gamma} \in
U({\mathfrak n}^-)^{- \gamma}$ is {\it in good position for $w$} if  for all
$\gl \in w\cdot  \Lambda$ we have
\begin{equation} \label{12os}
\theta^\gl_{\gamma} = \sum_{\pi \in {\overline{\bf P}}(\gamma)} a_{\pi,\gl}
e_{-\pi},
\end{equation}
where
the coefficients $a_{\pi,\gl} \in \ttk $ depend polynomially on $\gl$, and \bi%
\itema $2\deg  \; a_{\pi,\gl} \leq 2\ell(w)  +\cg(\gc)  - \cg(\pi)$
\itemb $a_{\pi^{\gamma},\gl}$ is a polynomial function of
$\gl$ of degree $\ell(w)$ with highest term equal to $c\prod_{\ga
\in N(w^{-1})}(\gl,\ga)$ for a nonzero constant $c$.
\end{itemize}\ed
\noi 	Consider  two orders on $\Delta^+$, and for $\pi \in
\overline{\bf P}(\gc)$, set $e_{-\pi} = \prod_{\alpha \in \Delta^+} e^{\pi
(\alpha)}_{- \alpha}$, and $\overline{e}_{-\pi} = \prod_{\alpha \in
\Delta^+} e^{\pi (\alpha)}_{- \alpha}$, the product being taken with
respect to the given orders.
\bl \label{1ant}Use the order on ${\overline{\bf P}}(\gc)$ defined just before Lemma \ref{1pig1}  to induce orders on
two bases ${\bf B_1} = \{e_{-\pi}|\pi \in {\overline{\bf P}}(\gc)\}$  and
${\bf B_2} = \{\overline{e}_{-\pi}|\pi \in {\overline{\bf P}}(\gc)\}$  for $U(\fn^-)^{-\gc}$.
Write $x \in U({\mathfrak n}^-)^{-\gc}\otimes S({\mathfrak h})$ as
\[x \; = \; \sum_{\pi \in {\overline{\bf P}}(\gc)}e_{-\pi}f_\pi \; = \;
\sum_{\pi \in {\overline{\bf P}}(\gc)}\overline{e}_{-\pi}g_\pi.\]
as in Equation {\rm (\ref{1owl})}. If the coefficients $f_\pi$ satisfy
\bi%
\itema $2\deg  \; f_{\pi} \leq 2\ell(w)  +\cg(\gc)  - \cg(\pi)$
\itemb $f_{\pi^{\gamma}}$ is a polynomial  of degree $\ell(w)$ with highest term equal to $c\prod_{\ga
\in N(w^{-1})}h_\ga$ for a nonzero constant $c$.
\end{itemize}
then the  coefficients  $g_\pi$ satisfy the same conditions.\el \bpf First we claim that the analog of Lemma \ref{1delorder1} holds, that is changing the  order on $\Gd^+$ in the definition of the $e_{-\pi}$
requires only the introduction of terms $e_{-\gs}$ where $\pi$ precedes $\gs$.
It is  enough to check this when the two orders differ only in that two neighboring roots are switched.  By definition of the Clifford filtration any even root vector is central modulo lower degree terms. Also for the distinguished or anti-distinguished  Borel subalgebra the commutator of two root vectors corresponding to
roots in $\overline{\Delta}^+_{1}$ is zero.
 So it is enough to check the case where the two roots are $\gd_i$ and $\gd_j$ with $i\neq j$. Because an even power of
$e_{-\gd_i}$ or $ e_{-\gd_j}$  is central modulo lower degree terms it is enough to check the case where
$\pi({\gd_i}) = \pi({\gd_j})=1.$ Here we have $[e_{-\gd_i}, e_{-\gd_j}] = e_{-\gd_i-\gd_j}$ up to a non-zero scalar multiple.  Then
\[{e}_{-\pi} =-\overline{e}_{-\pi} \pm\overline{e}_{-\gs}  \mbox{ modulo terms of lower degree in the Clifford filtration,} \] where $\gs$ is defined by
\be \label{2yam}
 \gs(\ga)=\left\{ \begin{array}
  {ccl}\pi(\ga) -1 &\mbox{if} \;\; \ga = \gd_i \mbox{ or } \gd_j,
\\ \pi(\ga) + 1 &\mbox{if} \;\; \ga = \gd_i + \gd_j,\\
\pi(\ga)& \mbox{ otherwise.}
\end{array} \right.\nn \ee
\noi Now the claim follows since $\cg(\pi) = \cg(\gs)$, but $|\pi| > |\gs|$.
The rest of the  proof is the same as the proof of Lemma {\rm \ref{1delorder}}.\epf
\noi If $\gc$ is a simple root, then $\theta_{\gamma} = e_{-\gc}$ satisfies the conditions of Theorem \ref{1aShap}. We recall for convenience Equations \eqref{1wuga} and \eqref{1wuga1} with minor modifications. Write $\gc = w\gb$ for some
$w \in W_{\rm nonisotropic}, w\neq 1$. Assume that
\be w =s_{\ga}u, \quad   \gamma' = u\gb, \quad \gamma = w\gb = s_{\alpha} \gamma',\nn
\ee
 with $\alpha \in \Pi_{\rm nonisotropic}$ and $\ell(w) =\ell(u) + 1.$
Suppose $\nu \in \Lambda$ and set  \be
\mu =
u\cdot \nu, \quad \lambda = w\cdot \nu
= s_{\alpha}\cdot \mu.\nn
 \ee
As before we assume that the positive roots are ordered so that for any partition $\pi$ we have $e_{-\pi} = e_{-\gs}e_{-\ga}^k$ for a non-negative integer $k$ and
 a partition  $\gs$ with $\gs(\ga)=0$.

\bl \label{1wpfg1} Suppose that $p, q$ and $\ga$ are as in Lemma \ref{11768} and
\begin{equation} \label{13nd} e^{p+q}_{-
\alpha}\theta^{\mu}_{\gamma'} = \theta^\lambda_{\gamma}e^p_{-
\alpha}.
\end{equation} 
Then the family $\theta^\gl_{\gamma} $ is in good position for $w$ if
the family $\theta^\mu_{\gamma'}$ is in good position for $u$. \el \bpf
If $\ga$ is odd, then  $p = 2\ell -1$ is odd,  and by Lemma \ref{1gag}, $q =2$.
Write $\theta^{\mu}_{\gamma'}$
as in (\ref{14cha7}) and then define
the $e_{- \pi'}^{(j)}$ as in (\ref{1qwer}).
Set $\gve(\gc') = 1$ if $\gc'$ is an even root and $\gve(\gc') = -1$ if $\gc'$ is odd. Then instead of (\ref{149}) we have, by (\ref{1f11})
\be \label{1f12}  e_{-\ga}^{2\ell+1}e_{- \pi'} = \sum^\ell _{i=0} \left( \begin{array}{c}
                \ell \\
                i \end{array}\right) e_{- \pi'}^{[2i+1]}e_{-\ga}^{2(\ell - i)}
 + \gve(\gc')\sum^\ell _{i=0} \left( \begin{array}{c}
                \ell \\
i \end{array}\right) e_{- \pi'}^{[2i]} e_{-\ga}^{2\ell - 2i +1} .\ee
\noi Parallel to the definition of  the $b_{j, \zeta}^{\pi'}$ in (\ref{199}), we set for sufficiently large $N$
\[e_{- \pi'}^{[j]}e^{N-j}_{- \alpha}= \sum_{\gz \in {\overline{\bf P}}(\gamma'+N\ga)} b_{j, \zeta}^{\pi'}e_{- \gz}.\]
For $x \in \mathbb{R}$ we denote the largest integer not greater than $x$ by $\left\lfloor x \right\rfloor$. Then if  $ b_{j,
\zeta}^{\pi'} \neq 0,$ we have
\begin{eqnarray} \label{198o}
\cg(\zeta) \leq \cg(\pi') + N - 2\left\lfloor \frac{j}{2} \right\rfloor.
\end{eqnarray}
Indeed this holds because by Lemma \ref{1it}, we have for such $j$
\be \label{1la1} e_{- \pi'}^{[j]}e^{N-j}_{- \alpha} \in  U_{\cg(\pi') + N - 2\left\lfloor \frac{j}{2} \right\rfloor}.\ee
 Define coefficients $a'_{\pi',\mu}$ 
by \begin{equation} \theta^{\mu}_{\gamma',m} =
\sum_{\pi' \in {\overline{\bf P}}(m\gamma')} a'_{\pi',\mu} e_{-\pi'},
\end{equation} 
as in \eqref{14cha7}. Then  
replacing (\ref{159}) we set,
\be\label{159o} c_{ \zeta,\lambda}
= \sum_{\pi' \in {\overline{\bf P}}(\gamma'),\;i \geq 0} \left(
\begin{array}{c}
\ell \\i\end{array}
\right) a'_{\pi',\mu} b_{2i+1, \zeta}^{\pi'}
+ \gve(\gc')\sum_{\pi' \in {\overline{\bf P}}(\gamma'),\;i \geq 0} \left(
\begin{array}{c}
\ell \\i\end{array}
\right) a'_{\pi',\mu} b_{2i, \zeta}^{\pi'}. \ee
Then we obtain the following variant of Equation (\ref{1lab})
\begin{eqnarray} \label{1labos} e^{2\ell+1}_{- \alpha} \theta^{\mu}_{\gamma'} &
= & \sum_{\zeta \in {\overline{\bf P}}(m\gc + N\ga )}
c_{ \zeta, \lambda}e_{-\zeta}e^{2\ell+1-N}_{- \alpha}.
\end{eqnarray}
In the cancellation step we find that $c_{ \zeta,\lambda}=0$ unless $\gz(\ga)\ge N-2$, and the bijection
\[f:{\overline{\bf P}}(\gamma) \longrightarrow \{\zeta \in {\overline{\bf P}}(\gamma' +N\ga)\;|\;\zeta(\ga) \geq N-2\},\]
is defined as in Equation (\ref{1yam}) with $m=1$ and $q=2.$ Then the coefficients  $a_{\pi,\lambda}$ are defined as in (\ref{1cwp2}).
Instead of Equations (\ref{1119})  and  (\ref{1fox}) we have, when $f\pi = \zeta,$ 
\begin{eqnarray} \label{1119o}\cg(\zeta) = \cg(\pi) + N -2,
\end{eqnarray}
and \be \label{1foxa} \; \mbox{deg} \;\; a_{\pi,\gl} \;
 \leq  \max\{\lfloor j/2 \rfloor + \; \mbox{deg} \;
a'_{\pi',\mu}\;|\; b_{j, \zeta}^{\pi'}
\neq 0\}. \ee
Hence using (\ref{198o}) in place of (\ref{198}), and then (\ref{1119o}) we obtain,
\begin{eqnarray}  \label{1scat} 2\mbox{deg} \;\; a_{\pi,\gl}
&\leq& 2\mbox{deg} \;\; a'_{\pi',\mu}+ \cg(\pi') + N - \cg(\gz)\\
&=& 2\mbox{deg} \;\; a'_{\pi',\mu}+ \cg(\pi')  - \cg(\pi) +2.\nonumber
\end{eqnarray}
Therefore by induction, and since $\cg(\gc') =\cg(\gc) $, 
\begin{eqnarray}  \label{1scat1} 2\mbox{deg} \;\; a_{\pi,\gl}
&\leq&
2\ell(u)+\cg(\gc')  - \cg(\pi) + 2\nonumber\\
&=&
2\ell(w)+\cg(\gc)  - \cg(\pi).\nonumber
\end{eqnarray}
giving  condition  (a) in Definition \ref{1sog}.
\\ \\
The proof in the case where $\ga$ is an even root is the same as in Section \ref{1s.5} apart from the inequalities.  If  $ b_{j,
\zeta}^{\pi'} \neq 0,$ then instead of  (\ref{198}), we have
\begin{eqnarray} \label{198e}
\cg(\zeta) \leq \cg(\pi') + 2(N - j).
\end{eqnarray}
Define 
a bijection
\[f:{\overline{\bf P}}(\gamma) \longrightarrow \{\zeta \in {\overline{\bf P}}(\gamma' +N\ga)\;|\;\zeta(\ga) \geq N-1\},\]
as in \eqref{1yam} with $m=q=1$.
 Then condition (a) follows since in place of (\ref{1119}) we have, if $f\pi =\gz$
\begin{eqnarray} \label{11199}\cg(\zeta) = \cg(\pi) + 2(N -1).
\end{eqnarray}
We leave the proof that (b) holds in  Definition \ref{1sog} 
to the reader.
\epf 
\noi Theorem \ref{1aShap} follows from Lemma \ref{1ant} in the same way that Theorem \ref{1Shap} follows from Lemma \ref{1wpfg}.

\br \label{cap}{\rm 
\v Sapovalov elements and their evaluations can also be constructed inside 
a Kostant $\Z-$form of the algebras $U(\fb^-)$ and $U(\fn^-)$. 
This is done for Lie algebras in  \cite{KL}. We briefly indicate how to do it for Lie superalgebras using the $\Z-$form given by Fioresi and Gavarini \cite{FG}.  \ff{The point of doing this is to 
allow for a change in the base ring, and in particular to enable passage to positive characteristic.} 
For $e \in \fn_0^-$, $h\in \fh\oplus \ttk$ and $b, c\in \N$, define 
\[e^{(b)}=\frac{e^b}{b !}\in U(\fn^-) \mbox{ and } 
\left(\begin{array}{c} h \\ c \end{array}\right) = \frac{h(h-1)\ldots (h-c+1)}{c!}\in U(\fh).\] 
Let $\{x_1, \ldots x_p\}$, $\{e_1, \ldots e_q\}$ and $\{h_1, \ldots h_n\}$ be bases for $\fn_1^-, \fn_0^-$ and $\fh$ respectively. Then for $A\in \Z_2^p, B \in \N^q$ and $C \in \N^n$,  
consider  the elements $x_A, e_B, h_C$ of $U(\fb^-)$ given by 
\[x_A = x_1^{a_1}\ldots x_p^{a_p}, 
\quad e_B = \frac{e_1}{b_1 !}\ldots \frac{e_q}{b_q !}, 
\quad h_C = \left(\begin{array}{c} h_1 \\ c_1 \end{array}\right)\ldots \left(\begin{array}{c} h_n \\ c_n \end{array}\right).\] 
By \cite{FG} Theorem 4.7, the products 
$x_A e_B h_C$ (resp. $x_A e_B $, $h_C$) form a $\ttk$-basis of 
$U(\fb^-)$ (resp. $U(\fn^-)$, $S(\fh)$)
 and the $\Z$-span of the set of all such products form a 
$\Z$-algebra $U_\Z(\fb^-)$ (resp. $U_\Z(\fn^-)$, $S_\Z (\fh)$). 
If $(\ad e)^{(b)}=\frac{(\ad e)^b}{b !},$ then by the proof of \cite{H} Proposition   25.5, the operator $(\ad e)^{(b)}$  leaves $U_\Z(\fn^-)$ and $U_\Z(\fb^-)$ invariant.  Now \eqref{1cow} may be written in the form 
\be 
 e^{(r)}a = \sum^r_{i=0} ((\ad e)^{(i)} a)e^{(r - i)}.\nn\ee
and \eqref{121nd} is equivalent to

\begin{equation} 
e^{(p)}_{- \alpha}e^{mq}_{- \alpha}\theta_{\gc',m} = \theta_{\gc,m} e^{(p)}_{- \alpha}.\nn
\end{equation}
Turning to \eqref{1f11} and using the notation introduced immediately before that equation, we define
\by x^{(2i)}  = e^{(i)}= e^i/i!= x^{2i+1}/i!
&&  x^{(2i+1)} = x^{(2i+1)}/i!\nn
\\ z^{^{(2i)}}  =\frac{(\ad x)^{2i}z}{i!}
&&z^{(2i+1)} = \frac{(\ad x)^{2i+1}z}{i!}\nn 
\ey
Then \eqref{1f11} can be written in the form 
\be 
x^{(2\ell+1)}z = \sum^\ell _{i=0} z^{(2i+1)}x^{(2\ell - 2i)}
 + (-1)^{\overline z}\sum^\ell _{i=0} z^{(2i)} x^{(2\ell - 2i +1)}.\nn
\ee
Finally we note that if $e= e_{-\ga}$  is a root vector and $h \in \fh,$
we have 
\[e\left(\begin{array}{c} h \\ c \end{array}\right) = \left(\begin{array}{c} h +\ga(h)\\ c \end{array}\right)e,\] 
and by \cite{H} Lemma 26.1, $\left(\begin{array}{c} h +\ga(h)\\ c \end{array}\right) \in S_\Z (\fh)$.
Using these equations and the inductive construction, 
compare \eqref{121nd} and \eqref{13nd}, it is easy to construct a \v{S}apovalov element 
$\theta_{\gamma,m}$ inside $U_\Z(\fb^-).$ 
}\er

\section{Changing the Borel subalgebra.}\label{1sscbs}
\subsection{Adjacent Borel subalgebras.}
 We consider the behavior of \v Sapovalov elements when the Borel subalgebra is changed.
Let ${\fb'},{\fb}''$ be arbitrary adjacent Borel subalgebras, 
 and suppose
\be \label{1fggb} \fg^\ga \subset {\fb}', \quad \fg^{-\ga} \subset {\fb''}\ee for some isotropic root $\ga$.
Let $S$ be the intersection of the sets of roots of  $\fb'$ and $\fb''$, $\fp =\fb' + \fb''$ and $\fr  = \bigoplus_{\gb \in S} \ttk e_{-\gb}.$ Then $\fr, \fp$ are subalgebras of $\fg$ with $\fg = \fp \oplus \fr$.  Furthermore $\fr$ is stable under $\ad e_{\pm \ga},$ and consequently, so is $U(\fr)$.
 Note that  \be \label{1rgb}\rho({\fb}'') = \gr({\fb'}) + \ga.\ee
Suppose $v_\mu$ is a highest weight vector for $\fb'$ with weight  $\mu=\mu(\fb')$, and let $N$ be the module generated by $v_\mu$.  Then \be \label{1ebl} e_\ga e_{-\ga} v_\mu = h_\ga v_\mu =(\mu, \ga)v_\mu . \ee
As is well known (see for example Corollary 8.6.3 in \cite{M}),  if $(\mu +\rho({\fb'}), \ga) \neq 0$, then
\be \label{1mon}\mu({\fb''}) + \rho({\fb''}) = \mu({\fb}') + \rho({\fb}').\ee
In this situation we call the change of Borel subalgebras from ${\fb'}$ to ${\fb}''$ (or vice-versa) a {\it typical change of Borels for $N$}.
Suppose that $\gc$ is a positive root of both ${\fb'}$ and ${\fb}''$, and that
$\theta_{\gamma,m}$ is a \v Sapovalov element
 corresponding to the pair $(\gc, m)$ using the negatives of the roots of ${\fb}''$. 
Consider the Zariski dense subset $\gL_{\gc, m}$ of ${\mathcal H}_{\gc, m}$ given by
\[\gL_{\gc, m} =  \{ \mu \in
{\mathcal H}_{\gc, m}|(\mu + \rho, \ga) \notin \bbZ  \mbox{ for all positive roots } \ga \neq \gc\}.\]
Since the coefficients of $\theta_{\gamma,m}$ are polynomials,
$\theta_{\gamma,m}$ is determined (as usual modulo a left ideal) by the values of
$\theta_{\gamma,m}v_\mu$ for
$\mu \in \gL_{\gc, m}$.
\\ \\
Assume $\mu \in \gL_{\gc, m}$, and for brevity set $\gth =\theta_{\gamma,m}(\mu).$ Let $v_{\mu'}$ be a highest weight  vector in a Verma  module $M_{\sfb'}(\mu')$ for ${\fb'}$ with highest weight  ${\mu'}\in \gL_{\gc, m}.$
Then $v_\mu =e_{-\ga} v_{\mu'}$ is  a highest weight  vector for ${\fb}''$
which also generates $M_{{ \sfb'}}(\mu')$.  Thus we can write
\[M_{{ \sfb'}}(\mu')=M_{{ \sfb''}}(\mu).\]
Next note that $u = \theta e_{-\ga} v_{\mu'}$ is  a highest weight  vector for ${\fb}''$, and $e_{\ga} \theta e_{-\ga} v_{\mu'} $ is a highest weight vector for ${\fb'}$ of weight $\mu' -m\gamma$ that generates the same submodule of $M_{{ \sfb'}}(\mu')$ as $u$.
We can write $\gth$ in a unique way as $\theta=e_\ga \theta_{1} + \theta_{2}$  with $\gth_i \in U(\fr)$. Then
\begin{eqnarray} \label{1slab}
 e_{\ga}\theta e_{-\ga} v_{\mu'} &=& e_{\ga}\theta_2 e_{-\ga} v_{\mu'} \nonumber\\
 &=& \gth'_1 e_{-\ga} v_{\mu'} \pm  \theta'_2 v_{\mu'} \nonumber
 \end{eqnarray}
where $\gth'_1 = [e_{\ga},\theta_2], \; \; \gth'_2 = (\mu',\ga)\theta_2 \in U(\fr).$
 Note that the term $e_{-m\pi^\gc}$ cannot occur in $e_\ga \gth_1$ or $\gth'_1 e_{-\ga}. $ 
  Allowing for possible re-order of positive roots used to define the $e_{-\pi}$
	(compare Lemma \ref{1delorder}) we conclude that modulo terms of lower degree,
	the coefficient of $e_{-m\pi^\gc}$ in $e_{\ga}\theta e_{-\ga} v_{\mu'}$ is equal to $\pm(\mu',\ga)$ times the coefficient of $e_{-m\pi^\gc}$ in  $\theta_2 e_\ga v_{\mu'}$.

\subsection{Chains of Borel subalgebras.}
Using adjacent Borel subalgebras it is possible to give an alternative construction of \v Sapovalov elements corresponding to an isotropic root $\gc$ which is a  simple root for some Borel subalgebra.  This condition always holds in type A, but for other types,  it is quite restrictive:  if $\fg = \osp(2m,2n+1)$ the assumption only holds
for roots of the form $\pm(\gep_i-\gd_j)$, while if $\fg = \osp(2m,2n)$ it holds only
for these roots and the root $\gep_m+\gd_n$. (Theorems  \ref{1Shap} and  \ref{1aShap} on the other hand apply to any positive isotropic root, provided we choose the appropriate Borel subalgebra satisfying Hypothesis (\ref{i}).)
\\ \\
Suppose that  
$\fb$ is a distinguished or anti-distinguished  Borel subalgebra, and let $\mathfrak{b}'\in \cB.$
Consider  a  sequence 
 \be \label{distm} \mathfrak{b} = \mathfrak{b}^{(0)}, \mathfrak{b}^{(1)}, \ldots, \mathfrak{b}^{(r)}=\fb' \ee
  of Borel subalgebras such that $\mathfrak{b}^{(i-1)}$ and $\mathfrak{b}^{(i)}$ are adjacent for $1 \leq i \leq r$.
  It follows from \cite{M} Theorem 3.1.3, that such a chain always exists. If there is no  chain of adjacent Borel subalgebras connecting $\fb$ and $\fb'$ of shorter length than (\ref{distm}), we set $d(\fb,\fb')=r.$
  In this case there are isotropic roots $\gb_i$ such that 
	\be \label{bar}\fg^{\gb_i} \subset \fb^{(i-1)}, \quad \fg^{-\gb_i} \subset \fb^{(i)}\ee
	for $i\in [r]$, and  $\gb_1,\ldots,\gb_r$ are distinct positive roots of $\fb.$ 
	Now suppose that $\gc$ is a positive isotropic root for $\fb$ and a simple root in some Borel subalgebra. Then define 
\be 
d(\gc)= \min\left\{\;r \;\vline  \mbox{ for some chain as in  \eqref{distm}, } \gc \mbox{ is a simple root in } \fb^{(r)}\right\}.  \ee \noi 
Suppose that $\gc$ is  a simple isotropic root of $\fb^{(r)}$, 
where $d(\gc)=r,$ and suppose that $\gl$ is a {\it general} element of $\cH_\gc$.  Specifically we take this to mean that each change of Borels in the chain \eqref{distm} is typical for the Verma modules  $M(\gl), M(\gl-\gc)$, induced from $\fb.$ 
The set of all such $\gl$ is Zariski dense in $\cH_\gc$. 
	Then set  

	\be \label{tom}v_0 = v_\gl, \;\;\gl_{0} = \gl, \;\; v_i= e_{-\gb_i}v_{i-1},\;\; \gl_{i} = \gl_{i-1}-\gb_i \;\mbox{ for }\; i \in [r],\ee
and 
\be \label{rom} u_r =e_{-\gc}v_r,\;\;  u_i = e_{\gb_{i+1}}\ldots e_{\gb_r}u_r \;\;\mbox{ for } \;\;0\le i \le r-1.\ee
 Then $u_i$ and $v_i$ are  highest weight vectors for the Borel subalgebra $\fb^{(i)}$ with weights $\gl_{i}$ and $\gl_{i}-\gc$ respectively. Since each change of Borels is typical, it follows that $(\gl_i+\gr_i,\gc)=0$ for all $i$, where $\gr_i$ is the analog of $\gr$ for 
$\mathfrak{b}^{(i)}.$ 
The observations of the previous Subsection give the following.
\bl\label{1sea} The leading term of the coefficient of $e_{-\gc}v_\gl$ in $$e_{\gb_{1}}\ldots e_{\gb_r}e_{-\gc}e_{-\gb_{r}}\ldots e_{-\gb_1}v_\gl$$ is, up to a constant multiple, equal to $\prod_{i=1}^r (\gl,\gb_i).$ \el
\noi 
\subsection{Relation to the Weyl group.}
When $\gc$ is a simple root of some Borel subalgebra, we relate the approach to \v Sapovalov elements by change of Borel to the approach using the Weyl group. 
The analysis reveals a set of roots $N(\gc)$ which has properties analogous to those of the set $N(w^{-1})$ from \eqref{nu}. 
\\ \\
Suppose that  
$\fb$ is the distinguished or anti-distinguished  Borel subalgebra, and let $\mathfrak{b}'\in \cB.$
Consider  the chain of Borels from \eqref{distm} and the roots $\gb_1,\ldots,\gb_r$ from \eqref{bar}.
\noi Recall the definition of $q(w,\ga)$ before Theorem  \ref{1Shap}.  The main result of this Subsection is the following. 
\bp \label{1alp}Suppose $\gc = w\gb$ with
$w \in W_{\rm nonisotropic}$ and $\gb$ a simple isotropic root, and that $\gc$ is a simple root of the Borel subalgebra $\mathfrak{b}^{(r)}$, where $d(\gc)=r$. Then
\[\{\gc-\gb_i| i=1, \ldots, r, (\gc,\gb_i) \neq 0\} = N(w^{-1}).
\] 
Moreover  $q(w,\ga) =1$ for all $\ga \in N(w^{-1}).$ 
\ep
\noi A corollary will be used in the next Subsection. Before proving the Proposition we need several preparatory results. The first is familiar from the Lie algebra case. 
\bl \label{H84} If $\gc, \ga$ are roots with $\ga$ non-isotropic, and the $\ga$-string through $\gc$ is 
\[ \gc -p\ga, \ldots, \gc, \ldots, \gc + q\ga\] then $(\gc, \ga^\vee) = p-q$. \el
\bpf Let $\fs$ be the subalgebra of $\fg$ generated by the root vectors $e_\ga$ and $e_{-\ga}$. Then $\fs$ is isomorphic to $\fsl(2)$ or $\osp(1,2)$.  In the first case this is proved as in the Lie algebra case, \cite{H} Prop. 8.4 (e) using the representation theory of $\fsl(2)$.  Essentially the same proof also works when $\fs \cong \osp(1,2)$.\epf

\bl \label{new}Suppose $d(\gc) =r$ and $\gc$ is a simple root of the Borel subalgebra $\mathfrak{b}^{(r)}$. Then
$\ga:=  \gc-\gb_r$ is a simple even root of $\mathfrak{b}^{(r-1)}$ and $\gc':= s_\ga \gc = \gc-\ga = \gb_r$.
\el
\bpf We show that the diagram below on the left (resp. right) is part of the Dynkin-Kac diagram for $\fb^{(r-1)}$ (resp. $\fb^{(r)}$).

\begin{picture}(60,40)(-53,-20)
\thinlines
\put(44,3){\line(1,0){39.8}}
\put(39.7,2.3){\circle{8}}
\put(83,-15){$\gb_r$}
\put(38,-15){$\ga$}
\put(205,-15){$-\gb_r$}
\put(170,-15){$\gc$}
\put(169,-0.3){$\otimes$}
\put(177,3){\line(1,0){36}}
\put(212.7,-0.3){$\otimes$}
\put(83.7,-0.3){$\otimes$}
\end{picture}

\noi Since $d(\gc)=r, \gc$ is not a simple root of $\fb^{(r-1)}$ and thus $\gc$ is connected to $-\gb_r$ as in the second diagram. By comparing the diagrams to those in Table 3.4.1 and the first row of Table 3.5.1 in \cite{M}, we see that the first diagram is the only possibility, 
and the $\ga$-string through $\gb_r$ consists of $\gb_r$ and $\gc$. Therefore by Lemma \ref{H84} $(\gc,\ga^\vee)=1,$ and the assertions about $\gc'$ follow easily.\epf

\noi Now suppose $\fb$ is any Borel in $\cB$ and set
\be Q(\fb)=\{\ga|\ga \in \Gd_1^+(\fb^\dist)\cap -\Gd_1^+(\fb)\}.\ee
and $Q(\gc) = Q(\fb)$ where $d(\gc)= d(\fb^\dist,\fb).$ 
Next define 
\be N(\gc) = \{\ga \in Q(\gc) |(\gc, \ga)\neq 0\}.\ee
This set is analogous to the set  $N(w^{-1})$ from  \eqref{nu}.

\bl \label{nef} In the situation of Lemma \ref{new}, we have a disjoint union
\[N(\gc) = s_\ga N(\gc')\ds \{\gc'\}.\] 
\el
\bpf

\noi Without loss of generality, we can assume that $\fg$ has type A, $\ga= \gep_{i}-\gep_{i+1}$, $\gc'= \gep_{i+1}-\gd_{j'}$ and $\gc = \gep_{i}-\gd_{j'}$. Set
\[R(\gc) = \{(k, \ell')| \gep_{k}-\gd_{\ell'} \in N(\gc)\}.\]
Then we need to show
\be \label{def}R(\gc) = \gt_{i,i+1}R(\gc')\ds  \{(i+1,j')\},\ee where $\gt_{i,i+1}$ 
is the transposition switching $i$ and $i+1$.
We can write the shuffles corresponding to $\gc$ and $\gc'$ as the concatenations of
\be \label{hon}A,A',i,j',i+1, B, B'\;\; \mbox{ and }\;\; A,i,A',i+1, j',B,B'\ee
respectively, where 
\[A = (1,\ldots,i-1), \;\;B= (i+2,\ldots, m), \]
\[ A' = (1',\ldots,j'-1),\quad B'=(j'+1,\ldots,n').\]
Thus 
\by 
R(\gc) &=& \{i\}\ti A' \;\;\bigcup \;\;(\{i+1\}\ds B) \;\;\ti\{j'\}, \nn\\
 R(\gc') &=& \{i+1\}\ti A' \;\;\bigcup B \;\;\ti\{j'\}\nn,\ey 
\noi and this clearly gives \eqref{def}.
\epf 

\noi {\it Proof of Proposition \ref{1alp}}. Suppose $\gc = w\gb$ with
$w \in W_{\rm nonisotropic}$ and $\gb$  simple isotropic. 
The first statement is clearly true for $\gc$ simple.  As in Lemma \ref{new} we assume that $\gc'= s_\ga \gc = \gc-\ga = \gb_r$. 
Set $A_w = \{\gc-\gb_i| i=1, \ldots, r, (\gc,\gb_i) \neq 0\} $ and $u= s_\ga w$.   
We show by induction on $|A_w|$ that $A_w =  s_\ga A_u  \ds \{\ga\},$  a disjoint union. 
The first equality below comes from the definition of $A_w$, the second from Lemma \ref{nef}, the third  is a simple rearrangement, and the fourth holds by induction. 
\by \label{uni} A_w &=& \{\gc-\gs|\gs\in N(\gc)\}  \nn\\
&=&  \{\gc-\gs|\gs\in s_\ga N(\gc')\} \ds  \{\gc-\gc'\}  \nn\\
&=& s_\ga\{\gc'-\gt|\gt\in N(\gc')\}  \ds \{\ga\} \nn\\
&=& s_\ga A_u  \ds \{\ga\}.\ey
This shows that $A_w$ satisfies the same recurrence as $N(w^{-1})$ from Equation (\ref{Nw}), and this gives the first statement. (Note also that since the union in  \eqref{uni} is disjoint, we have $\ell(w) = \ell(u)+1$.) 
By induction $q(u,\gs)=1$ for all $\ga \in N(u^{-1})$.  Thus $q(w, s_\ga\gs)=1$ for all such $\gs$. Since $q(w,\ga) =(w\gb, \ga^\vee)=(\gc,\ga^\vee)=1$ by Lemma \ref{new}, the second statement holds.
\hfill  $\Box$
\bc \label{ax}The leading term of the coefficient of $e_{-\gc}v_\gl$ in $\gth_\gc v_\gl $ is given by
$\prod^r_{i=1, (\gc,\gb_i) \neq 0} (\gl,\gb_i)$
  up to a scalar multiple, independently of $\gl\in \cH_\gc$.\ec
\bpf We work always up to a scalar multiple. Under the hypothesis of Theorem \ref{1aShap}  the leading term  equals 
$\prod_{\gb \in N(w^{-1})}(\gl,\gb)$.  
If instead we have the hypothesis of Theorem \ref{1Shap}, then since $q(w,\ga) =1$ for all $\ga \in N(w^{-1})$ 
by Proposition \ref{1alp}, we obtain the same leading term. Furthermore 
\[\prod_{\gb \in N(w^{-1})}(\gl,\gb)=\prod_{i: (\gc,\gb_i) \neq 0} (\gl,\gc-\gb_i).\]
\noi Since $(\gl+\gr,\gc) =0$ this gives the result.\epf

\subsection{\v Sapovalov elements and change of Borel. }
There exists a unique (modulo a suitable left ideal in $U(\fg)$) \v Sapovalov element $\gth_\gc^{(i)}$ for the Borel subalgebra $\fb^{(i)}$ and polynomials $g_i(\gl), h_i(\gl)$ such that
\be \label{1tlc}
e_{-\gb_i}\gth_\gc^{(i-1)}
e_{\gb_i}v_{i} = g_i(\gl)\gth_\gc^{(i)}v_{i},
\ee and
\be \label{1tlc1}
e_{\gb_i}\gth_\gc^{(i)}e_{-\gb_i}v_{i-1} = h_i(\gl)\gth_\gc^{(i-1)}v_{i-1}.
\ee
\bl \label{1112} For general $\gl\in \cH_\gc,$
\be \label{zxy}g_i(\gl)h_i(\gl)=
(\gl +\gr,\gb_{i})(\gl+\gr-\gc,\gb_{i}).\ee
\el
\bpf We have using first \eqref{1ebl} and the definition of $v_i$, then \eqref{1tlc} and \eqref{1tlc1},
\begin{eqnarray}
\label{1fff}
(\gl_{i-1}, \gb_{i} )(\gl_{i-1}-\gc,\gb_{i})\gth^{(i-1)}_\gc v_{i-1}
&=& e_{\gb_{i}}e_{-\gb_{i}} \gth_\gc^{(i-1)} e_{\gb_{i}}e_{-\gb_{i}}v_{i-1}
\nonumber \\
&=& e_{\gb_{i}}e_{-\gb_{i}} \gth_\gc^{(i-1)} e_{\gb_{i}}v_{i}
\nonumber \\
&=&
g_i(\gl)e_{\gb_{i}} \gth_\gc^{(i)} e_{-\gb_{i}} v_{i-1}
\nonumber \\
&=&
g_i(\gl)h_i(\gl)\gth_\gc^{(i-1)}v_{i-1}
.\nonumber \end{eqnarray}
\noi 
Now  each change of Borels is typical and  $\gb_{i}$ is simple isotropic for $\fb^{(i-1)}$, so
we have  $(\gl_{i-1}, \gb_{i} )$=$(\gl +\gr,\gb_{i}).$ The result follows.
\epf
 \bt \label{1rot} Set $F(\gc) = \{i\in [r]|(\gc,\gb_i)=0  \}$. There is a nonzero $c\in \ttk$ such that for all $\gl \in \cH_\gc$,

\[e_{\gb_{1}}\ldots e_{\gb_r}e_{-\gc}e_{-\gb_{r}}\ldots e_{-\gb_1}v_\gl= c\prod_{i \in F(\gc)}(\gl+\gr,\gb_i)
\gth_\gc v_\gl.\]
 \et
\bpf
It is enough to show this for all $\gl$ in a Zariski
dense subset of $\cH_\gc$. Thus we may assume that each change of Borels in \eqref{distm} is typical for $M(\gl)$.
Since $v_i= e_{-\gb_i}v_{i-1},$ and $u_{i-1}=
e_{\gb_{i}}u_{i},$ Equation (\ref{1tlc1}) and reverse induction on $i$ yield
  \be \label{1pt} u_{i-1}=\prod_{j=i}^r h_j(\gl)\gth_\gc^{(i-1)}v_{i-1}. \ee Hence
  \be \label{1nat} u_0 = \prod_{j=1}^r h_j(\gl) \gth_\gc v_{0}. \ee
Therefore by comparing the coefficient of $e_{-\gc} v_\gl$ on both sides of (\ref{1nat}), and using Lemma \ref{1sea} and Corollary \ref{ax}, we have modulo terms of lower degree, that
\be \label{1tat} \prod_{j=1}^r h_j(\gl)  = \prod_{i\in F(\gc)} (\gl,\gb_i). \ee
Now none of the functions $\gl \lra (\gl,\gb_i)$, for $i\in [r] $, on
$\cH_\gc$ is a multiple of another.
It follows that, up to a constant multiple,   $h_j(\gl) = (\gl,\gb_{j}) + \mbox{ a constant}$
    if $j \in F(\gc)$, and $h_j(\gl)$ is constant if
    $j \notin F(\gc).$
However we know from Lemma \ref{1112} that if
$(\gl,\gb_{j}) = 0$ then  $h_j(\gl)$ divides  $(\gl+\gr,\gb_{j}).$ Thus the result follows.
\epf \noi
Next suppose that $\fb'', \fb'$ are adjacent Borel subalgebras as in Equation {\rm (\ref{1fggb})}, and that $d(\fb,\fb'')= d(\fb,\fb')+1$.  We can find a sequence of Borel subalgebras as in (\ref{distm}) such that $\fb'=\mathfrak{b}^{(i-1)}$ and
  $\fb''=\mathfrak{b}^{(i)}$ .  Adopting the notation of Equations (\ref{1tlc}) and (\ref{1tlc1}),
  we can now clarify the relationship between the \v Sapovalov elements
$\gth_\gc^{(i-1)}$ and $\gth_\gc^{(i)}$.
\bc With the above notation, we have up to constant multiples,
 \bi \itema If $(\gc,\gb_i) = 0$, then $h_i(\gl) =  g_i(\gl) = (\gl +\gr,\gb_{i})$. \itemb
If $(\gc,\gb_{i}) \neq 0$, then $h_i(\gl) = 1$ and $g_i(\gl) = (\gl +\gr,\gb_{i})(\gl +\gr-\gc,\gb_{i}).$
\ei\ec
\bpf From the last paragraph of the proof it follows  that, up to a constant multiple, $h_j(\gl) = (\gl+\gr,\gb_{j})$
    if $j \in F(\gc)$. We obtain the result by looking at the degrees of both sides in 
	\eqref{1tat}.\epf

\section{Relations between \v Sapovalov elements.} \label{RS}
\subsection{Powers of \v Sapovalov elements.} \label{1zzprod}

\subsubsection{Isotropic Roots.}\noi First we record  an elementary but important property of the \v Sapovalov element $\theta_{\gamma}$ corresponding to an isotropic root ${\gamma}$.
\bt \label{1zprod} If $\gl \in \cH_\gc,$ then $ \theta^2_{\gamma} v_\gl = 0.$ Equivalently,
  $\theta_{\gamma}(\gl-\gc)\theta_{\gamma}(\gl) =0.$\et
\noindent \bpf It is enough to show this for all $\gl$ in a Zariski
dense subset of $\cH_\gc$. We assume $\gc = w\gb$ for $\gb \in \Pi$ where $\gb$ is isotropic, and $w\in  W_{\rm nonisotropic}$.
 Let $\gL$ be the subset of  ${\mathcal H}_{\beta}$ defined by
Equation  (\ref{1tar}), and suppose $\gl \in
w\cdot \Lambda$.
The proof is by induction on the length of $w.$ We can assume that
$w \neq 1.$ Replace $\mu$ with $\mu - \gc'$ and $\gl = s_\ga\cdot \mu$
with $s_\ga\cdot (\mu - \gc') = \gl - \gc$ in Equation (\ref{12nd}) or (\ref{13nd}).
Then $p$ is replaced by $p + q$ and we obtain
\[ e^{p + 2q}_{- \alpha}\theta^{\mu - \gc'}_{\gamma'} = \theta^{\lambda -
\gc}_{\gamma}e^{p + q}_{- \alpha}.\]
 Combining this with Equation
(\ref{12nd}) and using induction we have
\[0 = e^{p+ 2q}_{-\ga}\theta^{\mu - \gc'}_{\gamma'}\theta^{\mu}_{\gamma'} =
\theta^{\gl - \gc}_{\gamma} \theta^\gl_{\gamma}e^{p}_{-\ga} .\] The
result follows since $e_{-\ga}$ is not a zero divisor in
$U({\mathfrak n}^-).$\epf
\noi \br {\rm From the Theorem, if  $\gl \in \cH_\gc,$ there is  a sequence of maps
\be
\label{1let}  \ldots M(\gl-\gc)
 \stackrel{\psi_{\gl,\gc}}{\longrightarrow} M(\gl)
 \stackrel{\psi_{\gl+\gc,\gc}}{\longrightarrow} M(\gl+\gc) \ldots  \ee
such that the composite of two successive maps is zero.  
The map $\psi_{\gl,\gc}$ in (\ref{1let}) is defined by $\psi_{\gl,\gc}(x v_{\gl-\gc}) = x\gth_{\gc}v_\gl$. 
The complex \eqref{1let} can have  non-zero homology see Remark \ref{dig} (d).}\er
\subsubsection{Evaluation of \v{S}apovalov elements.} \label{nir} 
Up to this point we have only evaluated \v{S}apovalov element
$\theta_{\gamma,m}$ at points $\gl \in \cH_{m,\gc}$.
However to study the behavior of powers of  \v Sapovalov elements for non-isotropic roots,   we need to evaluate at arbitrary points, $\gl \in \fh^*$.  
Another situation where it is useful to do this is in the work of Carter \cite{Car} on the construction of orthogonal bases for non-integral Verma modules and simple modules in type A. 
However some care must be taken since the \v{S}apovalov element $\theta_{\gamma,m}$ is only defined modulo the ideal 
$U({\mathfrak b}^{-} ){\mathcal I}(\mathcal H_{\gc, m})$. 
\\ \\
\noi Suppose $e$ is locally $\ad $-nilpotent.  Equation \eqref{1cow} leads to a formula for conjugation by $e^r$ which extends to the 
$1$-parameter family of automorphism $\Gt_r, r\in \ttk$  of $A_S,$ given by 
\be \label{cows} \Gt_{r}(a) = \sum_{i\ge0} \left( \begin{array}{c}
                r \\
                i \end{array}\right) ((\ad e)^i a)e^{- i} .\ee
For a generalization, see \cite{Ma} Lemma 4.3. 
\\ \\
Fix a non-isotropic root $\ga$ and let $S$ be the multiplicative subset $S=\{e_{-\ga}^n\}$ in $U=U(\fn^-)$, 
and let $U_S$ be the corresponding Ore localization.
We use the generalized conjugation automorphisms $\Gt_x$ (for $x\in \ttk$) of $U_S$ defined in \eqref{cows} with $e =  e_{-\ga}$. For $x\in\N$ we have $\Gt_x(a) = e_{-\ga}^x a e_{-\ga}^{-x}.$ 
Thus  $\Gt_x(e_{-\ga}) = {e_{-\ga}}$, and 
\be\label{yam} \Gt_x(e_{-\ga}ue_{-\ga}^{-1}) = \Gt_{x+1}(u). \ee
As noted earlier, the value of $\Gt_x(u)$ for $x\in \ttk$ is determined by its values for $x\in\N$. To stress the dependence of $\Gt_x(u)$ on the root $\ga$ we sometimes denote it by $\Gt^\ga_x(u).$ 
\bl 
Suppose that  $\mu \in \cH_{\gc'}, $ $\gl = s_\alpha\cdot \mu,\;\;\gc = s_\alpha\gc',\;\;q = (\gc, \alpha^\vee).$
Assume that $ q  \in \mathbb{N}\backslash \{0\}$, and
if  ${\ga\in \Delta^+_{1} \backslash \overline{\Delta}^+_{1} }$, that $q = 2$.
If $x = (\mu + \rho, \alpha^\vee)$ we have
\be \label{bam} \theta_{\gamma}(\gl)= \Gt_{x}(e^{q}_{- \alpha}\theta_{\gamma'}(\mu)).\ee\el
\bpf If $\alpha \in\Pi\cap A(\mu)$ this follows from Lemma \ref{11768}.\epf
The issue is that we need to show that using  this definition, we need to show that
for  arbitrary points, $\gl \in \fh^*$, we have 
$\theta_{\gamma,m}(\gl)\in U$, and not just that $
\theta_{\gamma,m}(\gl)\in U_S,$ see  Corollary \ref{c4}. Fortunately there is an easy way to ensure this in all non-exceptional cases.\\ \\
  For $\gb \in \Pi$, let  $W_\gb$ be the subgroup of $W$ generated by   all simple reflections $s_\ga$ where $\ga$ is non-isotropic and  $\ga\neq\gb$. We consider a stronger version of hypotheses \eqref{iii} and \eqref{ii}.
\by \label{iv} && \mbox{The root }\gc \in \Gd^+ \mbox{ is such that for some }\\&\;\;& \gb \in \Pi, \mbox{  and } w\in W_\gb \mbox{ we have } \gc=w\gb. 
\nn\ey
\bl \label{sim} Hypothesis \eqref{iv} holds provided 
\bi \itema $\fg$ is a simple Lie algebra of type A-D, and $\gc$ is a positive root. 
\itemb $\fg$ is a contragredient Lie superalgebra, the basis $\Pi$ of positive roots satisfies \eqref{i} and $\gc$ is a positive isotropic root such that \eqref{ii} holds. 
\itemc $\fg$ is a contragredient Lie superalgebra of type A-D, the basis $\Pi$ of positive roots satisfies \eqref{i} and $\gc$ is a positive root such that \eqref{iii} or \eqref{ii} holds. 
\ei
\el
\bpf To prove (a), it suffices to exhibit a chain of roots of the form 
\be \label{jim} \gc=\gc_0, \gc_1, \ldots , \gc_m =\gb\ee
such that $\gb$ is simple, and that for $i \in [m]$, 
we have $\gc_i = s_{\gb_i} \gc_{i-1}$ for some $\gb_i \in \Pi \backslash \{\gb\}$ with $(\gc_i,\gb_i^\vee) > 0$.  
We use the same notation as Bourbaki \cite{Bo} Chapter 6.  
We treat the case where $\fg$ has type $B_n$ or $C_n$ first.  
Here there are two root lengths and all but one of the simple roots has the same length, the exception being the rightmost root $\ga_n$ of the Dynkin diagram.  
If $\gc$ and $\ga_n$ have the same length, then by \cite{H} Lemma 10.2, the result holds with $\gb= \ga_ n$. 
If this is not the case we have $\gc= \gep_i \pm \gep_j$ for some $i, j \in [n]$ with $i<j$.  By induction we can assume that $i=1$. If $\gc= \gep_1 + \gep_j$, the sequence of roots is 
\by \gc&=& \gep_1 + \gep_j,\; \gep_1 + \gep_{j+1},\; \ldots,\; \gep_1 + \gep_n, \nn\\
&&\gep_1 - \gep_n,\; \gep_1 - \gep_{n-1},\; \ldots,\; \gep_1 - \gep_2 =\gb. \nn\ey
It is easy to see that each term in the above sequence is obtained from its predecessor by a simple reflection, and that none of these reflections fixes the hyperplane orthogonal to $\gb$. If $\gc = \gep_1 - \gep_j$ we need only a terminal subsequence of the roots in the second row.  The same applies if $\fg$ has type $A_n$. If $\fg$ has type $D_n$, one such sequence is
\by \label{gym} \gc&=& \gep_1 + \gep_j,\; \gep_1 + \gep_{j+1},\; \ldots,\; \gep_1 + \gep_{n-1},\; \gep_1 + \gep_n, \\
&&\gep_1 - \gep_{n-1},\; \gep_1 - \gep_{n-2},\; \ldots,\; \gep_1 - \gep_2 =\gb. \nn\ey
Next note that (b) follows immediately since for $\gb$ isotropic we have $W_\gb = W_{\rm nonisotropic}.$ 
Finally to prove (c) we can assume that $\gc$ is non-isotropic, and then under the additional assumptions in (c) the we reduce to the statement in (a).\epf
\noi 
Now suppose that $\gc=w\gb$ is non-isotropic and  satisfies \eqref{iii} or \eqref{ii}. Suppose also that \eqref{iv} holds. We observe that in the inductive construction of the \v{S}apovalov element
$\theta_{\gamma,m}$, the assumption that $\gl \in \cH_{m,\gc}$ is only used for  the base case of the induction.
It follows that 
\bl If $\Pi_{\rm nonisotropic} \backslash \{\ga\} =\{\ga_1,\ldots, \ga_r\}$, and $\gl,\bar\gl \in \fh^*$ satisfy 
\be\label{yob}(\gl+\gr,\ga_i) =(\bar\gl+\gr,\ga_i) \mbox{ for } i\in[r],\ee
 then $\theta_{\gamma,m}(\gl) = \theta_{\gamma,m}(\bar\gl).$
\el
\bc \label{c4} For all $\gl \in \fh^*,$ and $j>0,$ we have $  \theta_{\gamma,j}(\gl)\in U$.\ec
\bpf 
Suppose $\bar\gl_j \in \cH_{\gamma,j}$ satisfies \eqref{yob} and note that $\theta_{\gamma,j}(\bar\gl) \in U$ by the cancellation step in the proof of Theorem \ref{1Shap}.\epf
\br {\rm We end this Subsection with a discussion of the uniqueness of the evaluations of  \v Sapovalov elements constructed above.  We first observe that even if $\fg=\fsl(3)$, and $\gc$ is the non-simple positive root, then unless $\gl\in \cH_{\gc,m}$ the construction depends on which of the simple roots we call $\gb$. We adopt the notation of \cite{M} Exercise 9.5.2. Thus the simple the positive roots are $\ga, \gb$ and $\gc = \ga +\gb.$ Choose negative root vectors $e_{-\ga} = e_{32}, e_{-\gb} = e_{21}$ and $e_{-\gc} = e_{31}.$ Suppose that 
$\mu =s_\ga \cdot\gl$, $\nu = s_\gb \cdot\gl$ and that $p=(\mu+\gr,\ga^\vee) =-(\gl+\gr,\ga^\vee) $, $q=(\nu+\gr,\gb^\vee) =-(\mu+\gr,\gb^\vee)$.  Then using $\gb,$ (resp. $\ga$) as the simple root in question, we obtain evaluations  $\theta_{\gc, 1}^{(\gb)}(\gl)$ (resp.
$\theta_{\gc, 1}^{(\ga)}(\gl)$)   
$$\theta_{\gc, 1}^{(\gb)}(\gl) = e_{-\gb}e_{-\ga}+(p+1)e_{-\gc} $$
$$\theta_{\gc, 1}^{(\ga)}(\gl) = e_{-\ga} e_{-\gb}-(q+1)e_{-\gc}.  $$
These are equal iff $p+q+1=0$, that is iff $\gl \in \cH_{\gc,1}.$ Coming back to the general case, we now fix the simple root  $\gb$. If  $\fg$ does not has type $A, B$ or $C$, then the sequence of roots in
\eqref{jim} is uniquely determined by $\gc$.  For $\fg = D_n$ the sequence is unique except that the subsequence 
\[ \gc_1=\gep_1 + \gep_{n-1},\; \gc_2=\gep_1 + \gep_n,\;\gc_3=\gep_1 - \gep_{n-1}\]
of
\eqref{gym} can be replaced by
\by 
\gc_1=\gep_1 + \gep_{n-1},\; \gc_2'=\gep_1 - \gep_n,\;
\gc_3=\gep_1 - \gep_{n-1}.\nn\ey  This is related to the fact that for $\fg= D_n$, the natural module is not uniserial as a module for a Borel subalgebra of $\fg$, \cite{M} Exercise 3.7.6. Let $\ga' = \gep_{n-1}- \gep_n$ and $\ga =   \gep_{n-1} + \gep_n,$ be the simple roots corresponding to the two rightmost nodes of the Dynkin diagram. We have 

\[ \gc_1= \gc_2 + \ga',\quad  \gc_2=\gc_3 + \ga, \] 
\[ \gc_1= \gc_2' + \ga,\quad  \gc_2'=\gc_3 + \ga'.\] 
Suppose that $\nu = s_\ga s_{\ga'} \cdot \mu.$   The construction of $\theta_{\gamma_1,m}(\nu) $ from  $
 \theta_{\gamma_3,m}(\mu) $ in 2 steps can be carried out in 2 different ways. We claim each way gives the same result. 
Suppose that $\gl = s_{\ga} \mu,$\;$ \gl' = s_{\ga'} \mu,$\; $ x = (\mu+\gr, \ga^\vee)$ and $y = (\mu+\gr, (\ga')^\vee).$ Then we need to show that 
\be \label{num} \Gt^\ga_x (e_{-\ga}^m\Gt^{\ga'}_y(e_{-\ga'}^m\gth_{\gc_3,m}(\mu)))
=
\Gt^{\ga'}_y (e_{-\ga'}^m\Gt^{\ga}_x(e_{-\ga}^m\gth_{\gc_3,m}(\mu))).\ee
 To show this set $\Gt^{\ga,\ga'}_{x,y} = \Gt^\ga_x \circ \Gt^{\ga'}_y.$
Since the roots $\ga$ and $\ga'$ are orthogonal, we have $\Gt^{\ga,\ga'}_{x,y} = \Gt^{\ga',\ga}_{y,x}$. Also $\Gt_x^\ga(e_{-\ga'}) = {e_{-\ga'}}$ and $\Gt_y^{\ga'}(e_{-\ga}) = {e_{-\ga}}$. 
Therefore both sides of \eqref{num} are equal to

\[\Gt^{\ga,\ga'}_{x,y} (e_{-\ga}^m e_{-\ga'}^m\gth_{\gc_3,m}(\mu)).\]
}\er

\subsubsection{Non-Isotropic Roots.} 
\bt \label{1calu} If $\gl \in \cH_{\gc,m}$, then 
\be\label{1CL}
 \theta_{\gamma,m}(\gl) =\theta_{\gamma,1}(\gl-(m-1)\gc)\ldots \theta_{\gamma,1}(\gl-\gc)\theta_{\gamma,1}(\gl).
\ee
Equivalently if $v_\gl$ is a highest weight vector of weight $\gl$ in the Verma module $M(\gl)$ we have $\theta_{\gamma,m}v_\gl = \theta_{\gamma,1}^m v_\gl$. If $\gl \in \cH_{\gc,m}$ this is a highest weight vector which is independent of the choice of $\gb$ in the inductive construction.
\et
\bpf 
Clearly \eqref{1CL} holds if $\gc$ is a simple root. Suppose that \eqref{121c} holds, and assume that $\gl\in w\cdot \gL$ where $\gL$ is defined in \eqref{1tar}.  
For $i=0,\ldots, m-1$ we have $$(\mu + \rho-i\gc', \alpha^\vee)= p+iq,$$ so by the inductive definition \eqref{yam} and \eqref{bam},  
\by \label{1inp} 
\theta_{\gamma,1}(\gl-i\gc)&=& \Gt_{ p + iq}(e^q _{- \alpha}\theta_{\gamma',1}(\mu-i\gc'))\nn\\
&=& \Gt_{ p }(e^{(i+1)q}_{- \alpha}\theta_{\gamma',1}(\mu-i\gc')e^{-iq}_{- \alpha}).\ey
Now using the corresponding result for $\theta_{\gamma',m}(\mu)$ we have
\by  \label{tea}e^{mq}_{- \alpha}\theta_{\gamma',m}(\mu) &=& e^{mq}_{- \alpha}\theta_{\gc',1}(\mu-(m-1)\gc')\ldots \theta_{\gc',1}(\mu-\gc')\theta_{\gc',1}(\mu)\nn\\
 &=&
e^{(m-1)q}_{- \alpha}(e^{q}_{- \alpha}\theta_{\gc',1}(\mu-(m-1)\gc'))e^{-(m-1)q}_{- \alpha}\cdot\nn\\
&&\cdot e^{(m-2)q}_{- \alpha}(e^{q}_{- \alpha}\theta_{\gc',1}(\mu-(m-2)\gc'))e^{-(m-2)q}_{- \alpha}\cdot\nn\\
&\cdots &\nn\\
& &  \cdot \;e^{kq}_{- \alpha}(e^{q}_{- \alpha}\theta_{\gamma',1}(\mu-k\gc') )e^{-kq}_{- \alpha}\cdot\nn\\ 
&\cdots & \nn\\
& &\cdot \;e^{q}_{- \alpha}(e^{q}_{- \alpha}\theta_{\gamma',1}(\mu-\gc'))e^{-q}_{- \alpha}\cdot e^{q}_{- \alpha}\theta_{\gamma',1}(\mu).\nn\ey
The result follows by applying the automorphism $\Gt_p$ to both sides and using \eqref{bam} and \eqref{1inp}.\epf

\subsection{Pairs of Roots.}
Next we consider relations between \v Sapovalov elements coming from different isotropic roots $\gc$ and $\gc'$. There are two cases depending on whether or not 
$(\gc,\gc') =0.$ 
\subsubsection{The non-orthogonal case.}
 Let $\Pi_{\rm nonisotropic}$ be the set of nonisotropic simple roots,  and ${Q}^+_{0} = \sum_{\ga \in \Pi_{\rm nonisotropic}} \mathbb{N}
 \ga.$
If $\gc,\gc'$ are  positive non-orthogonal isotropic roots, then
$\gc'=s_\ga \gc$ for some non-isotropic
$\ga\in {Q}^+_{0}$.  In this situation
the next result relates \v Sapovalov elements for $\gc, \gc'$ and $\ga$. We set $\gth_{\ga,0} =1.$

\bt  \label{man}Let  $\gc$ be a positive isotropic root and
$\alpha$ a non-isotropic root contained in ${Q}^+_{0}.$ Let  $v_{\gl}$ be a highest weight vector
in a Verma module with highest weight $\gl,$ and set $\gc' = s_\alpha\gc.$
Let $p = (\gl+\gr, \alpha^\vee)
$, and  assume $q = (\gc, \alpha^\vee) \in \mathbb{N}
\backslash \{0\}$. Suppose 
$\alpha \in A(\gl)$ and $p = (\gl + \rho, \alpha^\vee)$. 
 Then
\bi
\itema If $(\gl+\gr,\gc') = 0$
we have
\be \label{pin}\theta_{\gc}\theta_{\ga,p} v_{\gl}= \theta_{\ga,p+q} \theta_{\gc'} v_{\gl}.\ee
\itemb
If $(\gl+\gr,\gc) = 0$, and $p-q\ge 0,$
we have

\be \label{pun}\theta_{\gc'}\theta_{\ga,p} v_{\gl}= \theta_{\ga,p-q} \theta_{\gc} v_{\gl}.\ee\ei
\et

\bpf It suffices to prove (a) for all $\gl$ in the Zariski dense subset $\gL$ of
$\cH_{\gc'} \cap \cH_{{\ga,p}}$ given by

\[\gL=\{ \gl \in \cH_{\gc'} \cap \cH_{{\ga,p}}|A(\gl) = \{\ga\},\; B(\gl) = \{\gc'\}\}.\]
However for $\gl \in \gL$, $M(\gl)$ contains a unique highest weight vector of weight $s_\ga \cdot \gl -\gc.$ Therefore, $\theta_{\gc}\theta_{\ga,p} v_{\gl}$ and $ \theta_{\ga,p+q} \theta_{\gc'} v_{\gl}$ are equal up to a  scalar multiple. Now if $\pi^0$ is the partition of $\gc+p\ga$ with $\pi^0(\gs)=0$ for all non-simple roots $\gs,$ then it follows easily from the definition of \v Sapovalov elements, that $e_{-\pi^0}v_{\gl}$ occurs with coefficient equal to one in both $\theta_{\gc}\theta_{\ga,p} v_{\gl}$ and $ \theta_{\ga,p+q} \theta_{\gc'} v_{\gl}$,  and from this we obtain the desired conclusion.  The proof of (b) is similar.
\epf

\brs \rm{ (i) Perhaps the most interesting case of Equation (\ref{pin}) arises when $p=0$, since then we have an inclusion between submodules of a Verma module obtained by multiplying the highest weight vector $v_{\gl}$ by $\theta_{\gc}$ and $\theta_{\gc'}$.  Similarly the most interesting case of Equation (\ref{pun}) is when $p=q.$
\\
(ii) In the case that $\ga$ is a simple root, \eqref{pin} reduces to Equation (\ref{12nd}).}\ers
\subsubsection{The orthogonal case for special pairs.}\label{oc1}
 \noi  \noi Now we consider two isotropic roots $\gc_{1}, {\gc_{2}}$ such that $(\gc_{1},{\gc_{2}}) =0.$
In Subsection  \ref{oc} we show that for $\gl \in \cH_{\gc_{1},\gc_{2}}:=\cH_{\gc_{1}} \cap \cH_{{\gc_{2}}}$  the highest weight vectors
$\gth_{\gc_{2}}\gth_{\gc_{1}} v_\gl$ and $\gth_{\gc_{1}}\gth_{\gc_{2}}v_\gl$ are equal up to a constant multiple, the constant being a ratio of linear polynomials in $\gl$ differing only in their constant terms, see Equation (\ref{kat}) for the exact statement.  Here we consider some properties of such pairs of roots.\\ \\
If $\fg$ has two orthogonal isotropic roots then $\fg$ has defect at least two, see \cite{KaWa}. In particular $\fg$ cannot be exceptional.  So we assume that $\fg = \fgl(m,n), \osp(2m,2n)$ or $\osp(2m+1,2n)$ with $m, n \ge 2$.  As has become quite standard we  express the roots of $\fg$ in terms of linear forms $\gep_i, \gd_i \in \fh^*$, see \cite{K} 2.5.4 or \cite{M} Chapter 2.
The odd isotropic roots have the form

\be \label{or} \overline{\Delta}_1 = \{ \pm (\epsilon_i -\delta_j)|i\in [m], \;  j \in  [n]\},\ee if $\fg = \fgl(m,n)$ or

\be \label{nor} \overline{\Delta}_1 = \{ \pm \epsilon_i \pm \delta_j|i \in [m], \; j \in [n]\},\ee if $\fg$ is orthosymplectic.  We assume the  bilinear form $(\;,\;)$ on $\fh^*$ satisfies
\be
\label{edform}(\epsilon_i,\epsilon_j) = \delta_{i,j} = - (\delta_i,
\delta_j).\ee We say that $\Pi$ has $\osp$ type ($\spo$ type)
if it contains the simple root $\gep_m-\gd_1$ (resp. $\gd_n-\gep_1)$.  We also say that an isotropic  root is of type $\spo$ if it does not have the form $\gep_i-\gd_j$,  and type $\osp$ if it does not have the form
$\gd_i- \gep_j$.  Roots of the form $\gep_i+\gd_j$
are of both $\osp$ and $\spo$ type.
The following remarks and Lemma can easily be checked using \cite{K} Table VI, see  also \cite{M} Table 3.4.3.
The distinguished set of simple roots for types $B(m,n)$ with $m\ge 1$ and $D(m,n)$ have type $\spo$, and those of type $A(m,n)$ and $C(n)$ have type $\osp.$ The anti-distinguished set of simple roots, when it exists, has the opposite type. Note that the Lie superalgebra $C(n)$ has no anti-distinguished set of simple roots. We have

   \bl If $\Pi$ is either distinguished or anti-distinguished and $X$ is an orthogonal  set of positive isotropic roots,
	then all roots in $X$ have the same type.
	\el
\noi We say that an orthogonal  pair of isotropic roots $\gc_1, \gc_2$ is {\it special} if both roots are simple for the same Borel subalgebra. If this is the case, we can if necessary perform an odd reflection (see \cite{PS}, \cite{S1}, Section 3) and replace one of the roots by its negative, to assume that both $\gc_1$ and $ \gc_2$ have the same type.  From \eqref{or} or \eqref{nor}, there are positive even roots $\ga_1, \ga_2$ such that
\be\label{et2}\gi_1:= \gc_1+ \ga_1 = \gc_2+ \ga_2\ee is a  root.  It follows that
\be\label{et1}\gi_2  := \gc_1- \ga_2= \gc_2- \ga_1\ee is also a root, and that $\gi_1, \gi_2$ are positive orthogonal isotropic roots. 
Note that $\gc_1+ \gc_2=\gi_1+ \gi_2.$ Thus $0=(\gi_2,\gi_1+ \gi_2)= (\gi_2,\gc_1+ \gc_2),$ and we assume that
\be \label{rrr}(\gc_2, \gi_2)= 1 = -(\gc_1,\gi_2).\ee

\noi Now we can  state the main result on \v Sapovalov elements for special pairs of isotropic roots.
If the functions $a(\gl)$ and $b(\gl)$ are proportional we write $a(\gl)\doteq b(\gl)$.  
\bt \label{fiix} If $\gc_{1}, {\gc_{2}}$ is a special pair, and $\gi=\gi_2$ is as in \eqref{et1}, then for all $\gl \in \cH_{\gc_{1},\gc_{2}}$ we have
\be \label{kat1}[(\gl+\gr,\gi)-1]\gth_{\gc_{2}}(\gl - \gc_{1})\gth_{\gc_{1}}(\gl)\
\doteq
[(\gl+\gr,\gi)+1]\gth_{\gc_{1}}(\gl - {\gc_{2}})  \gth_ {\gc_{2}}(\gl).\ee
 \et
\noi 
We introduce some notation that is needed for the proof. For the remainder of Section \ref{RS}, we assume that $\fb$ is either the distinguished or anti-distinguished Borel subalgebra.  Consider a graph with vertices the set of Borel subalgebras having the same even part as $\fb$, with an edge connecting two Borels if they are connected by an  odd reflection.  If $\gc$ is an isotropic root, set $d(\gc)=r$ if $r$ is the shortest length of a path connecting $\fb$ to a Borel containing $\gc$. Such a path will be called a {\it path leading to $\gc.$}
Thus if $d(\gc_{1}) = r, d(\gc_{2})= s$, there are chains of adjacent Borel subalgebras, compare \eqref{distm}
 \be \label{dir} \mathfrak{b} = \mathfrak{b}^{(0)}, \mathfrak{b}^{(1)}, \ldots,
 \mathfrak{b}^{(r)}, \ee
and
 \be \label{dis} \mathfrak{b} = \mathfrak{b}^{[0]}, \mathfrak{b}^{[1]}, \ldots, \mathfrak{b}^{[s]} \ee
such that  $\gc_{1}$ and ${\gc_{2}}$ are simple roots of $\mathfrak{b}^{(r)}$ and $ \mathfrak{b}^{[s]}$ respectively. 
 There are odd roots $\gb_i$,  for $i \in [r]$ and $\gb_{[i]}$, for $ i \in [s]$,
such that for all $i$,
$$\fg^{\gb_i} \subset \fb^{(i-1)}, \quad \fg^{-\gb_i} \subset \fb^{(i)}, \quad \quad \fg^{\gb_{[i]}} \subset \fb^{[i-1]}, \quad \fg^{-\gb_{[i]}} \subset \fb^{[i]}.$$
Set
\[ F(1)= \{i \in [r] |(\gc_{1},\gb_i)=0\}, \quad  F(2)= \{i \in [s]| (\gc_{2},\gb_{[i]})=0\}.\]
We can arrange that the paths leading to $\gc_1$ and $\gc_2$ share an initial segment which is as long as possible.  This means that for some $t$, we will have $\gb_i= \gb_{[i]},$ for $i \in [t+1]$. 
\\ \\
\noi {\it Proof of Theorem \ref{fiix}}.
 Consider the Zariski dense subset
$\gL_{\gc_{1},\gc_{2}}$ of $\cH_{\gc_{1},\gc_{2}}$
given by
\[\gL_{\gc_{1},\gc_{2}}= \{\mu \in \cH_{\gc_{1},\gc_{2}}|(\mu + \rho, \ga) \notin \bbZ  \mbox{ for all positive isotropic roots } \ga \neq \gc_{1}, \gc_{2}\}.\]
It is enough to prove the  result for
$\gl \in \gL_{\gc_{1},\gc_{2}}.$ This ensures that each change of Borels in what follows is typical for $M(\gl)$,  except where one of the roots $\gc_{1}, \gc_{2}$ is replaced by its negative.
Next we define the following products of root vectors
\[e_T = e_{\gb_{1}} \dots e_{\gb_{t+1}}, \quad e_{-T} = e_{-\gb_{t+1}} \dots e_{-\gb_{1}},\]
 \[e_R = e_{\gb_{t+2}} \dots e_{\gb_r}, \quad e_{-R} = e_{-\gb_{r}} \dots e_{-\gb_{t+2}}, \]
\[e_S = e_{\gb_{[t+2]}} \dots e_{\gb_{[s]}},\quad e_{-S} = e_{-\gb_{[s]}} \dots e_{-\gb_{[t+2]}}. \]
The root $\gi := \gb_{t+1}$ is a simple root for the
Borel subalgebra
$\mathfrak{b}^{(t)}$,  and so $\gi$ corresponds to a node of the Dynkin-Kac diagram for
$ \mathfrak{b}^{(t)}$.
Let $\fk$ (resp. $\fl$) be the subalgebra of $\fg$ generated by root vectors (positive and negative) corresponding to the nodes to the left (resp. right) of this node.
We have $[\fk,\fl]=0$, and hence
\be \label{bats}[e_{\pm S},e_{\pm R}] = [e_{\pm S},e_{-\gc_{1}}]=[e_{\pm R},e_{-\gc_{2}}]  = 0.\ee
Set $$\Phi_S =e_Se_{-{\gc_{2}}}e_{-S} ,\quad \;\Phi_R=e_Re_{-\gc_{1}}e_{-R}.$$
We claim that
\be\label{god}\Phi_S \Phi_R=\pm  \Phi_R \Phi_S.\ee
Indeed by (\ref{bats}),
\begin{eqnarray}
\label{ma1}
\Phi_R \Phi_S
&=& \pm e_Re_{-{\gc_{1}}}e_{-R}e_Se_{-{\gc_{2}}}e_{-S}
\nonumber \\
&=& \pm e_Re_{-{\gc_{1}}}e_Se_{-R}e_{-{\gc_{2}}}e_{-S}
\nonumber \\
&=&  \pm e_Re_Se_{-{\gc_{1}}}e_{-\gc_{2}}e_{-R}e_{-S}
\nonumber \\
&=&  \pm e_Re_Se_{-{\gc_{2}}}e_{-\gc_{1}}e_{-R}e_{-S}
\nonumber\\
&=&  \pm e_Se_{-{\gc_{2}}}e_Re_{-S}e_{-\gc_{1}}e_{-R}
\nonumber\\
&=&\pm \Phi_S \Phi_R
.\nonumber
\end{eqnarray}
\noi From Theorem \ref{1rot} we obtain
\be \label{nod} e_{T}\Phi_Re_{-T}v_\gl \doteq \prod_{i \in F(1)}(\gl+\gr,\gb_i)\gth_{\gc_{1}} v_\gl. \ee
Similarly  using the fact that the expression in (\ref{nod}) is a highest weight vector of weight $\gl-\gc_{1}$, we have
\be \label{odd1}e_T\Phi_S e_{-T}e_{T}\Phi_Re_{-T}v_\gl \doteq  a_{\gc_{1}}\theta_{\gc_{2}}\theta_{\gamma_1}v_\gl, \ee
where
\[a_{\gc_{1}}=\prod_{i \in F(1)}(\gl+\gr,\gb_i)
\prod_{i \in F(2)}(\gl+\gr-\gc_{1},\gb_{[i]}),\]
Next  set
\be \label{cod} b_{\gc_{1}}=[(\gl+\gr,\gi)+1]\prod_{i=1}^t (\gl+\gr-\gc_{1},\gb_i),\ee
and $c_{\gc_{1}} = a_{\gc_{1}}/b_{\gc_{1}}.$
Switching the roles of $\gc_1$ and $\gc_2$, define similarly,
\[a_{\gc_{2}}=\prod_{i \in F(2)}(\gl+\gr,\gb_{[i]})
\prod_{i \in F(1)}(\gl+\gr-\gc_{2},\gb_i),
\quad b_{\gc_{2}}=[(\gl+\gr,\gi)-1]\prod_{i=1}^t (\gl+\gr-\gc_{2},\gb_i),\]
and $c_{\gc_{2}} = a_{\gc_{2}}/b_{\gc_{2}}.$ 
Since the leftmost factor in $e_T$ is $e_{\gb_{1}}$, we can write
 $e_{T}\Phi_Re_{-T}v_\gl =e_{\gb_{1}}w$ for some $w$ which is a highest weight vector for $\mathfrak{b}^{(1)}$ of weight $\gl-\gc_1-\gb_1$.  
On the other hand, the rightmost factor of $e_{-T}$ is $e_{-\gb_{1}}$,  and from (\ref{ebl})
 we have $e_{-\gb_{1}}e_{\gb_{1}}w =(\gl+\gr-\gc_{1},\gb_1)w$. Continuing like this with the other
 factors of $e_T$ and using \eqref{rrr}, gives
$e_{-T}e_{T}\Phi_Re_{-T}v_\gl \doteq  b_{\gc_{1}}\Phi_Re_{-T}v_\gl. $ \ff{ The factor $(\gl+\gr,\gi)+1$ in $b_{\gc_{1}}$, which plays a crucial role in \eqref{kat1}, arises at this point using \eqref{rrr}, since $e_{\gb_{t+1}}=e_\gi$ is a factor of $e_T$. The differences between the first factors in the definitions
of $b_{\gc_{1}}$ and $b_{\gc_{2}}$ are due to \eqref{rrr}.}   Therefore
 \be \label{odd}e_T\Phi_S e_{-T}e_{T}\Phi_Re_{-T}v_\gl \doteq  b_{\gc_{1}} e_T\Phi_S \Phi_Re_{-T}v_\gl. \ee
 Combining this with Equation (\ref{odd1}) we have

\be \label{rod1} e_T\Phi_S\Phi_Re_{-T}v_\gl \doteq  c_{\gc_{1}} \theta_{\gc_{2}}\theta_{\gamma_1}v_\gl.\ee
Similarly

 \be \label{sod}e_T\Phi_R\Phi_Se_{-T}v_\gl \doteq  c_{\gc_{2}} \theta_{\gc_{1}}\theta_{\gamma_2}v_\gl.\ee
Together with Equation (\ref{god}) this gives

\be \label{pod}c_{\gc_{1}} \theta_{\gc_{2}}\theta_{\gamma_1}v_\gl \doteq c_{\gc_{2}} \theta_{\gc_{1}}\theta_{\gamma_2}v_\gl.\ee
The proof of  
Theorem \ref{fiix} is completed by the Lemma below. Since the proof is rather technical, it is followed by an example illustrating some of the notation.\hfill  $\Box$

\bl \label{cop} We have $[(\gl+\gr,\gi)+1]c_{\gc_{1}} = [(\gl+\gr,\gi)-1]c_{\gc_{2}}$. \el 
\noi Set $G(j)= F(j)\cap [t+1]$ and $H(j) =F(j)\backslash G(j)$ for $j=1,2$. First we show
\bsul
If $k\in H(1),$ then $(\gc_2,\gb_k)=0$ and similarly if $k\in H(2),$ then $(\gc_1,\gb_{[k]})=0.$
\esul
\bpf
To show this, we assume  that $\gc_1 = \gep_i -\gd_{i'}$ and $ \gc_2 = \gep_j -\gd_{j'}$ with $i<j$ and $i'<j'$, compare also the first bullet in Theorem \ref{hog}.
(Note that \eqref{et2}-\eqref{rrr} hold with
$\gb_1 = \gd_{i'} -\gd_{j'}, \gb_2 = \gep_i-\gep_j,$ $\gi_1 = \gep_i -\gd_{j'}$
and $\gi_2 = \gep_j -\gd_{i'}.$) Set $\gi = \gi_2$.
\\ \\
In \cite{M} 3.3, the Borel subalgebras with the same even part as $\fb^\dist$ are described in terms of shuffles.  Here the notation is slightly different. We write a permutation $\gs$ of the set $\{1,2, \ldots,m, 1',2',\ldots,n'\}$ in one-line notation as
\[ \underline{\gs} = (\gs(1), \gs(2), \ldots ,\gs(m), \gs(1'),\gs(2'),\ldots,\gs(n')).\] 
Then we say that $\gs$ is a {\it shuffle} if  $1,2, \ldots,m$ and $ 1',2',\ldots,n'$ are subsequences of $\underline{\gs}$.
We can write the shuffles corresponding to $\fb^{(t+1)}$ and $\fb^{(r)}$ as the concatenations of
\[A,i,B,A',i',j,C,B'\;\; \mbox{ and }\;\; A,A',i,i',B,j,C,B'\]
respectively, where
\[A = (1,\ldots,i-1), \;\;B= (i+1,\ldots, j-1), \;\;C = (j+1,\ldots,m),\]\[ A' = (1',\ldots,i'-1),\quad B'=(i'+1,\ldots,n').\]
Then $\{\gb_k|k\in H(1)\}$ consists of all odd roots $\ga$ which are roots of $\mathfrak{b}^{(r)}$ but not roots of $\mathfrak{b}^{(t)}$, such that $(\gc_1,\ga)=0$. Looking at pairs of entries which occur in opposite orders in
the two shuffles, it follows that any such root $\ga$ is contained in the set of roots
\[\{\gep_i-\gd_{a'}, \gep_b-\gd_{a'}, \gep_b-\gd_{i'} | a' \in A', b\in B\},\]
and clearly all roots in this set are orthogonal to $\gc_2.$
\epf

\noi{\it Proof of Lemma \ref{cop}.}
By the Sublemma, if
$$z = \prod_{i \in H(1)}(\gl+\gr,\gb_i) \prod_{i \in H(2)}(\gl+\gr,\gb_{[i]}),$$
then
\[a_{\gc_{1}}=z\prod_{i \in G(1)}(\gl+\gr,\gb_i)
\prod_{i \in G(2)}(\gl+\gr-\gc_{1},\gb_{[i]}),\] and
\[a_{\gc_{2}}=z\prod_{i \in G(2)}(\gl+\gr,\gb_{[i]})
\prod_{i \in G(1)}(\gl+\gr-\gc_{2},\gb_i).\]
Now denote the complements of
$G(1), {G(2)}$ in $[1\ldots t+1]$ by $\overline{G(1)},
\overline{G(2)}$ respectively, and set

\[I_1 = {G(1)} \cap {G(2)}, \; I_2 = {G(1)} \cap \overline{G(2)}, \; I_3 = \overline{G(1)} \cap {G(2)}.\]
Note also that
$\overline{G(1)} \cap \overline{G(2)} =\{t+1\},$ but since $\ga_{t+1} = \gi$ is not orthogonal to $\gc_1$ it does not contribute to the product defining $a_{\gc_{1}} $.  Hence

\begin{eqnarray}
a_{\gc_{1}} &=&z\prod_{i\in I_1}(\gl+\gr,\gb_i)^2 \prod_{i\in I_2} (\gl+\gr,\gb_i) \prod_{i\in I_3}(\gl+\gr-\gc_{1},\gb_i).
\nonumber
\end{eqnarray}
and
$b_{\gc_{1}}=[(\gl+\gr,\gi)+1]\prod_{j = 1}^3 b_{\gc_{1}}^{(j)}$,
where  for $1\le j \le 3,$

\be \label{cads} b^{(j)}_{\gc_{1}}=\prod_{i \in I_j} (\gl+\gr-\gc_1,\gb_i).\ee
Canceling common factors of $a_{\gc_{1}}$ and $b_{\gc_{1}}$, it follows that 
$$[(\gl+\gr,\gi)+1]c_{\gc_{1}}  = z\prod_{i\in I_1}(\gl+\gr,\gb_i),$$ and similarly this is equal to $[(\gl+\gr,\gi)-1]c_{\gc_{2}}.$
\hfill  $\Box$\\ 
\vspace{0.1cm}

\noi This concludes the proof of Theorem \ref{fiix}.
\bexa {\rm Let $\fg = \fgl(4,4), \gc_1 = \gep_1 -\gd_{1}, \gc_2 = \gep_3 -\gd_{3}.$  Then $\gi = \gep_3 -\gd_{1}.$
Below we give four Dynkin-Kac diagrams for $\fg$.  The first is the anti-distinguished diagram, and the single grey node corresponds to the
simple root $\gep_4 -\gd_{1}.$  In the notation of Theorem \ref{fiix} the second, third and fourth diagrams correspond to the Borel subalgebras $\mathfrak{b}^{(t)}, \mathfrak{b}^{(r)},\mathfrak{b}^{[s]}$, and so $\gi, \gc, \gc'$  are simple roots of these subalgebras  respectively.
The node corresponding to the root $\gi$ is indicated by a square.
If $\gl \in \gL$, and $v_\gl$ is a highest weight vector in a Verma module with highest weight $\gl$, then $e_{-T}v_\gl, e_{-R}e_{-T}v_\gl$ and $e_{-S}e_{-T}v_\gl$ are highest weight vectors for the Borel subalgebras
$\mathfrak{b}^{(t)}, \mathfrak{b}^{(r)},\mathfrak{b}^{[s]}$ respectively.

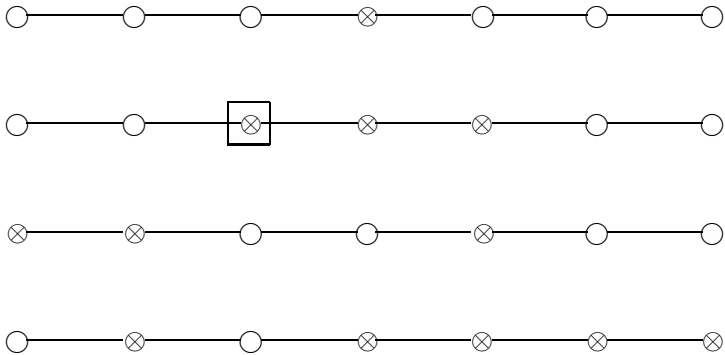
\begin{picture}(60,40)(-53,-20)
\thinlines
\put(1,3){\line(1,0){36.2}}
\put(46,3){\line(1,0){35.8}}
\put(41.7,2.3){\circle{8}}
\put(46,3){\line(1,0){35.8}}
\put(85.7,2.3){\circle{8}}
\put(90,3){\line(1,0){36}}
\put(125.7,-0.3){$\otimes$}
\put(133,3){\line(1,0){36}}
\put(173.7,2.3){\circle{8}}
\put(177,3){\line(1,0){36}}
\put(-2.7,2.3){\circle{8}}
\put(216.7,2.3){\circle{8}}
\put(220.5,3){\line(1,0){36}}
\put(260.3,2.3){\circle{8}}
\end{picture}

\begin{picture}(60,40)(-53,-20)
\thinlines
\put(-2.7,2.3){\circle{8}}
\put(1,3){\line(1,0){36.2}}
\put(46,3){\line(1,0){35.8}}
\put(41.7,2.3){\circle{8}}
\put(81.7,-0.3){$\otimes$}
\put(90,3){\line(1,0){36}}
\put(125.7,-0.3){$\otimes$}
\put(133,3){\line(1,0){36}}
\put(169,-0.3){$\otimes$}
\put(177,3){\line(1,0){36}}
\put(216.7,2.3){\circle{8}}
\put(77,11){\line(1,0){16}}
\put(77,-5){\line(1,0){16}}
\put(77,-5){\line(0,1){16}}
\put(93,-5){\line(0,1){16}}
\put(220.5,3){\line(1,0){36}}
\put(260.3,2.3){\circle{8}}
\end{picture}

\begin{picture}(60,40)(-53,-20)
\thinlines 
\put(-6.7,-0.3){$\otimes$}
\put(1,3){\line(1,0){36.2}}
\put(46,3){\line(1,0){35.8}}
\put(37.7,-0.3){$\otimes$}
\put(46,3){\line(1,0){35.8}}
\put(85.7,2.3){\circle{8}}
\put(90,3){\line(1,0){36}}
\put(129.7,2.3){\circle{8}}
\put(133,3){\line(1,0){36}}
\put(169.7,-0.3){$\otimes$}
\put(177,3){\line(1,0){36}}
\put(216.7,2.3){\circle{8}}
\put(220.5,3){\line(1,0){36}}
\put(260.3,2.3){\circle{8}}
\end{picture}

\begin{picture}(60,40)(-53,-20)
\thinlines
\put(1,3){\line(1,0){36.2}}
\put(46,3){\line(1,0){35.8}}
\put(46,3){\line(1,0){35.8}}
\put(85.7,2.3){\circle{8}}
\put(90,3){\line(1,0){36}}
\put(125.7,-0.3){$\otimes$}
\put(133,3){\line(1,0){36}}
\put(169,-0.3){$\otimes$}
\put(177,3){\line(1,0){36}}
\put(212.7,-0.3){$\otimes$}
\put(-2.7,2.3){\circle{8}}
\put(37.7,-0.3){$\otimes$}
\put(220.5,3){\line(1,0){36}}
\put(256.3,-0.3){$\otimes$}
\end{picture}

\noi
We have $t=1, r=3, s= 5$,
\[\gb_1 = \gep_4-\gd_1, \quad  \gb_2 = \gep_3-\gd_1, \quad  \gb_3 = \gep_2-\gd_1,\]\[
\gb_{[3]} = \gep_4-\gd_2, \quad  \gb_{[4]} = \gep_4-\gd_3, \quad  \gb_{[5]} = \gep_3-\gd_2,\]
\[F(1)=\emptyset, \quad  G(2) = \{1\}, \quad  H(2) = \{3\}\]
and
\[e_T = e_{45}e_{35},\quad  e_{-T}=e_{53}e_{54}, \]
\[e_R = e_{25},\quad e_{-R}=e_{52} ,\]
\[e_S = e_{46} e_{47} e_{36} ,\quad e_{-S}=e_{63} e_{74} e_{64}.\]
} \eexa

\subsubsection{On a result of Duflo and Serganova.}\label{cry}
To extend Theorem \ref{fiix}, we need a variant of a result of Duflo and Serganova, \cite{DS} Lemma 4.4 (2).
\noi First some definitions. Suppose that $\Pi$  is a basis of simple roots for $\fg$ satisfying (\ref{i}) and recall the group $W_{\rm nonisotropic}$ from Equation (\ref{ii}), defined using $\Pi$.
Let $W_1$ (resp. $W_2$) be the Weyl group of $\fo(\ell)$ with $\ell=2m$ or $2m+1$
(resp. the Weyl group of $\fsp(2n)$).   Then we have
$W_{\rm nonisotropic}=W_1\times S_n$ or $S_m \times W_2$ if $\Pi$ has type
$\spo$ or  $\osp$ respectively.
\bt \label{hog}Suppose $X$ is a set of orthogonal isotropic roots all  having the same type as $\Pi$. Then there exists $w\in W_{\rm nonisotropic},$  such that $wX$ is contained in the set of simple roots for some Borel subalgebra. In fact if $X=|k|$, there exists $w\in W$ such that $wX$ has one if the following forms.
\bi \item        $\{\epsilon_i -\delta_i\}_{i=1}^k $ if $\fg$ is of type A, or $\Pi$ is of $\osp$ type
\item                $\{\gd_{i}- \gep_i\}_{i=1}^k$ if  $\Pi$ is of $\spo$ type, and either $\fg= \osp(2m,2n)$\\         with $ k<m$ or $\fg = \osp(2m+1,2n)$
                            \item $\{\gd_{n-k+i} -\gep_{i}, \gd_{n} \pm \gep_{m}\}_{i=1}^{k-1}$ if $\Pi$ is of $\spo$ type and $\fg=\osp(2k,2n)$.
\ei
\et

\bpf
 First suppose
\be \label{clr} X = \{\gep_{f(i)} - a_{i}\gd_{h(i)}\}_{i=1}^k \ee
 where $a_{i} =\pm 1$ and $f:[k]\lra [m],$ and $h:[k]\lra [n]$.  Since $X$ is orthogonal, $f$ and $h$ are injective, so reordering $X$, we may assume that $f$ is increasing. Then using the $W$-action  we can assume $f(i)= h(i)=i$ for all $i\in [k]$.  If $ \fg$ has type A, then  $a_{i}=1$ for all $i$, and we have shown $X$ is conjugate to $\{\epsilon_i -\delta_i\}_{i=1}^k$. By changing the signs of the $\gd_i$ we also have the result  for $\fg$ is orthosymplectic, and $\Pi$ of $\osp$ type.  If $\Pi$ is of $\spo$ type, we start instead with
$X = \{\gd_{h(i)}- a_{i}\gep_{f(i)}\}_{i=1}^k$ and argue similarly unless $\fg=\osp(2m,2n),$ where we are only allowed to change an even number of signs of the $\gep_i$.  If $k<m$ there is still enough room for the argument to go through.  Otherwise $n \ge m = k$ and we see that $X$ is conjugate to
\[ \{\gd_{n-k+1} -\gep_{1}, \ldots \gd_{n-1} -\gep_{m-1}, \gd_{n} + (-1)^b \gep_{m}\},\]
where $b$ is the number of $a_{i}$ that are negative.  \epf

\subsubsection{The orthogonal case in general.}\label{oc}
\noi 
\noi Let $X=\{\gk_{1},\gk_{2}\}$ be a set of two orthogonal isotropic roots both having the same type as $\Pi$, and
let $w$ be as in Theorem \ref{hog}. Then set
  $\gc_{1}=w\gk_1,\gc_{2}=w\gc_2$ and let $\gi =\gi_2$ be as in \eqref{et1}.
We obtain \eqref{pee} in a more precise form.
\bt \label{fiix2} For $\gl \in \cH_{\gk_{1},\gk_{2}}$ we have
\be \label{kat}[(\gl+\gr,w^{-1}\gi)-1]\gth_{\gk_{2}}(\lambda - \gk_{1})\gth_{\gk_{1}}(\gl)\
\doteq
[(\gl+\gr,w^{-1}\gi)+1]\gth_{\gk_{1}}(\gl - {\gk_{2}})  \gth_ {\gk_{2}}(\gl).\ee
 \et
\bpf We use
induction on the length $\ell(w)$ of  $w\in W_{\rm nonisotropic}$, Theorem \ref{fiix} giving the result when $w=1$. Write  $w = us_\ga$ where $\ell(u)=\ell(w)-1$ and $\ga$ is a simple root. Assume that  $(\gk_{1},\ga^\vee)= q,$ $(\gk_{2},\ga^\vee)= q',$ and define $s_\ga \gk_{1} = \gt_1 =\gk_{1}-q\ga$, $s_\ga \gk_{2} = \gt_2 =\gk_{2}-q'\ga$.
If  $(\mu+\gr,\ga^\vee)= p,$ we assume that $p,q$ and $q'$  are non-negative integers. Then set $\gl = s_\ga\cdot\mu = \mu -p\ga$, $v_\gl=e_{-\ga}^p v_\mu$. By induction
\be \label{bad}[(\mu+\gr,u^{-1}\gi)-1]\gth_{\gt_2}(\mu - \gt_1)\gth_{\gt_1}(\mu)=[(\mu+\gr,u^{-1}\gi)+1]\gth_{\gt_1}(\mu - {\gt_2})  \gth_ {\gt_2}(\mu).\ee
 Then by Equation (\ref{12nd})

\be \label{nut} \gth_{\gk_{1}}e_{-\ga}^{p}v_\mu=e_{-\ga}^{p+q}\gth_{\gt_1}v_\mu,\ee
and

\be \label{not} \gth_{\gk_{2}}e_{-\ga}^{p+q}\gth_{\gt_1}v_\mu=e_{-\ga}^{p+q+q'}\gth_{\gt_2}\gth_{\gt_1}v_\mu. \ee
Hence
\begin{eqnarray}
\label{seat}
\gth_{\gk_{2}}\gth_{\gk_{1}}v_\gl
&=&  \gth_{\gk_{2}}\gth_{\gk_{1}}e_{-\ga}^{p}v_\mu
\nonumber\\
&=& \gth_{\gk_{2}}e_{-\ga}^{p+q}\gth_{\gt_1}v_\mu \quad \mbox{by (\ref{nut})}
\nonumber\\
&=&  e_{-\ga}^{p+q+q'}
\gth_{\gt_2}\gth_{\gt_1}v_\mu \quad \;\mbox{by (\ref{not})}.
\end{eqnarray}
Similarly  by first interchanging the pairs $(\gt_1,\gk_{1})$ and $(\gt_2,\gk_{2})$ in Equations (\ref{not}) and (\ref{nut}) we obtain
\be\label{yea}\gth_{\gk_{1}}\gth_{\gk_{2}}v_\gl=e_{-\ga}^{p+q+q'}\gth_{\gt_1}\gth_{\gt_2}v_\mu.\ee
Since $(\gl+\gr,w^{-1}\gi)= (\mu+\gr,u^{-1}\gi)$, we obtain the result from  Equations (\ref{bad}), (\ref{seat}) and  (\ref{yea}).\epf


\section{Highest weight modules with prescribed characters.}\label{jaf}
In this Subsection we assume that
$\fg$ is a basic classical simple  Lie superalgebra  and
$\Pi$ is a basis of simple roots satisfying hypothesis (\ref{i})  and either of (\ref{iii}) or (\ref{ii}).
Our goal is to prove Theorem \ref{newmodgen}.
Let $X$ be an orthogonal set of positive isotropic roots. When $(\gl + \gr, \gc) =
0$ for all $\gc \in X,$
we construct some highest weight modules $M^X(\lambda)$ with highest weight $\gl$
and character $ \tte^{\gl} p_X$.
\\ \\ We begin with a sketch of the construction. 
Define $\xi$ and $\widetilde{\lambda}$ as in \eqref{wld}. The main step in the proof of Theorem \ref{newmodgen} is the proof that the module 
$M^X({\widetilde{\lambda}})_B$ defined below in \eqref{tuv} has character $\tte^{\widetilde{\lambda}}p_X$.
First we show in Corollary \ref{ratss} that
the weight spaces of  $M^X(\widetilde{\lambda})$ satisfy  $\dim_B M^X(\widetilde{\lambda})^{\widetilde{\lambda}- \eta} \leq {\bf p}_{X}(\eta)$ for  all $\eta \in Q^+$.
Then if $X = Y \cup \{\gc\}$ where $\gc\notin Y,$ we show that there is an exact sequence
\be \label{ses}0 \lra L\lra {M^Y}(\widetilde{\gl})_B \lra N\lra 0,\ee 
where $L$, $N$ are homomorphic images of 
${M^X}(\widetilde{\lambda}-\gc)_B$ and ${M^X}(\widetilde{\lambda})_B$
respectively, see Equations \eqref{rod} and \eqref{yim}. The claim about the character of 
$M^{X}({\widetilde{\lambda} })_{B}$ then follows  by induction on $|X|$. Now let $M^{X}({\widetilde{\lambda} })_{A}$ be the submodule of 
$M^{X}({\widetilde{\lambda} })_{B}$ generated by the highest weight vector. We show that Theorem \ref{newmodgen} holds with
\be  \label{hut} M^{X}(\lambda) = M^{X}({\widetilde{\lambda} })_{A}/TM^{X}({\widetilde{\lambda} })_{A}.\ee

\subsection{An upper bound for the  dimension of certain weight spaces.}\label{jf1}
Let $X$ be an orthogonal set of isotropic roots.  In this Subsection,   $M=U(\fg)_Bv$ is a module with highest weight $\widetilde{\lambda}$ and highest weight vector $v$.  We assume $\gth_\gc v = 0$ for all $\gc \in X.$
In this Subsection we use the usual notation for partitions, see Subsection \ref{sss7.1}.
Let
$L$ be the subspace  of $M^{\widetilde{\lambda} -\eta}$  generated by all products $e_{-\pi} v$ with
$\pi\in {\bf {P}}_{X}(\eta).$
We fix an order on the positive roots such that for any partition $\pi$ the factors of the form $e_{-\gc}$, with
$\gc \in \Gd_1^+$  occur to the  right of the other factors in Equation
(\ref{negpar}) and among these factors, those with $\gc\in X$ occur farthest to the right.

\bp \label{rats} For all  $\pi \in {\bf {P}}(\eta)$ we have $e_{-\pi}v \in L$.
\ep
\noi For $S\subseteq \Gd^+_1$ we set $e_{-S} =\prod_{\gc \in S}e_{-\gc}$, and $\|S\| = \sum_{\gc \in S}\gc$. For any partition $\pi  \in {\bf {P}}(\eta)
$ we have a unique decomposition
\be \label{unique decomposition} e_{-\pi} =e_{-\gs}e_{-S},\ee
where $S\subseteq X$ and $\gs \in {\bf {P}}_X(\eta-\|S\|)$. Because of the way we have ordered root vectors, $e_{-S}v \in L$ implies
$e_{-\pi}v \in L$.  So it is enough to prove the result when $\eta = \|S\|$, equivalently $e_{-\gs}=1$. 
If this is the case and \eqref{unique decomposition} holds, we set $S_\pi = S.$
\\ \\
Suppose $a_1,\ldots, a_r\in \fn_1$ and $x_1,\ldots,x_s\in \fn_0$.
If $J=\{j_1<\ldots<j_t\}$ is a subset of $[s]$, we set
\[x_J = x_{j_1}\ldots x_{j_t}, \quad [a_i,x]_J = [[\ldots[a_i, x_{j_1}]\ldots ]x_{j_t}].\]
\bl \label{0.1} We have
\by\label{fla}  a_1\ldots a_r x_1\ldots x_s &=& \sum x_{J(0)} [a_1,x]_{J(1)}\ldots [a_r,x]_{J(r)}
\ey
where the sum is over all partitions $[s]= J(0)\cup J(1)\cup \ldots \cup J(r)$ of $[s]$.
\ff{We admit the possibility that some of the $J(i)$ are empty.}
 \el
\bpf An easy induction.\epf

\bl \label{dfg} If $R$  and $S$ are  subsets of $X$, and $\|R\|=\|S\|$ then $R=S.$

\el \bpf This is clear if $|X|=1$. Otherwise we can assume that $\fg$ is not exceptional, and then the result follows using an explicit description of the roots, compare \eqref{clr}. It is crucial that $X$ is an orthogonal set of roots.\epf

\noi {\it Proof of Proposition \ref{rats}.}
By \eqref{unique decomposition} and induction on $\eta$, it suffices to prove the result
when $\eta=\|S_\pi\|$.
We use induction on $s_\pi=|S_\pi|$.
If $s_\pi = 0,$ then $\pi \in {\bf {P}}_X(\eta)$ and the result holds by the definition of $L$.
Assume the result holds whenever $s_\gs<s_\pi=s$. Suppose $S_\pi=\{\gc_1,\ldots, \gc_s\}$.
Then

\be \label{fp} 0=\gth_{\gc_{1}}\ldots \gth_{\gc_{s}}v
=\sum_{\langle \gs\rangle}
p_{\langle \gs\rangle}e_{-\gs(1)} \ldots e_{-\gs(s)}v,\ee
where the sum is over all $s$-tuples $\langle \gs\rangle=(\gs(1), \ldots , {\gs(s)})$ with $\gs(i) \in {\bf {P}}(\gc_i)$ for $i\in[s]$, and $p_{\langle \gs\rangle}=p_{\langle \gs\rangle}(T) \in A.$
We write the term corresponding to ${\langle \gs\rangle} $ as
\be \label{lst} e_{\langle \gs\rangle}v=e_{-\gs(1)} \ldots e_{-\gs(s)}v.\ee
Now one term in \eqref{fp} is $e_{-\pi}v=e_{-\gc_1}\ldots e_{-\gc_s}v.$ 
We write this term $e_{-\pi}v$ also as $e_{\langle \pi\rangle}v$.
If ${\langle \gs\rangle} \neq {\langle \pi\rangle}$ we show that
$ e_{\langle \gs\rangle}v$
is, modulo terms that can be treated by induction, a $\ttk$-multiple of
$e_{\langle \pi\rangle}v$, see \eqref{ban}.
Since $\deg p_{\langle \gs\rangle}<\deg p_{\langle \pi\rangle}$
by Theorem \ref{1Shap} or Theorem \ref{1aShap}, this will give the result.
\\ \\
We have $\gs(i)(\gb_i)\neq0$ for a unique odd root $\gb_i$, and $\gb_i\le\gc_i$.
Also  $e_{-\gs(i)}=e_{-\gk_i}e_{-\gb_i}$ for some partition $\gk_i$ of ${\gc_i}-\gb_i$.
Thus
\be \label{m29} e_{\langle \gs\rangle}v= e_{-\gk_1}e_{-\gb_1}\ldots e_{-\gk_s}e_{-\gb_s}v,\ee
and by repeated use of Lemma \ref{0.1}, we move all odd root vectors $e_{-\gb_i}$ to the right in \eqref{m29}.
We do not however need a formula as explicit as \eqref{fla}.  Instead consider the multiset
\[\cA = \bigcup_{i=1}^s \{\ga^{\gk_i(\ga)}|\ga\in \Gd_0^+\}.\]
({The notation means that the root $\ga$ appears $\gk_i(\ga)$ times.})
Then  we see that modulo terms already known to satisfy the conclusion of the Proposition,
$e_{\langle \gs\rangle}v$ is a $\ttk$-linear combination of products of the form
\be \label{expr1} e_{-\go_1} \ldots e_{-\go_s}v,\ee
 where the $e_{-\go_i} \in (\ad U(\fn_0))e_{-\gb_i}$  are odd root vectors.  Thus there  is a multiset partition $\bigcup_{i=1}^s \cA_i$ of $\cA$ such that
$\go_i=\gb_i + \sum_{\ga\in \cA_i}\ga$. Note that since the bracket of two odd root vectors is even, and $[\fn^-_0,\fn_1^-]\subseteq \fn_1^-$, we can assume  by induction on $\eta$, that any terms $x_{J(0)}$
arising from \eqref{fla} are in fact constant.
We can further require that there is a partition $\go$ of $\eta$ such that the expression in \eqref{expr1} is equal to $e_{-\go}$ up to a permutation of the factors.
Then by induction on $s_\pi$, we can assume that the  $\gO =\{\go_1, \ldots, \go_s\}$ is a subset of $X$, and so by Lemma
\ref{dfg}, $\gO=\{\gc_1, \ldots, \gc_s\}$.  Now if $i\neq j$, then $[e_{-\gc_i},e_{-\gc_j}]=0$, so it follows that after all the rewriting we have
\be \label{ban} e_{\langle \gs\rangle}v = a_\gs e_{-\pi}v+b_\gs\ee
 with $a_\gs \in \ttk$ and $b_\gs \in L$.
By \eqref{fp}, $p(T)e_{\langle \pi\rangle}\in L$ for some polynomial $p$ with
 $\deg p = p_{\langle \pi\rangle}$. The result follows.
\hfill  $\Box$

\bexa {\rm 
We examine the proof in the case $\fg=\fgl(2,2)$ using the notation for roots and root vectors from Section \ref{Ch8}.
Let $X=\{\gb,\ga+\gb+\gc\}$, $\eta = \ga+2\gb+\gc$ and $M={M^X}(\widetilde{\lambda})_B^{\widetilde{\lambda}-\eta}$. 
Then Proposition \ref{rats} claims that $M$ is spanned over $B$ by all 
$e_{-\pi}v$ with ${\pi \in {{\bf P}_X(\eta)}}$.  
The most interesting stage of the proof arises when in the notation of Equation \eqref{unique decomposition} 
we have $e_{-\pi} =e_{-S}$ with $S\subseteq X$ and  we have already shown by induction that the result holds when $|S| =1$. 
Then we consider the product $0= e_{-\gb}\gth_{\ga+\gb+\gc} v.$ 
By Theorem \ref{1Shap} the coefficient of $e_{-\gb}e_{-\ga-\gb-\gc}v$ in $\gth_{\ga+\gb+\gc}v$ is quadratic in $T$. 
However the above product also contains the term $e_{-\gb} e_{-\ga} e_{-\gc} e_{-\gb}v$ with constant coefficient.  
Reordering creates the term $e_{-\gb}e_{-\ga-\gb-\gc}v$, 
and adding  this term does not change the quadratic nature of the coefficient of $e_{-\gb}\gth_{\ga+\gb+\gc}v$.  
The product also contains the terms $e_{-\gb} e_{-\ga} e_{-\gb-\gc} v$ and $e_{-\gb} e_{-\gc} e_{-\ga-\gb}v$ 
both with linear coefficient in $T$. 
Reordering these terms creates the new terms $e_{-\gb-\gc} e_{-\ga-\gb}v$, $e_{-\ga} e_{-\gb} e_{-\gb-\gc} v$ 
and $e_{-\gc} e_{-\gb} e_{-\ga-\gb}v$.  
The first of these has the form $e_{-\pi}v$ with $\pi \in {{\bf{{P}}_{X}(\eta)}}$ while the other two are contained in $M$ by induction.  
It follows that $e_{-\gb}e_{-\ga-\gb-\gc}v\in M.$ }
\eexa 

\bc \label{ratss}
With the same notation as the Proposition. 
 \bi
\itema The weight space $M^{\widetilde{\lambda} -\eta}_B$ is spanned over $B$ by all $e_{-\pi}v$ with $\pi \in {{\bf{{P}}}}_{X}(\eta)$.
\itemb
$\dim_B M^{\widetilde{\lambda}- \eta} \leq {\bf p}_{X}(\eta)$.\ei\ec
\bpf  Immediate.
\epf
We can in fact deduce more from the proof of the Proposition.  Let $L_A$ be the subspace of ${M_A}$ spanned over $A$ by all products $e_{-\pi} v$ where $\pi\in {\bf{P}}_X =\cup_{\eta \in Q^+} {\bf{P}}_{X}(\eta)$.
We are interested in situations with $L_A=M_A$.  This equality does not always hold,  see Theorem \ref{AA2} (f).  However we show that the condition holds if $\gl$ is replaced by $\gl+c\xi$ for all but finitely many values of $c$.
\bc \label{nob} Given $\gl, \xi$ as in Proposition \ref{rats}, set $\gl_c = \gl +c\xi$.  Then for all but finitely many $c \in \ttk$,
the weight space $M^{\widetilde{\lambda} -\eta}_A$ is spanned over $A$ by all $e_{-\pi}v$ with $\pi \in {\bf {P}}_{X}(\eta)$,
for all $\eta \in Q^+$. \ec

\bpf   For simplicity we assume the hypotheses of Theorem \ref{1aShap} hold.
By the choice of $\xi$ we have $(\ga^\vee, \xi) = a_\ga \neq 0$ for all $\ga\in \Gd^+_0$.
The leading term of $\gth_\gc$ evaluated at $\gl_c$ is up to a non-zero scalar multiple equal to
\be \label{grt} H_{\pi^{\gamma}}(\gl_c) = \prod_{\ga \in N(w^{-1})}((\ga^\vee,\gl) + ca_\ga),\ee
and each term in the product is zero for exactly one value of $c$.
For any other coefficient  of $H_\pi$ of $\gth_\gc$, $H_\pi(\gl_c)$ is a polynomial
in $c$ which  has degree strictly lower than the polynomial in \eqref{grt}. With these remarks the proof
is essentially the same as the proof of Theorem \ref{rats}.  The role of the indeterminate $T$ is played by $c$. \epf

\subsection{The Modules $M^{X}(\lambda).$}
Suppose $\gl \in \cH_X$,  and with $\widetilde{\lambda}$ as in \eqref{wld} define
\be \label{tuv} M^X({\widetilde{\lambda}})_B = M({\widetilde{\lambda}})_B/ \sum_{\gc \in X} U(\fg)_B \gth_\gc v_{{\widetilde{\lambda}}}.\ee
Then $M^X({\widetilde{\lambda}})_B$ is  a $U(\fg)_B$-module generated by a highest weight vector
$v^X_{{\widetilde{\lambda}}}$
(the image of $v_{{\widetilde{\lambda}}}$)
with weight  $\widetilde{\lambda}$. Set
$M^X({\widetilde{\lambda}})_A = U(\fg)_Av^X_{{\widetilde{\lambda}}} \subset
M^X({\widetilde{\lambda}})_B.$ Then

\[M^X({\widetilde{\lambda}})_A\ot_A B = M^X({\widetilde{\lambda}})_B.\]
\bt\label{zoo} 
The set \be \label{adb} \{e_{-\pi} v^X_{\widetilde{\lambda}}|\pi \in {{\bf{\overline{P}}}}_{X}(\eta)\},\ee is a
$B$-basis for
the weight space $M^X({\widetilde{\lambda}})_B^{\widetilde{\lambda} -\eta}$.
\et
\bpf By Proposition \ref{rats} the listed elements span $M^X({\widetilde{\lambda}})_B^{\widetilde{\lambda} -\eta}$, so it suffices to show
$\dim_B M^X(\widetilde{\lambda})_B^{\widetilde{\lambda}  - \eta}
={\bf p}_{X}(\eta ).$  Suppose that $X = Y \cup \{\gc\}$ where $\gc\notin Y.$
We show there is an exact sequence of $U(\fg)_B$-modules

\be\label{rod} 0 \lra L \lra {M^Y}(\widetilde{\gl})_B \lra N\lra 0\ee and surjective maps

\be\label{yim} {M^X}(\widetilde{\lambda}-\gc)_B \lra L, \quad {M^X}(\widetilde{\lambda})_B \lra N\ee
Indeed if $L=U(\fg)_Bw$ where
$w= \gth_\gc v^Y_{{\widetilde{\lambda}}}$, then $L$  is a highest weight module of weight
${\widetilde{\lambda}}-\gc.$ From 
Theorem \ref{fiix2}, and Theorem \ref{1zprod} we have $\gth_{\gc'}w = 0$ for all $\gc'\in X.$  Thus $L$ is an image of ${M^X}(\widetilde{\lambda}-\gc)_B$.  On the other hand, if $N$ is the cokernel of the inclusion of $L$ into ${M^Y}(\widetilde{\lambda})_B$,  it is clear that $N$ is    an image of
${M^X}(\widetilde{\lambda})_B.$  Hence using induction on $|X|$ for the first equality below, and then Equations \eqref{rod}, \eqref{yim} and Proposition \ref{rats} we obtain
\by {\bf p}_{Y}(\eta)
&=& \dim_B {M^Y}(\widetilde{\gl})_B^{\widetilde{\lambda}-\eta} =
\dim_B L^{\widetilde{\lambda}  - \eta} +\dim_B N^{\widetilde{\lambda}  - \eta}\nn\\
&\le& \dim_B {M^X}(\widetilde{\lambda}-\gc)_B^{\widetilde{\lambda}  - \eta} +
\dim_B M^X(\widetilde{\lambda})_B^{\widetilde{\lambda}  - \eta} \nn\\
&\le& {\bf p}_{X}(\eta - \gc)+{\bf p}_{X}(\eta ) ={\bf p}_{Y}(\eta). \nn \ey
It follows that equality holds throughout, and that the
maps in (\ref{yim}) are isomorphisms. Thus the
dimension of $M^X(\widetilde{\lambda})_B^{\widetilde{\lambda}  - \eta}$ is as claimed.
 \epf
\noi {\it Proof of Theorem \ref{newmodgen}.}
\noi The module $M^{X}(\lambda)$ defined in \eqref{hut} 
is generated by the image of $v^X_{{{\widetilde{\lambda}}}}$
which is a highest weight vector of weight $\lambda$.
 Finally the claim about the character of
$M^{X}(\lambda)$  follows from Theorem \ref{zoo} and the following considerations applied to the weight spaces $K, L$ of the modules
$M^{X}({\widetilde{\lambda} })_{R}$ for $R = A, B$ respectively.   If $K$ is
an ${A}$-submodule of a finite dimensional $B$-module $L$ such that $K\otimes_{A} B =
L$, then $\dim_\ttk K/TK = \dim_B L.$ \hfill  $\Box$
\bc \label{hco}
Suppose that  $\gamma$ is an odd isotropic root and $\lambda \in {\cH}_{\gc}$. Then the  kernel of the natural map
$M(\gl)\lra M^{\gamma}(\lambda)$ contains $U(\mathfrak{g})\theta_{\gamma}v_{{{\widetilde{\lambda}}}}$.
\ec
\bpf This follows from \eqref{tuv} since $\theta_{\gamma}v^X_{{{\widetilde{\lambda}}}}
\in U(\mathfrak{g})_{B}\gth_\gc v_{\widetilde{\lambda} }^X \cap U(\fg)_{A}v_{\widetilde{\lambda} }^X.$
\epf

\br\label{bod} {\rm We note some variations on \eqref{1let}. Suppose $\gl\in \cH_\gc$ for $\gc$ an isotropic root. By \eqref{rod} and \eqref{yim} with $Y=\emptyset$, we have an exact sequence
\be \label{duh} 0\lra M^{\gamma}({\widetilde{\lambda} -\gc})_{B}\lra {M}({\widetilde{\lambda} })_{B} \lra M^{\gamma}(\widetilde{\lambda})_{B}\lra 0\ee
where the first map sends $x v_{{\widetilde{\lambda} }-\gc}$ to $ x\gth_{\gc}v_{\widetilde{\lambda} }$.
\\
\\
The first map in \eqref{duh} also induces a map $M^{\gamma}({\widetilde{\lambda} -\gc})_{A}\lra {M}({\widetilde{\lambda} })_{A}$.  Clearly the kernel $N\cap T{M}({\widetilde{\lambda} })_{A}$ of the combined map
\be \label{bd}
N=M^{\gamma}({\widetilde{\lambda} -\gc})_{A}\lra {M}({\widetilde{\lambda} })_{A}
\lra {M}({\widetilde{\lambda} })_{A}/T{M}({\widetilde{\lambda} })_{A} \cong
 M(\lambda)
\ee
contains $TN$, but the containment can be strict, see Theorem \ref{A2} (c).  If this is the case the highest weight module $M^{\gamma}({{\lambda} -\gc})=N/TN$ will not embed in the Verma module
 $M(\lambda)$.
}\er

\subsection{Behavior in the most general cases.}\label{gb}
In the most general case, for $\gl \in \cH_X$ the modules $M^{X}(\gl)  $ are simple. Beyond this case we are interested in the behavior of $M^{X}(\gl) $  when $\gl$ lies on certain hyperplanes in $\cH_X$.  For example
in Proposition \ref{need}, we describe the general behavior when $\gl\in\cH_Y$ for an  orthogonal set of roots $Y$ containing $X$ such that $|Y|=|X|+1$.  

\bl \label{wry}\bi \itema Any orthogonal set of isotropic roots is linearly independent.

\itemb If $Y$ is an orthogonal set of isotropic roots and $\gb\in \Gd^+\backslash Y,$ then $\gb \notin \ttk Y.$
\ei \el

\bpf Part (a) is left to the reader.  It is similar to the proof of (b). Suppose $\gb\in \ttk Y.$ Then  $\gb$ is isotropic since $\ttk Y$ is an isotropic subspace of $\fh^*$.
Also $(\gb,Y)\neq 0$ by (a). The result is clear if $|Y|=1$, so we assume that $|Y|>1$.  This implies that $\fg$ is not exceptional.  Thus we have  $Y = \{\pm(\gep_{p_j}+a_j\gd_{q_j})\}_{j=1}^k,$ $\gb=(\gep_{p}\pm\gd_{q})$  where $a_j=\pm 1,$ and $ p, p_j \in [m],  q, q_j \in [n]$ for some $m, n$.  Suppose we have a relation
\[(\gep_{p}\pm\gd_{q}) +\sum_{j=1}^k b_j(\gep_{p_j}+a_j\gd_{q_j}) =0.\]
Since $\gep_{1},\ldots ,\gep_{m}$ are linearly independent we have $p=p_j$ for some $j$,
and $b_j=-1$, $b_\ell = 0$ for $\ell \neq j$.
Thus the relation is equivalent to  $\gep_{p}\pm\gd_{q}= \gep_{p_j}+a_j\gd_{q_j}.$ But this implies $q=q_j$ and $\gb \in Y$, a contradiction. \epf

\noi Suppose that $Y$ is an orthogonal set of (positive) isotropic roots. We introduce a suitable Zariski dense subset $\gL_Y$ of $\cH_Y$. 
For $S \subseteq {\overline{Y}}= \Gd^+_1\backslash Y,$ we have, for any highest weight vector $v_\gl$,  that $e_Se_{-S}v_\gl=p_S(\gl)v_\gl$ 
for some  $p_S\in S(\fh)$ with leading term $\prod_{\gb\in S} h_\gb$.  
It follows from Lemma \ref{wry} that if $V(p_S)$ is the zero locus of $p_S$, then $V(p_S) \cap \cH_Y$ is a proper closed subset of $\cH_Y$. (We remark that if $R$ is a subset of $S$ it need not be the case that $p_R$ divides $p_S$.)
Thus
\[ \gL_Y' = \cH_Y \backslash \bigcup_{S \subseteq {\overline{Y}}} V(p_S)\]
is a non-empty open subset of $\cH_Y.$
Now set
\[\gL_Y=\{\gl \in \gL_Y'|(\gl+\gr_0,\ga^\vee) \notin \bbZ \mbox{ for all non-isotropic roots } \ga \}.\]
Since $\gL_Y$ is obtained from $\gL_Y'$ by deleting a discrete countable union of hyperplanes, $\gL_Y$ is Zariski dense in $\cH_Y$. 
\\ \\
If $M$ is a $U(\fg)$-module, a {\it $\fg_0$-Verma flag} on $M$ is a
filtration by $\fg_0$-submodules
\be \label{end} 0 = M_0 \subset M_1 \subset \ldots \subset M_k = M\ee
such that for $i = 1,\ldots, k$
$M_i/M_{i-1}$ is a Verma module $M^0(\mu_i)$ for $U(\fg_0)$ with highest weight $\mu_i$. 
\bp \label{need} Suppose $Y$ is an isotropic and $Y = X\cup \{\gc\}$, where $\gc \notin X.$
\bi
\itema If  $\gl\in \gL_Y$, then $M^{Y}(\gl)$ is simple.
\itemb If  $\gl,\gl-\gc \in \gL_Y$,  there is a non-split exact sequence
\[0 \lra  M^{Y}(\gl-\gc) \lra  M^{X}(\gl) \lra M^{Y}(\gl)\lra 0\]
\ei \ep

\bpf
 By \cite{M} Theorem 10.4.5, $M(\gl)$ has a $\fg_0$-Verma flag.  The assumption on $\gl$ implies that if $\mu$ is the highest weight of any factor in this series, then $(\mu+\gr_0,\ga^\vee)$ is not an integer for any non-isotropic root $\ga$.  Thus the $U(\fg_0)$-Verma module
$M^0(\mu)$ is simple.  It follows that all $U(\fg_0)$-module composition factors of $M^Y(\gl)$ are Verma modules.   In addition we can order the subsets $S$ of ${\overline{Y}}$ so that the $U(\fg_0)$-submodules
$N_R(\gl) = \sum_{R\le S} U(\fg_0)e_{-S}v_\gl$ form a Verma flag in $M^Y(\gl)$. If $L$ is a non-zero submodule of $M^Y(\gl)$ we can choose $R$ so that $L\cap N_R(\gl) \neq 0,$ but $L\cap N_S(\gl) =0$ for any $S$ which properly contains $R.$ This implies that $e_{-R}v_\gl^Y$ is a $\fg_0$-highest weight vector in $L$. Hence (a) holds because $e_Re_{-R}v_\gl^Y$  is a non-zero multiple of $v_\gl^Y$. \\ \\
To prove (b), note first that $\gth_\gc v^X_{{{\lambda}}}$ is a highest weight vector in $M^{X}(\gl)$ so $U(\fg)\gth_\gc v^X_{{{\lambda}}}$ has a factor module which is isomorphic to $L(\gl-\gc)$.  However $M^{Y}(\gl-\gc)$ is a highest weight module with highest weight $\gl-\gc$. Hence $M^{Y}(\gl-\gc) \cong U(\fg)\gth_\gc v^X_{{{\lambda}}}$.  A similar argument shows that $M^{X}(\gl)/U(\fg)\gth_\gc v^X_{{{\lambda}}} \cong L(\gl).$ This gives the  sequence in (b). It does not split since $M^{X}(\gl)$ has a unique maximal submodule.
\epf

\section{The submodule structure of Verma modules.} \label{SV}\
In this section we apply Theorem \ref{Jansum} to determine the structure of Verma modules in the simplest cases.
First however we make some remarks concerning the Grothendieck group $K(\cO)$ and characters.
For $\lambda \in \mathfrak{h}^*$, set $D(\lambda) = \lambda - Q^+$ and
let $\mathcal{E}$ be the set of functions on $\mathfrak{h}^*$ which are zero
outside of a finite union of sets of the form $D(\lambda)$. Elements
of $\mathcal{E}$ can be written as formal linear combinations
$\sum_{\lambda \in \mathfrak{h}^*} c_{\lambda} \tte^{\lambda}$ where
$\tte^{\lambda}(\mu) = \delta_{\lambda\mu}$. We can make
$\mathcal{E}$ into an algebra using the convolution product, see \cite{M} Section 8.4 for details. 
 Let $C(\cO)$ be the additive subgroup of $\cE$ generated by the characters $\ch L(\gl)$ for $\gl \in \fh^*$.  It well known that  
there is an isomorphism from the group $K(\cO)$ to $C(\cO)$ sending $[M]$ to $\ch M$ for all modules $M \in \cO,$ see \cite{J1} Satz 1.11, \cite{M} Theorem 8.4.6. Hence we can work either in $K(\cO)$ or in $\cC(\cO)$ as it suits us.  A reason for doing the former was mentioned in Section C of the introduction.  However if we wish to carry out computations involving 
partition functions for example, then it is natural and easier to work with characters.

\bl \label{sat} Suppose $\lambda \in \mathfrak{h}^*$.
\bi \itema
If $A(\gl) = \{\ga\}$ and $B(\gl) = \emptyset$, then ${M}_{1}(\lambda) =
{M}(s_{\alpha}\cdot \lambda)$ and
${M}_{2}(\lambda) = 0,$
\itemb
If $B(\gl) = \{\gc\}$ and $A(\gl) = \emptyset$, then
${M}_{1}(\lambda) \cong M^{\gc}(\lambda -
\gc)$ and
${M}_{2}(\lambda) = 0.$
\ei
In both cases the submodule ${M}_{1}(\lambda)$ is simple. \el

\bpf Parts (a) and (b) follow easily from (\ref{lb}). The last
statement holds because ${M}_{1}(\lambda)$ is a self dual highest weight module. \epf
\noi If $M$ is a finitely generated $U(\fg)$-module, then  $M$ is finitely generated as a $U(\fg_0)$-module.  We can give $M$ a good filtration and then define the Gelfand-Kirillov
  dimension $d(M)$ and  Bernstein number  $e(M)$ of $M$ using the Hilbert polynomial of the associated graded $S(\fg_0)$-module.
  For a Verma module $M= M(\gl)$,  we set  $\ttd =d(M),$ and $\tte = e(M)$.
	Then $\ttd =|\Delta^+_{0}|,$  and $\tte = 2^{|{\Delta^+_{1}}|}$. For a more general result on induced modules, see
  \cite{M} Lemma 7.3.12. Also
  if the module $N$ has character $\tte^{\mu} {p}_{\gamma}$, we have $d(N)=\ttd$ and $e(N)=\tte/2.$

\bl \label{nl} If
$   0 = N_0 \subset N_1 \subset \ldots \subset N_k = {M}(\lambda)$ is a series of submodules of $M$,
then $\sum_i e(N_i/N_{i-1}) = \tte$ where the sum is over all factors $N_i/N_{i-1}$ with $d(N_i/N_{i-1})= \ttd$.\el \bpf See \cite{KrLe} Theorem 7.7.\epf
\noi We say that a $U(\fg)$-module $M$ is  {\it homogeneous} if $d(N) =d(M)$ for any non-zero submodule $N$ of $M.$ Any $\fg_0$ Verma module contains a unique minimal submodule, which is itself a Verma  module, it follows that a $\fg_0$ Verma module is homogeneous.
\bl \label{hom} If the  $U(\fg)$-module $M$ has a $\fg_0$-Verma flag, then $M$ is homogeneous. \el
\bpf The argument is well-known, but we outline the proof for convenience.
Consider  a filtration as in \eqref{end}.
 Let $N$ be a nonzero submodule of $M$
and choose $i$ minimal such that $N \cap M_i \neq 0$.  Then $N
\cap M_i$ is isomorphic to a nonzero submodule of $M^0(\mu_i)$
which is a homogeneous $U(\fg_0)$-module as observed above.
Hence
\begin{eqnarray*}
d(M^0(\mu_i)) &=& d(N \cap M_i) \\
                  & \leq & d(N) \leq d(M) .
\end{eqnarray*}
The result follows since $d(M^0(\mu_i)) = d(M) =
\ttd$.
\epf

\bc \label{veho} Any Verma module $M= M(\lambda)$ for $U(\fg)$
is homogeneous. \ec
\bpf
By \cite{M} Theorem 10.4.5, $M$ has a $\fg_0$-Verma flag. \epf

\bl \label{you}Suppose $A(\gl) = \{\ga\}$ and $B(\gl) = \{\gk\}$, and set $\mu =s_\ga \cdot \gl.$ Then

\bi \itema $A(\mu) = \emptyset$ and  $B(\mu) = \{s_\ga \gk\}$.
\itemb The unique  maximal ${}$ submodule of $\;M(\mu)\;$ is simple and  isomorphic to $\;{M}^{s_\ga \gk}(s_\alpha\cdot(\lambda -\gk)).$
\itemc If $N$ is a submodule of ${M}(\lambda)$ whose character is equal to $ \ch  M^\gk({\lambda -\gk})$ then
 ${M}(\mu)\cap N$ is the unique proper submodule of ${M}(\mu).$
\ei
 \el

 \bpf First $B(\mu) = \{s_\ga \gk\}$ since the bilinear form $(\;,\;)$, is $W$-invariant.
We prove $A(\mu) = \emptyset,$ in the case that $\Gd$ does not contain a non-isotropic odd root.
Let $\Gd_\gl = \{ \gc \in \Gd_0| (\gl+\gr,\gc^\vee) \in \Z\}$ be the integral subroot system determined by $\gl$. Then
$\Gd_\gl = \Gd_\mu$, and to show  $A(\mu) = \emptyset$ we may assume that  $\Gd_\gl = \Gd_0$.
Let $\mathtt P$ be the set of indecomposable roots in
$\Gd_0^+$.  Then $\mathtt P = \{\go_1, \go_2, \ldots, \go_r\}$ is a basis for the root system $\Gd_0$.
Since $A(\gl) = \{\ga\}$ it follows that $(\gl+\gr, \go_i^\vee)>0$ for some $i$ and then that $\ga=\go_i$.
But then $s_\ga$ permutes the positive roots other than $\ga$ and $A(\mu) = \emptyset$ follows.
 Part (b) follows from Lemma \ref{sat} (b) applied to $M(\mu)$. To prove (c) use the sum formula (\ref{lb})  in the form
\be \label{lib} \sum_{i > 0} [{M}_{i}(\lambda)] =[{M}(\mu)] + [M^{\gk}(\lambda -\gk)],\ee
Set  $N' = {M}(\mu)\cap N.$ If $N'= 0$, then
$M_1$ contains ${M}(\mu) \oplus N$. Combined with (\ref{lib}) and the hypothesis on $\ch N$ this implies that $M_2=0$ and
${M}(\mu) \oplus N=M_1$.  However $M_1/M_2$ is self-dual, which is impossible since it has
${M}(\mu)$ as a direct summand, and by (b) ${M}(\mu)$ is a non-simple highest weight module. Now $\tte(N) = \tte/2,$ so $N$ cannot contain a Verma submodule. Since
 $M(\mu)$ has length two by (b) the result follows.
 \epf

\noi Next for $q = (\gc, \alpha^\vee) \in \mathbb{N}
\backslash \{0\}$ where $\gc$ is isotropic, set $\gc'=s_\ga \gc = \gc-q\ga.$  We consider Verma modules $M(\gl)$ such that \be \label{but}
A(\gl) =\emptyset, \quad B(\gl) =\{\gc, \gc'\} .\ee Note that $B(\gl) =\{\gc, \gc'\}$ implies that $(\gl + \rho, \alpha^\vee)=0$. Now consider the set
\be \label{Gld}\gL=\{\gl| (\ref{but}) \mbox{ holds and }
A(\gl-\gc) =\emptyset, B(\gl-\gc) =\{\gc\}, A(\gl+\gc') =\emptyset, B(\gl+\gc') =\{\gc'\}
\}.\ee
 Since $(\gl+\gr-\gc, \alpha^\vee) = (\gl+\gr+\gc', \alpha^\vee)=-q$,   $\gL$ is a Zariski dense subset in $\mathcal H_{\gc} \cap \mathcal H_{\gc'}.$

\bl \label{nid} Suppose $\gl \in \gL$,  and that $K$ is the kernel of the map
 ${M}(\lambda)\lra {M}^{\gc'}(\lambda).$ Let $L$ be the socle of $M(\gl)$. Then the Jantzen filtration
on  ${M}(\lambda)$ satisfies
 $${M}_{3}(\lambda)=0\subset L={M}_{2}(\lambda)\subset K={M}_{1}(\lambda)\subset{M}(\lambda).$$
This is the unique composition series for $M(\gl)$.
Furthermore
$K=U(\fg)\gth_{\gc'} v_\gl,\;$
$L=U(\fg)\gth_\gc v_\gl\cong M^{\gc}(\gl- \gc)
\cong L(\gl- \gc)$, $K/L \cong L(\gl-\gc')$ and  $d(K/L)<\ttd$.
\el
\noi \bpf The hypotheses and Lemma  \ref{sat} imply that
$M(\gl-\gc)$ has length two with simple top
$M^\gc(\gl-\gc)$ and $M(\gl+\gc')$ has length two with simple socle
$M^{\gc'}(\gl)$. Thus $M^{\gc'}(\gl)$ is the unique simple image of $M(\gl)$, so $M_1(\gl)=K$.
Note $K$ has the same character as $M^{\gc'}(\lambda -\gc')$.  Thus $d(K) =\ttd$ and $e(K) = \tte/2.$
Now the socle of $M(\gl)$  contains a copy of $M^\gc(\gl-\gc)$, and we have
  $d(M^\gc(\gl-\gc)) =\ttd$ and $e(M^\gc(\gl-\gc)) = \tte/2.$
From Lemma \ref{hom}
we see that $L$ is isomorphic to $M^\gc(\gl-\gc)\cong L(\gl-\gc)$. Since $\gth_\gc v_\gl$ is a highest weight vector with weight $\gl-\gc$ we have $L = U(\fg)\gth_\gc v_\gl$.
It follows from this and  (\ref{lb}) that
in the Grothendieck group $K(\cO)$

\by \label{lob} \sum_{i > 0} [ {M}_{i}(\lambda)] &=& [ M^{\gc}(\lambda -\gc)]+ [ M^{\gc'}(\lambda -\gc')]\nn\\
&=& [K]+[L].\nn\ey
Therefore  $\sum_{i >1} [ {M}_{i}(\lambda)] = [L]$,
and this gives the statements about ${M}_2(\lambda)$ and ${M}_3(\lambda)$.  Finally $\ttd(K/L) <\ttd$ by Lemma \ref{nl}.
\epf
\bp \label{fad} Suppose that $A(\nu) = \{\ga\}$, $B(\nu) = \{\gc\}$ with
$p=(\nu + \rho, \alpha^\vee) \in {\mathbb N} \backslash\{0\}$, and that $\gl = \nu-\gc\in \gL$, as defined in \eqref{Gld}.
Set $\mu = s_\ga \cdot \nu$.  If $p =(\gc, \ga^\vee),$
then the lattice of submodules of $M={M}(\nu)$ is as in Figure \ref{fig1}, where
\be \label{cid} V_1  = U(\fg)\gth_{\ga,p}v_{\nu} \cong M(\mu), \quad V_2 =\Ker M(\nu)\lra M^\gc(\nu), \quad V_3 = U(\fg) \gth_{s_\ga\gc}\gth_{\ga,p}v_{\nu}.\ee
The unique maximal submodule of $V_i$ is  $L_i$ where
\be \label{cfcrs} L_1=L(\mu)\cong {M}^{s_\ga \gc}(\mu),\quad L_2 = L(\gl-\gc'),\quad L_3=L(\gl)\cong {M}^{s_\ga \gc}(\gl).\ee
The
Jantzen filtration is given by
\be \label{lid}{M}_{3}(\nu) = 0, \quad {M}_{2}(\nu) = V_3 =V_1 \cap V_2, \quad {M}_{1}(\nu) = V_1+ V_2.\ee
\ep

\bpf Define the $V_i$ and $L_i$ by \eqref{cid} and \eqref{cfcrs}.
Note that $(\gl+ \rho, \alpha^\vee) =0$, so $s_\ga \cdot \gl=\gl.$
Since  $(\nu + \rho, \alpha^\vee) \in \N \backslash\{0\},$
$M(s_\ga \cdot \nu)$ embeds in $M(\nu)$,
and by Lemma  \ref{sat}, $V_1= {M}(\mu)$  has length two with socle $V_3$.
 Now $\tte(V_3)= \tte(M(\nu)/V_3) = \tte/2$, so the same argument
as in Lemma \ref{nid} yields that $V_3$ is the socle of $M(\nu).$
Hence since $V_2$ has length two, $V_3 \subseteq V_1 \cap V_2 \subseteq V_1$
and  $V_2$ cannot contain $V_1$, by looking at characters. Thus $V_3 =V_1 \cap V_2$.
Now 
by (\ref{lb}) we have in the Grothendieck group $K(\cO)$

\be \label{sob} \sum_{i > 0} [ {M}_{i}(\nu)] = [ M^{\gc}(\nu -\gc)]+ [ M(\mu)].\ee
Since  $d(L_3)=d(L_1)= \ttd$,  and $e(L_3)=e(L_1)= \tte/2,$ we have

\be \label{zob} \tte\ge e (M_1(\nu)) \ge |{M}_{1}(\nu):L_3| + |{M}_{1}(\nu):L_1| \ge \tte.\ee
Thus $|{M}_{1}(\nu):L_3|\le 1$, 
but by (\ref{sob}),
$\sum_{i > 0} |M_{i}(\nu):L_3| =2,$
so $V_3\subseteq {M}_{2}(\nu).$
Hence

\begin{eqnarray}  \label{rob}
[ {M}_{1}(\nu)]+[V_3] +\sum_{i \ge 3} [ {M}_{i}(\nu)]
&\le&
\sum_{i > 0} [ {M}_{i}(\nu)]
=[ M(\mu)]+ [ V_2]\nonumber\\
&=&
[ M(\mu)+  V_2]
+[ M(\mu)\cap V_2]
\\
&\le& [ {M}_{1}(\nu)]+[V_3].
\nonumber\end{eqnarray}
Therefore ${M}_{3}(\nu)=0,
{M}_{2}(\nu) = V_3$ and ${M}_{1}(\nu) = V_2+V_1.$  The only thing left to show is that
$V_2/V_3\cong L(\gl-\gc').$
However we have $\ch V_2 = \tte^{\nu-\gc} {p}_\gc$, and by Lemma \ref{nid} with $\gl=\nu-\gc$, this is also the character of the length two module $M(\gl)/M_2(\gl)$.  Thus $V_2$ has  length two and $V_2/V_3$ is simple.
We have $\ch V_2/V_3= \ch V_2-\ch V_3$ and this can easily be calculated, compare
Lemma \ref{ink}, and this gives the correct highest weight.
\epf
\brs {\rm \label{2.6}
\noi \bi \itema In the situation of Proposition \ref{fad}, the modules $V_2$ and  $M(\gl-\gc)/M_2(\gl-\gc)$ need not be isomorphic, see Theorem \ref{A2} with $n=1$. The same example shows that  $V_2$ need not be a highest weight module.
\itemb The hypothesis in Lemma  \ref{nid} holds in Theorem \ref{AA2} when $n=0,$ with  $\gc'=\gb.$
\ei}\ers
\noi  In the Theorem below, both cases arise when $\fg =\fsl(2,1)$, (or $\fgl(2,1)$) using the anti-distinguished Borel subalgebra. In the notation of Subsection \ref{C8},  take $\gk =\gc$ if $q = 1,$ and $\gk = \gb$ if $q=-1.$

\bt Suppose  $q= (\gk, \ga^\vee) = \pm1$ and $p  = (\gl + \rho, \alpha^\vee)\in \mathbb{N} \backslash \{0\}$
where $\gk$, $\ga$ are isotropic and non-isotropic respectively.
Then  for general $\gl$ such that $A(\gl) = \{\ga\}$, $B(\gl) = \{\gk\}$,
 the lattice of submodules of $M={M}(\lambda)$ is as in Figure \ref{fig1} with the submodules $V_i$ as in $(\ref{cid}).$  Moreover the 
Jantzen filtration is given by $(\ref{lid})$ and if $p>q=1$ or $q=-1$, then $V_2 = U(\fg) \gth_{\gk}v_\gl$.
\et

\bpf Define $V_1, V_2, V_3$ by Equation  (\ref{cid}) with $\nu$ replaced by $\gl$, and $\gc$ by $\gk$. Set $\mu = s_\ga \cdot \gl.$  Then by Lemma \ref{you} $V_1 = M(\mu)$ has length 2 with socle $V_3$.
We use induction on $p$. 
If  $\nu =  \gl-q\gk,$ then $(\nu + \rho, \alpha^\vee)=p-1$. Different proof strategies are necessary
depending on the sign of $q$. If $q=1$ we use the map from $M(\nu)$ to $M(\gl)$ sending $v_\nu$ to $\gth_\gk  v_{\gl} $, and if $q=-1$,  we use the map from $M(\gl)$ to $M(\nu)$ sending $v_\gl$ to $\gth_\gk  v_{\nu} $.
\\ \\
Suppose that $q=1.$ If $p=1$, then all the assertions hold by Proposition \ref{fad},
so assume  $p\ge 2$. By Proposition \ref{fad} if $p=2$, or induction if $p>2$, $M(\nu)$ has a length two factor module with character
$\tte^\nu {p}_\gk = \ch V_2$.
Hence $V_2$ has length two with socle isomorphic to $V_3$.
Also $U(\fg) \gth_\gk  v_{\gl} \subseteq V_2$ by Corollary \ref{hco}.  Since
$\gth_\gk  v_{\gl}$ is a highest weight vector with weight $\gl-\gk$, and the only highest weight of $V_3$ is
$s_\ga\cdot(\gl-\gk) \neq \gl-\gk$, it follows that
$U(\fg) \gth_\gk  v_{\gl} =V_2$.
Because $e(M(s_\ga\cdot\gl)) =\tte$,
$L(s_\ga\cdot(\gl-\gk))=\tte/2$
 and $M(s_\ga\cdot\gl)$ contains
$L(s_\ga\cdot(\gl-\gk))$  with multiplicity one, we have $|M(\gl):L(s_\ga\cdot(\gl-\gk))|=|{M}_{2}(\gl):L(s_\ga\cdot(\gl-\gk))|=1.$
Now the sum formula \eqref{lb} takes the form
\[
\sum_{i\ge1}[M_i(\gl)] = [V_1]+ [V_2] =
[L(s_\ga\cdot\gl)]+
[L(\gl-\gk)]+
2[L(s_\ga\cdot(\gl-\gk))].
\] This easily gives
$M_1(\gl) = V_1+V_2$, $M_2(\gl) = V_1\cap V_2=V_3$  and $M_3(\gl) = 0.$
\\ \\
Now suppose
$q=-1.$  If $p=1,$ the conditions of Lemma \ref{nid} hold in the general case with
$\nu, \gk'=s_\ga\cdot\gk$ and $\gk$ in place of $\gl$ and $\gc$ and $\gc'$ respectively.
Thus
$\gth_{\gk'}  v_{\nu}$ generates
$M_2(\nu)$, and the map $xv_\gl \lra x\gth_{\gk} v_\nu$ induces an isomorphism
$M(\gl)/V_2\lra M_1(\nu)=U(\fg) \gth_\gk  v_{\nu}\cong M^\gk(\nu-\gk)$. 
 Therefore
 \be \label{rid} \ch (M(\gl)/V_2) =\ch M^\gk(\nu-\gk) = \tte^\gl {p}_\gk. \ee
If $p>1$ then (\ref{rid}) holds by induction and a similar argument.
Now as in the proof of Proposition \ref{fad}, (see Equation \eqref{zob}) we see that $|{M}_{1}(\lambda):V_3|\le 1$ and (\ref{rob}) holds.
\\ \\
Also from (\ref{rid}) $M(\gl)/V_2$ has length two.
We know $M_1(\gl)/M_2(\gl) = V_1/V_3 \oplus V_2/V_3$ where  $V_1/V_3\cong L(s_\ga\cdot \gl)$ and $V_2/V_3$ is a highest weight module with
highest weight $\gl-\gk$ such that $|V_2/V_3: L(s_\ga\cdot \gl)| =0.$ Hence $V_2/V_3$ is a self-dual highest weight module and so is simple.
\epf

\noi Next we consider the structure of $M(\gl)$ when
$B(\lambda)= \{\gc, {\gc'}\}$ for orthogonal roots $\gc$ and $\gc'$.

\bl \label{fix} If $\gl\in \cH_{\gc} \cap \cH_{\gc'},$ then for all but only finitely many $c \in \ttk$ we have
\be \label{we1}\gth_\gc(\lambda +c\xi - \gc)\gth_{{\gc'}}(\gl+c\xi)\neq 0\ee and
\be \label{we2}\gth_{{\gc'}}(\lambda +c\xi - \gc)\gth_\gc(\gl+c\xi)\neq 0.\ee\el
\bpf Set $\widetilde{\lambda}= \gl +T\xi$.  It follows  from Corollary \ref{rmd} that when $\gth_\gc(\widetilde{\lambda}- \gc)
\gth_{{\gc'}}(\widetilde{\lambda})v_{\widetilde{\lambda}}$
 is written as a $A$-linear combination of terms $e_{-\pi}v_{\widetilde{\lambda}}$, the coefficient of $e_{-\gc}e_{\gc'} v_{\widetilde{\lambda}}$ is a polynomial in $T$ of degree $d_\gc +d_{{\gc'}}$. Hence  (\ref{we1}) holds for all but  finitely many $c$, and a similar argument applies to (\ref{we2}).\epf
\noi \bl \label{tri} For general $\gl \in \cH_{\gc} \cap \cH_{{\gc'}}$ we have
$$[M(\gl):L(\gl-\gc-\gc')]=1.$$ \el

\bpf \ff{This proof is due to Vera Serganova.}
 We  require
$\gth_\gc(\lambda  - \gc)\gth_{{\gc'}}(\gl)\neq 0$ or $\gth_{{\gc'}}(\lambda - \gc)\gth_\gc(\gl)\neq 0$,
as well as some further conditions that arise in the proof. This implies $[M(\gl):L(\gl-\gc-\gc')]\ge 1.$ By \cite{M}, Theorem 10.4.5
$M(\gl)$ has a series with factors  which are Verma modules $M^0(\mu)$  for  $\fg_0$, and $M^0(\gl-2\gr_1)$ occurs exactly once as a factor in this series. Furthermore $M^0(\gl-2\gr_1) =L^0(\gl-2\gr_1)$ for general $\gl$, so as a $\fg_0$-module
$$[M(\gl):L^0(\gl-2\gr_1)]=1.$$
So it is enough to show that  \be \label{tag}[L(\gl-\gc-\gc'):L^0(\gl-2\gr_1)]\ge1.\ee
or equivalently
$[L({{{\overline{\gl}}}}):L^0({\overline{\gl}}-\eta)]\ge1$
 where $\overline{\gl} = \gl-\gc-\gc',$ and $-\eta =\gc+\gc' -2\gr_1$. Let $X$ be the set of all positive isotropic roots $\gs$ different from $\gc$ and $ \gc'$, and $e_{-X}$ be the ordered product of all root vectors $e_{-\gs}$ where $\gs\in X$.  The weight space $M(\overline{\gl})^{\overline{\gl}-\eta}$ has a basis consisting of vectors of the form $e_{-\pi}v_{\overline{\gl}}$, with $\pi$ a partition of $\eta$, and
 $e_{-X}v_{\overline{\gl}}$ is one such basis element.
We claim that for  general ${\overline{\gl}}$,

\be \label{tug}e_{-X} v_{\overline{\gl}} \notin \fn_0^-M({\overline{\gl}}) + I_{\overline{\gl}}\ee
where $I_{\overline{\gl}}$ is the maximal submodule of $M({\overline{\gl}}).$  
 Let $(\;,\;)$ be a contravariant form on $M({\overline{\gl}})^{{\overline{\gl}}-\eta}$ with radical equal to $I({\overline{\gl}})^{{\overline{\gl}}-\eta}$. By the proof of \cite{M} Lemma 10.1.2, we see that as a polynomial in ${\overline{\gl}}$, the degree of
$g({\overline{\gl}}) =(e_{-X}v_{\overline{\gl}}, e_{-X}v_{\overline{\gl}})$ is greater than  the degree of $(wv_{\overline{\gl}},e_{-X}v_{\overline{\gl}})$ for any $w \in (\fn_0^-U(\fn^-)v_{\overline{\gl}})^{{\overline{\gl}}-\eta}$. Moreover the leading  term of $g({\overline{\gl}})$ is $\prod_{\gs \in X}({\overline{\gl}},\gs)$ which is non-zero at general elements of $\cH_{\gc} \cap \cH_{\gc'}$.  This implies (\ref{tug}). We deduce from this that the image of $e_{-X}v_{\overline{\gl}}$ in $H_0(\fn_0^-,L({\overline{\gl}}))^{\overline{\gl}-\eta}$ is non-zero.  This says that
$L({\overline{\gl}})$ contains a $\fg_0$ highest weight vector of weight ${\overline{\gl}-\eta}$,
 and (\ref{tag}) follows.
\epf

\noi This gives a representation theoretic proof of \eqref{pee}.
\bc \label{kt1} There is a rational function $a$ of $\gl \in \cH_{\gc} \cap \cH_{{\gc'}}$ such that
\be \label{kt}
\gth_{\gc'}(\lambda - \gc)
\gth_\gc(\gl)=a(\gl)\gth_ \gc(\gl - {\gc'})  \gth_ {\gc'}(\gl).\ee
\ec

\bpf
Up to a scalar multiple, $M(\gl)$ can contain at most one highest weight vector with weight $\gl-\gc-{\gc'}$.  However both
$\gth_{\gc'}(\lambda - \gc)
\gth_\gc(\gl)v_\gl$ and $\gth_ \gc(\gl - {\gc'})  \gth_ {\gc'}(\gl)v_\gl$ are both highest weight vectors with this weight, so (\ref{kt}) holds.
\epf

\section{An (ortho) symplectic example.} \label{1cosp}
\bexa \label{1ex2.4}
{\rm  A crucial step in the construction of \v Sapovalov elements was the observation in the proofs of Lemmas \ref{1wpfg}  and \ref{1wpfg1} that the term $c_{ \zeta,\lambda}$ defined in Equation (\ref{159}) are zero unless $\zeta(\ga) \geq
N-mq,$ (using the notation of the Lemmas).  We give an example where the individual terms on the right of Equation (\ref{159}) are not identically zero, and verify directly that  the sum itself is zero. This cannot happen in Type A.  The key difference in the examples below seems to be that it is necessary to apply Equation (\ref{121nd}) more than once with the same simple root $\ga$. Consider the Dynkin-Kac diagram below for the Lie superalgebra $\fg = \osp(2,4)$.}\eexa

\vspace{0.4cm}

\vspace{0.4cm}
\begingroup
\setlength{\unitlength}{0.10in}
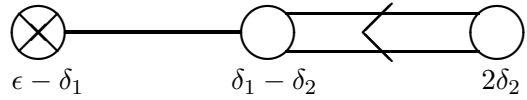
\begin{picture}(-30,-10)
\thicklines
\put(14.414,0.0){\line(1,0){9.23}}
\put(25.95,1.0){\line(1,0){10.1}}
\put(25.95,-1.0){\line(1,0){10.1}}
\put(30,0){\line(1,-1){1.5}}
\put(30,0){\line(1,1){1.5}}
\put(12,1){\line(1,-1){2}}%
\put(12,-1){\line(1,1){2}}
\put(13,0){\circle{2.828}}
\put(25,0){\circle{2.828}}
\put(37,0){\circle{2.828}}
\put(11.6,-3.0){$\gep-\gd_1$}
\put(23.1,-3.0){$\gd_1-\gd_2$}
\put(36.2,-3.0){$2\gd_2$}
\end{picture}

\vspace{1.4cm}

\endgroup

\noi Let $\gb = \epsilon - \gd_1, \ga_1 = \gd_1 - \gd_2, \ga_2 = 2\gd_2$, be the corresponding simple roots.  If we change the grey node to a white node we obtain the Dynkin diagram for $\fsp(6)$.  In this case the simple roots are
$\gb=\delta_{0} - \delta_{1}, \ga_1 = \delta_{1} - \delta_{2} $ and $\ga_2=2 \delta_{2}$. Let
$e_{-\gb}, e_{-\ga_1}, e_{-\ga_2}  $ be the negative simple root vectors.
The computation of  the \v Sapovalov elements
$\gth_1, \gth_2, \gth_3$ for the roots $\gb + \ga_1, \gb + \ga_1+ \ga_2$ and $\gb + 2\ga_1+ \ga_2$ respectively,
is the same for $\osp(2,4)$ and for $\fsp(6)$.
Then define the other negative root vectors by
$$e_{- \ga_1 - \ga_2} = [e_{-\ga_1},e_{- \ga_2 }  ],  \quad \quad
e_{- 2\ga_1 - \ga_2} = [e_{-\ga_1}, e_{- \ga_1 - \ga_2}],  $$
$$e_{-\gb- \ga_1} =[e_{-\ga_1} ,e_{-\gb}
], \quad e_{-\gb- \ga_1-\ga_2} = [e_{- \ga_2 } ,e_{-\gb- \ga_1}], \quad e_{-\gb- 2\ga_1 - \ga_2} = [e_{-\ga_1},e_{-\gb- \ga_1-\ga_2}]. $$
It follows from the Jacobi identity that
$$[e_{-\gb}, e_{- \ga_1 - \ga_2}] = e_{-\gb- \ga_1-\ga_2},
\quad \quad
[e_{- \ga_1 - \ga_2},e_{-\gb- \ga_1}] = e_{-\gb- 2\ga_1 - \ga_2},$$ and
$$[e_{-\gb},e_{- 2\ga_1 - \ga_2}] = 2e_{-\gb- 2\ga_1-\ga_2}.$$
We order the set of positive roots so that for any  partition $\pi$, $e_{-\ga_1}$ occurs first if at all in $e_{-\pi}$,  and any root vector $e_{-\gs}$ with $\gs$ an odd root occurs last.\\ \\
Let $s_1, s_2$ be the reflections corresponding to the simple roots $\ga_1, \ga_2.$ 
Then for $\gl \in \fh^*$ define $\gl_1 =s_1\cdot\gl, \;\gl_2 = s_2\cdot\gl, \;\mu = s_1\cdot\gl_2$. Let
$$(\gl + \gr,\ga^\vee_1) = p = -(\mu + \gr,(\ga_1+\ga_2)^\vee) $$ and  \[(\gl_1 + \gr,\ga^\vee_2) = (\gl + \gr,(2\ga_1  + \ga_2)^\vee) = q= -(\mu + \gr,(2\ga_1+\ga_2)^\vee) ,\] $$(\gl_2 + \gr,\ga^\vee_1) = (\gl + \gr,(\ga_1+\ga_2)^\vee) = r = -(\mu + \gr,\ga_1^\vee) .$$
Then $r = 2q - p.$ Let $\gc$ be any positive root that involves $\gb$ with non-zero coefficient when expressed as a linear combination of simple roots. We compute the
\v Sapovalov elements
$\gth_{\gc,1}$ for $\fsp(6)$ and $\gth_{\gc}$ for $\osp(2,4).$
To do this we use Equation (\ref{121nd}).   We can assume $\gc \neq \gb.$ Suppose that $p, q, r$ are nonnegative integers. Then
$$e_{-\ga_1}^{p+1}e_{-\gb} = \gth_1e_{-\ga_1}^{p}$$
$$e_{- \ga_2 }^{q+1}\gth_1 = \gth_2e_{- \ga_2 }^{q}$$    $$e_{- \ga_1 }^{r+1}\gth_2 = \gth_3e_{- \ga_1 }^{r}.$$ 
In the computations below we write  $e_{-\pi}$, for $\pi$ a partition (resp. $\gth_i$ for $i=1,2,3$) in place of  $e_{-\pi}v_\gl$ (resp. $\gth_i v_\gl$).
First note that
$$[e_{-\ga_1}^{p+1} ,e_{-\gb}
] = (p+1)e_{-\gb- \ga_1} e^p_{-\ga_1}$$
$$[e_{- \ga_2 }^{q+1} ,e_{-\gb- \ga_1}] =
(q+1)e_{-\gb- \ga_1-\ga_2} e_{- \ga_2 }^q$$
$$[e_{- \ga_2 }^{q+1} ,e_{-\ga_1}]  =-(q+1)e_{- \ga_1 - \ga_2}
e_{- \ga_2 }^q.$$
This easily gives
 $$\gth_1 = (p+1)e_{-\gb- \ga_1} + e_{-\gb}e_{-\ga_1} = pe_{-\gb- \ga_1} + e_{-\ga_1}e_{-\gb}.$$
We order the set of positive roots so that for any  partition $\pi$, $e_{-\ga_2}$ occurs last if at all in $e_{-\pi}$,  and any root vector $e_{-\gs}$ with $\gs$ an odd root occurs first.
\begin{eqnarray}\label{1la2}
\gth_2
&=& (p+1)[(q+1)e_{- \gb - \ga_1 - \ga_2} +e_{-\gb- \ga_1}e_{-\ga_2}] +e_{-\gb}[ e_{-\ga_1}e_{- \ga_2} - (q+1)e_{- \ga_1 - \ga_2} ].
\nonumber
\end{eqnarray}
Next order the set of positive roots so that for any  partition $\pi$, $e_{-\ga_1}$ occurs last if at all in $e_{-\pi}$,  and any root vector $e_{-\gs}$ with $\gs$ an odd root occurs first. To find $\gth_3$  we use
$$[e_{-\ga_1}^{r+1} ,e_{-\gb-\ga_1-\ga_2}
] = (r+1)e_{-\gb- 2\ga_1 -\ga_2} e^r_{-\ga_1},$$
$$[e_{-\ga_1}^{r+1} ,e_{-\gb-\ga_1}e_{-\ga_2}
] = (r+1)e_{-\gb- \ga_1}e_{-\ga_1-\ga_2}e^
{r}_{-\ga_1}
+\left( \begin{array}{c}
                r+1 \\
                2 \end{array}\right)e_{-\gb- \ga_1}e_{-2\ga_1 -\ga_2} e^{r-1}_{-\ga_1},$$
$$[e_{-\ga_1}^{r+1} ,e_{-\gb}e_{-\ga_1-\ga_2}
] = (r+1)[e_{-\gb}e_{- 2\ga_1 -\ga_2} e^r_{-\ga_1}
 +e_{-\gb- \ga_1}e_{-\ga_1 -\ga_2} e^r_{-\ga_1} +re_{-\gb-\ga_1}e_{-2\ga_1-\ga_2}
e_{\ga_1}^{r-1}],$$
$$e_{-\ga_1}^{r+1} e_{-\gb}e_{-\ga_2}e_{-\ga_1}
 =
e_{-\gb}[e_{-\ga_2}
e_{- \ga_1}^{2} +
(r+1)e_{- \ga_1-\ga_2}e_{-\ga_1}
+\left( \begin{array}{c}
                r+1 \\
                2 \end{array}\right)e_{- 2\ga_1 -\ga_2}] e^{r}_{-\ga_1}$$
$$+ (r+1)e_{-\gb- \ga_1}[e_{-\ga_2}
e_{- \ga_1} +
re_{- \ga_1-\ga_2}]e^{r}_{-\ga_1}
+(r-1)\left( \begin{array}{c}
                r+1 \\
                2 \end{array}\right)e_{-\gb}e_{- 2\ga_1 -\ga_2} e^{r-1}_{-\ga_1}.$$
The above equations allow us to write $e_{- \ga_1 }^{r+1}\gth_2$  in terms of elements $e_{-\pi}$ with $\pi$ a partition of $\gb + (r+2)\ga_1 + \ga_2$.
We see that the term $e_{-\gb- \ga_1}e_{- 2\ga_1 -\ga_2} e^{r-1}_{-\ga_1}$  occurs in $e_{- \ga_1 }^{r+1}\gth_2$ with coefficient
\[\left( \begin{array}{c}
                r+1 \\
                2 \end{array}\right)[(p+1) -2q + (r-1)] = 0.\]
This is predicted by the cancellation step in the proof of Lemma \ref{1wpfg}. In the remaining terms, $e_{- \ga_1 }^{r}$ can be factored on the right, and this yields
\begin{eqnarray}\label{1th3}
\gth_3 &=&
(p+1)(q+1)(r+1)e_{-\gb- 2\ga_1 -\ga_2} +(p+1)(q+1)e_{-\gb- \ga_1-\ga_2}e_{-\ga_1} \\
&+& (q+1)(r+1)e_{-\gb- \ga_1}e_{-\ga_1 -\ga_2}-(p/2)(r+1)e_{-\gb}e_{- 2\ga_1 -\ga_2}\nonumber \\
&+& 2(q+1)e_{-\gb- \ga_1}e_{-\ga_2}e_{-\ga_1}+ (r-q+1)e_{-\gb}e_{-\ga_1-\ga_2}e_{-\ga_1}+e_{-\gb}e_{-\ga_2}
e_{- \ga_1}^{2}. \nonumber
\end{eqnarray}
Using the opposite orders on positive roots to those used above to define the $e_{-\pi}$ we obtain
\begin{eqnarray}\label{1la3}
\gth_2 &=& p[qe_{- \gb - \ga_1 - \ga_2} +
e_{- \ga_2 }
e_{-\gb- \ga_1}]
+[ e_{-\ga_2}
e_{- \ga_1 } - qe_{- \ga_1 - \ga_2} ]e_{-\gb},
\nonumber
\end{eqnarray}
and
\begin{eqnarray}\label{1tho}
\gth_3 &=&
pqr e_{-\gb- 2\ga_1 -\ga_2} +pq
e_{-\ga_1}
e_{-\gb- \ga_1-\ga_2}\\
&+& qr e_{-\ga_1 -\ga_2}e_{-\gb- \ga_1}
-(r/2)(p+1)e_{- 2\ga_1 -\ga_2}e_{-\gb}\nonumber \\
&+& 2q e_{-\ga_1}e_{-\ga_2}e_{-\gb- \ga_1} +
 (r-q-1)e_{-\ga_1}e_{-\ga_1-\ga_2}e_{-\gb}
+e_{- \ga_1}^{2}
e_{-\ga_2}e_{-\gb}.\nonumber
\end{eqnarray}
\br {\rm It seems remarkable that all the coefficients of $\gth_3$ in (\ref{1th3}) and (\ref{1tho}) are products of linear factors.  This is also true in the Type A case, see Equations (\ref{1shtpa}) and (\ref{1shtpb}). A partial explanation of this phenomenon is given by specializing these coefficients to zero.
Vanishing of these coefficients gives rise to factorizations of $\gth_3$ as in the examples below.
} \er
\bi \itema If $p=(\gl + \gr,\ga_1^\vee) =0$, then $r = 2q$ and we have \begin{eqnarray}\label{1la4}
\gth_3 &=&
 2q^2 e_{-\ga_1 -\ga_2}e_{-\gb- \ga_1}
-qe_{- 2\ga_1 -\ga_2}e_{-\gb}\nonumber \\
&+& 2q e_{-\ga_1}e_{-\ga_2}e_{-\gb- \ga_1} +
 (q-1)e_{-\ga_1}e_{-\ga_1-\ga_2}e_{-\gb}
+e_{- \ga_1}^{2}
e_{-\ga_2}e_{-\gb} \nonumber \\&=& \gth_{\ga_1+\ga_2}\gth_{\gb+\ga_1}.\nonumber
\end{eqnarray}
\itemb
If $q=(\gl + \gr,(2\ga_1  + \ga_2)^\vee) =0$ then $p = -r$, $\gth_2=
\gth_{\gb +\ga_1+\ga_2} = \gth_{\ga_2}\gth_{\gb+\ga_1}$, and we have \begin{eqnarray}\label{1la5}
\gth_3 &=&
[(p/2)(p+1)e_{- 2\ga_1 -\ga_2}
 -(p+1)e_{-\ga_1}e_{-\ga_1-\ga_2}
+e_{- \ga_1}^{2}
e_{-\ga_2}]e_{-\gb}\nonumber \\&=& \gth_{2\ga_1+\ga_2}\gth_{\gb}.\nonumber
\end{eqnarray}
\itemc If $r=(\gl + \gr,(\ga_1  + \ga_2)^\vee) = 0$, then $p = 2q$,
 and
we have \begin{eqnarray}\label{1la6}
\gth_3 &=&
2q^2
e_{-\ga_1}
e_{-\gb- \ga_1-\ga_2}\nonumber \\
&+& 2q e_{-\ga_1}e_{-\ga_2}e_{-\gb- \ga_1} -
 (q+1)e_{-\ga_1}e_{-\ga_1-\ga_2}e_{-\gb}
+e_{- \ga_1}^{2}
e_{-\ga_2}e_{-\gb}\nonumber \\
 &=&
e_{-\ga_1}[2q^2
e_{-\gb- \ga_1-\ga_2}
 2q e_{-\ga_2}e_{-\gb- \ga_1} -
 (q+1)e_{-\ga_1-\ga_2}e_{-\gb}
+e_{- \ga_1}
e_{-\ga_2}e_{-\gb}]\nonumber \\&=& \gth_{\ga_1}\gth_{\gb+\ga_1+\ga_2}.\nonumber
\end{eqnarray}
\ei
Similarly if $p=-1,$ (resp. $q=-1,$ $r=-1$) then (\ref{1tho}) yields the factorizations $\gth_3=\gth_{\gb+\ga_1}\gth_{\ga_1+\ga_2}$, (resp. $\gth_2=
 \gth_{\gb+\ga_1}
\gth_{\ga_2},$ $\gth_3=\gth_{\gb}\gth_{2\ga_1+\ga_2}$, and  $\gth_3=\gth_{\gb+\ga_1+\ga_2}\gth_{\ga_1}$). On the other hand we see that $p$ divides the coefficients of $e_{-\gb- 2\ga_1 -\ga_2}$ and $e_{-\ga_1}
e_{-\gb- \ga_1-\ga_2}$ in (\ref{1tho}) since when $p = 0$, $\gth_3=\gth_{\ga_1+\ga_2}\gth_{\gb+\ga_1}$ can be written as a linear combination of different $e_{-\pi}$. In this way we obtain explanations for all the linear factors in (\ref{1th3}) and (\ref{1tho}) with the exception of the coefficients $r-q\pm 1$ of
$e_{-\gb}e_{-\ga_1-\ga_2}e_{-\ga_1}$ and $e_{-\ga_1}e_{-\ga_1-\ga_2}e_{-\gb}.$
At this point it may be worthwhile mentioning that $r-q =(\gl + \gr,\ga_2^\vee)$.  In addition equality holds in the upper bounds given in Theorem \ref{1Shap} for the degrees of all the coefficients in (\ref{1th3}) and (\ref{1tho}).

\section{The Type A Case.}\label{1s.8}
\subsection{Lie Superalgebras.}\label{1far} We construct the elements $\theta_{\gamma}$ in Theorem \ref{1Shap}
explicitly when ${\mathfrak g} = \fgl(m,n).$  Suppose that $\gc = \gep_r
- \gd_{s}.$ For $1 \leq  i<j \leq  m$ and $1 \leq  k<\ell \leq  n$
define roots $\gs_{i,j}, \gt_{k,\ell}$ by
\[ \gs_{i,j} = \gep_{i} - \gep_{j}, \quad \gt_{k,\ell} = \gd_k - \gd_\ell. \]
Suppose   $B = (b_{i,j})$ is a $k \times \ell$ matrix with entries
in $U({\mathfrak n}^{-})$, $I \subseteq \{1, \ldots , k \}=[k]$ and $J
\subseteq \{1, \ldots , \ell \}.$  We denote the submatrix of $B$
obtained by deleting the $ith$ row for $i \in I,$ and the $jth$
column for $j \in J$ by $_IB_J.$ If either set is empty, we omit the
corresponding subscript. When $I = \{ i\},$ we write $_iB$ in place
of $_IB$  and likewise when $|J| = 1.$ If $I$ or $J$ equals $[p]$
we write $_{[p]}B$ or $B_{[p]}$.
If   $k = \ell$ we define two noncommutative determinants of $B,$
the first working from left to right, and the second working from
right to left.
\be \label{lr}{\stackrel{\longrightarrow }{{\rm det}}}(B) =  \sum_{w \in \mathcal{S}_k} sign(w) b _{w(1),1} \ldots b _{w(k),k}, \ee
\be \label{rl}{\stackrel{\longleftarrow}{{\rm det}}}(B) =  \sum_{w \in \mathcal{S}_k} sign(w) b _{w(k),k} \ldots
b _{w(1),1}.\ee If $k = 0$ we make the convention that
${\stackrel{\longleftarrow}{{\rm det}}}(B) =
{\stackrel{\longrightarrow}{{\rm det}}}(B) = 1.$ 
We call ${\stackrel{\longrightarrow}{{\rm det}}}(B) $ and ${\stackrel{\longleftarrow}{{\rm det}}}(B)$ the {\it LR and RL determinants} of $B$ respectively. We note the following cofactor expansions down the first and last columns
\by {\stackrel{\longrightarrow }{{\rm det}}}\;B
&=& \sum_{j=1}^k (-1)^{j+1} b_{j1} {\stackrel{\longrightarrow }{{\rm det}}}({}_{j}B_{1})\nn\\ 
&=& \sum_{j=1}^k (-1)^{k+j+1}  {\stackrel{\longrightarrow }{{\rm det}}}({}_{j}B_{k})b_{jk}.\nn
\ey   
These are easily derived from \eqref{lr} by grouping  terms. 
There is a similar formula for cofactor expansion down the last column.
 \\ \\
Consider the
following matrices with entries in $U(\fn^-)$
 \be\label{wok}
A^+(\gl, r) =  \left[ {\begin{array}{ccccc}
e_{r+1,r} &  e_{r+2,r} & \hdots & e_{m,r}\\
-a_1 & e_{r+2,r+1} & \hdots & e_{m,r+1} \\
0&-a_2 & \hdots & e_{m,r+2} \\
 \vdots  & \vdots & \ddots & \vdots              \\
0&  \hdots & \ldots&
e_{m,m-1}   \\
0 & \hdots & 0& -a_{m-r}\\
\end{array}}
 \right],\ee

\be \label{gov} A^-(\gl, s) =  \left[ {\begin{array}{ccccc}
e_{m+s,m+s-1}&e_{m+s,m+s-2} & \hdots & e_{m+s,m+2}  & e_{m+s,m+1} \\
 b_{s-1} &e_{m+s-1,m+s-2} & \hdots & e_{m+s-1,m+2} &  e_{m+s-1,m+1} \\
0& b_{s-2} & \hdots & e_{m+s-2,m+2} &   e_{m+s-2,m+1} \\
 \vdots  & \vdots & \ddots & \vdots& \vdots \\ 
0&0&  \hdots & b_2& e_{m+2,m+1} \\
0& 0& \hdots & \hdots 0& b_1  
\end{array}}
 \right],
\ee 
where for $\gl\in \cH_\gc$ we set
\be \label{eat} a_i=(\gl + \gr,\gs_{r,r+i}^\vee ) \mbox{ and } b_i = (\gl + \gr,\gt_{i,s}^\vee ).\ee
Also  let $B^+(\gl, r)$ (resp. $B^-(\gl, s)$) be the matrices obtained from $A^+(\gl, r)$ (resp. $A^-(\gl, s)$) by replacing each $a_i$ by $a_i +1$ (resp. replacing each $b_i$ by $b_i +1$).
\\ \\

Observe that in $A^+(\gl, r)$ and $A^-(\gl, s)$ the number of rows exceeds the number of columns by one. We also consider two degenerate cases.  In general $A^+(\gl, r)$ and $B^+(\gl, r)$ are $m-r+1\ti m-r$ matrices, so
if $r=m$ then $A^+(\gl, r)$
and $B^+(\gl, r)$ are ``matrices with zero columns"
  In this case we ignore the summation over $j$ in the following formulas,  replacing
${\stackrel{\longrightarrow}{{\rm det}}}(_jA^+(\gl, r))$ by 1, $e_{m+i,j+r-1}$ by $e_{m+i,m}$  and $i+j+r+m$ by $i+1+r+m$. Similar remarks apply to the case where $s = 1.$ 
\\ \\ 
Below we present two determinantal formulas for the \v Sapovalov element $\gth_\gc$ evaluated at $\gl\in \cH_\gc$.  
The key  differences are the placement of the odd root vectors $e_{m+i,j+r-1}$ and the types of the determinants used.
 \bt \label{1shgl} For $\gl \in \cH_\gc$, we have 
\be \label{1shtpa} \theta_{\gamma}(\gl)
= {\sum_{j=1}^{m-r+1} \sum_{i= 1}^s
(-1)^{i+1}{\stackrel{\longrightarrow}{{\rm det}}}(_jA^+(\gl, r))
\;{\stackrel{\longrightarrow}{{\rm det}}}( _{s+1-i}A^-(\gl, s))}e_{m+i,j+r-1}.\quad
\quad \ee
\be \label{1shtpb}= {\sum_{j=1}^{m-r+1} \sum_{i= 1}^s e_{m+i,j+r-1}\;(-1)^{i+1}{\stackrel{\longleftarrow}{{\rm det}}}(_jB^+(\gl, r))\;{\stackrel{\longleftarrow}{{\rm det}}}( _{s+1-i}B^-(\gl, s))}.\quad \quad \ee
\et
\bpf We prove (\ref{1shtpa}). The proof of  (\ref{1shtpb}) is similar. 
 For the isotropic simple root $\gb = \gep_m-\gd_1$, (\ref{1shtpa}) reduces to $\gth_\gb(\gl)  = e_{m+1, m}.$
Note that the overall sign in \eqref{1shtpa} is determined by the condition that the coefficient of $e_{-\pi^0}$ in $\theta_{\gamma}(\gl)$ is equal to 1, and that this term arises when the last rows of  $A^-(\gl, s)$ and 
$A^+(\gl, r))$ are deleted before taking determinants.
Suppose that $\ga = \gd_{s} -\gd_{s+1}, \gc =
\gep_r - \gd_{s}$ and $\gc' = \gep_r - \gd_{s+1} = s_\ga \gc.$
Assuming the result for $\gc$ we prove it for $\gc'$.  The result
for $\gep_{r -1} - \gd_{s}$ can be deduced in a similar way. Set
$e_{-\ga}
= e_{m+s+1,m+s}.$ For $1 \leq i \leq s-1$ we have 
\[ \gt_{i,s+1} =  \gd_i - \gd_{s+1} = s_\ga \gt_{i,s}.\]  
Consider the matrix
 $A^-(\gl',s+1)$  below, which replaces the matrix $A^-(\gl, s)$ in
Equation (\ref{1shtpa}) in the analogous expression for
$\theta_{\gc'}(\gl')$. 
$$
\left[ {\begin{array}{ccccc}
e_{m+s+1,m+s}&e_{m+s+1,m+s-1} & \hdots & e_{m+s+1,m+2}  & e_{m+s+1,m+1} \\
 (\gl' + \gr,\gt_{s,s+1}^\vee ) &e_{m+s,m+s-1} & \hdots & e_{m+s,m+2} &  e_{m+s,m+1} \\
0& (\gl' + \gr,\gt_{s-1,s+1}^\vee ) & \hdots & e_{m+s-1,m+2} &   e_{m+s-1,m+1} \\
 \vdots  & \vdots & \ddots & \vdots& \vdots \\ 
0&0&  \hdots & (\gl' + \gr,\gt_{2,s+1}^\vee ) & e_{m+2,m+2} \\
0& 0& \hdots & 0 & (\gl' + \gr,\gt_{1,s+1}^\vee ) 
\end{array}}
 \right].
$$
Suppose that $(\gl + \gr,\ga^\vee )= p$ and
let $\gl' = s_\ga\cdot\gl.$  Then
\[ (\gl' + \gr,\gt_{i,s+1}^\vee ) = (\gl + \gr,\gt_{i,s}^\vee ),\] for $1 \leq i \leq s-1$
, and this means that the last $s-1$ subdiagonal entries of
$A^-(\gl',s+1)$ and $A^-(\gl, s)$ are equal.  Also
the entry in second row first column  of
$A^-(\gl',s+1)$ is equal to $-p.$
If we remove the first row from $A^-(\gl',s+1)$  the first column of the resulting matrix will have only one non-zero entry $-p$.
If in addition we remove this column, we obtain the matrix $A^-(\gl, s)$. Therefore 
 \be \label{1stra2} -p\;{\stackrel{\longrightarrow}{{\rm det}}}(_1A^{-}(\gl,s)) = 
{{\stackrel{\longrightarrow}{{\rm det}}}}( {}_1A^-(\gl',s+1)).\ee
Similarly by removing the second row from $A^-(\gl',s+1)$
and noting that $e_{m+s+1,m+s} $ commutes with all entries in $
{}_1A^{-}(\gl,s)$,
we see that
\be \label{1stra1}
e_{m+s+1,m+s} {\stackrel{\longrightarrow}{{\rm det}}}(
{}_1A^{-}(\gl,s))=
{\stackrel{\longrightarrow}{{\rm det}}}(
{}_1A^{-}(\gl,s)) e_{m+s+1,m+s} \quad  =
{\stackrel{\longrightarrow}{{\rm det}}}( _2A^-(\gl',s+1)). \ee
Equation (\ref{121a})  in this situation takes the form
 \begin{equation} \label{122nd}
e^{p + 1}_{- \alpha}\theta_\gc(\gl) = \theta_{\gc'}(\gl') e^p_{-
\alpha}.
\end{equation}
If   $r \leq  k \leq m +s -1$ we have
 \begin{equation} \label{123nd}
e^{p + 1}_{- \alpha} e_{m+s,k} = (pe_{m+s+1,k} +
 e_{m+s+1,m+s} e_{m+s,k})e^p_{- \alpha}.
\end{equation}
We now consider  two cases: in the first entries
in ${\stackrel{\longrightarrow}{{\rm
det}}}( _{s+1-i}A^-(\gl, s))$ are
replaced by entries in
${\stackrel{\longrightarrow}{{\rm det}}}(
_{\{1,s+2-i\}}A^-(\gl',s+1)_1)$.
Suppose $1 \leq  i \leq  s-1$, and $r \leq  k \leq m$.
 Then $
e_{- \alpha}$ commutes with $e_{m+i,k}$ and all entries in the
matrix $_{s+1-i}A^-(\gl, s)$ except for those in the 
first row. Replacing
$e_{m+s,k}$ in the matrix
$_{s+1-i}A^-(\gl, s)$ by $e_{m+s+1,k}$
yields the matrix  $_{\{1,s+2-i\}}A^-(\gl',s+1)_1.$ Hence Equation
(\ref{123nd}) gives the first equality below.  For the second we
use cofactor expansion down the first column
\begin{eqnarray} \label{124nd}
e^{p + 1}_{- \alpha}\;{\stackrel{\longrightarrow}{{\rm
det}}}( _{s+1-i}A^-(\gl, s))e_{m+i,k} &=&
[p\;{\stackrel{\longrightarrow}{{\rm det}}}(_{\{1,s+2-i\}}A^-(\gl',s+1)_1)\nonumber \\
& +& e_{m+s+1,m+s}{\stackrel{\longrightarrow}{{\rm det}}}(_{s+1-i}A^-(\gl, s)) ]e_{m+i,k}e^p_{- \alpha}\nonumber\\
&=& \;{\stackrel{\longrightarrow}{{\rm det}}}(
_{s+2-i}A^-(\gl',s+1))e_{m+i,k} e^p_{- \alpha}.
\end{eqnarray}
In the second case 
entries in
${\stackrel{\longrightarrow}{{\rm
det}}}( _1A^-(\gl, s))$
 are unchanged but the factor $e_{m+s,k}$ is replaced.
If $r \leq  k \leq m$, then all
entries in the matrix $_1A^-(\gl, s),$ commute with $e_{m+s,k}$ and
$e_{- \alpha}, $ so by Equation (\ref{123nd}) we get the first
equality below, and the  second equality comes from Equations
(\ref{1stra2}) and (\ref{1stra1})
\be\label{12576}
e^{p + 1}_{- \alpha} \;{\stackrel{\longrightarrow}{{\rm
det}}}( _1A^-(\gl, s))e_{m+s,k}e^{-p}_{- \alpha}=
\;{\stackrel{\longrightarrow}{{\rm det}}}( {}_1A^{-}(\gl,s))
[pe_{m+s+1,k} +
 e_{m+s+1,m+s}e_{m+s,k}] \nonumber \ee
\be\quad\quad\quad\quad\quad\quad\quad\quad\quad = -{\stackrel{\longrightarrow}{{\rm det}}}( _{1}A^-(\gl',s+1))e_{m+s+1,k}  + {\stackrel{\longrightarrow}{{\rm det}}}(_{2} A^-(\gl',s+1))e_{m+s,k}.\ee
 Since ${\stackrel{\longrightarrow}{{\rm det}}}(_jA^+(\gl, r))$ commutes with $e_{- \alpha}$ and $e_{m+j,i}$ for all $i,j$, 
it follows from the induction assumption and Equations (\ref{122nd}), (\ref{124nd}) and (\ref{12576}) that 
\be \label{1shtpab} \theta_{\gc'}(\gl')= \sum_{j=1}^{m-r+1} \sum_{i= 1}^{s+1} 
(-1)^{i+1}{\stackrel{\longrightarrow}{{\rm det}}}(_jA^+(\gl, r)) \;
{\stackrel{\longrightarrow}{{\rm det}}}( _{s+2-i}A^-(\gl',s+1))e_{m+i,j+r-1},\quad \quad \nonumber 
\ee as desired.  \epf


\noi By inserting odd root vectors as an extra column in either
$A^+(\gl, r)$  or $A^-(\gl, s)$, we obtain a variation of Equation \eqref{1shtpa} 
 where the odd root vectors are inserted into one of the determinants. Note that the two determinants in \eqref{1shtpa} commute.  We give the details only for $A^+(\gl, r)$.
For $i\in [s]$, let $C_{(i)}(\gl, r)$ be the matrix obtained from $A^+(\gl, r)$ by adjoining the vector
\be\label{hug} (e_{m+i,r}, e_{m+i,r+1}, \ldots, e_{m+i,m} )^\trp\ee
as the last column.
 \bt \label{1thgl}With the above notation, suppose that $\gc= \gep_r
- \gd_{s}$ and $\gl\in\cH_\gc$.  Then  we have
 \be \label{1thtpa} \theta_{\gamma}(\gl)
=
\sum_{i= 1}^s
(-1)^{i+1}{\stackrel{\longrightarrow}{{\rm det}}}(_{s+1-i}
A^-(\gl, s))\;{\stackrel{\longrightarrow}{{\rm det}}}(C_{(i)}(\gl, r)).
\quad
\quad \ee
\et
\bpf This follows by cofactor expansion of
${\stackrel{\longrightarrow}{{\rm det}}}(C_{(i)}(\gl, r))$ down the last column. If $r=m$, ${\stackrel{\longrightarrow}{{\rm det}}}(C_{(i)}(\gl, r))$ should be interpreted as $e_{m+i,m}$.
\epf
\br \label{did} {\rm For convenience we record the facts.
\bi \itema If 
$A''$ is the matrix obtained from $A^-(\gl, s)$ by deleting the row and column containing $b_k$, then there are no entries in $A''$ 
of the form $e_{m+k,*}$ or $e_{*,m+k}$.
\itemb If 
 $C_{(i)}''$ is the matrix obtained from  $C_{(i)}(\gl, r)$ by deleting the row and column containing $a_\ell$, 
then there are no entries in $C_{(i)}''$ of the form $e_{\ell,*}$ or $e_{*,\ell}$.
\ei}\er
\subsection{Lie Algebras.}
Let $\fg= \fgl(m)$, and $\ga = \gep_r -\gep_\ell $.  We give a determinantal formula for the  \v{S}apovalov element $\gth_{\ga,1}$. Set $\gs_{i,j} = \gep_{i} - \gep_{j}$. Consider the
following matrices with entries in $U(\fn^-)$.

$$
C_\gl =  \left[ {\begin{array}{ccccc}
e_{r+1,r} &  e_{r+2,r} & \hdots & e_{\ell,r}\\
-a_1 & e_{r+2,r+1} & \hdots & e_{\ell,r+1} \\
0&-a_2 & \hdots & e_{\ell,r+2} \\
 \vdots  & \vdots & \ddots & \vdots              \\
0&  \hdots & -a_{j-1}&
e_{\ell,\ell-1}   
\end{array}}
 \right],\; D_\gl=  \left[ {\begin{array}{cccc}
e_{\ell,\ell-1}&e_{\ell,\ell-2} & \hdots   & e_{\ell,r} \\
 1-a_{1} &e_{\ell -1,\ell-2} & \hdots & e_{\ell -1,r} \\
0& 1-a_{2} & \hdots &    e_{\ell -2,r} \\
 \vdots  & \vdots & \ddots & \vdots \\ 
0&0&  1-a_{j-1} & e_{r+1,r} 
\end{array}}
 \right].
$$ 
where $j=\ell-r$ and $a_k=(\gl + \gr,\gs_{r,r+k}^\vee ) $ for $k\in [j-1]$.

\bt \label{1shgl2}
The \v{S}apovalov element $\gth_{\ga,1}$ is given by
\[ \gth_{\ga,1}(\gl) =
{\stackrel{\longrightarrow}{{\rm det}}}\;C_\gl = {\stackrel{\longrightarrow}{{\rm det}}}\;D_\gl, \]
for $\gl$ such that $ (\gl + \gr,\gs_{r,\ell}^\vee )=1$.\et
\bpf Similar to the proof of Theorem \ref{1shgl}. \epf
\noi Note that if $ (\gl + \gr,\gs_{r,\ell}^\vee )=1$ and $ (\gl + \gr,\gs_{r,s}^\vee )=0$, then $ (\gl + \gr,\gs_{s,\ell}^\vee )=1$ and $ (\gl + \gr-
\gs_{s,\ell},\gs_{r,s}^\vee )=1$.
\bc \label{1sapc}  In the above situation, for any highest weight vector $v_\gl$ of weight $\gl$, we have
$\gth_{\gep_{r,1} -\gep_{\ell,1}}v_\gl = 
\gth_{\gep_{r,1} -\gep_{s,1}} \gth_{\gep_{s,1} -\gep_{\ell,1}}v_\gl$.
\ec
\bpf Under the given hypothesis the matrix $C_\gl$ in Theorem \ref{1shgl2}  is  block upper triangular.\epf
\noi
\bc \label{11.5} For $p\ge 1$ and $\gl\in \cH_{\ga,p}$, we have
\[\gth_{\ga,p}(\gl) = {\stackrel{\longrightarrow}{{\rm det}}}\;C_{\gl-(p-1)\ga}\ldots {\stackrel{\longrightarrow}{{\rm det}}}\;C_{\gl-\ga}\;{\stackrel{\longrightarrow}{{\rm det}}}\;C_{\gl}.\]
\ec
\bpf Combine Theorems \ref{1calu} and \ref{1shgl2}. \epf
\noi  The above result may be viewed as a version of \cite{CL} Theorem 2.7. 

\subsection{Expansions of the Determinantal Formulas.}\label{expa}
Returning to $\fg= \fgl(m,n),$ we obtain two expansions of the determinantal formula from \eqref{1thtpa}. 
In the first case we assume that $\gc, \gc' \in B(\gl)$, where $\gc' = s_\alpha\gc = \gc -\ga$.
This means that $(\gl + \gr,\ga^\vee ) =0$ and by  Equation  \eqref{pin} we have, in this situation $\theta_{\gc}v_{\gl}= \theta_{\ga,1} \theta_{\gc'} v_{\gl}$.
Theorem \ref{stc} below could be considered as a refinement of this equation when $\gl$ is replaced by
${{\widetilde{\lambda}}}=\gl+T \xi$.
It will be used to study the modules ${M^X}({\lambda})$ in Subsection \ref{pea}. The idea is to express the extra term that arises as a \v Sapovalov element for a proper subalgebra of $\fg$. There is a second case, which  will be used in a similar way in Subsection \ref{pie}.
\subsubsection{Hessenberg Matrices.}
A matrix $B=(b_{ij})$ is {\it  $($upper$)$ Hessenberg} if $b_{ij}=0$ unless $i\le j+1$. 
If $B$ is $n\ti n$ Hessenberg, we say that  $B=(b_{ij})$ is {\it   Hessenberg of order $n$}.
 The determinant of such a matrix is a signed sum of $2^{n-1}$ (suitably ordered) terms  $ \prod_{i=1}^n b_{\nu(i),i}$ for certain permutations $\nu$.  
\bl \label{hes}
Suppose that $B$ is Hessenberg of order $n$.
\bi \itema For a fixed $q\in [n-1],$ let $\ttT = b_{q+1q}.$ Then 

\be \label{sam}{\stackrel{\longrightarrow}{{\rm det}}} B = 
-\ttT  {\stackrel{\longrightarrow}{{\rm det}}}B'' + {\stackrel{\longrightarrow}{{\rm det}}}B',\ee
where  $B'$ and  $B''$ are obtained from  $B$ by setting $\ttT =0,$ 
and by deleting the row and column containing $\ttT$ respectively. 
\itemb The matrix $B'$ is block upper triangular, with two diagonal
blocks which are upper Hessenberg of order $q$ and $n-q$. 
\itemc The matrix 
$B''$ is upper Hessenberg of order $n-1$. Also any term in  the expression \eqref{lr} for ${\stackrel{\longrightarrow }{{\rm det}}}(B)$ which contains a factor  of the form $b_{iq}$ or $b_{q+1j}$ 
 cannot occur in ${\stackrel{\longrightarrow}{{\rm det}}} B''$.
\itemd
If $B$ is 
 $(n+1)\ti n$  Hessenberg with $b_{q+1,q}=0$, and $i\in [q]$, then 
the submatrix  obtained from $B$ by deleting row $i$ is singular $($both determinants are zero$)$. 
\ei
\el 
\bpf Part (a) follows by separating the products in \eqref{lr} that contain $\ttT$ from those that do not. Note that $\ttT$ commutes with all entries in $B$, and that the order of all other factors of the products is unchanged.  The rest is easy. \epf \noi 
\subsubsection{Two Factorizations.}
At first we do not make any assumptions on $\gl\in \fh^*$. However as in Subsection \ref{nir}, we need to evaluate \v Sapovalov elements at arbitrary points, $\gl \in \fh^*$.  To do this we take Equation   \eqref{1thtpa} as the definition. 
 We write $C_{(i)}'$ and $C_{(i)}''$ for  the matrices obtained from  $C_{(i)}(\gl, r)$ by setting $\ttT =0,$ 
and by deleting the row and column containing $\ttT$ respectively. Similarly 
we write $A'$ and $A''$ for  the matrices obtained from  $A=A^-(\gl, s)$ 
by setting $\ttS=0$, and by deleting the row and column containing $\ttS$ respectively.  
The same notation is used for submatrices of 
$C_{(i)}$ and $A$.
Because $-\ttT$ is on the subdiagonal, we have by Lemma  \ref{hes} 

\by \label{2thtpa} {\stackrel{\longrightarrow}{{\rm det}}}C_{(i)}
&=&\ttT \;{\stackrel{\longrightarrow}{{\rm det}}}C_{(i)}''
+ {\stackrel{\longrightarrow}{{\rm det}}}C_{(i)}'.
\ey 

\noi For $1\le i\le s$, let  $E_{(i)}= {}_{[j]}{C_{(i)}}{}_{[j]}$ be the matrix 
obtained from $C_{(i)}$ by deleting the first $j=\ell-r$ rows and columns. 
Also let $F=\; _{[s-k]}A^-(\gl, s)_{[s-k]} $. 
For an even root $\ga$, we abbreviate $\gth_{\ga,1}$ to $\gth_{\ga}$. 
\bl \label{haha}
\bi 
\itema
If $k+1\le i\le s$, then $_{s+1-i}A'$ is singular. 
\itemb
If $i\in [k]$, the matrix $_{s+1-i}A'$ is block upper triangular with upper triangular block having LR determinant $\gth_{\ga_2}$
 and lower triangular block $_{k+1-i}F$. Hence
\be \label{hoho}{\stackrel{\lra}{{\rm det}}} _{s+1-i}A' = \gth_{\ga_2}\ti{\stackrel{\lra}{{\rm det}}}_{k+1-i}F.\ee
\itemc
For all $i\in [s]$,
 the matrix  $C_{(i)}'$ is block upper triangular with with upper triangular block having LR determinant $\gth_{\ga_1}$
and  lower triangular block $E_{(i)}$. Hence 
\be\label{hihi}{\stackrel{\longrightarrow}{{\rm det}}}C_{(i)}'=\gth_{\ga_1}\ti{\stackrel{\longrightarrow}{{\rm det}}}E_{(i)}.\ee
\ei \el
\bpf Part (a) follows  from Lemma \ref{hes} (d). By Lemma \ref{hes} (b) $_{s+1-i}A'$ is block upper triangular, and using Theorem \ref{1shgl2}, the upper triangular block has determinant equal to $\gth_{\ga_2}$. This gives (b) and the proof of (c) is similar.\epf
\noi We apply Lemma \ref{hes} in the case that $B=  {}_{s+1-i}A$.

\bl \label{sky}
\bi 
\itema  If $i\leq k-1$, then 
\be \label{cue}{\stackrel{\lra}{{\rm det}}}  _{s+1-i}A = 
-\ttS  {\stackrel{\lra}{{\rm det}}} _{s-i}A'' + {\stackrel{\lra}{{\rm det}}} _{s+1-i}A'.\ee
\itemb  If $i\geq k+1$, then  
\be \label{men}{\stackrel{\lra}{{\rm det}}}  _{s+1-i}A = 
\ttS  {\stackrel{\lra}{{\rm det}}} _{s+1-i}A''.\ee 
\itemc 
\be \label{kip}{\stackrel{\lra}{{\rm det}}}  _{s+1-k}A = 
{\stackrel{\lra}{{\rm det}}} _{s+1-k}A'.\ee 
\ei\el
\bpf Note that if $p<q$, then row $q$ of $A$ is row $q-1$ of $_p A$.  
If $i\leq k-1$ it follows that $_{s-i}(A'')$ is obtained from $_{s+1-i}A$ by deleting the row and column containing $\ttS$, 
that is $_{s-i}(A'') = (_{s+1-i}A)''$.  
Similarly if $i\geq k+1$ then $_{s+1-i}(A'') = (_{s+1-i}A)''$. 
In both cases we obtain $_{s+1-i}(A') = (_{s+1-i}A)'  $, but by Lemma \ref{haha} (a) $_{s+1-i}A'$ is singular if $i\geq k+1$.
The different signs in \eqref{cue} and \eqref{men} arise since  $\ttS$ is on the subdiagonal of $_{s+1-i}A$ in case (a), and on the diagonal of $_{s+1-i}A$ 
in case (b).
\epf

\subsubsection{Some General Linear Subalgebras.}
To state the results, we require a bit more notation.  For the rest of this Subsection, the element $\gr$ defined for $\fgl(m,n)$ in Equation \eqref{rde} will be denoted by $\gr_{m,n}.$
In addition let $(\;,\;)_{m,n}$ be the bilinear form on $\fh^*$ defined by 
\[(\gep_i,\gep_j)_{m,n} = - (\gd_i,\gd_j)_{m,n} = \gd_{i,j}\]
for all relevant indices $i, j$. 
Set $({\lambda} + \gr,\ga^\vee ) =T$.
The extra term in \eqref{nice} below comes from a  element for a subalgebra of $\fg$ isomorphic to $\fgl(m-1,n)$.
Suppose that $V=\ttk^{m|n}$ is a super vector space of dimension $(m|n)$ and identify
$\fgl(m,n)$ and $\fgl(V)$ by means of the standard basis $e_1,\ldots,e_{m+n}.$ Do the same for
$\fgl(m-1,n)=\fgl(\ttk^{m-1|n})$ using the standard basis $\overline{e}_1,\ldots,\overline{e}_{m+n-1}$ for $\ttk^{m-1|n}$. Fix $k, \ell $ with $1\le \ell \le m$, and $1\le k\le n$. Then let 
${\overline\psi}^\ell :\overline{\fg} = \fgl(m-1,n)\lra \fg=\fgl(m, n)$ 
be the embedding induced by the map $V=\ttk^{m-1|n}\lra \ttk^{m|n}$ 
sending $\overline{e}_i$ to $e_i$ for $i<\ell$ and   $\overline{e}_i$ to $e_{i+1}$ for $i\ge \ell$ 
(so that $e_\ell$ is not in the image of this map).
\\ \\
The embeddings of general linear subalgebras defined in the above paragraph are all we need for Theorem \ref{stc}, but 
for Theorem \ref{nex} we need two variations. 
Define maps $\ttk^{m|n-1}\lra \ttk^{m|n}$, (resp. $\ttk^{m-1|n-1}\lra \ttk^{m|n}$) of superspaces using basis elements in a similar way, with  $e_{m+k}$ not in the image ( resp. $e_\ell$ and $e_{m+k}$  not in the image). Then let $\underline{\psi}_k:\underline{\fg} =\fgl(m,n-1)\lra \fg$ and
${\underline{\overline{\psi}}_k^\ell }:
\underline{\overline{\fg}} = \fgl(m-1,n-1)\lra \fg$
 be the embeddings of Lie superalgebras induced by these maps. Finally, let ${\fh}, {\overline{\fh}}$ and $\underline{\overline{\fh}}$ be the diagonal Cartan subalgebras of ${\fg}, {\overline{\fg}}$ and $\underline{\overline{\fg}}$,  and let \[{\overline\phi^\ell }:\fh^*\lra
{\overline{\fh}^*}, \quad {\underline\phi_k}:\fh^*\lra
{\underline{\fh}^*} \;
\mbox{ and }\; \underline{\overline{\phi}}_k^\ell :\fh^*\lra \underline{\overline{\fh}}^*\]
 be the maps dual to the restriction of ${\overline\psi}^\ell, \underline{\psi}_k$
and
$\underline{\overline{\psi}}_k^\ell $ to ${\overline{\fh}}, \underline{\fh}$ and $\underline{\overline{\fh}}$
respectively.
\\ \\
We fix $\ell \in [m]$ and $k \in [n]$, and then use the shorthand 
\[ \overline{\psi} =\overline{\psi}^\ell,\quad    \underline{\psi} = \underline{{\psi}}_k,\quad  \underline{\overline{\psi}} =\underline{\overline{\psi}}_k^\ell, \quad  \overline{\phi} =\overline{\phi}^\ell,\quad    \underline{\phi} = \underline{{\phi}}_k,\quad  \underline{\overline{\phi}} =\underline{\overline{\phi}}_k^\ell.\]
For $\ga\in \fh^*$ we define
\[\overline{\ga} =\overline{\phi}(\ga),\quad    \underline{\ga} = \underline{{\phi}}(\ga),\quad  \underline{\overline{\ga}} =\underline{\overline{\phi}}(\ga).\]
Observe that  $\;{\overline\phi}({\gep_{\ell}}) =0$ and the restriction of ${\overline\phi}$ to 
$\fh'=\span\{ \gep_p,\gd_q|p\neq \ell\}$ 
is an isomorphism onto ${\overline\fh}^*.$  
Given $\gl\in \fh^*$, we define ${\lambda}_1\in {\overline{\fh}}^*$, by  
\be \label{flx} ({\lambda}_1+\gr_{m-1,n},  {\overline\gb})_{m-1,n} = ({\lambda}+\gr_{m,n},\gb)_{m,n}\ee
 for $\gb \in \fh'$.  
 Similarly we define  ${\lambda}_2\in \underline{{\fh}}^*, \; {\lambda}_3\in \underline{\overline{\fh}}^*$ by  
\be \label{fax} ({\lambda}_2+\gr_{m,n-1},{\underline\gb})_{m,n-1} = ({\lambda}+\gr_{m,n},\gb)_{m,n},\ee 
\be \label{fox} ({\lambda}_3+\gr_{m-1,n-1},\underline{\overline\gb})_{m-1,n-1} = ({\lambda}+\gr_{m,n},\gb)_{m,n}.\ee
In \eqref{fax} and \eqref{fox} we assume that 
$$\gb \in \span\{ \gep_p,\gd_q|q\neq k\} \mbox{ and } \gb \in \span\{ \gep_p,\gd_q|p\neq \ell, q\neq k\}$$ respectively.
\subsubsection{Pullbacks.}
From now on $A= A^-(\gl, s) $. Let $A''= {}_{s+1-k}A_{s-k} $.  This is the matrix obtained from $A$ by deleting the row and column containing $\ttS$. 
By Remark \ref{did}, all entries from $\fg$ in $A''$ are in the image of the map
$\;{\underline\psi}_{k}:\underline{\fg} = \fgl(m,n-1)\lra \fg=\fgl(m, n)$. 
This allows us to define a matrix $A[k]$ with 
entries in $U({\underline{\fg}})$ called the {\it pullback of $A''$ under ${\underline{{\psi}}}_k$}
in the following way.  
\bi \itemi The matrix 
$A[k]$ is Hessenberg of the same size as $A''$.
\itemii All subdiagonal entries of $A[k]$ and $A''$ are equal.
\itemiii If $j \ge i$, then the entry in row $i$ and column $j$ of $A[k]$ is the element
of ${\underline{\fg}}$ which is mapped by ${\underline\psi}$ to the entry in row $i$ and column $j$ of $A''$.
\ei
We sometimes write $\underline{\psi}(A[k])= A''$ in this situation.  We have $\underline{\psi}(_{s-i}A[k]) =_{s-i}A''.$
 Note that this is consistent with 
\eqref{fax}.  In fact the latter just says that (ii) holds. The point of Equations \eqref{flx}-\eqref{fox} is that they 
make it clear that $\gc\in B(\gl)$ implies 
$\overline{\gc} \in B(\lambda_1)$ etc, which we need for the applications. 
\\ \\
Next we list the other pullbacks we require.
Then we explain why these pullbacks are defined. 
Let 
\be \label{lax} 
f(i) =
\left\{ \begin{array}
  {ccc} i & \mbox{ if } & i\le k-1\\
i+1 & \mbox{ if }  & i\ge k.
\end{array} \right.\ee 
We define 
\be \label{sax} 
C_{(i)}[k] \mbox{ (resp. } E_{(i)}[k]) \mbox{  to be the pullback  } \mbox{ of } C_{(f(i))}  \mbox{ (resp. } E_{(f(i))})  \mbox{ under  } {\underline{{\psi}}} 
,\ee
\be \label{six} F[{\ell}] \mbox{ (resp. } C_{(i)}[{\ell}]) \mbox{ to be the pullback of } 
F \mbox{ (resp. } C_{(i)}'') 
\mbox{ under } \overline{\psi},\ee
\be \label{sox}
A''[k,{\ell}] \mbox{ (resp. } C_{(i)}[k,{\ell}])  \mbox{ to be the pullback of } 
A'' \mbox{ (resp. } C''_{(f(i))})
\mbox{ under } \underline{\overline{\psi}}. \ee
The pullbacks of $C_{(i)}''$ and both matrices in \eqref{sox} exist by Remark \ref{did}.
For \eqref{sax} note that all entries in $A^+(\gl, r)$ are elements of $ \fgl(m)\op 0 \subseteq \fg_0$.  
Also if $i\neq k$, the vector in \eqref{hug} contains no element of the form $e_{m+k,*}$ or $e_{*,m+k}$, 
so the same is true for $C_{(i)}$, and $E_{(i)}, i\neq k$.  Note that $C'_{(k)}$ is not considered in \eqref{sax}. 
The argument for  $F$ is similar but easier, all entries in $F={}_{[s-k]}A^-(\gl, s)_{[s-k]} $ 
are elements of  $0\op \fgl(n) \subseteq \fg_0$. Thus $F$ contains no entries of the form $e_{\ell,*}$ or $e_{*,\ell}$.

\subsubsection{\v Sapovalov elements.} \label{dry}
We can use Theorem \ref{1thgl} to  express the {\v Sapovalov elements} for

\be   \gc_1, \quad		{\gb_1} = {\underline{\gc_1+\ga_2}}, \quad {\gb_2} = {\overline{\ga_1+\gc_1}}, \quad 
  {\underline{\overline{\gc}}} = {{\underline{\overline\phi} (\gc)}}.
\ee  
in terms of the matrices introduced above. Note that these are roots for different Lie superalgebras. 
\bl Given $\gl\in \fh^*$, define $\gl_1-\gl_3$ by \eqref{flx}-\eqref{fox}. Then the \v Sapovalov elements for the roots 
$ \gc_1, \gb_1,  \gb_2$ and 
${\underline{\overline{\gc}}}$ are given by 
\be \label{moe} \theta_{\gc_1}(\gl) =\sum_{i= 1}^k
(-1)^{i+1}{\stackrel{\lra}{{\rm det}}}_{k+1-i}F\ti 
{\stackrel{\longrightarrow}{{\rm det}}}E_{(i)},
 \ee 
\be\label{mad} \gth_{\gb_1}({{\lambda}_1})=
\sum_{i= 1}^{s-1} (-1)^{i}{\stackrel{\lra}{{\rm det}}}_{s-i}A[k] \ti 
{\stackrel{\longrightarrow}{{\rm det}}}\;E_{(i)}[k],\ee

\be\label{mud}\gth_{{\gb_2}}({{\lambda}_2})=\sum_{i= 1}^k (-1)^{i+1}{\stackrel{\lra}{{\rm det}}}_{k+1-i}F[\ell] 
\ti{\stackrel{\longrightarrow}{{\rm det}}}{C}''_{(i)}[\ell],\ee

\be \label{mix}\gth_{\underline{{\overline{\gc}}}}({{\lambda}_3})= 
\sum_{i= 1}^{s-1}(-1)^{i+1}\;{\stackrel{\lra}{{\rm det}}}_{s-i}A''[k,{\ell}]\ti 
 {\stackrel{\longrightarrow}{{\rm det}}}C_{(i)}''[k,{\ell}].
\ee
\el
\bpf This follows from Theorem \ref{1thgl} with $F$ 
(resp. $A[k]$,  $F[\ell], A''[k,{\ell}]$) playing the role of $A^-(\gl, s)$,  
and $E_{(i)},$ (resp. $E_{(i)}[k],$ ${C}''_{(i)}[\ell], C_{(i)}''[k,{\ell}])$ playing the role of 
$C_{(i)}(\gl, r)$.
\epf 
\subsubsection{Determinantal Expansions: Main Results.} \label{det}
Assume that $\gc = \gep_r- \gd_{s},$ $\ga = \gep_r-  \gep_{\ell }$  and
$\gc'=s_\ga\gc=\gep_\ell- \gd_{s},$ where $r<\ell=r+j\le m.$
\ff{ The case were $\ga$ is a root of $\fgl(n)$ can be handled similarly.} 
Note that  $\;{\overline\phi}({\gep_{\ell}}) =0$ and the restriction of ${\overline\phi}$ to $\fh' =\span\{ \gep_i,\gd_j|i\neq \ell\}$ is an isomorphism onto ${\overline\fh}.$   Define $a_i, b_i$ as in \eqref{eat}. 
\bt \label{stc} With the above notation, set $\ttT=({\lambda}+\gr,\ga^\vee)=a_j.$  Then  
if $\overline{\gc} = {\overline\phi} (\gc)$, we have 
\be\label{nice}\gth_\gc v_{{\lambda}} = [\overline{\psi}(\gth_{\overline{\gc}}(
{{\lambda_1}}))\ttT
+\gth_{\ga,1}\gth_{\gc'}]v_{{\lambda}}.
\ee\et
 \noi We remark that in the definition of $A^+({{\lambda} }, r)$, and hence also in $C_{(i)}({{\lambda} }, r)$ we have 
We give the proof after proving Theorem \ref{nex}.
\\ \\
\noi Next we suppose $\ga_1, \ga_2$ (resp. $\gc_1, \gc_2$) are distinct positive even (resp. odd) roots	
such that $\ga_1 +\gc_1 +\ga_2 \;=\gc_2$. 
Then set $\;\gc=\gc_2$. 
Without loss of generality we may assume that $\ga_1 = \gep_r- \gep_{\ell }$, $\ga_2 = \gd_k- \gd_{s},$ 
$\gc_1 =\gep_{\ell }- \gd_k$, where $r<\ell=r+j\le m.$ 
In the result below, we treat $\ttT= a_j,$ and $ \ttS= b_k$ as indeterminates, and obtain expansions of Equation \eqref{1thtpa}.
\bt \label{nex} We have 
\be\label{mice}
\gth_{\gc_2}v_{{\lambda}} = 
[\gth_{\ga_2}\gth_{\ga_1}\gth_{\gc_1}
-\gth_{{{\ga_1}}}\overline{\psi}(\gth_{{\gb_1}}({{\lambda_1}}))\ttS +
\gth_{{\ga_2}}\underline{\psi}(\gth_{{\gb_2}}({{\lambda_2}}))\ttT
-\underline{\overline\psi}(\gth_{\underline{{\overline{\gc}}}}({{\lambda_3}}))\ttS\ttT]
v_{{\lambda}}.\ee
\et
\br \label{dan}
{\rm We need to know that certain products in $U(\fn^-)$ commute, and this can be shown based on a consideration of the
weights of  their factors.  
Write $\gep_{m+i} =\gd_i$ for $i\in [n]$. Then we say that a product $u=u_1 \ldots u_t$ has {\it top weight}
$\gep_p -\gep_q,$ where $1\le p< q\le m+n,$  if for $i\in [t]$, $u_i$ has weight $\gep_{p_i} - \gep_{q_i}$ (with $p_i<q_i$) and 
$\sum_{i=1}^{t}\gep_{p_i} - \gep_{q_i} =\gep_{p} - \gep_{q}$. 
Clearly $v, w$ commute if they are respectively linear combinations of elements having  top weights 
$\gep_p -\gep_q$ and  $\gep_a -\gep_b$ with $q<a$. }\er
\noi 
{\it Proof of Theorem \ref{nex}.}
We substitute \eqref{2thtpa} and the expressions for $_{s+1-i}A^-({\gl}, s)$ from Lemma \ref{sky} into \eqref{1thtpa} to obtain 
\by \label{skn} \theta_{\gamma}({\gl})
&=&\sum_{i= 1}^{k-1}(-1)^{i+1}
(-\ttS  {\stackrel{\lra}{{\rm det}}} _{s-i}A'' + {\stackrel{\lra}{{\rm det}}} _{s+1-i}A')
\ti(\ttT \;{\stackrel{\longrightarrow}{{\rm det}}}C_{(i)}''
+ {\stackrel{\longrightarrow}{{\rm det}}}C_{(i)}')\nn\\
&+&\sum_{i= k+1}^s(-1)^{i+1} \ttS  {\stackrel{\lra}{{\rm det}}} _{s+1-i}A''\ti(\ttT \;{\stackrel{\longrightarrow}{{\rm det}}}C_{(i)}''
+ {\stackrel{\longrightarrow}{{\rm det}}}C_{(i)}')\nn\\
&+&(-1)^{k+1} {\stackrel{\lra}{{\rm det}}} _{s+1-k}A'\ti(\ttT \;{\stackrel{\longrightarrow}{{\rm det}}}C_{(k)}''
+ {\stackrel{\longrightarrow}{{\rm det}}}C_{(k)}').\ey
To complete the proof we compute the coefficients of $\ttS\ttT, \ttS, \ttT$ and the constant term in  \eqref{skn}.
\\
\\
 \noi (i) 
The constant coefficient 
in \eqref{skn} is the first expression below. We use  Lemma \ref{haha} (a), then  \eqref{hoho}, \eqref{hihi} and then the fact that $\gth_{\ga_1}$ commutes with all entries in ${\stackrel{\lra}{{\rm det}}}_{k+1-i}F$ (Remark \ref{dan}).  The final equality holds by \eqref{moe}  
\by \label{2thtpc} 
\sum_{i= 1}^s (-1)^{i+1}{\stackrel{\lra}{{\rm det}}}_{s+1-i}A'\ti {\stackrel{\longrightarrow}{{\rm det}}}C_{(i)}
&=&
\sum_{i= 1}^k(-1)^{i+1}{\stackrel{\lra}{{\rm det}}}_{s+1-i}A'\ti {\stackrel{\longrightarrow}{{\rm det}}}C_{(i)}'\nn\\
&=&
\sum_{i= 1}^k (-1)^{i+1} \gth_{\ga_2}\ti{\stackrel{\lra}{{\rm det}}}_{k+1-i}F\ti
\gth_{\ga_1}\ti{\stackrel{\longrightarrow}{{\rm det}}}E_{(i)}\nn\\
&=&
\sum_{i= 1}^k (-1)^{i+1}\gth_{\ga_2}\gth_{\ga_1}\ti{\stackrel{\lra}{{\rm det}}}_{k+1-i}F\ti
{\stackrel{\longrightarrow}{{\rm det}}}E_{(i)}\nn\\
&=&  \gth_{\ga_2}\gth_{\ga_1}\gth_{\gc_1}v_{{\lambda}}.
\ey
\\ \\
\noi (ii) 
Using first (a pullback of) \eqref{hihi}, then the fact that $\gth_{{\underline\ga_1}}$ commutes with all entries in $A[k]$, and then \eqref{mad} we have
\by \label{4thtpa} 
\sum_{i= 1}^{s-1}(-1)^{i}{\stackrel{\lra}{{\rm det}}}_{s-i}A[k]\ti {\stackrel{\longrightarrow}{{\rm det}}}C_{(i)}'[k] 
&=&
\sum_{i= 1}^{s-1}
(-1)^{i}{\stackrel{\lra}{{\rm det}}}_{s-i}A[k]\ti \gth_{{\underline\ga_1}}\ti{\stackrel{\longrightarrow}{{\rm det}}}\;E_{(i)}[k]\nn\\
&=&
\sum_{i= 1}^{s-1} (-1)^{i}\gth_{{\underline\ga_1}}\ti{\stackrel{\lra}{{\rm det}}}_{s-i}A[k] \ti {\stackrel{\longrightarrow}{{\rm det}}}\;E_{(i)}[k]\nn\\
&=&  
-\gth_{{\underline\ga_1}}\gth_{{\gb_1}}({{\lambda_1}}).
\ey   
By \eqref{sax} and the fact that $\underline{{\psi}}$ maps  ${\;_{s-i}}A[k]\;$ to 
	$\;_{s-i}A''$, 
$\underline{{\psi}}$ maps \eqref{4thtpa} to
\be
\sum_{i= 1}^{k-1}(-1)^{i}\;{\stackrel{\lra}{{\rm det}}} _{s-i}A'' \ti {\stackrel{\longrightarrow}{{\rm det}}}C_{(i)}'
+\sum_{i= k}^{s-1}(-1)^{i}\;{\stackrel{\lra}{{\rm det}}} _{s-i}A'' \ti {\stackrel{\longrightarrow}{{\rm det}}}C_{(i+1)}'\nn
\ee and this shows that the coefficient of $\ttS$ in as stated in \eqref{skn}.
\\ \\
\noi (iii) 
Next, by \eqref{hoho} and \eqref{mud} we have 

\by 
\sum_{i= 1}^{k}(-1)^{i+1}{\stackrel{\lra}{{\rm det}}}_{s+1-i}A'[\ell]\ti {\stackrel{\longrightarrow}{{\rm det}}}C''_{(i)}[\ell] 
&=&\sum_{i= 1}^k (-1)^{i+1}\gth_{\overline{\ga_2}} {\stackrel{\lra}{{\rm det}}}_{k+1-i}F[\ell] \ti{\stackrel{\longrightarrow}{{\rm det}}}{C}''_{(i)}[\ell]\nn\\
& =& \gth_{{\overline\ga_2}}\gth_{{\gb_2}}({{\lambda_2}}),\nn
\ey
and $\overline{{\psi}}$ maps this to 
\be\sum_{i= 1}^k
(-1)^{i+1}
{\stackrel{\lra}{{\rm det}}}_{s+1-i}
A' \ti {\stackrel{\longrightarrow}{{\rm det}}}C''_{(i)},\nn\\
\ee which is the coefficient of $\ttT$ in \eqref{skn}.
\\ \\
(iv) Finally, let \[A''= {}_{s+1-k}A_{s-k} = \underline{\overline{\psi}}(A[k,\ell]). \]
By \eqref{mix}, 
\be -\gth_{\underline{{\overline{\gc}}}}({{\lambda}_3})= 
\sum_{i= 1}^{s-1}(-1)^{i}\;{\stackrel{\lra}{{\rm det}}}_{s-i}A''[k,{\ell}]\ti 
 {\stackrel{\longrightarrow}{{\rm det}}}C_{(i)}''[k,{\ell}],\nn
\ee
 and $\underline{\overline{\psi}}$ maps this to

\be\sum_{i= 1}^{k-1}(-1)^{i}\;{\stackrel{\lra}{{\rm det}}}_{s-i}A'' \ti {\stackrel{\longrightarrow}{{\rm det}}}C_{(i)}''
+\sum_{i= k}^{s-1}(-1)^{i}\;{\stackrel{\lra}{{\rm det}}}_{s-i}A'' \ti {\stackrel{\longrightarrow}{{\rm det}}}C_{(i+1)}''
\nn
\ee which is the coefficient of $\ttS\ttT$ in \eqref{skn}.
\hfill  $\Box$
\\ \\
{\it Proof of Theorem \ref{stc}}
The proof is similar to the proof of Theorem \ref{nex}, but easier.
Instead of \eqref{skn} we have 

\be \label{sin} \theta_{\gamma}(\widetilde{\gl})
=
\sum_{i= 1}^s(-1)^{i+1}{\stackrel{\lra}{{\rm det}}}(_{s+1-i}A^-(\widetilde{\gl}, s))
\ti(-\ttT \;{\stackrel{\longrightarrow}{{\rm det}}}C_{(i)}''
+ {\stackrel{\longrightarrow}{{\rm det}}}C_{(i)}')\nn.\ee
The proof is completed by computing coefficients as before.
\hfill  $\Box$

\br{\rm 
Suppose that $(\gl+\gr, \alpha^\vee_1) = (\gl+\gr, \alpha^\vee_2)=1$ and assume that  
$(\xi,\gc_1)= (\xi,\gc_2)=0$, and $(\xi,\ga^\vee_1) = (\xi,\ga^\vee_2) = 1.$ 
Then set ${{\widetilde{\lambda}}}=\gl+T \xi$. 
Suppose that $\ga_1 = \gep_r- \gep_{\ell }$, $\ga_2 = \gd_k- \gd_{s},$ 
$\gc_1 =\gep_{\ell }- \gd_k$, where $r<\ell=r+j\le m.$ 
Then $(\widetilde{\lambda}+\gr, \ga_1^\vee) = (\widetilde{\lambda}+\gr, \ga_2^\vee) = T+1$.  
Hence  in the notation of Theorem \ref{nex} we have $\ttS=\ttT = T+1$. 
Thus if  $\gth_{\gc_1}v_{\widetilde{\lambda}}=\gth_{\gc_2}v_{\widetilde{\lambda}} =0$, then \eqref{mice} yields, over $B$
\be\label{mice2}  
0=[\underline{\overline\psi}(\gth_{\underline{{\overline{\gc}}}}({{{\widetilde{\lambda}}_3}}))(T+1)
+\gth_{{\ga_1}}\underline{\psi}(\gth_{{\gb_1}}({{{\widetilde{\lambda}}_2}}))
-\gth_{{{\ga_2}}}\overline{\psi}(\gth_{{\gb_2}}({{{\widetilde{\lambda}}_1}}))]
v_{\widetilde{\lambda}}.\ee
This gives a generalization of \eqref{rat2}. 
}\er

\section{{Survival of \v{S}}apovalov elements in factor modules.}\label{1surv}
 Let  $v_{\gl}$ be a highest weight vector
in a Verma module $M(\gl)$ with highest weight $\gl,$ and suppose $\gc$ is an odd root with
$(\gl+\gr, \gc) =0$. We are interested in the condition that the image of
$\theta_{\gc}v_{\gl}$ is non-zero in various factor modules of $M(\gl)$.
\subsection{Independence of \v{S}apovalov elements.}
Given $\gl \in \fh^*$ recall the set $B(\gl)$ defined in Section \ref{1s.1}, and define a ``Bruhat order" $\le$ on $B(\gl)$ by $\gc' \le \gc$ if $\gc-\gc'$ is a sum of positive even roots. Then introduce a relation $\downarrow$ on $B(\gl)$  by $\gc' \da \gc$ if $\gc' \le \gc$ and $(\gc,\gc') \neq 0.$ If $\gc \in B(\gl)$, we say that $\gc$ is $\gl$-{\it minimal} if $\gc' \da \gc$  with $\gc' \in B(\gl)$ implies that $\gc' = \gc$.
For $\gc \in B(\gl)$ set $B(\gl)^{-\gc} = B(\gl)\backslash \{\gc\}$. We say $\gc$ is {\it independent
at } $\gl$ if $$\theta_{\gc}v_{\gl}\notin \sum_{\gc' \in B(\gl)^{-\gc} }U(\fg)\theta_{\gc'} v_{\gl}.$$
\bp \label{1pot}If $\gc' \da \gc$ with $\gc' \in B(\gl)$ and $\gc' < \gc$, then $\theta_{\gc}v_{\gl}\in U(\fg)\theta_{\gc'} v_{\gl}.$ \ep
\bpf The hypothesis implies that
$(\gc,\ga^\vee)>0$
and $\gc =s_\ga \gc'$. Thus the   result follows from Theorem \ref{man}.\epf
\noi By the  Proposition,  if we are interested in the independence of the \v{S}apovalov elements $\gth_\gc$ for distinct isotropic roots, it suffices to  study only $\gl$-minimal roots $\gc$.
\\ \\
For the rest of this section we assume that $\fg = \fgl(m,n)$. We use 
Equation
(\ref{1shtpa}), and  order the positive roots of $\fg$ so that each summand in this equation is a constant multiple of $e_{-\pi}$ for some $\pi \in {\overline{\bf P}}(\gamma)$. For such $\pi$  the odd root vector is the rightmost factor of $e_{-\pi}$, that is we have $e_{-\pi}\in U(\fn^-_0)\fn^-_1$.
\bl \label{1car}If If $\gc$ is $\gl$-minimal, then $e_{-\gc} v_\gl$ occurs with non-zero coefficient in $\gth_\gc v_\gl$.
\el
\bpf Assume $\gc = \gep_r
- \gd_{s}.$ Then if $\ga = \gep_r - \gep_i$ with $r<i$,  or $\ga = \gd_j - \gd_s$ with $j<s$ we have $(\gl + \gr, \ga^\vee) \neq 0$, since $\gc$ is $\gl$-minimal.  Thus the entries on the superdiagonals of $A^+(\gl, r)$ and $A^-(\gl, s)$ are non-zero.  Thus the result follows from Theorem \ref{1shgl}.\epf
\bt \label{1boy} The isotropic root $\gc$ is independent at  $\gl$ if and only if $\gc$ is $\gl$-{minimal}. \et
\bpf Set $B=B(\gl)^{-\gc} $. If
$\gc$ is not $\gl$-{minimal}
then $\gc$ is not independent at  $\gl$
by Proposition \ref{1pot}. Suppose that $\gc$ is $\gl$-{minimal}
and
$$\theta_{\gc}v_{\gl}\in \sum_{\gc' \in B}U(\fg)\theta_{\gc'} v_{\gl} = \sum_{\gc' \in B}U(\fn^-)\theta_{\gc'} v_{\gl},$$ then by comparing weights
\be \label{1la}\theta_{\gc}v_{\gl}\in \sum_{\gc' \in B}U(\fn^-_0)e_{-\gc'} v_{\gl},\ee
But
Lemma \ref{1car} implies that
$$\theta_{\gc}v_{\gl}\equiv ce_{-\gc}v_{\gl}\mod \sum_{\gc' \in B}U(\fn^-_0)e_{-\gc'} v_{\gl}$$ for some non-zero constant $c.$  
By  (\ref{1la}) this contradicts the PBW Theorem.
\epf
\subsection{Survival of \v{S}apovalov elements in Kac modules.}
For $\fg = \fgl(m,n)$ we have $\mathfrak{g}_1 = \mathfrak{g}_1^+ \oplus
\mathfrak{g}_1^-,$ where $\mathfrak{g}_1^+ $ (resp. $\mathfrak{g}_1^-$) is the set of block  upper (resp. lower) triangular matrices.
Let $\fh$ be the Cartan subalgebra of $\fg$ consisting of diagonal matrices, and set $\mathfrak{p} =  \mathfrak{g}_0  \oplus \mathfrak{g}_1^{+}.$ Next let
 \begin{eqnarray*}
  P^+ &=&
   \{ \lambda \in \mathfrak{h}^*|(\lambda, \alpha^\vee) \in \mathbb{Z}, (\lambda, \alpha^\vee) \geq 0 \quad
 \mbox{for all}\quad \alpha \in \Delta^+_0 \}
 \end{eqnarray*}
For $\gl \in P^+,$ let $L^0(\lambda)$ be the (finite dimensional) simple  $\mathfrak{g}_0$-module with highest weight $\gl$.  Then
$L^0(\lambda)$ is naturally a $\fp$-module and we define the Kac module $K(\gl)$ by
\[K(\lambda) =  U(\mathfrak{g}) \otimes_{U(\fp)}
L^0(\lambda).\]
 Note that as a $\mathfrak{g}_0$-module
\[K(\lambda) =  \Lambda(\mathfrak{g}^-_1) \otimes
L^0(\lambda).\]
The next result is well-known.  Indeed two
methods of  proof are given in  Theorem 4.37 of \cite{Br}.  The
second of these is based on Theorem 5.5 in \cite{S2}.  We give a short
proof using  Theorem \ref{1shgl}. We now assume that the roots are ordered as in Equation (\ref{1shtpb}), that is  with the odd root vector first.
\bt \label{1shapel} If $\lambda$ and $\lambda - \epsilon_r+\delta_s$
belong to $P^+$ and $(\lambda+\rho,\epsilon_r - \delta_s) = 0$, then
$$[K(\lambda):L(\lambda-\epsilon_r+\delta_s)] \neq 0.$$ \et
 \bpf Set $\gamma = \epsilon_r -
\delta_s$.
 Let $\theta_\gamma(\lambda)$ be as in Theorem \ref{1shgl}.  Then
$w = \theta_\gamma(\lambda)v_\lambda$ is a highest weight vector in
the Verma module $M(\lambda)$ with weight $\lambda - \gamma$.  It
suffices to show that the image of $w$ in the Kac module
$K(\lambda)$ is nonzero.  We have an embedding of $\fg_0$-modules
\[ \fg^-_1 \otimes L^0(\lambda) \subseteq \Lambda \fg^-_1 \otimes L^0(\lambda). \]
The elements $e_{m+j,i+r-1}$ in Equation (\ref{1shtpb}) form part of
a basis for $\fg^-_1$. Furthermore the coefficient of
$e_{m+j,i+r-1}$ belong to $U(\fn^-_0)$. Therefore it suffices to
show that the coefficient of $e_{m+s,r}$ in this equation is
nonzero. This coefficient is found by deleting the first column  of the matrix $B^+(\gl, r)$ and the last row of
$B^-(\gl, s)$ and taking determinants of the resulting matrices, which have only zero entries above the main diagonal. We find that the coefficient of $e_{m+s,r}$ is
\[  \pm \prod_{k=1}^{m-r} (1 - (\lambda + \rho, \gs_{r,r+k}^\vee))
\prod^{s-1}_{k=1} (1 - (\lambda + \rho, \gt_{k,s}^\vee )).    \]
Since $\lambda \in P^+, (\lambda + \rho, \gs_{r,r+k}^\vee) \geq 1$
with equality if and only if $k = 1$ and $(\lambda, \epsilon_r -
\epsilon_{r+1}) = 0$.    This cannot happen if $\lambda - \gamma \in
P^+$, so the first product above is nonzero, and similarly so is the
second. \epf

\section{The Sum Formula.}\label{sf} 
\subsection{An analog of the \v Sapovalov  determinant.} \label{7.4}
Throughout this Section we assume that $\fg= \fgl(m,n).$ Although ${\bf{\overline{P}}}_X(\eta) = {\bf P}_X(\eta)$ for all $\eta,$ we continue to use the former notation because some arguments hold outside of Type A. 
The goal is to give a Jantzen sum formula for the modules ${M^X}({\lambda})$ introduced in Section
\ref{jaf}. This is done by first computing  a \v Sapovalov  determinant for these modules. 
We define an $A$-valued bilinear form on ${M^X}(\widetilde{\lambda})_A$ as in \cite{M} Corollary 8.2.11 and Equation (8.2.14). 
Let $F_\eta^X(\widetilde{\lambda})$ be the restriction of this form to the weight space 
$M^X(\widetilde{\lambda})_A^{\widetilde{\lambda}-\eta}.$ 
(Equation \eqref{dff} below can be taken as the definition of $F_\eta^X(\widetilde{\lambda}))$. 
The determinant of this form has the important property, see
\cite{J1} Lemma 5.1  or \cite{M} Lemma 10.2.1,
\be \label{yet}v_{T}(\det F^X_\eta(\widetilde{\lambda})) = \sum_{i>0} \dim M_i^{X}(\lambda)^{\lambda - \eta},\ee
where $\{\dim M_i^{X}(\lambda)\}$ is the Jantzen filtration. 
However there is a related determinant $\det G^X_\eta(\widetilde{\lambda})$ whose leading term is easier to compute.
The elements $e_{-\pi}v_{\widetilde{\lambda}}^X$ with $\pi \in {\bf{\overline{P}}}_{X}(\eta)$
belong to ${M^X}(\widetilde{\lambda})_A^{\widetilde{\lambda}  - \eta}$,
and  form a $B$-basis for ${M^X}(\widetilde{\lambda})_A^{\widetilde{\lambda}  - \eta}\ot_A B$,
but they do not in general form a basis for ${M^X}(\widetilde{\lambda})_A^{\widetilde{\lambda}  - \eta}$
as an $A$-module.
We define $G^X_\eta(\widetilde{\lambda})$ to be the $A$-bilinear  form on $M^X(\widetilde{\lambda})_A^{\widetilde{\lambda}  - \eta}$ such that  for ${\pi,\gs \in {\bf{\overline{P}}}_{X}(\eta)}$, 

\be \label{shy} G^X_\eta(\widetilde{\lambda})(e_{-\gs}v_{\widetilde{\lambda}}^X,e_{-\pi}v_{\widetilde{\lambda}}^X)=[\gz_A(e_{\gs}e_{-\pi})(\widetilde{\lambda})],\ee
where $\gz_A:U(\fg)_A \lra S(\fh)_A$ is the Harish-Chandra projection, \cite{M} (8.2.13). We note that $\det G_\eta^X$ depends on the ordering of the basis,
as can be seen already in the case of $\fgl(2,1)$. 
However its leading term, which we denote by $\LT \det G_\eta^X$ is well-defined up to a scalar multiple.\\ \\
Our goal  in this section is to compute the determinants $\det F_\eta^X$ and $ \LT \det G_\eta^X$.
We point out at the outset some complications  that arise which are not present in the classical case \cite{H2} Theorem 5.8, \cite{M} Theorem 10.2.5. The first is a rather minor point: these determinants should really be
considered as elements of $\cO(\cH_X)=S(\fh)/\cI(\cH_X),$ but we shall express them as elements of $S(\fh)$ which map to the corresponding elements of $\cO(\cH_X)$. Remarkably the determinant $F_\eta^X$ factors a product of linear terms with leading terms of the form $h_\ga$ with $\ga$ a root, see Theorem \ref{shapdet}.
Apart from having to deal with two determinants, the first real complication arises since it is possible to have distinct non-isotropic positive roots
$\ga_1,$  ${\ga_2}$ such that $h_{\ga_1}$ and $h_{\ga_2}$ are proportional mod $\cI(\cH_X).$
 Indeed suppose that there are distinct orthogonal isotropic roots
 $\gc_1,$ ${\gc_2}\in X$ and  consider the hypotheses

\be \label{had} \ga_1^\vee \equiv -\ga_2^\vee \mod \Q \gc_1 +\Q \gc_2. \ee
or
\be \label{cad} \ga_1^\vee \equiv \ga_2^\vee \mod \Q \gc_1 +\Q \gc_2 \ee
If either \eqref{cad} or \eqref{had} holds, no third non-isotropic positive root $\ga_3$ satisfies
$\ga_1^\vee \equiv \pm \ga_3^\vee \mod \Q X$.
\\ \\
\noi
For example in type A, suppose that $i<j, k <\ell,$

\be\label{wan} \ga_1 = \gep_i -\gep_j, \;\ga_2 =\gd_k-\gd_\ell\ee and either
\be \label{won} \gc_1 = \gep_j -\gd_k, \; \gc_2 = \gep_i -\gd_\ell
                                         \ee or
\be\label{wun}
\gc_1 = \gep_j-\gd_\ell, \;\gc_2 = \gep_i -\gd_k.\ee
If we have \eqref{won}, then $\ga_1 +\gc_1 +\ga_2 =\gc_2$ so \eqref{had} holds.
If we have \eqref{wun} instead, then $\ga_1+ \gc_1 = \ga_2 +\gc_2$, so \eqref{cad} holds.
Note that for $\gl\in \cH_X$ we have,
$(\gl+\gr,\ga_1^\vee)= -(\gl+\gr,\ga_2^\vee)$ if \eqref{had} holds, and
$(\gl+\gr,\ga_1^\vee)= (\gl+\gr,\ga_2^\vee)$ if \eqref{cad} holds.  
The former case is easily dealt with.  To deal with the latter, 
let $E_{X}$ be the set of pairs $(\ga_1,\ga_2)$ such that \eqref{cad} holds.
In this case we assume $\ga_1 +\gc_1 +\ga_2 =\gc_2$, and
we frequently let $\ga$ denote $\ga_1$, especially to avoid double subscripts.
Similarly for compactness we write $[\ga] =(\ga_1,\ga_2)$.
A subalgebra of $\fg$ which is a direct sum of root spaces and is isomorphic to $\fgl(2,2)$ will be called a 
$\fgl(2,2)$-{\it subalgebra}. When \eqref{cad} holds we let $\fk[\ga]$ denote the $\fgl(2,2)$-subalgebra whose positive part is generated by the root vectors for the roots $\ga_1, \ga_2$ and $ \gc_1$. 
Let  $W'(\ga)$  (resp. $V'(\ga)$) be the set of all  (resp. all odd) positive roots of this subalgebra, and set $W(\ga)=W'(\ga)\cup X$,
$V(\ga)=V'(\ga)\cup X$.
For $\gl \in \cH_X$ define
\be\label{efg} E_X(\gl) = \{(\ga_1,\ga_2) \in E_{X}|(\gl+\gr,\ga_1^\vee) =1\}.\ee
There are two more situations to deal with. First set
$$B_X= \{\gc \in (\Gd^+_1\; \backslash \; X)|  \gc  \mbox{ is isotropic and } (\gc ,X) =0\}.$$ 
Finally let $C_X$ be the set of positive non-isotropic roots $\ga$, such that there is a unique
isotropic root $\gc\in X$ with
$(\gc,\ga^\vee)\neq0,$ and set $\gc = \Gc(\ga)$.
In this circumstance, since $X$ is an orthogonal set of roots, it follows that  $s_\ga \gc \notin X.$
Note that if $\ga \in C_X$ and $\gc=\Gc(\ga)$ we have $h_{{s_\ga \gc} } \equiv \pm h_\ga$
mod $\cI(\cH_X)$.
Let $C^+_X$ be the subset of $C_X$ consisting of those $\ga$ for which $(\Gc(\ga),\ga^\vee) >0.$ Thus
for  $\ga \in C^+_X,$ we have  $s_\ga \gc =\gc-\ga.$
For $\ga \in C_X^+$, set $Z(\ga) = X \cup s_\ga X= X\cup \{s_\ga\gc\}$. Then for $\gl\in \fh^*$, define
\[C_X(\gl) = \{\ga \in C_X^+|Z(\ga) \subseteq B(\gl)\}.\]
Note that if $\ga\in  C_X(\gl)$ we have $(\gl+\gr,\ga)=0$.
 \\ \\
\noi Next fix $\eta$ and  consider the following products
\begin{eqnarray}
D'_0&=& \prod_{\ga \in C^+_{X},
 \gc = \Gc(\ga)}h_{\ga}^{2 {\bf p}_{Z(\ga) }(\eta-\gc)}
\label{yaws}\\
D''_0&=&
\prod_{[\ga] \in E_{X}}h_{\ga}^{2({\bf p}_{V(\ga)}(\eta -\gc_2)+{\bf p}_{V(\ga)}(\eta -\gc_1-\gc_2)-{\bf p}_{W(\ga)}(\eta -\gc_1-\gc_2))},\label{yews}\\
D_{1} & = & \prod_{\alpha \in \ovd^+_{0}} \;\;
\prod^\infty_{r=1}
            (h_{\alpha} + (\rho, \alpha) - r(\alpha,\alpha)/2)^
            {{\bf p}_{s_\ga X}(\eta - r \alpha)}, \label{yew1}\\
D_2 & = &
\prod_{\gc  \in B_X } (h_{\gc } +(\rho, \gc ))^{{\bf p}_{X\cup\{\gc \}}(\eta - \gc )},\label{yew3}\\
D_3 &=&\prod_{\ga \in C^+_X, \gc = \Gc(\ga)} (h_{\ga}+ (\rho, \alpha) )^{{\bf p}_X(\eta)-{\bf p}_{s_{\ga} X}(\eta)},\label{yew4}\\
D_4&=& \prod_{[\ga] \in E_{X}}(h_{\ga}+ (\rho, \alpha_1) - (\alpha_1,\alpha_1)/2) ^{2{\bf p}_{W(\ga)}(\eta -\gc_1-\gc_2)}.\label{yew5}
\ey
Then
set $D_0=D'_0D''_0$ and $D^X_{\eta}=D_{1}D_{2}D_{3}D_4.$ Our main result on the \v Sapovalov  determinant is as follows. 
\begin{theorem} \label{shapdet}

\bi \itema Modulo the ideal defining $\cH_X$,  and up to a nonzero constant factor
\be \label{la8}
\det F^X_{\eta} =D^X_{\eta}.\ee
\itemb Up to a nonzero constant factor, $\det G^X_{\eta}$ has the same leading term as $D_0\det F^X_{\eta}.$\ei
\et

\brs \bi \itema{\rm
The factors in $D_0$ come from comparing the bilinear forms, see Theorems \ref{son} and \ref{epin}, also
\eqref{hid} (resp.  \eqref{rat3}-\eqref{cat1}) for the case of $\fgl(2,1)$
(resp. $\fgl(2,2)$). 
The other factors come from representation theory. The factors in $D_1$ are analogs of the classical ones, those in $D_2$ come from isotropic roots orthogonal to $X$ and those in $D_3$ come from Lemma  \ref{eon}. For $D_4$ see Theorem \ref{name}.
\itemb We note that the contributions from $D''_0$ and $D_4$ to the leading term simplify when combined, since $W(\ga)$ is eliminated.  We obtain
\begin{eqnarray}\label{yew}
\LT\;D''_0D_4&=&
\prod_{[\ga] \in E_{X}}h_{\ga}^{2({\bf p}_{V(\ga)}(\eta -\gc_2)+{\bf p}_{V(\ga)}(\eta -\gc_1-\gc_2))}.\;\;\;\;\;\;\;\ey
\itemc
\noi To compare the  factors in the leading term of
$\det G^X_{\eta}$ and $D_0D^X_{\eta}$,
we use some shorthand. Thus if $H\in S(\fh)$ and $h\in \fh,$ let $|H:h|$ denote the multiplicity of $h$ in the leading term  of $H$. Similarly $|H/H':h|$ means $|H:h|-|H':h|$. Note that
if $\ga \in \Gd^+$, then $\ga$ can belong to at most one of the sets $B_X, C_X, E_X$. Hence to find
$|D_\eta^X:h_\ga|$, it is enough to find  $|D_1D_i:h_\ga|$ for $i=2,3,4$, see also Remark \ref{dab}.
}\ei
 \ers

\bexa {\rm Suppose $\fg =\fgl(3,2)$ and set $\gc_1=\gep_2-\gd_1, \gc_2 =\gep_1-\gd_2, \ga_1 = \gep_1-\gep_2,$ and $\ga_2=\gd_1-\gd_2.$ 
Then $\ga_1 +\gc_1 +\ga_2 =\gc_2$.
There are three $W$-orbits on the set of sets of orthogonal isotropic roots of size two and we consider a representative from each orbit. If $X =\{\gc_1, \gc_2\}$, then  
$E_X=\{(\ga_1,\ga_2)\}$ and $s_{\ga_1}X = s_{\ga_2}X  =  \{\gep_1-\gd_1, \gep_2-\gd_2 \}$. 
We have 
\[ C_X =\{\gep_1-\gep_3, \gep_2-\gep_3 \} = C_X^+, \;\Gc(\gep_1-\gep_3) = \gep_1-\gd_1, \mbox{and } \Gc(\gep_2-\gep_3) = \gep_2-\gd_2.\]
 If $X' =\{\gep_1-\gd_1, \gep_3-\gd_2 \}$ and    
$X'' =\{\gep_2-\gd_1, \gep_3-\gd_2 \}$, then $C_{X'}^+ =\{\gep_1-\gep_2 \}$  and $C_{X''}^+$ is empty.
We have $|C_{X'}| =|C_{X''}|=2$ and  $E_{X'}, E_{X''}$ are singletons.
}\eexa

\subsection{Comparison of the Bilinear Forms.}\label{com}
\noi Our aim is to compare the $T$-adic valuation of the determinants  of the bilinear forms
$F^X_\eta(\widetilde{\lambda})$ and $G^X_\eta(\widetilde{\lambda})$ from Subsection  \ref{7.4}.
The main results, Theorems \ref{son} and \ref{epin} explain the presence of the terms $D_0'$ and $D_0''$ respectively in Theorem \ref{shapdet}.
Theorem \ref{son} generalizes the behavior in the case of $\fgl(2,1)$, see \eqref{hid} with $n=0$,
while Theorem \ref{epin} generalizes the behavior in the case of $\fgl(2,2)$, see \eqref{rat3}.
The proof of both results follows the same pattern, which we outline first.
In this subsection we denote the highest weight vector $v_{\widetilde{\lambda}}^X$ in ${M^X}(\widetilde{\lambda})_B$ simply by $v_{\widetilde{\lambda}}$.
\subsubsection{Strategy of the Proofs.}
  \noi As before, $A=\ttk[T]$ and $B=\ttk(T)$.   \noi Let $C$ be the local ring $C=A_{(T)}$ with maximal ideal $(T)=TC$.
	As far as the $T$-adic valuation is concerned, we may work over $C$ rather than over $A$.
	If $\eta\in Q^+$, then by Theorem \ref{zoo}, a non-zero $p \in {M^X}(\widetilde{\lambda})_B^{\widetilde{\lambda}-\eta}$ can be written uniquely in the form

\be \label{sea} p=\sum_{\gt \in {{\bf{\overline{P}}}_X(\eta)}}a_\gt e_{-\gt}v_{\widetilde{\lambda}}^X \in
{M^X}(\widetilde{\lambda})_B^{\widetilde{\lambda}-\eta},\ee
with $a_\gt \in B$, and we set $\Supp\; p =\{\gt \in {\bf{\overline{P}}}_X(\eta)|a_\gt \neq 0\}$.
\\ \\
Also if $\ga\in \Gd_0^+$, an {\it $\ga$-principal part} of $p$ is a partition $\gs\in \Supp\; p$
such that $\gs(\ga)>\gt(\ga)$ for all $\gt\in \Supp\; p, \gt \neq \gs$.
Clearly the $\ga$-principal part of $p$ is unique if it exists.
\\ \\
Fix  $\gl\in \cH_X$, $\eta\in Q^+,$ and set $M_C = {M^X}(\widetilde{\lambda})_C^{\widetilde{\lambda}-\eta}$.
Let $N_C$ the $C$-submodule of $M_C$ spanned by the set $\{e_{-\pi}v_{\widetilde{\lambda}}| \pi \in {\bf{\overline{P}}}_{X}(\eta)\}$.
By the elementary divisor theorem, there is a $C$-basis
$v_1,\ldots, v_k$ for $M_C$ and integers $0\le z_1 \le  \ldots \le z_k$ such that

\be \label{rex} w_1=T^{z_1}v_1,\ldots, w_k=T^{z_k}v_k\ee
 is a $C$-basis for $N_C$.  Let $z_{\gl,\eta} = 2\sum_{i=1}^k z_i$. 
Suppose we use these bases to find $\det G^X_{\eta}(\widetilde\gl)$ and $\det F^X_{\eta}(\widetilde\gl)$.  
Then if we factor $T^{\sum_{i=1}^k z_i}$ from both the rows and the columns of the Gram matrix of $G^X_{\eta}(\widetilde\gl)$ we obtain 
the Gram matrix of $F^X_{\eta}(\widetilde\gl)$. Therefore

\be \label{ron}   z_{\gl,\eta} = v_T(\det G^X_{\eta}(\widetilde\gl))-v_T(\det F^X_{\eta}(\widetilde\gl)). \ee
We need to consider more general bilinear forms. Suppose that ${\bf V}$ is a basis for $M_B= M_C\ot B$ contained in $M_C$,
and define

$$G^X_{\eta,{\bf V}}(\widetilde\gl) =\det (\gz(b^t b'))_{b,b'\in {\bf V}}. $$
Thus if {\bf W} is a $C$-basis for $N_C$, then $v_T(\det G^X_{\eta,{\bf W}}(\widetilde\gl))=v_T(\det G^X_{\eta}(\widetilde\gl))$.
At the other extreme, suppose that ${\bf U}$ is a basis of $M_C$ as a $C$-module,
then 
\be\label{dff} v_T(\det G^X_{\eta,{\bf U}}(\widetilde\gl))=v_T(\det F^X_{\eta}(\widetilde\gl)).\ee
Consider bases ${\bf W}= \{w_1,\ldots, w_k\}$ and ${\bf V}= \{v_1,\ldots, v_k\}$  such that \eqref{rex} holds and

\be \label{nub}  N_C=\CS{\bf W}\subseteq \CS
 {\bf V}  \subseteq \CS {\bf U}.\ee
Then we have

\be \label{uuu}
z_{\gl,\eta} \ge v_T(\det G^X_{\eta}(\widetilde\gl))-v_T(\det G^X_{\eta, {\bf V}}(\widetilde\gl)),\ee
and we show that when the hypothesis of Theorem \ref{son} (resp. Theorem \ref{epin})  hold for an even root $\ga$, the right side of
\eqref{epin} is bounded below by $|\LT D'_0:h_\ga|$ (resp. $|\LT D''_0:h_\ga| $).
\\ \\
\noi At this point the proofs of the two results begin to diverge.  For Theorem \ref{son}, to construct  {\bf V, W} we use the explicit expressions for $\gth_\gc v_{\widetilde{\lambda}}$ given in Theorem \ref{stc}.
This yields non-zero elements $p, q\in {M^X}(\widetilde{\lambda})_C^{\widetilde{\lambda}-\gc}$
such that $p+qT=\gth_\gc v_{\widetilde{\lambda}}=0$.
Then for certain partitions $\gs$ we set $p_\gs= e_{-\gs}p$ and $q_\gs= e_{-\gs}q$.
The basis {\bf  W} contains all of the $p_\gs$, and  ${\bf V}$ is obtained from ${\bf W} $ by replacing each $p_\gs$ by $q_\gs$.
 Then to prove Theorem \ref{son}, we simply have to count the number of partitions affected.
The situation for Theorem \ref{epin} is a bit more complicated.  We begin with Theorem \ref{nex}, but then there are two different cases involved and certain partitions need to be excluded.

\subsubsection{The Factor $D_0'$.}\label{pea}
We start with the factor $D_0'$ from Equation \eqref{yaws}.
Suppose that $\ga\in C^+_X$, $\gc = \Gc(\ga)$ and set $Z=Z(\ga)=X \cup s_\ga X$. We adopt the notation of Theorem  \ref{stc}, with $\ttT=T$.
Thus $\gc = \gep_r- \gd_{s},$ $\ga = \gep_r-  \gep_{\ell }$  and
$\gc'=s_\ga\gc=\gep_\ell- \gd_{s}.$ Set $j=\ell-r$. 

\bt \label{son}
For  general $\gl \in \cH_X \cap \cH_{\ga, 0}$, we have

\be \label{laid}z_{\gl,\eta} \ge 2 \;{\bf p}_{Z(\ga)}(\eta-\gc) = |\LT D'_0:h_\ga| .\ee
\et

\noi   The equality  in \eqref{laid} holds by the definition of $D_0'$, see \eqref{yaws}. To prove the inequality, we need some preparation.  
The hypothesis in the Theorem implies that $\gc, \gc'=s_\ga \gc =\gc-\ga \in B(\gl)$ for some odd root $\gc$.
If $p= \gth_{\ga,1}\gth_{\gc'}v_{\widetilde{\lambda}}$ and
$q=\overline{\psi}(\gth_{\overline{\gc}}(\overline{\phi}({\widetilde{\lambda}})))v_{\widetilde{\lambda}}$,
then by \eqref{nice} $\gth_\gc v_{\widetilde{\lambda}}= p+qT$.
Suppose $\gs \in {\bf{\overline{P}}}_{Z}(\eta-\gc)$, define $p_\gs, q_\gs$ as in the previous subsection and set 
$\gs'=\gs+\pi^\ga+\pi^{\gc'}$. Thus 
\[\gs'(\xi)
=\left\{ \begin{array}
  {cl}
	\gs(\ga)+1 &\mbox{if} \;\; \xi=\ga
\\
	1 &\mbox{if} \;\; \xi=\gc'
\\
\gs(\xi) &\mbox{otherwise}.
\end{array} \right. \]
Note that $\gs'\in {\bf{\overline{P}}}_{X}(\eta),$ since  $\gc =\gc'+\ga$, $Z=X \cup \{\gc'\}$  and $\gs\in {\bf{\overline{P}}}_{Z}(\eta-\gc)$.
\bl \label{app} For general $\gl \in \cH_X \cap \cH_{\ga, 0}$  
\bi
\itema 
The coefficient of
$e_{-\gs'}v_{\widetilde{\lambda}}$ in $\gth_{\ga,1}\gth_{\gc'}v_{\widetilde{\lambda}}$ is invertible in $C$.
\itemb The $\ga$-principal part of $p_\gs$ is  $\gs'$.
\ei\el

\bpf 
Define $a_i$ and $b_j$ as in \eqref{eat}.  Then
\by  \widetilde{a_i}&=&(\widetilde{\lambda} + \gr,\gs_{r,r+i}^\vee ) ={a_i} + T(\xi,\gs_{r,r+i}^\vee )\\
\widetilde{b_i} &=& (\widetilde{\lambda} + \gr,\gt_{i,s}^\vee )\;\;\;\; = 
b_i +T(\xi,\gt_{i,s}^\vee ).\nn 
\ey 
To find the coefficients of $p$ and $q$ we evaluate the determinants appearing in the determinants \eqref{wok} and \eqref{gov} with  $\gl$ replaced by $\widetilde{\lambda}$, that is with $a_i,b_i$ replaced by $\widetilde{a_i}, \widetilde{b_i}$ respectively. 
Thus if the $b_i,$ and the $a_i$  are non-zero, then the 
$\widetilde{b_i},$ and the $\widetilde{a_i}$  (with $i\neq j$) are units in $C$. Thus for general ${\lambda}$ the coefficients of $p$ and $q$ are units in $C$.
\\ \\
Now to prove (b), observe first that $\pi^\ga +\pi^{\gc'} \in \Supp\; \gth_{\ga,1}\gth_{\gc'}v_{\widetilde{\lambda}}$,
and for general $\gl$, $e_{-\gs}e_{-\ga}e_{-\gc'}v_{\widetilde{\lambda}}$ occurs in
$e_{-\gs}\gth_{\ga,1}\gth_{\gc'}v_{\widetilde{\lambda}}$ with a coefficient  that is a unit in $C$.
Recall that odd root vectors occur last in the product $e_{-\gs}$, and that the root vectors for odd positive roots anti-commute since we are working in type $A.$
Thus the only factor in $e_{-\gs}e_{-\ga}e_{-\gc'}v_{\widetilde{\lambda}}$ that
could be out of order is the middle root vector  $e_{-\ga}$.
Also for any term $e_{-\go}$ that arises by taking commutators of $e_{-\ga}$ with root vectors for roots that are involved in $\gs$,
the partition $\go$ satisfies $\go(\ga)<\gs'(\ga)$.\\ \\
We consider in more detail the odd root vectors $e_{-\go}$ that
arise in the way described in the above paragraph, using commutators with $e_{-\ga}$. We need to check that for such $\go$ we have 
$\go(\gk)=0$  for all $\gk\in X$. 
(Clearly this also holds if  $e_{-\go}$ arises from commutators of $e_{-\ga}$ with other even roots.)
Since $\ga\in C^+_X$, we have $(\gc,\ga^\vee)=1.$ 
Thus without loss of generality, we can assume that $\ga=\gep_j-\gep_k, \gc' = \gep_k-\gd_r $ and $\gc =\gc'+\ga=\gep_j-\gd_r$.  Then set
\[I(\gs) =\{s\in [n]\;|\; \gs(\gep_k-\gd_s) = 1, \gs(\gep_j-\gd_s) = 0, s\neq r\}\]
and for $\gs \in {\bf{\overline{P}}}_{Z}(\eta-\gc)$ and $s\in I(\gs)$, define  a partition $ \gs^{(s)}$ of $\eta$ by
\[\gs^{(s)}(\xi)
=\left\{ \begin{array}
  {cl}
	0 &\mbox{if} \;\; \xi=\gep_k-\gd_s
\\
	1 &\mbox{if} \;\; \xi=\gep_j-\gd_s
\\
1 &\mbox{if} \;\; \xi=\gc'
\\
\gs(\xi) &\mbox{otherwise}.
\end{array} \right. \] Note that $\gep_j-\gd_s\notin X$ because $(\gep_j-\gd_s,\gc)\neq 0$ and $X$ is an orthogonal set of roots containing $\gc$. Since $\gs \in {\bf{\overline{P}}}_{Z}(\eta-\gc)$ it follows that $\gs^{(s)} \in {\bf{\overline{P}}}_{X}(\eta)$.
Now we have, modulo the span of terms obtained from commutators of $e_{-\ga}$ with even root vectors,

\be \label{007} e_{-\gs}e_{-\ga}e_{-\gc'} \equiv e_{-\gs'} + \sum_{s\in I(\gs)}  c_\gs e_{-\gs^{(s)}}, \ee

\noi with non-zero $c_\gs\in \ttk$, and   $\gs^{(s)}(\ga) = \gs(\ga)=\gs'(\ga)-1$  for all $s\in I(\gs)$.
Thus $\gs'$ is the $\ga$-principal part of $p_\gs$.
\epf

\bexa {\rm Suppose that $\fg=\fgl(3,3)$, $\ga=\gep_1-\gep_3, \gc'=\gep_3 -\gd_2, \gc=\gep_3 -\gd_2,$ $e_{-\gs} = e_{63}e_{43}$ 
and $X=\{\gc\}$.  

Suppose we want to collect all odd root vectors on the right. (Note that all root vectors for positive odd roots commute). 
Then $e_{-\ga}e_{-\gc'}v_{\widetilde{\lambda}}$ occurs with non-zero coefficient 
in $\gth_{\ga,1}\gth_{\gc'}v_{\widetilde{\lambda}}$. 
However, since the root vectors $e_{63}$ and $e_{43}$ are odd, $e_{-\ga}$ is 
out of order in $e_{-\gs}e_{-\ga}e_{-\gc'} v_{\widetilde{\lambda}}=e_{63}e_{43}e_{31}e_{53}v_{\widetilde{\lambda}}$.  
In the notation of the Lemma, $j=1, k = 3, r =2$ and $I(\gs)= \{1, 3\}$, and we have 
\by e_{-\gs}e_{-\ga}e_{-\gc'} 
&=& e_{31}e_{63}e_{43}e_{53}\nn\\
&+&  e_{61}e_{43}e_{53}\nn\\
 &+& e_{63}e_{41}e_{53}.\nn\ey
Now all terms are in the correct order, and $e_{-\ga}$ occurs only in the first term.} 
\eexa

\noi Next recall that $N_C= \CS \{e_{-\pi}v_{\widetilde{\lambda}}| \pi \in {\bf{\overline{P}}}_{X}(\eta)\}$, and set
\[{\bf W}_{1} = \{e_{-\pi}v_{\widetilde{\lambda}}|\pi\in{\bf \overline{P}}_{X}(\eta),
\pi\neq  \gs' \mbox{ for any } \gs\in{\bf \overline{P}}_{Z}(\eta-\gc)\}.\]

\bl \label{0,111} For general $\gl \in \cH_X \cap \cH_{\ga, 0}$
the set
\[{\bf W}= \{p_\gs v_{\widetilde{\lambda}}|\gs\in{\bf \overline{P}}_{Z}(\eta-\gc)\}\cup {\bf W}_1\]
is a $B$-basis for ${M^X}(\widetilde{\lambda})_B^{\widetilde{\lambda}-\eta}$ and $\CS {\bf W} = N_C$.
\el

\bpf Since $|{\bf W}| = \dim_B {M^X}(\widetilde{\lambda})_B^{\widetilde{\lambda}-\eta}= \rank_C  N_C$, it suffices to prove the last statement.
Obviously $L=\CS {\bf W}\subseteq N_C$.
If  $\pi \in {\bf \overline{P}}_{X}(\eta)$ we show by induction on $\pi(\ga)$ that $e_{-\pi}v_{\widetilde{\lambda}} \in L$.
By definition of ${\bf W}$
we only have to show this when $\pi=\gs'$ for some $\gs \in {\bf \overline{P}}_{Z}(\eta-\gc)$.
Write
\by \label{sos} p_\gs &=& u_\gs e_{-\gs}e_{-\ga}e_{-\gc'} v_{\widetilde{\lambda}}
+ \sum_{\go \in {\bf{\overline{P}}}_{X}(\eta):\go\neq \gs'} a_{\gs,\;\go} e_{-\go}v_{\widetilde{\lambda}}\nn\\
&=& u_\gs e_{-\gs'} v_{\widetilde{\lambda}} + \sum_{\go \in {\bf{\overline{P}}}_{X}(\eta):\go\neq \gs'}  b_{\gs,\;\go} e_{-\go}v_{\widetilde{\lambda}},
\ey

\noi where $a_{\gs,\;\go}, b_{\gs,\;\go} \in C,$  and $u_\gs $ is a unit in $C$.
Each term $e_{-\go}v_{\widetilde{\lambda}}$ in the sum \eqref{sos} is obtained by expanding
$e_{-\gs}\gth_{\ga,1}\gth_{\gc'}v_{\widetilde{\lambda}}$ and re-ordering (as in the proof of the previous lemma). 
Observe that when this is done, any $\go$ arising from a commutator of $e_{-\ga}$
with another root vector satisfies $\go(\ga)<\pi(\ga)$.
Note that for general  $\gl,$ $e_{-\gc'}v_{\widetilde{\lambda}}$ occurs in $\gth_{\gc'}v_{\widetilde{\lambda}}$ with a coefficient which is a unit in $C$.
Since $\gs(\gc')=0$,    any term $e_{-\go}v_{\widetilde{\lambda}}$ in the expansion with $\go =\gt'$ (where $\gt \in {\bf{\overline{P}}_{Z}}(\eta-\gc)$) must be a multiple of $e_{-\gc'}
v_{\widetilde{\lambda}} \in \Supp \;\gth_{\gc'}v_{\widetilde{\lambda}}$.
Using Lemma \ref{app}, it follows that for all $\go$ such that  $b_{\gs,\;\go}\neq 0$ in \eqref{sos} we have either $\go \neq \gt'$ or $\go(\ga)<\pi(\ga).$
In either case $e_{-\go}v_{\widetilde{\lambda}} \in L$, so since $u_\gs$ is a unit in $C$, we obtain the result from \eqref{sos}.
\epf

\noi {\it Proof of Theorem \ref{son}.}
We claim first that the set
$\{p_\gs v_{\widetilde{\lambda}}|\gs\in{\bf \overline{P}}_{Z}(\eta-\gc)\}$ is $B$-linearly independent.
Suppose that
\[
\sum_{\gs\in{\bf \overline{P}}_{Z}(\eta-\gc)} a_\gs p_\gs =0,
\]
and set $k=\max\{\gs(\ga)|a_\gs \neq 0\}.$ Then from \eqref{sos} we have
\[ \sum_{\gs\in{\bf \overline{P}}_{Z}(\eta-\gc):\gs(\ga)=k} c_\gs e_{-\gs'} =0,\]
where $c_\gs =a_\gs u_\gs,$ for $u_\gs$ as in \eqref{sos}. But by Theorem \ref{zoo} the elements
$e_{-\pi}v_{\widetilde{\lambda}}$ with $\pi\in{\bf \overline{P}}_{X}(\eta)$ are $B$-linearly independent, 
and $\gs'\in {\bf{\overline{P}}}_{X}(\eta)$ by Lemma \ref{app}.
Thus all $c_\gs$ are zero, a contradiction which proves the claim.
\\ \\
Next set
\[{\bf V}= \{q_\gs |\gs\in{\bf \overline{P}}_{Z}(\eta-\gc)\}\cup {\bf W}_1.\]
Since $p_\gs +q_\gs T =0,$  it follows from  Lemma \ref{0,111} that
 there are bases $ w_1,\ldots,  w_k$  and $ v_1,\ldots,  v_k$ for $N_C$ and $\CS{\bf V} $ respectively,
and integers $0\le z_1 \le  \ldots \le z_k$ such that \eqref{rex} holds
and at least ${\bf p}_{Z}(\eta-\gc)$ of the $z_i$ are positive.
Hence $$v_T(\det G^X_{\eta}(\widetilde\gl))-v_T(\det G^X_{\eta, {\bf V}}(\widetilde\gl))\ge
2 \;{\bf p}_{Z}(\eta-\gc)$$ and combined with \eqref{uuu}, this yields \eqref{laid}.
\hfill  $\Box$

\subsubsection{The Factor $D_0''$.}\label{pie}

We now turn our attention to the factor $D_0''$ from Equation \eqref{yews}.
If  $(\ga_1,\ga_2)\in E_X$ there are roots
$\gc_1,\gc_2\in X$ such that $\ga_1 +\gc_1 +\ga_2 =\gc_2$.
Set $\ga =\ga_1$. The main result is as follows.
\bt \label{epin}
If $\eta \ge \gc$ 
then for general $\gl \in \cH_X\cap \cH_{\ga,1}$ we have
\be \label{laid2}z_{\gl,\eta} \ge 2({\bf p}_{V(\ga)}(\eta -\gc_2)+{\bf p}_{V(\ga)}(\eta -\gc_1-\gc_2)-{\bf p}_{W(\ga)}(\eta -\gc_1-\gc_2)) = |\LT D''_0:h_\ga| .\nn\ee
\et
\noi
Again the equality follows from the definition of $D_0''$, see \eqref{yews}. For the rest we first introduce some notation. The proof is similar to the proof of Theorem \ref{son}, so we omit some of the details. First define
\[\Gt_\eta^{(1)} =\{\pi\in{\bf{\overline{P}}}_{V(\ga)}(\eta -\gc_1-\gc_2)|\pi \notin {\bf{\overline{P}}}_{W(\ga)}(\eta -\gc_1-\gc_2)\}, \quad  \Gt_\eta^{(2)}= {\bf{\overline{P}}}_{V(\ga)}(\eta -\gc_2),\]
 and set $\Gt_\eta= \Gt_\eta^{(1)} \cup \Gt_\eta^{(2)}$.
\\ \\
Since $(\gl+\gr, \ga^\vee) =1,$ we have $(\widetilde{\lambda}+\gr, \ga^\vee) =T+1.$ Adopting the notation of Subsection \ref{det} we define $\gl_1,\gl_2,$ and $ \gl_3$ as in \eqref{flx}, \eqref{fax} and \eqref{fox}. Then from Equation \eqref{mice} we see that if
$\gth_{\gc_2}v_{\widetilde{\lambda}} =\gth_{\gc_1}
v_{\widetilde{\lambda}}= 0$ we have, (compare  \eqref{rat3} for the case $\fgl(2,2)$)
\\
\be \label{ump} [\underline{\overline\psi}(\gth_{\underline{{\overline{\gc}}}}(\underline{\overline{\phi}}({\widetilde{\lambda}})))(T +1)+
\gth_{{\ga_1}}\underline{\psi}(\gth_{{\gb_1}}({{\lambda_2}})) -\gth_{{{\ga_2}}}\overline{\psi}(\gth_{{\gb_2}}({{\lambda_1}}))]v_{\widetilde{\lambda}}=0.\ee
Let
\[p_1 =\gth_{{\ga_1}}\underline{\psi}(\gth_{{\gb_1}}({{\lambda_2}}))v_{\widetilde{\lambda}},\quad 
p_2 =\gth_{{{\ga_2}}}\overline{\psi}(\gth_{{\gb_2}}({{\lambda_1}}))v_{\widetilde{\lambda}},\quad
q= \underline{\overline\psi}(\gth_{\underline{{\overline{\gc}}}}({{\lambda_3}}))v_{\widetilde{\lambda}}
\]
and $p=p_1-p_2+q.$ Then in $U(\fn^-)_A^{-\gc_2}v_{\widetilde{\lambda}}$ we can write \eqref{ump} as $p+qT=0,$
\bl \label{log} \bi \itema $|\Gt_\eta^{(1)}|= {\bf p}_{V(\ga)}(\eta -\gc_1-\gc_2)-{\bf p}_{W(\ga)}(\eta -\gc_1-\gc_2)$.
\itemb If $\gs \in \Gt_\eta^{(1)}$, then either $\gs(\ga_1)>0$  or $\gs(\ga_2)>0.$
\itemc If $\gs \in {\bf{\overline{P}}}_{V(\ga)}(\eta -\gc_1-\gc_2)$,  then
$e_{-\gs}p\neq 0$. 
\ei \el
\bpf The first two parts follow easily from the definitions, so we turn to (c). If $\pi_1= \pi^{\ga_1}+\pi^{(\gc_1 + \ga_2)}$ and $\pi_2=\pi^{\ga_2}+\pi^{(\ga_1+\gc_1)}$, 
\ff{ We can assume $\ga_1 = \gep_r -\gep_\ell, \;\ga_2 =\gd_k-\gd_s,
\gc_1 = \gep_\ell- \gd_k$ and $\gc_2 = \gep_r -\gd_s.$ Then putting odd root vectors last, we have
$e_{-\pi_1}= e_{\ell,r} e_{m+s,\ell}$ and $e_{-\pi_2}= e_{m+s,m+k}e_{m+k,r}$.} 
we have by Lemma \ref{hes} (c)

\be \label{p1p} \pi_1\in \Supp \; p_1,\;
\pi_1\notin \Supp \; p_2,\; \pi_1\notin \Supp\; q,\ee
and
\be \label{p2p} \pi_2\in \Supp \; p_2,\;
\pi_2\notin \Supp \; p_1,\; \pi_2\notin \Supp\; q.\ee
 \noi Next suppose $\gs \in {\bf{\overline{P}}}_{V(\ga)}(\eta -\gc_1-\gc_2)$.  Then
$\gs(\ga_1+\gc_1)= \gs(\gc_1 + \ga_2)=0.$
Combined with \eqref{p1p} and \eqref{p2p}, this shows that $e_{-\gs}p\neq 0$. \epf
\noi 
\bl \label{axe} Suppose $\gs\in {\bf{\overline{P}}}_{V(\ga)}(\eta -\gc_1-\gc_2)-{\bf{\overline{P}}}_{W(\ga)}(\eta -\gc_1-\gc_2)$,
and define \[\gs' = \gs +\pi^{\gc_1 + \ga_2}+\pi^{\ga_1+\gc_1} \;\in \;{\bf{\overline{P}}}_{X}(\eta).\] 
Then for general $\gl$ and for $i=1,$ or 2, $\gs'$ is the $\ga_1$ {\rm(}resp. $\ga_2${\rm)} principal part of
$e_{-\gs}e_{-(\gc_1 + \ga_2)}e_{-(\ga_1+\gc_1)}v_{\widetilde{\lambda}}$. 
\el
\bpf
By  Lemma  \ref{log} we have either $\gs(\ga_1)> 0$, or $\gs(\ga_2)> 0$, and we assume the former.
Define  the partition $\hat \gs$ by
\[\hat \gs(\xi)  =\left\{ \begin{array}
  {ccl}\gs(\ga_1)- 1 &\mbox{if} \;\; \xi=\ga_1
\\\gs(\xi) &\mbox{otherwise}
\end{array} \right. \]
Then set
$p_\gs =  e_{-\hat \gs} e_{-(\ga_1+\gc_1)}p$ and 
$q_\gs =  e_{-\hat \gs} e_{-(\ga_1+\gc_1)}q$.
\\ \\
We have, up to re-ordering terms
\be \label{ess1} e_{-\gs}e_{-(\gc_1 + \ga_2)}e_{-(\ga_1+\gc_1)}v_{\widetilde{\lambda}}= e_{-\hat \gs}
e_{-\ga_1}e_{-(\gc_1 + \ga_2)}e_{-(\ga_1+\gc_1)}v_{\widetilde{\lambda}}\ee
Up to re-ordering these elements
are in $\Supp\; p_\gs$ and for a general $\gl\in \cH_X\cap \cH_{\ga,1}$,  they occur in $p_\gs$ with coefficient that is a unit in $C$.
The result follows. (Compare \eqref{dog1} and \eqref{cat1} for the case of $\fgl(2,2).$)
\epf
\noi We need to consider additional relations arising in the following way.
For $\gs \in \Gt_\eta^{(2)}={\bf{\overline{P}}}_{V(\ga)}(\eta -\gc_2)$ we have the relation $e_{-\gs}(p+qT)=0.$  Set $p_\gs =e_{-\gs}p$  and $q_\gs = e_{-\gs}q.$ Set $\gs' = \gs+ \pi^{\ga_1}+\pi^{(\gc_1 + \ga_2)} \in{\bf \overline{P}}_{X}(\eta).$ As in Lemma \ref{axe}
it can be shown that for general $\gl$, $\gs'$ is the $\ga_1$ principal part of $p_\gs$. 
\bl \label{eee}
Set
\[{\bf W}_{1} = \{e_{-\pi}v_{\widetilde{\lambda}}|\pi\in{\bf \overline{P}}_{X}(\eta),
\pi\neq  \gs' \mbox{ for any } \gs\in\Gt_\eta\}.\]
Then Equation \eqref{nub} holds with
\[{\bf W}= \{p_\gs |\gs\in\Gt_\eta\}\cup {\bf W}_1, \quad
{\bf V} =\{q_\gs |\gs\in\Gt_\eta\}\cup {\bf W}_1.\]
Furthermore {\bf W} and {\bf V}
are $B$-bases for ${M^X}(\widetilde{\lambda})_B^{\widetilde{\lambda}-\eta}$.\el

\bpf Similar to the proof of Lemma \ref{0,111}.\epf

\noi {\it Proof of Theorem \ref{epin}.}
First we use Theorem \ref{zoo} to show the elements of {\bf W} are linearly independent and $N_C=\CS{\bf W}$. Suppose that
\[\sum_{\gs \in \Gt_\eta^{(1)}} a_\gs p_\gs + \sum_{\gs \in \Gt_\eta^{(2)}} b_\gs p_\gs =0,\]
with coefficients $a_\gs, b_\gs \in B.$
Then both sums must be zero because for $\gs$ in the  first sum  we have
\[\gs'(\ga_1+\gc_1) =\gs'( \gc_1 + \ga_2)=1,\]
but for $\gs$ in the second sum we have
\[\mbox{either } \gs'(\ga_1+\gc_1) =0 \mbox{ or } \gs'( \gc_1 + \ga_2)=0.\]
Arguing as in the proof of Theorem \ref{son}, it follows that all coefficients $a_\gs, b_\gs$ are zero.
Hence \eqref{rex} holds for the bases {\bf W} and {\bf V},
and at least \[{\bf p}_{V(\ga)}(\eta -\gc_2)+{\bf p}_{V(\ga)}(\eta -\gc_1-\gc_2)-{\bf p}_{W(\ga)}(\eta -\gc_1-\gc_2) \]
 of the $z_i$ are positive.
The proof is concluded as before.
\hfill  $\Box$
\bc
For all $\ga\in \Gd^+_0$,
\be \label{(c)}|\LT \det G^X_\eta:h_\ga| \ge |\LT D_0 :h_\ga| + |\LT  \det F^X_\eta :h_\ga|.\ee
\ec
\bpf Combine Theorems \ref{son} and \ref{epin}.\epf
\br{\rm In this Subsection we have been concerned with situations where $z_{\gl,\eta} $ is as large as possible.  However for general $\gl\in \cH_X$, we have $z_{\gl,\eta} = 0$. Indeed if set $\gl_c = \gl +c\xi$ for $c\in \ttk$, it  follows from Corollary \ref{nob} that  for all but finitely many $c$, $z_{\gl_c,\eta} = 0$ for all $\eta,$ that is
\be \label{poe}   u\det  G^X_\eta(\widetilde\gl_c) = \det F^X_\eta(\widetilde\gl_c),\ee for some unit $u$ in $C$.}
\er

\subsection{Partition Identities.} \label{Partitions.}

We refer to $p, p_X$ and $p_\ga$ as defined by \eqref{pfun} as {\it
partition functions}.  They are functions in the following sense.  Suppose $\nu, \mu\in Q^+$.
If the $\Z$-linear function  $\tte^{\nu}$ from $\Z Q$ to $\Z$ is defined 
by $\tte^{\nu}(\mu) = \gd_{\nu,\mu}$, then 
${p_X}(-\mu) = \sum_\eta {\bf p}_X(\eta)\tte^{-\eta}(-\mu) = {\bf p}_X(\mu)$. Note that 
$\tte^{-\ga}p_X(-\mu) = {\bf p}_X(\mu-\ga)$.
 We have the following relations between partition functions.

\bl \label{ink}
Suppose $\ga \in C_{X}, \gc = \Gc(\ga)\in X$ with $\gc' = s_\alpha\gc
 \notin X.$
Set $Y = s_\ga X,$ $ Z = X\cup Y$, and $\mu_0 = \mu-\gc$.  Then
\bi \itema If $\gc'= \gc -\ga,$ then
\[({p}_X- {p}_{s_\ga X})= {\tte^{-\gc'}(1 - \tte^{- \ga})}{p}_Z,\] and
\[{{\bf p}_{X}(\mu)}-{{\bf p}_{s_\ga X}(\mu)} ={\bf p}_Z{(\mu_0  +\alpha)}-{{\bf p}_Z(\mu_0 )}.\]
\itemb If  $\gc'= \gc +\ga,$ then
\[{{\bf p}_{X}(\mu)}-{{\bf p}_{s_\ga X}(\mu)} ={\bf p}_Z{(\mu_0 -\alpha)}-{{\bf p}_Z(\mu_0)}.\]
\ei
\noi Now suppose $(\ga_1,\ga_2) \in E_{X}$ and $Y = s_{\ga_1} X= s_{\ga_2} X$.  Write $\ga_1 +\gc_1 +\ga_2 =\gc_2$.  Then
if $V=X\cup Y$,
\bi \itemc \[{p}_Y-{p}_{X}=
[(1 + \tte^{-\gc_1})(1 + \tte^{-\gc_2})-(1 + \tte^{-\gc_1- \ga_1})(1 + \tte^{-\gc_1- \ga_2})]{p}_{V}.\]
\itemd
\[{\bf p}_Y{(\mu)}-{{\bf p}_{X}(\mu)}=
{{\bf p}_{V}(\mu-\gc_1)} +{{\bf p}_{V}(\mu-\gc_2)} -
{{\bf p}_{V}(\mu-\ga_1-\gc_1)}-{{\bf p}_{V}(\mu-\ga_2-\gc_1)}.\]
\ei\el
\bpf The first statement in (a) follows since
\by p_X - p_{s_\ga X} &=&[(1 + \tte^{-\gc'})-(1 + \tte^{-\gc})]p_Z\nn\\
&=&{\tte^{-\gc'}(1 - \tte^{- \ga})}{p}_Z.\nn\ey
The second follows by evaluation at $-\mu$. The proofs of  the other parts are  similar.\epf

\subsection{The Leading Term.}\label{lete}
We make a direct computation of the leading term of $ \det G^X_{\eta}$  by adapting the proof of \cite{M} Lemma 10.1.3, and show the result is consistent with Theorem \ref{shapdet}. We prove Theorem \ref{shapdet} by showing that for each factor of $D_\eta$ we have \be \label{por} |D_\eta:h_\ga|\ge |G_\eta:h_\ga|.\ee  Since
$D_\eta$ and $G_\eta$ have the same leading term, it follows that there can be no factors of $D_\eta$ other than those listed. We need to consider the cases where $\gl$ lies on a hyperplane inside $\cH_X$, and among such $\gl$ it suffices to look at the most general behavior.  In such cases we will have equality in \eqref{por}.
\\ \\
In the next result we use the notation $B_X, C_X, E_X, V(\ga)$ 
 as in the statement of Theorem \ref{shapdet}.

\bl \label{gro} The leading term of $ \det G^X_{\eta}$ has the form $\LT \det  G^X_\eta =G_{1}G_{2}G_{3}G_4$ where
\by \label{la1}
G_1&=&\prod_{\alpha \in \ovd^+_{0}}\prod^{\infty}_{r=1}h_\ga^{{\bf p}_{X}(\eta-r\alpha)},\quad
G_3 =\prod_{\ga \in C_{X}, \gc = \Gc(\ga)} \prod_{\pi \in {\bf \overline{P}}_{X}(\eta)}
 h_{{\ga} }^{\pi({s_\ga \gc} )},\nn\\
G_2&=&\prod_{\gc  \in B_X } h_{\gc }^{{\bf p}_{X\cup\{\gc \}}(\eta - \gc )}, \;
G_4 =\prod_{[\ga] \in E_{X}}
\prod_{\pi \in {\bf \overline{P}}_{X}(\eta)} h_{{\ga} }^{\pi({\ga_1+\gc_1} )+\pi({\ga_2+\gc_1} )}.\nn
\ey

\br \label{dab} {\rm
Note that 
\be \label{nat}G_3 =\prod_{\ga \in C_{X}}  h_{{\ga} }^{{\bf{\overline{p}}}_{X}(\eta)-{\bf{\overline{p}}}_{s_\ga X}(\eta)}\ee
and 
\be \label{nit}G_4 =\prod_{[\ga] \in E_{X}} h_{{\ga} }^{{\bf{\overline{p}}}_{V(\ga)}(\eta-\ga_1-\gc_1 )+ 
{{\bf{\overline{p}}}}_{V(\ga)}(\eta-\ga_2-\gc_1) +2{\bf{\overline{p}}}_{V(\ga)}(\eta-\gc_1-\gc_2)}.\ee
The exponent on $G_3$ here is the same as the exponent on $D_3$ in \eqref{yew4}, however the product is over $C_X$ in the former 
and $C_X^+$  in the latter. It is worth noting that $D_2$ and $G_2$ are the only terms in Theorem \ref{shapdet} and Lemma \ref{gro}
that involve odd roots.  Also the sets $C_X$ and $E_X$ are disjoint. This means that in checking the multiplicities of the factors of the product in Theorem \ref{shapdet} (as well as Lemmas \ref{gro} and \ref{xyz}), we can study 3 cases separately. Note that we need to consider $D_1$ and $G_1$ in both Cases 2 and 3. 
\bi \item Case 1. Factors of $D_2$ and $G_2$. 
\item Case 2. Factors of $D_0'', D_1, D_4$, $G_1$ and $ G_4$.
\item Case 3. Factors of $D_0', D_1, D_3$, $G_1$ and  $G_3$. \ei
}\er
\el
\noi {\it Proof of Lemma \ref{gro}.}
We adapt the proof of \cite{M} Lemma 10.1.3. The leading term of $\det G^X_{\eta}$  is the same as the leading term of the product of the diagonal entries, and this equals

\be \label{why} \prod_{\ga  \in(\Gd^+\; \backslash \; X)}\;
\prod_{\pi \in {\bf \overline{P}}_{X}(\eta)} h_{\ga }^{\pi(\ga )}.\ee
We have to consider various possibilities for roots
$\ga \in \Gd^+\; \backslash \; X$.
Suppose first that $\ga \in \ovd^+_{0}$.
If $\ga\notin C_X$ and neither \eqref{cad} or \eqref{had} hold. 
then
the multiplicity of $h_\alpha$ is $\sum^\infty_{r=1}{\bf p}_{X}(\eta - r\alpha)$ exactly as in \cite{M} Lemma 10.1.3,
while if  \eqref{had} holds the multiplicity is
 \be \label{add}\sum^\infty_{r=1}{\bf p}_{X}(\eta - r\alpha_1)+
\sum^\infty_{r=1} {\bf p}_{X}(\eta - r \alpha_2).\ee
This accounts for the exponent in the inner product in \eqref{why}, and gives the term $G_1$.
 However if $\ga\in C_X$ and $ \gc = \Gc(\ga)$, then since $h_\ga  \equiv \pm h_{{s_\ga \gc} }$
 we get the additional contribution coming from $G_3.$ This deals with Case 3, and we turn to Case  2. If
$(\ga_1,\ga_2) \in E_{X}$ we need to add
\[{\bf{\overline{p}}}_{V(\ga)}(\eta-\ga_1-\gc_1 )+ {{\bf{\overline{p}}}}_{V(\ga)}(\eta-\ga_2-\gc_1) 
+2{{\bf{\overline{p}}}}_{V(\ga)}(\eta-\gc_1-\gc_2)\] to \eqref{add} to get the multiplicity of $h_\ga$, whence the term $G_4$.
Finally Case 1 is straightforward: if
$\ga  \in B_X $ as in $G_2$, the multiplicity of $h_{\ga }$ in the leading term of  $\det G^X_{\eta}$
is
\[\sum_{\pi \in {\bf \overline{P}}_{X}(\eta)} \pi(\ga) =
{{\bf p}_{X\cup\{\ga \}}(\eta - \ga )}.
\]
\hfill  $\Box$
\bl \label{xyz} Modulo the ideal defining $\cH_X$,  we have
$|\det G^X_{\eta}:h| = |D_0D^X_{\eta}:h|$ for all $h \in\fh.$
 \el
\bpf If $\ga\in \ovd^+_{0}
$ and $s_\ga X = X$, it is easily checked that  $ h_{\ga}$ has the same multiplicity in $G_\eta$ and $ D_\eta$.
Next if ${\ga  \in B_X }$, that is $\ga$ is an isotropic root which is orthogonal to $X$,  then $ h_{\ga}$ has the same multiplicity in the leading terms of  $G_2$ and $ D_2$. Indeed $G_2$ {\it is} the leading term of $D_2$.
The multiplicity of $h_\ga$ in other terms is zero.
\\ \\
Next consider the case where $\ga\in  C_X $, that is Case 3 in  Remark \ref{dab}. Suppose that
$\gc = \Gc(\ga)\in X,$ but $\gc'=s_\ga \gc \notin X$.  Set $Z=X\cup \{\gc'\}$ and $\eta_0 =\eta-\gc$. Then
using  Lemma \ref{ink} we obtain telescoping sums
\by \label{stu}|G_1 :h_\ga| -| D_1 :h_\ga| &=&\sum_{r\ge 1}{{\bf p}_{X}(\eta - r \alpha)}-{{\bf p}_{s_\ga X}(\eta - r \alpha)}\nn\\
&=& \left\{
\begin{array}{cr}
{\bf p}_Z(\eta_0)& \mbox{if } {\ga \in C^+_X}, \\
-{\bf p}_Z(\eta_0 - \alpha) & \mbox{otherwise}.
\end{array}
\right.  \ey
Next note that if $\gc\in X$, and $s_\ga \gc\notin X$ then
\by  \label{G4}|G_3 :h_\ga| &=&
\sum_{\pi \in {\bf \overline{P}}_{X}(\eta)} {\pi({s_\ga \gc} )} \nn\\
&=&|{\pi \in {\bf \overline{P}}_{X}(\eta)}|{\pi({s_\ga \gc} )}=1|\;=\;{\bf p}_Z(\eta- {s_\ga \gc})\nn\\
&=& \left\{
\begin{array}{ll}
{\bf p}_Z(\eta_0 + \alpha)\;\;\; \mbox{if} & {\ga \in C^+_X, \gc = \Gc(\ga)}\\
{\bf p}_Z(\eta_0 - \alpha) & \mbox{otherwise}.
\end{array}
\right.  \ey
Also if $\ga \in C^+_X$ and $\gc = \Gc(\ga)$ then by  Lemma \ref{ink} and  \eqref{yaws},
\be \label{G5} | D_3 :h_\ga|=  {\bf p}_Z(\eta_0 + \alpha)-{\bf p}_Z(\eta_0) \mbox{ and  } | D_0' :h_\ga|= 2{\bf p}_Z(\eta_0). \ee
Therefore using \eqref{stu}, \eqref{G4} and \eqref{G5} we have
\by && |G_1 :h_\ga| -|D_1 :h_\ga| +|G_3 :h_\ga|-|D_0' :h_\ga| -|D_3 :h_\ga|\nn\\ &=& {\bf p}_Z(\eta_0) +{\bf p}_Z(\eta_0 + \alpha) -2{\bf p}_Z(\eta_0)+ {\bf p}_Z(\eta_0)-{\bf p}_Z(\eta_0 + \alpha)=0. \nn \ey
showing that  $| D_\eta:h_\ga|= |G_\eta:h_\ga|.$ The case where $\ga \in C_X \backslash C^+_X$ is  similar, but easier since the multiplicity of $h_\ga$ in $D_0'$ and $D_3$ is zero.  \\ \\
\noi Now consider the case of a pair $(\ga_1,\ga_2) \in E_{X}$, that is Case 2 in Remark \ref{dab}.
Set $\ga=\ga_1, Y = s_{\ga}X$ and $V=X\cup Y$.
We claim that 
\be \label{cri} | D_1/G_1 :h_\ga|  =  -2{{\bf p}_{V}(\eta - \gc_2)}+\sum_{i=1}^2{\bf p}_{V}(\eta - \alpha_i-\gc_1).\ee
Indeed by \eqref{yew1} and the expression for $G_1$ in Lemma \ref{gro},

\[ | D_1/G_1 :h_\ga|  =
\sum^\infty_{r=1} \sum_{i=1}^2{{\bf p}_{Y}(\eta - r \alpha_i)}-
            {{\bf p}_{X}(\eta - r \alpha_i)}\nn \]
Next using Lemma \ref{ink},
\by | D_1/G_1 :h_\ga| =\sum^\infty_{r=1} \sum_{i=1}^2 {\bf p}_{V}(\eta - r\alpha_i-\gc_1)
           &+& 
	\sum^\infty_{r=1} \sum_{i=1}^2 {\bf p}_{V}(\eta - r\alpha_i-\gc_2)\nn\\
-\sum^\infty_{r=1} \sum_{i=1}^2 {\bf p}_{V}(\eta - r\alpha_i- \ga_1 -\gc_1)&-& \sum^\infty_{r=1} \sum_{i=1}^2 {\bf p}_{V}(\eta - r\alpha_i-\ga_2-\gc_1).\nn\ey
Introducing $\gc_2 =\ga_1 +\gc_1 +\ga_2$ into the second line above, we obtain
the following, which evidently collapses to the right side of \eqref{cri}
\by | D_1/G_1 :h_\ga| =\sum^\infty_{r=1} \sum_{i=1}^2 {\bf p}_{V}(\eta - r\alpha_i-\gc_1)
           &+& 
	\sum^\infty_{r=1} \sum_{i=1}^2 {\bf p}_{V}(\eta - r\alpha_i-\gc_2)\nn\\
-\sum^\infty_{r=1} \sum_{i=1}^2 {\bf p}_{V}(\eta - (r+1)\alpha_i-\gc_1)&-& \sum^\infty_{r=1} \sum_{i=1}^2 {\bf p}_{V}(\eta - (r-1)\alpha_i-\gc_2).\nn\ey
Taking into account the contributions given by $D_0''D_4$
from \eqref{yew}, and $G_4$ from Lemma \ref{gro}
we obtain \be| D_1/G_1 :h_\ga|  +|D_0''D_4:h_\ga| -|G_4 :h_\ga|=0
\label{hope}\ee

\epf

\bc \label{tcp} To complete the proof of Theorem \ref{shapdet}, it suffices to show that $\det F^X_{\eta}$ is divisible by each of the factors $D_1,\ldots,D_4$ defined in Equations \eqref{yew1}-\eqref{yew5}.\ec

\bpf Immediate.  Since the leading terms match there can be no more factors.\epf
\bc \label{fara} For all
roots $\ga$, $|D^X_{\eta}:h_\ga|\ge |\det F^X_{\eta}:h_\ga|$.
\ec  \bpf Combine Lemma \ref{xyz} with  \eqref{(c)}. \epf

\subsection{Parabolic Induction and Twisting Functors.}\label{agm}
\subsubsection{Introduction.}
Let $X$ be an orthogonal set of isotropic roots. The aim of this Subsection is to construct certain modules that appear in Jantzen sum formula Theorem \ref{Jansum101}.  
In  the next Subsection we consider a
special case where the modules can be constructed by parabolic induction.  
In Subsection \ref{nmvtf} 
the general  result will be deduced using  twisting functors.
\\ \\
Our next goal is to show that if $[\ga] =(\ga_1,\ga_2)\in E_X(\gl)$ 
there is a module  in the category $\cO$ with character $\tte^\gl p_{W(\ga)}$. 
\ff{See \eqref{efg} for notation.}
\ff{This is the first and only time in this paper that we consider 
modules with a prescribed character of the form $\tte^\gl p_Y$ where $Y$ contains even roots.}
These modules play an important role in the Jantzen sum formula for $M^X(\gl)$. 
We briefly recall what this means.
The hypothesis $[\ga] \in E_X$ means that there exist odd roots $\gc_1, \gc_2\in X$ and
even roots $\ga_1, \ga_2$, such that  $\ga_1^\vee \equiv \ga_2^\vee \mod \Q \gc_1 +\Q \gc_2$.  This happens when 
$\gc_2 = \ga_1 +\gc_1+\ga_2$. Using the fact that $X$ is orthogonal, it is easy to see that if $\ga_3$ is any other positive even root
such that $\ga_1^\vee \equiv \ga_3^\vee \mod \Q X$, then $\ga_3=\ga_2$, that is either $\ga_1$ or $\ga_2$    in $[\ga] =(\ga_1,\ga_2)\in E_X(\gl)$ determines the other.  
The additional hypothesis on $\gl$ is that $(\gl+\gr,\ga_1^\vee) =1$. 
Let $\fk[\ga]$ be the $\fgl(2,2)$ subalgebra associated to $[\ga]$, and $W'(\ga)$ 
(resp. $V'(\ga)$) be the set of all  (resp. all odd) positive roots of $\fk[\ga]$. 
We assume that $X'$  is an orthogonal set of isotropic roots which is orthogonal to $W'(\ga)$ and there are disjoint unions  
$X=V'(\ga)\ds X' $ and  $W(\ga)=W'(\ga) \ds X'.$
\\ \\
In the special case where $\fk[\ga]=\fk$ 
 is the subalgebra $\fk$ from  \eqref{str} with $(k,\ell)=(2,2)$, and $X'=\{\gep_{k+i}- \gd_{\ell+i}\}_{i=1}^{r}$ the result follows from Proposition  \ref{propmain}. In general we need to
prove the  result when $[\ga]$ and $\gl$ are replaced by $v[\ga] =(v\ga_1,v\ga_2)\in E_{vX}(v\cdot\gl)$ for $v\in W$.  
\ff{ Note that in Type A, $v\gr_1=\gr_1$ and hence $v\cdot\gl= v\circ\gl$.}
 For this we use twisting functors. The result is as follows.

\bt 
\label{con1} Suppose $\fk[\ga]$ is as above, ${[\ga] \in E_{X}(\gl)}$, 
and $v\in W$ is such that $N(v)$ is disjoint from $\Gd^+(\fk)$. 
Then the module 
$T_v \Ind_{\fp}^{\fg} \;\ttk_\gl$ has  character $\tte^{v\cdot\gl} p_{vW(\ga)}$.
\et
\bpf This follows from the more general Theorem \ref{con2} below. \epf
\subsubsection{Parabolic Induction.}\label{pari}
In the next Subsection, we use twisting functors to construct some new modules that we need for the Jantzen sum formula. 
The underlying reason for the existence of these modules is the fact that $\fgl(1,1)$ has many one dimensional modules. Here we study a
special case where the modules can be constructed by parabolic induction.  
Suppose that $\fk,\fl$ are subalgebras of $\fg$, such that $\fk$ is isomorphic to $\fgl(k,\ell)$ and for some $r,t_1, t_2$
\be\label{str}\fl = \fk \op \fgl(1,1)^r\op\fgl(1|0)^{t_1}\op \fgl(0|1)^{t_2}.\ee
\noi 
 First we explain how $\fl$ is embedded in $\fg=\fgl(m,n)$.
To do this it is convenient to introduce a certain partial flag in $\ttk^{m|n}$. This will also allow us to define some 
important subalgebras of $\fg$.
Let $e_1, \ldots, e_m,$ $e_{1'}, \ldots ,e_{n'}$ be the standard basis for $\ttk^{m|n}$
and consider the diagram below. The $i^{th}$ set of numbers (for $i\ge 0$) separated by a pair of
vertical line is used as the index set for a $\Z_2$-graded subspace $\vw^{(i)}$ of $\ttk^{m|n}$.
\begingroup

\vspace{0.55cm}


\setlength{\unitlength}{0.14in}
\begin{picture}(-15,-10)(18.2,0) \label{fig2}
\put(17.6,-1){\line(0,1){3}}
\put(23.6,-1){\line(0,1){3}}
\put(23.6,-1){\line(0,1){3}}
\put(18,-1){$1',2',\ldots ,\ell'$}
\put(18,1){$1,\;2,\;\ldots ,k$}
\put(24,-1){$(\ell +1)'$}
\put(24,1){$\;\;k+1\;$}
\put(27.6,-1){\line(0,1){3}}
\put(28.5,0){$\cdots$}
\put(30.6,-1){\line(0,1){3}}
\put(31,-1){$(\ell +r)'$}
\put(31,1){$\;k+r\;$}
\put(34.6,-1){\line(0,1){3}}
\put(35.2,0){$k+r+1$}
\put(39.7,-0.4){\line(0,1){1.5}}
\put(40.5,0){$\cdots$}
\put(42.5,-0.4){\line(0,1){1.5}}
\put(43.1,0){$m$}
\put(44.3,-0.4){\line(0,1){1.5}}
\put(44.9,0){$(\ell +r+1)'$}
\put(51.2,0){$\cdots$}
\put(50.5,-0.4){\line(0,1){1.5}}
\put(53.8,0){$n'$}
\put(53.2,-0.4){\line(0,1){1.5}}
\put(55,-0.4){\line(0,1){1.5}}
\end{picture}

\vspace{0.6cm}

\endgroup

\noi We assume that $\fl$ contains $\fh$ (equivalently $m= k+r+t_1$ and $n=\ell+r+ t_2$).  This gives $t+1$ subspaces (where $t:=m -k + t_2$), which we number as
 \[\vw^{(0)} =\span\{e_1,\ldots, e_k,e_{1'},\ldots, e_{\ell'}\},\;
\vw^{(1)} =\span\{e_{k+1}, e_{(\ell+1)'}\}, \ldots, \vw^{(t)}=\ttk e_{n'}.\]
\noi Now set
$\vv^{(-1)}=0,$  and $\vv^{(i)} = \vw^{(0)}\op\ldots \op \vw^{(i)}$ for $0\le i\le t$.
In this notation, let
\[\fk =  \{x\in \fg|x\vw^{(0)}\subseteq \vw^{(0)},
\; x\vw^{(i)} =0 \mbox{ for } i>0\},\]
\[\fp = \{x\in \fg|x\vv^{(i)}\subseteq \vv^{(i)} \mbox{ for all } i\}, \quad
\fl = \{x\in \fg|x\vw^{(i)}\subseteq \vw^{(i)} \mbox{ for all } i\}.\] 
Next define
\[\fv_0 = \{x\in \fg_0|x\vw_j^{(i)}\subseteq \vv_j^{(i-1)},\;\;\mbox{ for all } i \mbox{ and } j=1,2\},\]
\[\fv_1^+ = \{x\in \fg|x\vw_1^{(i)}\subseteq \vv_0^{(i-1)},\;\;x\vw_0^{(i)} = 0 \mbox{ for all }  i\},\]
\[\fv_1^- = \{x\in \fg|x\vw_0^{(i)}\subseteq \vv_1^{(i-1)},\;\;x\vw_1^{(i)} = 0 \mbox{ for all } i\},\]
and then set $\fv_1= \fv_1^+\op\fv_1^-$, $\fv=\fv_0\op\fv_1$.  Note that $\fp$ is the subalgebra of $\fg$
stabilizing the partial flag
\[ 0 \subseteq \vv^{(0)} \subseteq \vv^{(1)} \subseteq \ldots \subseteq \vv^{(t) },\]
and $\fv$ is the ideal in $\fp$ consisting of all elements that act nilpotently on the flag.
We have $\fp=\fl\op \fv$. Also $\fk\cong \fgl(k,l)$. If $\fc$ is an ad-$\fh$ stable subspace of $\fg$, let $\Gd(\fc)$ be  the set of roots of $\fc$. 
The set of even and odd roots in $\Gd(\fc)$ is denoted by $\Gd_0(\fc)$ and $\Gd_1(\fc)$ respectively. 
Let $\fw$ be the ad-$\fh$ stable subspace of $\fg$ with $\Gd(\fw) = - \Gd(\fv)$.  Then 
$\fg= \fw \op \fp$. 
\bexa \label{54}{\rm The first diagram yields a flag in $\ttk^{5|4}$  The second diagram then shows the resulting subspaces of $\fg=\fgl(5,4)$, constructed  using the above recipes.}\eexa
\begingroup

\vspace{0.3cm}


\setlength{\unitlength}{0.14in}
\begin{picture}(-15,-10)(18.2,1) \label{fig101}
\put(27.6,-1){\line(0,1){3}}
\put(31.0,-1){\line(0,1){3}}
\put(33.5,-1){\line(0,1){3}}
\put(36.0,-0.4){\line(0,1){1.5}}
\put(38.5,-0.4){\line(0,1){1.5}}
\put(41.0,-0.4){\line(0,1){1.5}}
\put(29.2,-1){$1'$}
\put(28.5,1){$1,\;2$}
\put(32,-1){$2'$}
\put(32,1){$3$}
\put(34.6,-1){$3'$}
\put(34.6,1){$4$}
\put(37.1,0){$5$}
\put(39.6,0){$4'$}
\end{picture}

\vspace{0.7cm}
\endgroup

    \[ \begin{tabular}{|c|c|c|c|c||c|c|c|c|}\hline
$    \fk$ &$\fk$ & $\fv_0$ & $\fv_0$ & $\fv_0$ &$    \fk$ &$\fv_1^+$ & $\fv_1^+$ & $\fv_1^+$ \\ \hline 
$    \fk$ &$\fk$ & $\fv_0$ & $\fv_0$ & $\fv_0$  & $    \fk$ &$\fv_1^+$  & $\fv_1^+$ & $\fv_1^+$ \\ \hline  
$    0$ &0 & $\fgl(1,1)^{(1)}$ & $\fv_0$ & $\fv_0$  & 0&$\fgl(1,1)^{(1)}$ & $\fv_1^+$ & $\fv_1^+$ \\ \hline 
$    0$ &0& 0& $\fgl(1,1)^{(2)}$ & $\fv_0$& 0&0 & $\fgl(1,1)^{(2)}$  & $\fv_1^+$ \\ \hline 
$    0$ &0& 0& 0 & $\fgl(1,0)$  & 0&0 & 0 & $\fv_1^+$  \\ \hline \hline
$\fk$ &$\fk$ & $\fv_1^{-}$ & $\fv_1^{-}$ & $\fv_1^{-}$  & $\fk$&$\fv_0$ & $\fv_0$ & $\fv_0$ \\ \hline
0&0& $\fgl(1,1)^{(1)}$ & $\fv_1^{-}$ & $\fv_1^{-}$  &   0 &$\fgl(1,1)^{(1)}$ & $\fv_0$ & $\fv_0$ \\ \hline   
0&0& 0 & $\fgl(1,1)^{(2)}$  & $\fv_1^{-}$  & 0&0 & $\fgl(1,1)^{(2)}$ & $\fv_0$ \\ \hline 
0&0& 0& 0& 0& 0&0& 0& $\fgl(0,1)$ \\ \hline
    \end{tabular}  \]
Non-zero entries in the second diagram correspond to the non-zero matrix entries in $\fp$. Superscripts are used to distinguish the two copies of $\fgl(1,1).$
\\ \\		
\noi Now set  $\vw = \vw^{(0)}\op \vw^{(r+1)}\op\ldots \op \vw^{(t)}$, and 
\[\fu^+ = \{x\in \fg|x\vw_1^{(i)}\subseteq \vw_0^{(i)},\;\;x\vw_0^{(i)} = 0 \mbox{ for } i \in [r], \;\;x\vw=0 \},\]
\[\fu^- = \{x\in \fg|x\vw_0^{(i)}\subseteq \vw_1^{(i)},\;\;x\vw_1^{(i)} = 0 \mbox{ for } i \in [r], \;\;x\vw=0 \}.\]
If $\fu =\fu^+\op\fu^-$, we have  $\fl=\fk\op\fu+\fh$.
Note that the entries in $\fu^+$ (resp. $\fu^-$) occur in the copies of $\fgl(1,1)$ above (resp. below) the main diagonal. 
The set of weights of $\fu^+$ is 
\be \label{lip} \{\gep_{k+j} - \gd_{(l+j)'}|j \in [r]\}. \ee
Next define
 \[Z_0 = \Gd^+(\fv_0), \quad Z_1^\pm=\Gd(\fv_1^\pm), \quad Z_1 = Z_1^+\cup Z_1^-, \quad \widehat{Z}_1 = Z_1^+\cup -Z_1^- .\]

\bl \bi \itemo \itema If $L_0={\prod_{\alpha \in \Gd_0^+}} (1 -\tte^{ -\alpha})$, we have 
	\be  \label{abc} {\prod_{\alpha \in \Gd({\fl_0^+})} (1 - \tte^{ -\alpha})}{\prod_{\alpha \in Z_{0}} (1 - \tte^{ -\alpha})} = L_0.\ee
\noi  
\itemb Set 
$\gs = \sum_{\ga\in Z^-_1}\ga$. Then 
\be \label{lem} \prod_{\alpha \in Z_{1}} (1 +\tte^{- \alpha})=\tte^{-\gs}\prod_{\alpha \in \widehat{Z}_1}(1 +\tte^{- \alpha}).\ee
\itemc If $\gr(\fb')$ is the analog of $\gr$ for $\fb'$ we have
$\gr(\fb')=\gr + \gs. $
\itemd There is a unique Borel subalgebra $\fb'$ contained in $\fp$ such that 
\[ d(\fb,\fb') =\min\{d(\fb,\fb'')| \fb''  \mbox{ is a Borel subalgebras contained in } \fp\}\]
\iteme We have $\Gd^+(\fb') = Z_1\cup \{\gep_{k+j} - \gd_{(l+j)'}|j \in [r]\}.$ 
\ei
\el
\bpf  For (a) note we have a disjoint union $\Gd_0^+ = \Gd(\fl^+_0) \ds Z_0,$ from which
  \eqref{abc} follows. Similarly (b) holds since $Z_{1}^+$ is common to $Z_{1}$ and  
$\widehat{Z}_1$.  For (c) note that $Z_1^-=\Gd(\fv_1^-)$ is precisely the set of odd roots of $\fb'$ which are not roots of $\fb^\dist.$ We leave (d), (e) to the reader.
\epf
\noi Note that  
\be \label{rad} (\gl+\gr(\fb^\dist),X) = (\gl'+\gr(\fb'),X) = 0.\ee
 \noi The following lemma is straightforward.
\bl \label{oil} If $\gl\in \fh^*$, then $\gl$ defines a  one dimensional representation of $\fl$
if and only if $(\gl+\gr,\gc)=0$  for all odd $\gc \in X,$ 
  and  $(\gl,\ga)=0$ for each even simple root $\ga$ of $\fl$. \el
\bp \label{propmain} 
Assume that $\gl\in \fh^*$ defines a  one dimensional representation $\ttk_\gl$ of $\fl,$ and regard  $\ttk_\gl$ 
as a $\fp$-module by allowing $\fv$ to act trivially. Let $X=\Gd^+(\fl)$. Then 
\bi
\itema The module $\Ind_{\fp}^{\fg} \ttk_\gl$ has character 
\be \label{sot} \ch \Ind_{\fp}^{\fg} \;\ttk_\gl   =\frac{{}\prod_{\alpha \in Z_{1}} (1 +\tte^{- \alpha})}
{\prod_{\alpha \in Z_{0}} (1 -\tte^{ -\alpha})}\;\tte^\gl 
= \frac{{}\prod_{\alpha \in \widehat{Z}_1} (1 +\tte^{- \alpha})}{
\prod_{\alpha \in Z_{0}} (1 -\tte^{ -\alpha})}\;\tte^{\gl-\gs}=\tte^\gl p_X.\ee 
\itemb
The module $\Ind_{\fp_0}^{\fg_0}\;\ttk_\gl $ has character 
\by \label{wed} \ch \Ind_{\fp_0}^{\fg_0}\;\ttk_\gl 
&=&
{\prod_{\alpha \in Z_{0}} (1 -\tte^{ -\alpha})}^{-1}\;\tte^\gl  = 
{\prod_{\alpha \in \Gd({\fl_0^+})} (1 -\tte^{ -\alpha})}\tte^\gl/L_0. \ey 
\ei
 \ep

\bpf (a) Let $\fw$ be the $\ad \fh$-stable subalgebra of $\fg$ such that $\Gd(\fw)=-\Gd(\fv)$ and $\fg=\fw\op \fp$.
By \cite{M} Corollary 6.1.5, if $S(\fw)$ is the supersymmetric algebra on $\fw$, we can write $U(\fg) = S(\fw)\ot U(\fp)$ as right $U(\fp)$-modules. 
Now $\Gd(\fw)= Z_0 \cup Z_1$, Hence the first equality in 	(a) since $S(\fw)$ is the tensor product of the exterior algebra on $\fw_1$ with the symmetric algebra on $\fw_0$. 
The second equality comes from \eqref{lem}, and the third is an easy consequence.
The proof of the first equality in 	(b) is similar, and the second follows from \eqref{abc}.	
	\epf
\noi Let $X=\Gd^+_1(\fl)$. We note the following general behavior of $\Ind_{\fp}^{\fg} \;\ttk_\gl   $.
\bl \label{clo} If $B(\gl)=X$, then $M=\Ind_{\fp}^{\fg} \;\ttk_\gl$  is a highest weight module for the distinguished Borel subalgebra   of $\fg$ with highest weight $\gl'=\gl-\gs$. \el
 \bpf    Consider the Borel subalgebra $\fb$ with the same even part as
the distinguished Borel subalgebra $\fb^\dist$ of $\fg$, and with
$$\fb^\dist_1 \cap \fb_1=(\fb^\dist_\fk)_1\op\fu^+\op\fv_1,$$ where $\fb^\dist_\fk$ is the distinguished Borel subalgebra of $\fk$.
 Then $\fb \subseteq \fp$, so  $M$ is a highest weight module for $\fb$ with highest weight $\gl$.
Also there is a  sequence
 \be \label{1distm1} \mathfrak{b} = \mathfrak{b}^{(0)}, \mathfrak{b}^{(1)}, \ldots,
 \mathfrak{b}^{(u)}. \ee
  of Borel subalgebras of $\fg$, such that $\mathfrak{b}^{(i-1)}$ and
  $\mathfrak{b}^{(i)}$ are adjacent for $1 \leq i \leq u$, and $\mathfrak{b}^{(u)} = \mathfrak{b}^\dist$. The assumption $B(\gl)=X$ ensures that each change of Borels from $\mathfrak{b}$ to $ \mathfrak{b}^\dist$ is typical for $M$, see  \eqref{1mon}.
Now $M$ is  generated by a highest weight vector $v_\gl$ with respect to $\fb$, 
and there are odd root vectors $e_{-\ga_1},\ldots, e_{-\ga_u}$ such that $m=e_{-\ga_u} \ldots e_{-\ga_1}v_\gl$ is a highest weight vector 
for $\fb^\dist$. The assumption $B(\gl)=X$ also implies that $m$ generates $M$ and the result follows.
\epf \noi  
\subsubsection{New Modules via Twisting Functors.}\label{nmvtf}
Now consider any $\fh$-stable subalgebra $\fl'$ of $\fg=\fgl(m,n)$ such that $\fl' \cong \fk \op \fgl(1,1)^r$, with $\fk$ a $\fgl(k, \ell)$-subalgebra of $\fl$. By abuse of notation we write $\fgl(1,1)^r$ for the isomorphic copy of this subalgebra inside $\fl'$.  We denote the set of even (resp. odd) simple roots of $\fk$ by $\Pi_0(\fk)$, (resp. $\Pi_1(\fk)$) and do likewise for  
$\fl'$.
\bl We have 
\bi
\itema $$\Pi_0(\fk)=\{\gep_{s_i}- \gep_{s_{i+1}}\}_{i=1}^{k-1} \cup \{\gd_{t_i}- \gd_{t_{i+1}}\}_{i=1}^{\ell-1},$$
 $\Pi_1(\fk)=\{\gep_{s_k}- \gd_{t_{1}}\}$, $\Pi_0(\fl')=\Pi_0(\fk)$, and 
 $\Pi_1(\fl')=\Pi_1(\fk)\cup\{\gep_{c_i}- \gd_{d_{i}}\}_{i=1}^r,$ where $1\le s_1<s_2<\ldots s_k\le m$, $1\le t_1<t_2<\ldots t_\ell\le n$ 
and $c_i \notin \{{s_i}\}_{i=1}^k$, $d_i \notin \{{t_i}\}_{i=1}^\ell$ 
\itemb Let $i' = n+i$ for $i\in [\ell]$. Suppose $w\in W=S_m \ti S_n$ satisfies 
$w(k+i) = c_i$ for  $i\in[r]$, $w(\ell+i)' = d_i'$ for $i\in[r]$
$wi = s_i,$ for $i\in [k]$, $wi' = t_i$ $i\in [\ell]$. Then 
$w\Gd^+(\fl)= \Gd^+(\fl')$ 
and $N(w)\cap \Gd^+_0(\fl)=\emptyset.$ \ei\el
\bpf Left to the reader. \epf

\bt 
\label{con2} If $\gl$ defines a one-dimensional representation of $\fp$, then
 the module 
$T_w \Ind_{\fp}^{\fg} \;\ttk_\gl$ has  character $\tte^{w\cdot\gl} p_{\Gd^+(\fl')}$.
\et
\bpf See \cite{M101}. \epf
\subsection{An Equivalence of Categories and Generic Cases.}\label{lab}
\subsubsection{The Equivalence of Categories.}\label{gcvtf}
In \cite{CMW}, twisting functors are used to show that every block of the category $\cO$ for a Lie superalgebra in Type A is equivalent to an integral block of some (possibly smaller) Lie superalgebra of Type A. We show how this result can be used to reduce the study of the modules $M^X(\gl)$ in the most general cases  to the study of $\fgl(2,1)$ and $\fgl(2,2).$
The main result of \cite{CMW} states that
\begin{theorem}\label{thmmain}
Every block $\cO_\gl$, $\gl\in \fh^*$, is equivalent to an integral
block of the category $\cO$ for some direct sum of general linear Lie superalgebras.
\end{theorem}
\noi The weight $\gl$ determines a subroot system $\Gd_\gl$ of $\Gd$ with Weyl group
$W_\gl$ which is the subgroup of the Weyl group $W$ generated by
${s_\alpha}$ with $\alpha$ an even root of $\Gd_\gl$.
 The proof in \cite{CMW} shows that there is a Lie superalgebra
$\fg_{[\widehat{\gl}]}$ and an integral
weight  $\widehat{\gl}$ for $\fg_{[\widehat{\gl}]}$ such that $\mathcal O_\gl$
is equivalent to a block $\cO_{\widehat{\gl}}^\fl$
of the category of $\cO^{\fl}$ for $\fl = \fg_{[\widehat{\gl}]}$.
The block $\cO^{\fl}_{{\widehat{\gl}}}$ is the Serre subcategory of
$\cO^{\fl}$ generated by simple objects of the form $L_{\fl}(\mu)$
with $\mu=w\cdot(\gl-\sum_{j}k_j\alpha_j)$, for  $w\in W_{\gl}$ and
$\{\alpha_j\}$ a set of mutually orthogonal odd isotropic roots
satisfying $(\gl+\rho_{\fl},\alpha_j)=0$  and $k_j\in\Z$.
Under this equivalence,  the Verma module $M(w\cdot \gl)$ corresponds to
the Verma module $M_\fl(w\cdot \widehat{\gl})$ for $\fl$ with highest weight $w\cdot \widehat{\gl}$, and the simple module
$L(w\cdot\gl)$ to $L_\fl(w\cdot\widehat{\gl})$ for all $w\in W_\gl$.
In this set-up  $\fl$  is a Levi subalgebra of $\fg$ containing a Cartan subalgebra $\fh$ of $\fg$.

\bp \label{ec} Suppose that $B(\gl) = X$ and $\fl = \fg_{[\widehat{\gl}]}$ is as in \eqref{str}. 
\bi \itema There is an equivalence of categories between the block $\cO_{\widehat{\gl}}^\fl$ and the block $\cO_{\widehat{\gl}}^{\fk+\fh}$.
\itemb There is an  equivalence of categories  $\cO_{\widehat{\gl}}^{\fk+\fh}\lra \cO_{{\la}}$obtained by combining $(a)$ with  Theorem {\rm \ref{thmmain}}.  This equivalence  sends $N$ to
$\Ind^\fg_\fp N$.
\ei
\ep
 \bpf  The key points are as follows: any simple $\fl$-module is the tensor product of a simple $\fk$-module and simple modules for the copies of $\fgl(1,1), \fgl(1,0)$ and $\fgl(0,1)$ in \eqref{str}. From the definition of $\cO_{\widehat{\gl}}^\fl$ as a Serre subcategory and the fact that $X'' \subseteq B(\gl)$ we see that the only simple $\fgl(1,1)$-modules that can appear are one dimensional, and are annihilated by $\fu$. Such modules arise canonically from simple modules in   $\cO_{\widehat{\gl}}^{\fk+\fh}$. 
This proves (a) and (b) holds by construction.\epf 

\subsubsection{Lessons learned from the $\fgl(2,2)$ case: Roots in $E_X$.}
\noi Now suppose that $X$ is an orthogonal set of isotropic roots and
$(\ga_1,\ga_2)\in E_{X}$.  Suppose that $\ga_1 +\gc_1 +\ga_2 =\gc_2$ where $\gc_1, \gc_2\in X$.
Combining the equivalence of categories from Theorem \ref{thmmain}, with Proposition \ref{propmain} we deduce most of the following results from the $\fgl(2,2)$ case, see Theorems 
\ref{genst}, \ref{stz1} and \ref{mnt}. However, the fact that \v Sapovalov  elements give highest weight vectors is of course true in general. Also the statements about the Jantzen filtration hold because the modules in question are multiplicity free, 
and each quotient $M^X_i(\gl)/M^X_{i+1}(\gl)$ is contragredient with composition factors of multiplicity one, and therefore semisimple. 
The phrase ``for general $\gl$" means that, up to $W$-conjugacy, $\fl = \fg_{[\widehat{\gl}]}$ is as in \eqref{str} with $\fk =\fgl(2,2)$, and we assume (without any loss of generality) that
\[\Gd_0^+(\fk) =\{\ga_1,\ga_2\}, \; B(\gl) = X \mbox{ and } \{\gc_1,\gc_2\} \subseteq \Gd_1^+(\fk)\cap X.\]
Note that the number of $\fgl(1,1)$ terms in \eqref{str} is $r=|X|-2$.
\bt\label{wane} Let $n$ be an integer greater than 1. For general $\gl \in \cH_X\cap \cH_{\ga_1,n}$ where $(\ga_1,\ga_2) \in E_{X}$,
\bi \itema The socle $\cS$ of $M^X(\gl)$ is  $M^X(s_1s_2\cdot \gl) \cong L(s_1s_2\cdot \gl).$
\itemb We have
\[U(\fg)\gth_{\ga_1,n}v_\gl \cong M^Y(s_1\cdot \gl), \quad U(\fg)\gth_{\ga_2,n}v_\gl \cong M^Y(s_2\cdot \gl)\] 
\[M^X_1(\gl)=U(\fg)\gth_{\ga_1,n}v_\gl+U(\fg)\gth_{\ga_2,n}v_\gl\]
and \be \label{tan} \cS=U(\fg)\gth_{\ga_1,n}\gth_{\ga_2,n}v_\gl =U(\fg)\gth_{\ga_1,n}v_\gl\cap U(\fg)\gth_{\ga_2,n}v_\gl \subseteq M_2^X(\gl).\ee
\ei The lattice of submodules of $M^X(\gl)$ is as in Figure \ref{fig1} with $V_1=U(\fg)\gth_{\ga_1,n}v_\gl$ and $V_2 =U(\fg)\gth_{\ga_2,n}v_\gl.$
\et
\bt \label{name} For general $\gl \in \cH_X\cap \cH_{\ga_1,1}$ where $(\ga_1,\ga_2) \in E_{X}$,
\bi \itema the socle $\cS$ of $M^X(\gl)$ satisfies
\be \label{stz} \cS = L(\gl-\gc_1-\gc_2)\oplus L(s_1s_2\cdot \gl) =
U(\fg)\gth_{\ga_1,1}v_\gl\cap U(\fg)\gth_{\ga_2,1}v_\gl,  \ee
and \be \label{ton} \cS \subseteq M^X_2(\gl).\ee
\itemb We have
\[U(\fg)\gth_{\ga_1,1}v_\gl/L(\gl-\gc_1-\gc_2) \cong M^Y(s_1\cdot \gl), \quad U(\fg)\gth_{\ga_2,1}v_\gl/L(\gl-\gc_1-\gc_2) \cong M^Y(s_2\cdot \gl),\] 
\[L(s_1s_2\cdot \gl) =M^X(s_1 s_2\cdot \gl),\]
 and \[U(\fg)\gth_{\ga_1,1}v_\gl+U(\fg)\gth_{\ga_2,1}v_\gl=M^X_1(\gl).\]

\ei
\et
\br {\rm It follows from the sum formula, Theorem \ref{Jansum101}, that equality holds in \eqref{tan} and 
\eqref{ton}, that is $\cS = M^X_2(\gl)$ in both cases. In addition we have $M^X_3(\gl)=0$.}\er

\bt \label{bits} Suppose that  \eqref{had} holds, and that $m=(\gl+\gr,\ga_1^\vee)\in\Z.$ Set $Y=s_{\ga_1} X$.
 Then  for general $\gl$ in $\cH_X\cap \cH_{{\ga_1,m}}$
\bi
\itema if $m>0$ then ${M}^X(\lambda)$ has a unique proper submodule $M^{Y}(s_{\ga_1}\cdot\lambda)={M}^X_{1}(\lambda).$
\itemb If $m< 0$, then $M^X(\gl)$ has a unique proper submodule $M^{Y}(s_{\ga_2}\cdot\lambda)={M}^X_{1}(\lambda).$
\ei
\et
\subsubsection{Lessons learned from the $\fgl(2,1)$ case: Roots in $C_X$.}
\noi Next we  consider the case where $\ga\in C_X$.  Most of the following result can be deduced using twisting functors
from the $\fgl(2,1)$ case,  but
since parts of it seem to work outside type A, we give a direct proof. The highest weight vectors of 
$M^X(\gl)$ and $ M^{Y}(\gl)$ are denoted by $v^X_\gl$ and $v^{Y}_\gl$ respectively.

\bt \label{eon} Suppose $\ga \in C_{X}, \gc = \Gc(\ga)$ and set $\gc' = s_\alpha\gc, Y =s_\ga X$, ${Z(\ga)}=X \cup Y$.
\bi
\itema If $\ga\in C_X\backslash C_X^+$, then for general
$\gl$ with $ B(\gl) = X\cup Y$,  $M^X(\gl)$ is simple.
\itemb If $\ga \in C_{X}(\gl)$, we have an onto map $M^X(\gl)\lra M^{Y}(\gl)$ sending $v^X_\gl$ to $v^{Y}_\gl$.
\itemc Let ${K}^{X,\ga}({\lambda})$ be the kernel of the map $M^X(\gl)\lra M^{s_\ga X}(\gl)$.  Then
\be \label{cut} \ch {K}^{X,\ga}({\lambda}) = \tte^{\gl}({p}_X- {p}_{s_\ga X}) 
= {\tte^{\gl-\gc'}(1 - \tte^{- \ga})}{p}_{Z(\ga)}.\ee
 Equivalently
\[\dim {K}^{X,\ga}({\lambda})^{\gl-\eta} ={\bf p}_{Z(\ga)}(\eta-\gc') - {\bf p}_{Z(\ga)}(\eta-\gc).\]
\itemd For general $\gl$ such that $\ga\in C_{X}(\gl)$, ${K}^{X,\ga}({\lambda})=M_1^X(\gl)$, and this submodule is simple.
\ei
\et

\bpf We deduce (a) from the $\fgl(2,1)$ case. 
Let $\fk=\fgl(2,1)$ in \eqref{str}, and denote highest weight modules for $\fk$ by $M_\fk$.  Then switching to the notation of the Appendix, what has to be checked is that the $\fk$-module $M_\fk^\gb(-\gr)$ is simple.  However from Theorem \ref{AA2} (e), $M_\fk^\gb(-\gr)\cong M_\fk(-\gr)/V_1$ where $V_1$ is the maximal submodule of $M_\fk(-\gr)$. It is clear that the map in (b) is surjective if it is well-defined. 
The domain and target of the proposed map are both obtained as factors of $M(\gl)$ by imposing certain relations. 
To show that it is well-defined we need to show that the relations satisfied by $M^X(\gl)$ are also satisfied by  $M^{s_\ga X}(\gl)$. 
The relations needed to define $M^X(\widetilde{\lambda})_B$ and $M^{s_\ga X}(\widetilde{\lambda})_B$
are the same except that we require 
$\gth_\gc v_{\widetilde{\lambda}} = 0$ in  $M^X(\gl)_B$ and $\gth_{\gc'} v_{\widetilde{\lambda}} = 0$ in $M^{s_\ga X}(\gl)_B$. 
Now in any module in which $\gth_{\gc'} v_{\widetilde{\lambda}} = 0$ we have $\gth_{\ga,1}\gth_{\gc'}v_{\widetilde{\lambda}}=0$, 
and hence by Theorem \ref{stc}, $\gth_\gc v_{\widetilde{\lambda}} \in TU(\fg)_A v_{\widetilde{\lambda}}$. 
Reducing mod $T$ we obtain the statement in (b).
\\ \\
 The first equality in \eqref{cut} follows from (b), the second from Lemma \ref{ink}.
For (d), note that by  (a)
$M^Y(\gl)$ is simple for general $\gl$, hence ${K}^{X,\ga}({\lambda})=M_1^X(\gl)$.
Again we compare to the $\fgl(2,1)$ case. Then in the notation of Theorem \ref{AA2}, the map $M^X(\gl)\lra M^{Y}(\gl)$ is the map 
$M_\fk(\gl)/V_2 \lra M_\fk(\gl)/V_1$ and the kernel of this map is simple.
 \epf 



\subsection{Proof of Theorem  \ref{shapdet}.}\label{pfs}
We use Corollary \ref{tcp} and consider the factors $D_1,\ldots,D_4$ defined in Equations \eqref{yew1}-\eqref{yew5}. We show that the multiplicity $m(\ga,k,\eta)^X$ of $h_{\ga}-k$ in $\det F^X_\eta$ is exactly that predicted in the Theorem for each relevant root $\ga$.  To do this it is enough to show that the multiplicity of $h_{\ga}$ in the  
leading term of  $\det F^X_\eta$ is as predicted. First an easy Lemma, see \cite{KK}.
\bl \label{eek} Suppose  $M^X_1(\lambda)$ is simple and isomorphic to $M^Y(\lambda-\mu)$ for some $\mu \in Q^+$. Then
\bi \itema
  the Jantzen filtration
on $M^X(\lambda)$ takes the form
\be\label{wad} M^X_{1}(\lambda) = M^X_{2}(\lambda) \ldots= M^X_{c_k}(\lambda) \neq 0 =M^X_{c_k+1}(\lambda). \ee
\itemb the multiplicity of $h_\ga-k$ in $\det F^X_\eta$ is
$c_{k}{{\bf p}_{Y}(\eta - \mu)}.$
\ei
\el
\bpf Obviously (a) holds, and then (b) follows from  \eqref{yet}.\epf 
\noi {\it Proof of Theorem \ref{shapdet}.}  We consider 3 cases.\\ \\
\noi 
{\bf Case 1.}  First consider the factors of $D_2$. Suppose that $\gc$ is isotropic $(\gc,X)=0$. 
Suppose that  $\gl$ satisfies the conditions  of
Lemma \ref{need} (b), and use the notation of the Lemma.  Then $M^X_{1}(\lambda) = M^Y(\lambda-\gc)$ is simple.  Hence by \eqref{yet} and \eqref{wad}, the
multiplicity of $h_{\gc} + (\rho, \gc)$ in $\det F^X_\eta$
is $v_{T}(\det F^X_\eta(\widetilde{\lambda})) = a{\bf p}_{X}(\eta - \gc)$.  The computation of the leading term in Lemma \ref{gro} implies $a=1$ as claimed.
\\ \\
{\bf Case 2.} Now consider the case where \eqref{cad} holds.
Without loss we can assume that $\ga_1 +\gc_1 +\ga_2 =\gc_2$.  Suppose that $\gl_0$  is a general point in $\cH_X\cap \cH_{\ga_1,0}.$ Then if $\gl_k =\gl_0-k\gc_1,$ for $k\ge1$ we have
$(\gl_k+\gr,\ga^\vee)=k$. We use Theorems \ref{wane} and \ref{name}. 
The  socle $\cS_k$ of $M^X(\gl_k)$ satisfies $\cS_k\subseteq M^X_2(\gl_k)$.
If $k>1$, then
the lattice of submodules of $M^X(\gl_k)$
 is as in Figure \ref{fig1} with
\be \label{pup} V_1(\gl_k)=M^Y(s_1\cdot\gl_k), \;V_2(\gl_k)=M^Y(s_2\cdot\gl_k), \mbox{ and }  V_3(\gl_k)=M^X(s_1s_2\cdot\gl_k),\ee
 where $Y = s_1 X=s_2X$. Thus the composition factors of $M^X_1(\gl_k)$ are
\[L_1(\gl_k)=L(s_1\cdot\gl_k), \;L_2(\gl_k)=L(s_2\cdot\gl_k) \mbox{ and } L_3(\gl_k)=L(s_1s_2\cdot\gl_k),\]
Also $[L_1(\gl_k)], [L_2(\gl_k)],$ and $[L_3(\gl_k)]$ generate the same subgroup of the Grothendieck group $K(\cO)$ as $[V_1(\gl_k)], [V_2(\gl_k)],$ and $[V_3(\gl_k)]$. 
If $k=1$, then the
subgroup of the Grothendieck group $K(\cO)$ generated by the composition factors of $M_1^X(\gl)$ is the same as
that generated by $[L(\gl_1-\gc_1-\gc_2)]$ and $[V_1(\gl_1)], [V_2(\gl_1)],$  $[V_3(\gl_1)]$ (as in \eqref{pup}). 
Thus the composition factors of $M^X_1(\gl_k)$ are
\[L_1(\gl_k)=L(s_1\cdot\gl_k), \;L_2(\gl_k)=L(s_2\cdot\gl_k) \mbox{ and } L_3(\gl_k)=L(s_1s_2\cdot\gl_k),\]
Also $[L_1(\gl_k)], [L_2(\gl_k)],$ and $[L_3(\gl_k)]$. Note that 
$[M^X(s_1s_2\cdot\gl_k)]=[L_3(\gl_k)]$.  
\\ \\
In the classical case, the next step is to use the linear independence
of certain partition functions \cite{KK}, \cite{M} Theorem 10.2.5.  Here we use a generalization of this argument. 
We have
\be \label{mik}
\sum_{i>0} [M_i(\gl_k)]= a_k[V_1(\gl_k)]+b_k[V_2(\gl_k)]+c_k[V_3(\gl_k)] +\gd_{k,1}d[L(\gl-\gc_1-\gc_2)],
\ee for integers $a_k,b_k, c_k, d$.  Since $\cS_1 \subseteq M^X(\gl_1)$ we have $d\ge2$. Thus the multiplicity $m(\ga_1,k,\eta)^X$ of $h_{\ga_1}-k$ in $\det F^X_\eta$ is
\[ a_k{{\bf p}_{Y}(\eta - k\ga_1)} +b_k{{\bf p}_{Y}(\eta - k\ga_2)}
+c_k{{\bf p}_{X}(\eta - k\ga_1-k\ga_2)} +\gd_{k,1}d{\bf p}_{W(\ga)}(\eta-\gc_1-\gc_2).\]
Now the multiplicity of $h_\ga$ in
$\LT  \det F^X_{\eta}$ is at least $\sum_{k\ge 1}m(\ga_1,k,\eta)^X.$
We claim that
\be \label{expr}\sum_{k\ge 1}{{\bf p}_{Y}(\eta - k\ga_1)} +{{\bf p}_{Y}(\eta - k\ga_2)}+2{\bf p}_{W(\ga)}(\eta-\gc_1-\gc_2)\ge \sum_{k\ge 1}m(\ga_1,k,\eta)^X.\ee
By \eqref{yew1} and \eqref{yew5}, the multiplicity of $h_{\ga_1}$ in the leading term of $D^X_\eta$ equals the left side of \eqref{expr}. Thus the claim follows from Corollary \ref{fara}.
\\ \\
Now if $\eta = i\ga_1$ for $i\ge 1$,  then ${{\bf p}_{Y}(\eta - k\ga_2)}={\bf p}_{W(\ga)}(\eta-\gc_1-\gc_2) =0$ for all $k\ge 1$.
Thus we obtain a system of inequalities
\[\sum_{k\ge 1}{{\bf p}_{Y}((i- k)\ga_1)}\ge
\sum_{k\ge 1} a_k{{\bf p}_{Y}((i - k)\ga_1)}.\]
Taking $i=1, 2,\ldots$ we have
$${\bf p}_{Y}(0)\ge a_1{\bf p}_{Y}(0), \;\; {\bf p}_{Y}(0) +{\bf p}_{Y}(\ga_1) \ge a_2{\bf p}_{Y}(0) +a_1{\bf p}_{Y}(\ga_1) \mbox{ etc. }$$ 
Since $|M^X_1(\gl_k):L_j(\gl_k)|>0$ for $j=1,2$, it follows that $a_k$ and $b_k$ are positive. 
We conclude that $a_k=1$ and similarly $b_k=1$ for all $k\ge1$.
Since $|M^X_1(\gl_k):L_j(\gl_k)|=1$ and $\cS_k\subseteq M^X_2(\gl_k)$ for all $k,$ we have from \eqref{mik} that $c_k \ge 0$ and
\[2{\bf p}_{W(\ga)}(\eta-\gc_1-\gc_2)\ge d{\bf p}_{W(\ga)}(\eta-\gc_1-\gc_2) +\sum_{k\ge 1}c_k{{\bf p}_{X}(\eta - k\ga_1-k\ga_2)}. \]
Hence arguing as before, using $\eta =i(\ga_1+\ga_2)$ for $i\ge 1$, we find $c_k=0$ for all $k$,  and 
$d=2$.
\\ \\
{\bf Case 3.}
\noindent
  Suppose next that $\ga=\alpha_1 \in \ovd^+_{0}$ and \eqref{cad} {\it does not} hold. We can handle several cases simultaneously by introducing some more notation. First we define the index set $I$ by
\[
I= \left\{
\begin{array}{cl}
\Z & \mbox{if there is a root } \ga_2 \mbox{ such that  } \ga_1^\vee \equiv -\ga_2^\vee \mod \Q X, \\
\N& \mbox{ if } \ga\in C_X^+\\
\N \backslash \{0\} & \mbox{otherwise}.
\end{array}
\right.\]
Next for $r\in I$, define ${\bf q}_{Y}(\eta - r \alpha)$ by
\[
{\bf q}_{Y}(\eta - r \alpha) =\left\{
\begin{array}{crc}
{\bf p}_{Y}(\eta - r \alpha_1)& \mbox{if } r>0, \\
{\bf p}_{Y}(\eta + r \alpha_2)& \mbox{if } I =\Z &\mbox{ and } r<0\\
0 & \mbox{otherwise}.&
\end{array}
\right.\]	
Suppose $r\in I, r\neq 0$ and choose a general $\gl\in \cH_X$ such that $r=(\gl+\gr,\ga_1^\vee).$
Then ${M}^X_{1}(\lambda)$ is simple by Theorem \ref{bits}.  Hence by Lemma \ref{eek}, the Jantzen filtration has the form \eqref{wad}
and the multiplicity of $h_{\alpha} + \rho(h_{\alpha}) -
r(\alpha, \alpha)/2$ in $\det F^X_\eta$ is equal to $c_r {\bf q}_{Y}(\eta - r \alpha)$, for some positive integer $c_r$.
If $r=0$, and $\ga \in C^+_{X}$, set $\gc = \Gc(\ga)$ and $\gc' = s_\alpha\gc$. Then by Theorem \ref{eon},
 ${M}^X_{1}(\lambda) = K^{X,\ga}(\lambda),$ is simple, 
and for some positive integer $c_0$, the multiplicity of $h_{\alpha} + \rho(h_{\alpha})$ in $\det F^X_\eta$ is equal to
\[
c_0( {\bf p}_{Z(\ga)}(\eta-\gc') - {\bf p}_{Z(\ga)}(\eta-\gc)) .\]
  Thus using Corollary \ref{fara},  the multiplicity of $h_{\alpha}$ in the leading term of $\det F^X_\eta$ is at least 
\[c_0( {\bf p}_{Z(\ga)}(\eta-\gc') - {\bf p}_{Z(\ga)}(\eta-\gc)) +
\sum_{r\in I, r\neq 0}c_r{\bf q}_{Y}(\eta - r \alpha)=|\LT \det F^X_\eta:h_\ga| \le |D^X_\eta :h_\ga|\]
On the other hand the multiplicity of $h_\ga$ in the leading term of $D^X_\eta$ is
\[|D_1D_3:h_\ga|= {\bf p}_{Z(\ga)}(\eta-\gc') - {\bf p}_{Z(\ga)}(\eta-\gc) +
\sum_{r\in \Z, r\neq 0}{\bf q}_{Y}(\eta - r \alpha).\]
Reasoning as in Case 2, we deduce that $c_r=1$ for all $r\in I,$ $c_r =0$ otherwise. So the multiplicity of
 $(h_{\alpha} + (\rho, \alpha) - r(\alpha,\alpha)/2)$ is as claimed.
Taken together, Cases 1-3 show  that $D^X_{\eta}$ divides $\det F^X_{\eta}$.  Therefore the result follows from Corollary \ref{fara}.
\hfill  $\Box$
\subsection{The Jantzen sum formula. }\label{jasu} 
\noi {\rm We have the following analog of the Jantzen sum formula. 
\bt \label{Jansum101}
For all $\lambda \in \cH_X$, we have

\by \label{tin} \sum_{i>0} \ch {M}^X_i({\lambda}) &=&
\sum_{\alpha \in A(\lambda)} \ch {M}^{s_\ga X} (s_\alpha \cdot \lambda)
 +\sum_{\gc  \in
B_X(\lambda)}
\ch {M}^{X\cup\{\gc \}}({\lambda}-\gc ) \\
&+& \sum_{\ga  \in
C_X(\lambda)}
\ch {K}^{X,\ga}({\lambda})+\sum_{[\ga] \in E_{X}(\gl)} \ch M^{W(\ga)}(\gl-\gc_1-\gc_2).\nn\quad \quad\ey
\et


\noindent \bpf We have 
\begin{equation}\label{24}\sum_{i>0} \ch {M}^X_{i}
(\lambda) = \sum_{i>0} \sum_{\eta} \dim
{M}^{X}_{i}(\lambda)^{\lambda - \eta} \tte^{\lambda -
\eta}, \end{equation} and by \eqref{yet},
 \begin{equation}\label{25} \sum_{i>0}
\dim {M}^{X}_{i}(\lambda)^{\lambda - \eta} =
v_{T}(  \det F_\eta^X (\widetilde{\lambda})).
\end{equation}

 \noi Now use the factorization of $  \det F_\eta^X $ given by
Theorem \ref{shapdet}. For
$\alpha \in {\Delta}^+_0$

\[ v_{T}((h_{\alpha} + \rho(h_{\alpha}) - r(\alpha,\alpha)/2)
(\widetilde{\lambda})) = \left\{
\begin{array}{ll}
1 \;\;\; \mbox{if} & (\lambda + \rho, \alpha^\vee) = r \\
0 & \;\;\;\;\; \mbox{otherwise}.
\end{array} \right . \]

\noi
and for $\gc  \in B_X$

\[ v_{T}((h_{\gc } + \rho(h_{\gc })) (\widetilde{\lambda})) =
\left\{
\begin{array}{ll}
1 \;\;\; \mbox{if} & (\lambda + \rho, \gc ) = 0
\\ 0 &
\;\;\;\;\; \mbox{otherwise.}
\end{array}
\right . \]
We have similar expressions if $\ga\in C_X^+$ and $(\ga_1,\ga_2) \in E_{X}$.
\noi From this we obtain

\by\label{26}
v_{T}(  \det F_\eta^X (\widetilde{\lambda})) &=& \sum_{\alpha \in
A(\lambda)} {\bf p}_{s_\alpha X}(\eta - (\lambda + \rho, \alpha^\vee) \alpha )  +
\sum_{\gc  \in B_{X}(\lambda)} {\bf p}_{X\cup \{\gc \}}(\eta - \gc ).\nn\\
&+&
\sum_{\alpha \in C_X(\lambda)}\sum_{\eta} ({\bf p}_{Z(\ga)}(\eta-\gc') - {\bf p}_{Z(\ga)}(\eta-\gc) )\\
&+&\sum_{[\ga] \in E_{X}(\gl)} {\bf p}_{W(\ga)}(\eta-\gc_1-\gc_2).\nn
\ey
\noi We combine
Equations (\ref{24}) - (\ref{26}) to conclude

\begin{eqnarray}
 \sum_{i>0} \ch {M}^X_i({\lambda})
& = & \sum_{\alpha \in A(\lambda)} \sum_{\nu}
{\bf p}_{s_\alpha X}(\nu)\tte^{s_\alpha \cdot \lambda - \nu}
+ \sum_{\alpha \in C_X(\lambda)}\sum_{\eta} ({\bf p}_X(\eta)-{\bf p}_{s_{\ga} X}(\eta))
\tte^{\lambda - \eta}
\nonumber\\
& + & \sum_{\gc  \in B_{X}(\lambda)} \sum_{\nu} {\bf p}_{X\cup \{\gc \}}(\nu)\tte^{\lambda - \gc  - \nu}
+\sum_{[\ga] \in E_{X}(\gl)} {\bf p}_{W(\ga)}(\eta-\gc_1-\gc_2)\tte^{\lambda - \eta}
\nonumber
\end{eqnarray}
and this easily yields \eqref{tin}.
\epf


\pagebreak

\appendix
\Bc {\Large\bf  Appendices} \Ec
\section{Anti-distinguished Borel subalgebras.} \label{pip}
\noi In \cite{K} Table VI, Kac gave a particular diagram for each contragredient Lie superalgebra that we will call {\it distinguished.} 
\ff{The term distinguished Borel subalgebra was later introduced by Kac to refer to the Borel subalgebras corresponding to these diagrams, \cite{Kac2} Proposition 1.5.}
\ff{We remark that there are some omissions in this table.  Corrections appear in \cite{FSS} (and elsewhere).}
The distinguished diagram contains at most one grey node, attached to a simple isotropic root vector.  
Unless  $\fg = \osp(1,2n)$, $\osp(2,2n), D(2,1,\ga)$ or $F(4)$ there is exactly one other diagram with this property.  
This diagram  will be called {\it anti-distinguished.} 
The anti-distinguished diagram for $\fgl(m,n)$ is the same as the distinguished diagram for $\fgl(n,m)$.
If $\fg = \osp(1,2n)$  there is only diagram, and it contains no grey node,  while if 
$\fg = \osp(2,2n)$   the anti-distinguished diagram  contains exactly two grey nodes. 
If $\fg = F(4)$, then apart from the  distinguished diagram, there are two other diagrams with a unique grey node.  
We call all these diagrams {\it anti-distinguished}. 
\\ \\
Except in type A, the anti-distinguished diagrams are given in the table below. 
For $\fg=F(4)$ or $G(3)$ we follow the notation of \cite{FSS}. 
In the first two rows, the symbol 
$\begin{picture}(10,10)(50,0) \thinlines \put(55,2.8){\circle{6}} \put(52.7,.10){$\bullet$} \end{picture}$ 
represents a node which is either grey or white. 
Each diagram in these rows has a unique grey vertex. 
Each anti-distinguished diagram corresponds to  one or more Borel subalgebras, which we also call {\it anti-distinguished}. 
The number of such Borels is indicated in the table.  
For $\fg$ of type A, the grey node in the  distinguished (resp. anti-distinguished diagram)  
corresponds to the simple root $\gep_m-\gd_1$ (resp. $\gd_n-\gep_1$).  
The corresponding Borel subalgebras consist of the upper (resp. lower) triangular matrices. 
If $\fg=D(2,1,\ga)$ there are three Borel subalgebras that share the same diagram as the distinguished diagram for $\osp(4,2)$.  
We arbitrarily declare one of these to be distinguished and another anti-distinguished.   
Without doing this, some results would not apply  to the  Lie superalgebra $D(2,1,\ga)$.

\begin{tabular}{|c|c|c|} \hline {}
 algebra  &
 \# Borels & anti-distinguished diagram
\\ \hline
&               &           \\
&                &                \\
$\begin{array}{c}
                B(m,n)    \\
                m >0 \end{array}$ & 1& \begin{picture}(310,10)(26,0)
\thinlines \put(55,3){\circle{6}}
  \put(58,3){\line(1,0){46}}
	\put(162,3){\line(1,0){46}}
	\put(214,3){\line(1,0){46}}
	\put(265.2,5){\line(1,0){46}}
\put(265.2,1){\line(1,0){46}}
  \put(308,-2.3){\huge $\bullet$}
\put(130,3){$\ldots$}
\put(211,3){\circle{6}}
\put(263,3){\circle{6}}
\put(208.3,0.30){$\bullet$}
\put(260.6,0.30){$\bullet$}
\put(52.6,.30){$\bullet$}
\put(295,3){\line(-1,-1){10}}
  \put(295,3){\line(-1,1){10}}
  \put(335,20){}
\end{picture}\\&                                & \\ 
      \hline
			&               &           \\
&                &                \\
$C(n)$    & 1& \begin{picture}(310,10)(26,0)
\thinlines \put(55,3){\circle{6}}
  \put(58,3){\line(1,0){46}}
	\put(162,3){\line(1,0){46}}
	\put(214,3){\line(1,0){46}}
\put(130,3){$\ldots$}
\put(211,3){\circle{6}}
\put(263,3){\circle{6}}
\put(296.5,-19){$\otimes$}
\put(296,17){$\otimes$}
  \put(266.0,4){\line(2,1){31}}
  \put(265,0.5){\line(2,-1){31.5}}
	\put(300,-12.5){\line(0,1){29.0}} 
  \put(340,0){}
\end{picture}\\&                                & \\ &               &           \\
      \hline
			&               &           \\
&                &                \\
$D(m,1)$    & 2& \begin{picture}(310,10)(26,0)
\thinlines \put(55,3){\circle{6}}
  \put(58,3){\line(1,0){46}}
	\put(162,3){\line(1,0){46}}
	\put(214,3){\line(1,0){46}}
	\put(266.7,3){\line(1,0){45}}
  \put(315,3){\circle{6}}
\put(130,3){$\ldots$}
\put(211,3){\circle{6}}
\put(259,.3){$\otimes$}
\put(295,3){\line(-1,-1){10}}
  \put(295,3){\line(-1,1){10}}
  \put(231,3){\line(1,-1){10}} 
	\put(231,3){\line(1,1){10}}
  \put(335,20){}
\end{picture}\\&                                & \\
      \hline
&               &           \\
&                &                \\
$\begin{array}{c}
                D(m,n)    \\
                m, n >1 \end{array}$ & 2& \begin{picture}(310,10)(26,0)
\thinlines \put(55,3){\circle{6}}
  \put(58,3){\line(1,0){46}}
	\put(162,3){\line(1,0){46}}
	\put(214,3){\line(1,0){46}}
	\put(265.2,5){\line(1,0){46}}
\put(265.2,1){\line(1,0){46}}
  \put(314,3){\circle{6}}
\put(130,3){$\ldots$}
\put(211,3){\circle{6}}
\put(263,3){\circle{6}}
\put(208.6,0.30){$\bullet$}
\put(52.6,.30){$\bullet$}
  \put(285,3){\line(1,1){10}}
	\put(285,3){\line(1,-1){10}} 
 
 \put(335,20){}
\end{picture}\\&                                & \\ 
      \hline
&               &           \\
&                &                \\
$F(4)$
    & 1& \begin{picture}(310,10)(26,0)
\thinlines \put(55,3){\circle{8}}
  \put(59,3){\line(1,0){82}}
	\put(58,0.5){\line(1,0){84}}
	\put(58,5.5){\line(1,0){84}}
	\put(148.4,3){\line(1,0){76}}
	\put(232,5.2){\line(1,0){80}}
	\put(232,.8){\line(1,0){80}}
\put(315,3){\circle{8}}
	\put(55,3){\circle{8}}
\put(140.7,.3){$\otimes$}
\put(229,3){\circle{8}}
  \put(266,3){\line(1,-1){10}} 
\put(266,3){\line(1,1){10}}
\put(105,3){\line(-1,-1){10}}
  \put(105,3){\line(-1,1){10}}
  \put(335,20){}
\end{picture}\\&                                & \\
      \hline

		&               &           \\
&                &                \\
$F(4)$
    & 1& \begin{picture}(310,10)(26,0)
\thinlines \put(55,3){\circle{8}}
  \put(59,3){\line(1,0){82}}
	\put(58,0.5){\line(1,0){84}}
	\put(58,5.5){\line(1,0){84}}
	\put(233,3){\line(1,0){78}}
\put(147,5.3){\line(1,0){79}}
\put(147.5,.8){\line(1,0){78}}\put(315,3){\circle{8}}
	\put(55,3){\circle{8}}
\put(140.7,.3){$\otimes$}
\put(229,3){\circle{8}}
\put(105,3){\line(-1,-1){10}}
  \put(105,3){\line(-1,1){10}}
	\put(186,3){\line(1,-1){10}} 
\put(186,3){\line(1,1){10}}
\end{picture}\\&                                & \\
      \hline
					&               &           \\
&                &                \\
$G(3)$
    & 1& \begin{picture}(310,10)(26,0)
\thinlines 
  \put(99,5.2){\line(1,0){85.5}}
	\put(99,0.5){\line(1,0){85.5}}
	\put(190,3){\line(1,0){79}}
\put(189.5,5.3){\line(1,0){80}}
\put(189.5,.8){\line(1,0){80}}
\put(273,3){\circle{8}}
\put(89.7,-3.3){\Huge $\bullet$}
\put(182.5,.3){$\otimes$}
	\put(225,3){\line(1,-1){10}} 
\put(225,3){\line(1,1){10}}
\end{picture}\\&                                & \\
      \hline
					
\end{tabular}

\vspace{0.4cm}

\section{Low Dimensional Cases}\label{ldc}
This appendix serves two purposes.  First it illustrates the general theory developed in the rest of the paper.
Secondly and more fundamentally, for the factorization of the  \v{S}apovalov  determinant for $\fgl(m,n)$ we need to know the structure of the modules $M^X(\gl)$ in certain generic cases.  This is done by a reduction to the cases of $\fgl(2,1)$ and $\fgl(2,2)$.
 \subsection{The Cases $\fsl(2,1)$ and $\fgl(2,1)$.}\label{C8}
We suppose that $\mathfrak{g} = \fsl (2,1)$, though everything we say works with minor modifications for $\fgl(2,1)$. Let\[h = e_{1,1} - e_{2,2},\; \;z = e_{1,1} + e_{2,2} + 2e_{3,3}
\; \; \mbox{and} \; \;\mathfrak{h} = \mbox{ span } \{h,z\}.\]
The notation used here is the same as \cite{M} Exercises 10.5.3 and 10.5.4.
We use the distinguished  set of positive roots
${\Delta^+} = \{ \alpha, \beta, \gc = \alpha + \beta
\},$ with $\alpha$ even and $\beta$ odd  simple roots.  We
define  negative root vectors
\[e_{-\alpha} = e_{2,1}, \; \; e_{-\beta} = e_{3,2},\;\;
e_{-\gc} = e_{3,1}.\] \noi The \v{S}apovalov element $\gth_\gc$ is given by

\[\gth_\gc = e_{-\beta}e_{-\ga} +e_{-\gamma}h.\]
Note that $\fg$ has three Borel subalgebras whose even part is the standard upper triangular Borel subalgebra of $\fg_0$.
We label these as $\fb^{(1)}$-$\fb^{(3)}$ where $\fb^{(1)}, \fb^{(2)}, \fb^{(3)}$ respectively have set of simple roots
\[ \{\ga,\gb\}, \quad \{\gc,-\gb\}, \quad \{-\gc,\ga\}.\] We also consider  the parabolic subalgebras \[\fp = \fb^{(1)}+\fb^{(2)},\quad \fq = \fb^{(2)}+ \fb^{(3)}.\] Note that $\fg =\fp \oplus  \fm =\fq \oplus  \fn$ where $\fm$ and $\fn$ are abelian subalgebras with roots $-\ga, -\gc$ and $-\ga, \gb$ respectively. Thus $U(\fm) = \ttk[e_{-\ga}, e_{-\gc}]$ and $U(\fn) = \ttk[e_{-\ga}, e_{\gb}]$.  If $\mu\in \fh^*$ defines a one-dimensional $\fp$-module, we write $\Ind_\fp^\fg \;\ttk_\mu $ for the induced $\fg$-module.  Similar notation is used for one-dimensional  modules induced from $\fq.$
\\ \\
\noi The figure below
will be used repeatedly when we study the submodule structure of Verma modules and the modules $M^X(\gl)$ in certain cases.
\[
\xymatrix{
&M\ar@{-}[d]&\\
&V_1+V_2&\\
V_{1} \ar@{-}[dr] \ar@{-}[ur] &&
V_{2} \ar@{-}[dl] \ar@{-}[ul]&\\
& V_3= V_1\cap V_2\ar@{-}[d] & \\
&0&}
\]
\bg \label{fig1}\eg

\noi
Note that using the notation of the figure, we have in $K(\cO)$ that
\be \label{tod} [V_1+V_2] +[V_3]  = [V_1] +[V_2].\ee

\noi Suppose that  $\lambda = \gb+n\gc$ where $n$ is a non-negative integer. This is the most interesting case, because $\gl +\gr=n\gc$ a multiple of the non-simple odd positive root $\gc.$
\bt \label{A2} If $\lambda = \gb+n\gc$ with $n\ge 1$, we have
\bi \itema The  lattice
submodules of $M={M}(\lambda)$ is as in Figure \ref{fig1}, where $V_i  = U(\mathfrak{g})v_i$ for
\be \label{kid} v_1 = e_{-\ga}^{n}v_\gl,\quad  v_2 = e_{-\gc}e_{-\gb}v_\gl,\quad \mbox{ and } v_3 = e_{-\gb}v_1 = e_{-\ga}^{n-1}e_{\gb}v_2.\ee

\itemb The module $M/V_2$ has character $\tte^{\lambda} {p}_{\gamma}.$ Thus $M^\gc(\gl)=M/V_2$ is the factor module of $M$ obtained from the specialization process described in Theorem \ref{newmodgen}. 
\itemc The following are equivalent
\bi
\itemi $n>1$
\itemii $V_2 = U(\fg)\gth_\gc v_\gl$.
 \itemiii $V_2$  is a highest weight module.
\itemiv The kernel $N\cap T{M}({\widetilde{\lambda} })_{A}$ of the map from
$N= M^{\gamma}({\widetilde{\lambda} -\gc})_{A}$ to ${M}({{\lambda} })$
in \eqref{bd} equals $TN.$
\ei
\ei \et
\bpf For the statement about the lattice of submodules see \cite{M3} or \cite{M} Exercise 10.5.4.
 We consider
 highest weight modules obtained by specialization.
Let ${\widetilde{\lambda}} = \lambda +T\gc$ and consider
the $U(\fg)_B$-module
$M({\widetilde{\lambda}})_{B}$.  We have
 \by \label{hid} u &:=& \gth_\gc v_{{\widetilde{\lambda}}} = (e_{-\beta} e_{-\ga}+(T+n-1)e_{-\gamma} )v_{{\widetilde{\lambda}}}\\
&=& (e_{-\ga} e_{-\gb}+(T+n)e_{-\gamma} )v_{\widetilde{\lambda}}.\nn \ey
Let  $\bar{v}_{\widetilde{\lambda}}$ be the image of $v_{\widetilde{\lambda}}$ modulo $U(\fg)_B u$. Since $T+k$ is invertible in $B$ for any $k$, it follows that
$$U(\fg)_B u = B[e_{-\ga},e_{-\gb}]u \;\mbox{ and }\; U(\fg)_B \bar{v}_{\widetilde{\lambda}}
= B[e_{-\ga},e_{-\gb}]\bar{v}_{\widetilde{\lambda}}.$$
 Thus if $N=U(\fg)_Au$
and $N'=U(\fg)_A\bar{v}_{\widetilde{\lambda}},$
 the proof of Theorem \ref{newmodgen} shows that for any $n$, $N/TN$ and
$N'/TN'$ are highest weight modules with characters $\tte^{\lambda -\gc} {p}_{\gamma}$
and $\tte^{\lambda} {p}_{\gamma}$ respectively.
\\ \\
To prove (b) first note that $e_{-\gb}v_{{\lambda} } $ is a highest weight vector for $\fb^{(2)}$
which generates $M(\gl).$ However $h_\gc e_{-\gb}v_{\lambda}  = 0$, so $e_{-\gc} e_{-\gb}v_{{\lambda} } $ generates a proper submodule. It follows easily that
\be \label{inm}V_2\cong\Ind_\fq^\fg \;\ttk_{\gl-\gb-\gc} = U(\fn)\ttk_{\gl-\gb-\gc},\mbox{  and } M/V_2\cong\Ind_\fq^\fg \;\ttk_{\gl-\gb} = U(\fn)\ttk_{\gl-\gb},\ee
which yields the first statement in (b).  The second statement holds since $V_2$ is the
unique submodule with this character.
\\ \\
Now we prove (c). Suppose first that $n>1$. 
Since $h_\gb v_2 = (1-n) v_2$ it follows from \eqref{hid} that (ii) and (iii) hold, so $V_2 \cong M^\gc(\gl-\gc)$.
Thus $N/TN\cong M^\gc(\gl-\gc)$ embeds in $M(\gl)$, so  $N\cap T{M}({\widetilde{\lambda} })_{A}=TN.$
\\ \\
Now suppose $n = 1.$ Then the submodule $V_3$ generated by $\gth_\gc v_{{\lambda}} =v_3$ has character $\tte^{\gb} {p}_{\gb} \neq \tte^{\lambda -\gc} {p}_{\gamma}.$ Since $e_{\gb}v_2\neq 0$, $V_2$ is not a highest weight module. We have shown that (ii), (iii) do not hold. Note in addition that by Equation  \eqref{hid}, 
$e_{-\gb}uv_{{\widetilde{\lambda} }} = Te_{-\gb}e_{-\gc}v_{\widetilde{\lambda} } \in N\cap T{M}({\widetilde{\lambda} })_{A}$.
Since $e_{-\gb}e_{-\gc}v_{\widetilde{\lambda} } \notin N$, (iv) does not hold. \epf 
\br {\rm The element $v_2\in M(\gl)$ is not a highest weight vector for the distinguished Borel subalgebra.  
However, compare \eqref{hid},
 \be \label{suv}e_{\beta}v_2 = \gth_\gc v_{{{\lambda}}}=(e_{-\beta} e_{-\ga}+(n-1)e_{-\gamma} )v_{{{\lambda}}}=
(e_{-\ga} e_{-\gb}+ne_{-\gamma} )v_{{\lambda}}.\ee}
\er
\bt \label{AA2} Suppose that $n= 0$, that is $\gl=\gb =-\gr$.  Then
\bi
\itemd ${M}(\lambda)$ has a unique composition series of length 3
 $${M}(\lambda) \supset V_1
= U(\mathfrak{g}) e_{-\beta} v_{{\lambda}}
 \supset V_2 = U(\mathfrak{g})e_{-\ga} e_{-\gb} v_{{\lambda}}
 \supset 0,$$ and $V_1/V_2$ is isomorphic to the trivial module.
\iteme $M(\gl)/V_1$ and $M(\gl)/V_2$ have characters $\tte^{\lambda} {p}_{\beta}$ and $\tte^{\lambda} {p}_{\gamma}$ respectively.
\itemf The module $M^\gc(\gl)=M(\gl)/V_2$   has a basis consisting of the images of the elements $e_{-\ga}^mv_{
\gl}$,
$e_{-\ga}^me_{-\gc} v_{\gl}$ with $m \ge 0$
 and $e_{-\gb}v_{\gl}$.
\ei \et \bpf First (d) is shown in \cite{M3} or \cite{M} Exercise 10.5.4.
It is easy to see that $V_1\cong\Ind_\fp^\fg \;\ttk_{\gl-\gb} = U(\fm)\ttk_{\gl-\gb},$ and $M/V_1\cong\Ind_\fp^\fg \;\ttk_{\gl} = U(\fm)\ttk_{\gl}.$ This gives $\ch M(\gl)/V_1 = \tte^{\lambda} {p}_{\beta}$ so by Lemma \ref{ink}, $\ch M(\gl)/V_2 =
\tte^{\lambda} {p}_{\gamma}$.  Finally (f) holds  because
in the factor module $M^{\gamma}(\lambda)_A$ of $M({\widetilde{\lambda}})_A$ we have \[e_{-\ga}^{m+1} e_{-\gb}v_{\widetilde{\lambda}}= -Te_{-\ga}^m e_{-\gamma} v_{\widetilde{\lambda}},\]
so that $M^{\gamma}(\lambda)_A^{\widetilde{\lambda}-\gc-m\ga}$ has $A$-basis $e_{-\ga}^m e_{-\gamma}v_{\widetilde{\lambda}}$. \epf
 \brs \label{dig}{\rm
\bi \itema Many of the submodules and factor modules of $M(\gl)$  can be constructed using induced modules and their
duals. When $n=1$ in addition to \eqref{inm}, we have
$V_3\cong\Ind_\fp^\fg \;\ttk_{s_\ga \cdot\gl-\gb} = U(\fm)\ttk_{s_\ga\cdot\gl-\gb}$  and $V_1\cong M(s_\ga \cdot \gl)$.
\itemb When $n=0$ we have  $V_1\cong\Ind_\fp^\fg \;\ttk_{\gl-\gb} = U(\fm)\ttk_{\gl-\gb}$  and
$M/V_1\cong\Ind_\fp^\fg \;\ttk_{s_\ga \cdot\gl} = U(\fm)\ttk_{s_\ga\cdot\gl}$ and
$V_2\cong\Ind_\fq^\fg \;\ttk_{\gl-\gb-\gc} = U(\fn)\ttk_{\gl-\gb-\gc}.$
\itemc
 When $n=0$, the module $M^\gc(\gl)=M(\gl)/V_2$ has a one dimensional trivial submodule  This shows that unlike the case of Verma modules, the modules $M^\gc(\gl)$ need not have a filtration with factors which are Verma modules for $\fg_0$, compare \cite{M} Theorem 10.4.5. This example raises issues for the Jantzen filtration, see Lemma \ref{eon} and the factor $D_3$ in Theorem \ref{shapdet}.
Finally we note that $$(M(\gl)/V_2)^\vee \cong\Ind_\fq^\fg \;\ttk_{\gl-\gb} = U(\fn)\ttk_{\gl-\gb}.$$ This is not a highest weight module for the standard Borel.
\itemd We claim that the complex \eqref{1let} is exact at $M(\gl)$ iff $n\neq 0,1$. 
First if $n=0$, it follows easily from \eqref{suv} that $\Img {\psi_{\gl,\gc}}= V_3$ and $\Ker {\psi_{\gl+\gc,\gc}}=V_2.$ If $n=1$, then as noted in the proof of Theorem \ref{A2}, $\Img {\psi_{\gl,\gc}}= V_3$.  We have 
$${\psi_{\gl+\gc,\gc}}(v_2) =e_{-\gc }e_{-\gb}(e_{-\gb}e_{-\ga}+ e_{-\gc})v_{\gl+\gc}=0.$$
Thus $V_2\subseteq\Ker {\psi_{\gl+\gc,\gc}},$ and it is not hard to see that equality holds. If $n>1$, then by Theorems \ref{A2} and 
that $\Img {\psi_{\gl,\gc}}= V_2\subseteq \Ker {\psi_{\gl+\gc,\gc}}.$ 
Equality must hold since any submodule strictly containing $V_2$ has finite codimension.
Similar arguments can be used for the case $\lambda = \gb+n\gc$ with $n$ negative. 
\ei} \ers

\subsection{The Case of $\fgl(2,2)$.}\label{Ch8}
The Jantzen sum formula has some terms which are not Verma
modules.  This suggests that we look at some more modules. Suppose
that $X$ is a set of odd orthogonal roots, and $(\gl + \gr, \ga) =
0$ for all $\ga \in X.$
We consider the modules $M^X(\lambda)$ with
 character $ \tte^{\gl} p_X$ defined in Section \ref{jaf}.
For the remainder of this appendix, let $\mathfrak{g}$ be the Lie superalgebra
$\fgl(2,2)$. Then $\mathfrak{g} = \fn^- \oplus \mathfrak{h} \oplus \fn$ where
$\fn^-, \fh$ and $\fn$ are the subalgebras of lower triangular,
diagonal and upper triangular matrices.  Let
$$e_\ga = e_{1,2},\;\; e_\gb = e_{2,3},\;\;e_\gc = e_{3,4},$$ so
that the simple roots of $\fn$ are $\ga,\;\; \gb,\;\; \gc.$  Also
define
$$e_{\ga + \gb} = e_{1,3},\;\; e_{\gb +\gc} =
e_{2,4},\;\; e_{\ga + \gb + \gc} = e_{1,4}.$$ For each positive
root $\eta,$ let $e_{-\eta}$ be the transpose of $e_\eta.$
Let $\fh^*$ have basis $\epsilon_{1}, \epsilon_{2}, \delta_{1},
\delta_{2} $ as usual with $(\epsilon_{i},  \epsilon_{j}) =  -
(\delta_{i},  \delta_{j}) = \delta_{i,j}.$ The inner product of
the roots is given by the table
\[
 \begin{tabular} {|c|c|c|c|} \hline
 &$\ga$ & $\gb$ & $\gc$\\  \hline
$\ga$ & 2& -1& 0\\ \hline
  $\gb$ &-1 & 0 & 1\\  \hline
$\gc$ & 0 & 1 & -2\\ \hline
 \end{tabular}
  \]
Let $\fh$ be the Cartan subalgebra with basis \[h_{\ga} = e_{1,1}
- e_{2,2},\;\;\;h_{\gc} = e_{4,4} - e_{3,3}, \;\;\; h_{\gb} =
e_{3,3} + e_{2,2}.
\]
Note that with these definitions $<h_\gl,h_\mu> = (\gl, \mu)$
where $<\;,\;>$ is the supertrace,  so $h_\gl v_\mu = (\gl,
\mu)v_\mu.$
 Let
\[ \fb = \fb_1 = \left[ \begin{array}{cccc}
* & * & * & *  \\ 0 & * & * & *  \\ 0 & 0 & * &
* \\ 0 & 0 & 0 & *
\end{array} \right] \] be the distinguished Borel, and
\[  \fb_2 = \left[ \begin{array}{cccc}
* & * & * & *  \\
0 & * & $0$ & *  \\
0 & * & * & *\\
 0 & 0 & 0 & *
\end{array} \right] .\]
Note that $\fb, \fb_2$ are adjacent Borels.  We also  consider the
subalgebras
\[\fp = \left[ \begin{array}{cccc}
 * & * & *   & * \\
 0   & * & 0   & * \\
 * & * & * & * \\
 $0$ & * & 0   & *
\end{array} \right] ,
\fm = \left[ \begin{array}{cccc}
 $0$ & $0$ & 0   & 0 \\
 *   & $0$ & *   & 0 \\
 $0$ & $0$ & $0$ & $0$ \\
 * & $0$ & *   & $0$
\end{array} \right] \] and $\fq = \fb_1 + \fb_2$.
Note that the supertrace, $\ST$ is given by
\[\ST = \gep_1+\gep_2-\gd_1-\gd_2.\] We think of $\ST$ as a one dimensional representation of $\fg$  with kernel $\fsl(2,2)$.
Let
\be \label{per} X = \{\gb, \ga +\gb+\gc\}, \quad Y = \{\ga +\gb, \gb+\gc\}.\ee
Each of  $X,Y$ is an isotropic set of positive roots, but only $Y$ is the set of simple roots for some Borel subalgebra.  Note that $s_\ga(X) =s_\gc(X)= Y$. When $(\gl+\gr,\gb)=0$, $\gl$ defines a one-dimensional $\fq$-module,
$\ttk v_\gl$, and we define 
$M^\gb(\lambda)= {\mbox{ Ind}}^{\sfg}_{\sfq} \;\ttk v_\gl.$
\subsection{\v Sapovalov Elements.} \label{1s}
Suppose $\gl  \in \fh^*$ satisfies
\be \label{ustare} h_{\ga} v_{\gl}  = av_{\gl}, \quad \quad h_{\gb} v_{\gl} =
0, \quad \quad h_{\gc} v_{\gl} =  -cv_{\gl} \ee and
$$(\gl, \ga^\vee) = a,\;\;\; (\gl, \gc^\vee) = c.$$
Set $\gs = s_\ga \cdot \gl, \;\gt = s_\gc \cdot \gl, \;\nu = s_\ga s_\gc \cdot \gl.$ Then
\[(\gs+\gr, \ga+\gb)= (\gt+\gr, \gb+\gc)= (\nu+\gr, \ga+\gb+\gc)=0.\]
Also
\[(\gs+\gr, \ga^\vee) = -(a+1),\;\;\; (\gs+\gr, \gc^\vee) = c+1.
\]
\[(\gt+\gr, \ga^\vee) = a+1,\;\;\; (\gt+\gr, \gc^\vee) = -(c+1).
\]
\[(\nu+\gr, \ga^\vee) = -(a+1),\;\;\; (\nu+\gr, \gc^\vee) = -(c+1).
\]
First to find the \v Sapovalov element for the root $\ga + \gb,$ note that $$ e_{-\ga}^{p+1} e_{-\gb} = [e_{-\ga}e_{-\gb}- pe_{-\ga -\gb
}] e_{-\ga}^{p}.$$
Therefore using \eqref{121nd} with $p=(s_\ga\cdot\gl +\gr, \ga^\vee) = -(a+1),$ we obtain
\[\gth_{\ga +\gb}(\gs)= e_{-\gb}e_{-\ga}- (a+2)e_{-\ga -\gb} . \] So
\be \label{a+b}\gth_{\ga +\gb}= e_{-\gb}e_{-\ga}+ e_{-\ga -\gb
}h_\ga. \ee
Similarly we have
\be \label{c+b} \gth_{\gb+\gc}= e_{-\gb}e_{-\gc}+ e_{-\gb-\gc}h_\gc. \ee


\be \label{a+b+c} \gth_{\ga+\gb+\gc}= e_{-\gb}e_{-\gc}e_{-\ga}+ e_{-\ga -\gb}e_{-\gc}h_\ga
+e_{-\gb-\gc}e_{-\ga}h_\gc +e_{-\ga-\gb-\gc}h_\ga h_\gc. \ee
This gives
\bl With $a, c$ as in \eqref{ustare} the \v Sapovalov element for the root $\ga +\gb +\gc$ satisfies
\[\gth_{\ga +\gb +\gc}(\gl) = e_{-\ga}e_{-\gc}e_{-\gb} -(a+1)(c+1)e_{-\ga -\gb-\gc} + (a+1)e_{-\gc}e_{-\ga -\gb} -(c+1)e_{-\ga}e_{-\gb -\gc}.\]\el
We will be especially interested in the deformed case with $a=c$. Note that $(\gr, \ga^\vee) = (\gr, \gc^\vee) = 1$
and $(\gr, \ga+\gb+\gc) =  (\gr, \gb) = 0.$ Thus  for ${{\lambda}\in\cH_X }$ we can take ${\widetilde{\lambda} }=\gl+T\gr.$
Then  when $a=c$, we have
\[\gth_{\ga +\gb +\gc}(\widetilde{\gl}) = e_{-\ga}e_{-\gc}e_{-\gb} -(a+T+1)^2e_{-\ga -\gb-\gc} + (a+T+1)e_{-\gc}e_{-\ga -\gb} -(a+T+1)e_{-\ga}e_{-\gb -\gc}.\]

\noi In the deformed setting when
$e_{- \gb}v_{{\widetilde{\lambda} } }=0$,  and we  can always cancel the factor $a+T+1$ in   $M^\gb({{\widetilde{\lambda}}})_B$.  Hence
\by \label{ing}\widetilde{w}:&=&[(a+T+1)e_{-\ga -\gb-\gc} - e_{-\gc}e_{-\ga -\gb} +e_{-\ga}e_{-\gb -\gc}]
v_{\widetilde{\gl}}.\\
&=& [(a+T)e_{-\ga -\gb-\gc} - e_{-\gc}e_{-\ga -\gb }+e_{-\gb -\gc}e_{-\ga}]v_{\widetilde{\gl}} \nn \\
&=& [(a+T-1)e_{-\ga -\gb-\gc} - e_{-\ga -\gb }e_{-\gc}+e_{-\gb -\gc}e_{-\ga}]v_{\widetilde{\gl}} \in U(\fg)_B \gth_{\ga +\gb +\gc} v_{{\widetilde{\lambda}}}\nn 
\nn
\ey

\subsection{Change of Borel.}\label{cob}

In the table below we define 6 Borel subalgebras of $\fg$
by listing the corresponding sets of simple roots. Adjacent entries in the table correspond to adjacent Borel subalgebras.
\[
 \begin{tabular} {|c|c|c|} \hline
 $\fb^{(1)} \;\; \{\ga, \gb, \gc\}$& &\\   \hline
 $\fb^{(2)} \;\; \{ \ga + \gb, -\gb, \gb + \gc \}$ &
 $\fb^{(3)}  \;\; \{ -\ga -\gb, \ga, \gb + \gc\}$& \\ \hline
 $\fb^{(4)} \;\; \{\ga + \gb , \gc, -\gb - \gc\}$ &
 $\fb^{(5)} \;\; \{
  -\ga  - \gb,  \ga + \gb+ \gc,- \gb -\gc\}$&%
   $\fb^{(6)} \;\; \{\gc, -\ga - \gb -\gc, \ga\}$\\
   \hline
 \end{tabular}
  \]
\noi Next let 
\be \label{fff} f_1 = e_{-\ga - \gb},\quad f_2 = e_{- \gb - \gc},\quad f_3 = e_{-\ga -
\gb - \gc}.\ee
The negatives of the weights of these elements are
\[\gs_1 = {\ga +
\gb},\; \gs_2 = { \gb + \gc}, \; \gs_3 = {\ga +
\gb+ \gc}.\]
	Suppose that $(\gl+\gr,\gb)= (\gl+\gr,\ga+\gb+\gc)=0$, and let $v_\gl$ be a highest weight vector in
$M= M^\gb(\lambda)$ or $M^X(\lambda)$ with weight $\gl$. We  list some highest weight
vectors in  $M$
for various Borel subalgebras. These elements will be used later to analyze the structure of the modules $M^\gb(\lambda)$ and $M^X(\lambda)$. Starting from $\fb^{(1)}$ and $ \fb^{(2)}$ we perform odd reflections, until we arrive at a highest weight vector $v_\gk$ for $\fb^{(5)}.$
It turns out that $v_\gk$ is
a highest weight vector  for $ \fb^{(6)}$ in $M^\gb(\lambda)$, and also a $\fg$ highest weight vector
which maps to zero in $M^X(\lambda)$, see Theorem \ref{gnn}. Working now only in $M^\gb(\lambda)$
 we reverse the process until we return to a
highest weight vector $w$ for $\fb^{(1)}$ and $ \fb^{(2)}$.
The element $w$ arises from a  \v Sapovalov element, see \eqref{ing} and Lemma \ref{w}.
The elements in the  table below are highest weight vectors for $\fg_0$ since $\fb^{(1)}, \ldots,\fb^{(6)}$ all have the same even part.
These elements  will be used to study the decomposition $M^\gb(\lambda)$ and  $M^X(\lambda)$ as $\fg_0$-modules.
\[ \begin{tabular}{|c|c|c|c|c|c|c|} \hline
Highest weight vector& Borel subalgebra &Weight
\\ \hline
$v_\gl$&  $\fb^{(1)}, \fb^{(2)}$&
      $\gl$
\\ \hline
$e_{-\ga - \gb}v_\gl$&  $\fb^{(3)}$&     $\gl - \gs_1$
\\ \hline
$e_{- \gb - \gc}v_\gl$&  $\fb^{(4)}$&   $\gl - \gs_2$
\\ \hline
$e_{-\ga - \gb}e_{-\gb-\gc}v_\gl$&  $\fb^{(5)}$&  $\gl - \gs_1 - \gs_2$
\\ \hline
 $v_\kappa = e_{-\ga - \gb - \gc}e_{-\gb-\gc}e_{-\ga - \gb}v_\gl$&  $\fb^{(5)}$ and $\fb^{(6)}$&
$\gl - \gs_1 - \gs_2 -\gs_3$
 \\ \hline
 $x = e_{\ga + \gb}v_\kappa$&  $\fb^{(4)}$&   $\gl - \gs_2 -\gs_3$
 \\ \hline
$y = -e_{\gb+\gc}v_\kappa$&  $\fb^{(3)}$&    $\gl - \gs_1 -\gs_3$
\\ \hline
 $e_{\gb+\gc}e_{\ga + \gb }v_\kappa$&  $\fb^{(2)}$ and $\fb^{(1)}$&  $\gl - \gs_3,$ see Lemma \ref{quack}
\\ \hline
$w$&  $\fb^{(1)}, \fb^{(2)}$&        $\gl-\gs_3$
\\ \hline
 \end{tabular}
  \]

\subsection{General Computations in $M^\gb(\lambda)$.}\label{cms}

Suppose $\gl = a(\ga+\gb+\gc).$\ff{ We could consider more generally $\mu = a(\ga+\gb+\gc)+b\gb,$ but the analysis is very similar. This is because
\[M^Y(\gl) \ot \ST^{\ot c} \cong M^Y(\gl + c(\ga+2\gb+\gc) ).\]
}At first there
is no condition on $a$.
We will construct $M^X(\lambda)$ as a factor module of $M^\gb(\lambda)$. Before doing this
it will help to introduce some special elements of $M^\gb(\lambda)$.  Since $(\gl, \ga + \gb + \gc) = (\gl,
\gb) = 0,$ that is $(\gl,X) = 0,$   $\gl$
defines a one dimensional $\fq$-module $\ttk v_{\gl}.$  Now \eqref{ustare} becomes
\be \label{sstar} h_{\ga} v_{\gl}  =
av_{\gl}, \quad \quad h_{\gb} v_{\gl} =  0, \quad \quad h_{\gc}
v_{\gl} =  -av_{\gl}.\ee Also we have
$$(\gl, \ga^\vee) = (\gl, \gc^\vee)=a.$$
 Let $M = M^\gb(\gl)$ be the induced $\fg$-module:
$M=\Ind^\fg_\fq \;\ttk_\gl$ and
\be\label{vk}v_\kappa = e_{-\ga-\gb-\gc}e_{-\gb-\gc }
 e_{-\ga-\gb }v_\gl.\ee The module $M$ has character
$\tte^{\gl}p_\gb.$
\bl \label{w} The element
\by \label{ww} w &=& [(a+1)e_{-\ga -\gb-\gc} - e_{-\gc}e_{-\ga -\gb } + e_{-\ga}e_{-\gb -\gc}] v_\gl \nn\\
&=& [ae_{-\ga -\gb-\gc} - e_{-\gc}e_{-\ga -\gb }
+e_{-\gb -\gc}e_{-\ga}]v_\gl\\
&=& [(a-1)e_{-\ga -\gb-\gc} - e_{-\ga -\gb }e_{-\gc}+e_{-\gb -\gc}e_{-\ga}]v_\gl
\nn
\ey
 is a highest weight vector for $\fg,$ and \[w = [e_{-\ga -\gb-\gc}(h_\ga+1) - e_{-\gc}e_{-\ga -\gb } +
e_{-\ga}e_{-\gb -\gc}] v_\gl \] \el
\bpf This follows from \eqref{ing} or a direct computation.\epf 
\noi
Next using
(\ref{sstar})
\by \label{dog} x =e_{\ga + \gb }v_\gk &=&  e_{-\gc} e_{-\ga-\gb } e_{-\gb-\gc }v_\gl - ae_{-\ga-\gb-\gc} e_{-\gb-\gc
}v_\gl, \nn\\
&=& e_{-\ga-\gb } e_{-\gb-\gc }e_{-\gc}
v_\gl - (a-1)e_{-\ga-\gb-\gc} e_{-\gb-\gc
}v_\gl
\ey and
\by \label{cat} y = -e_{\gb+\gc }v_\gk &=& e_{-\ga} e_{-\ga-\gb } e_{-\gb-\gc }v_\gl + ae_{-\ga-\gb-\gc} e_{-\ga-\gb}v_\gl\nn \\
&=& e_{-\ga-\gb } e_{-\gb-\gc }e_{-\ga}  + (a-1)e_{-\ga-\gb-\gc} e_{-\ga-\gb}v_\gl.
\ey 

\br {\rm In the deformed case, when $a=0$ and $(T\xi,\ga^\vee) = T$, we obtain from \eqref{ing} that modulo $U(\fg)_B \gth_{\ga +\gb +\gc} v_{{\widetilde{\lambda}}}$
\by\label{rat2} Te_{-\ga -\gb-\gc}v_{\widetilde{\gl}} &\equiv& [e_{-\gc}e_{-\ga -\gb }-e_{-\gb -\gc}e_{-\ga}]v_{\widetilde{\gl}}\\
&\equiv & [e_{-\ga -\gb }e_{-\gc}-e_{-\ga}e_{-\gb -\gc}]v_{\widetilde{\gl}}\nn\ey
Equivalently
\be\label{rat3} (T+1)e_{-\ga -\gb-\gc}
v_{\widetilde{\gl}}
 \equiv [e_{-\ga -\gb }e_{-\gc}
-e_{-\gb -\gc}e_{-\ga}]v_{\widetilde{\gl}}\ee
Similarly when $a=0$ we obtain from \eqref{dog} and \eqref{cat} that

\be \label{dog1} e_{-\gc} e_{-\ga-\gb } e_{-\gb-\gc }v_{\widetilde{\gl}}   \equiv Te_{-\ga-\gb-\gc} e_{-\gb-\gc
}v_{\widetilde{\gl}}  \ee and
\be \label{cat1} e_{-\ga} e_{-\gb-\gc } e_{-\ga-\gb}v_{\widetilde{\gl}} \equiv T e_{-\ga-\gb} e_{-\ga-\gb-\gc}v_{\widetilde{\gl}}.\ee
For a generalization of \eqref{rat3} see \eqref{mice2}. These congruences are important when we compare the determinants of the two
bilinear forms $F^X_\eta$ and  $G^X_\eta$ in Subsection \ref{com}. 
(Note that we can obtain \eqref{dog1} and \eqref{cat1} by multiplying \eqref{rat2} by
$e_{-\gb -\gc}$
and $e_{-\ga -\gb }$ respectively.)}\er
\bl \label{wot}
We have \be \label{kwik}
e_{-\ga-\gb }w = y
\ee
\be
\label{quick}
e_{-\gb-\gc }w = x\ee
\be \label{kwak}
e_{-\ga-\gb }x = \pm (a-1) v_\gk
\ee
\be
\label{quock}
e_{-\gb-\gc }y = \pm (a-1) v_\gk.\ee
Therefore if $a\neq 1$, $v_\gk\in U(\fg)w$.\el
\bpf  To prove \eqref{kwik} note that by \eqref{ww} and \eqref{cat},
\[e_{-\ga-\gb }w = f_1[(1-a)f_3+f_2e_{-\ga}]v_\gl=y.\]
The proof of \eqref{quick} is similar.  The proofs of   \eqref{kwak} and \eqref{quock} are easier. When \eqref{dog}
is multiplied by $e_{-\ga - \gb}$ one of the terms becomes
zero and the other a multiple of $v_\gk$, giving \eqref{kwak}. The last statement follows immediately.
\epf
\bc \label{II} For any $a$ we have $x, y\in U(\fg)w \cap U(\fg)v_\gk$.
\ec
\bpf This follows from  \eqref{dog}, \eqref{cat}, \eqref{kwik} and \eqref{quick}.
\epf
\noi Next set  
\be \label{k11} z =  e_{-\gc}f_3f_1v_\gl -e_{-\ga} f_3f_2v_\gl.\ee 
\bl \label{yamo}
We have 
\bi
\itema $e_\gb v_\gk = \pm z.$
\itemb $ e_{\ga}z = \pm x,  \;\; e_{\gc}z = \pm y.$
\itemc $ e_{\gb}z = 0.$
\itemd $e_{-\gb} e_{-\ga}v_\gl = f_1 v_\gl$, and $e_{-\gb} e_{-\gc}v_\gl = -f_2 v_\gl$
\ei \el \bpf
For (a) we note that 
\[e_\gb v_\gk = \pm f_3[e_\gb, f_1f_2]v_\gl = \pm f_3(e_{-\ga}f_2 - e_{-\gc}f_1)v_\gl = \pm z.\]
Then (c) follows from (a), (b) is an easy verification, and (d) follows since $e_{-\gb}v_\gl= 0$.
\epf

\bl \label{quack}
\be\label{zin} e_{\gb+\gc }x =
(1-a)w.\ee
\el
\bpf
First note that
$$e_{\gb+\gc }e_{-\ga-\gb } e_{-\gb-\gc}v_\gl  =  -e_{-\ga - \gb}h_{\gb+\gc}v_\gl = ae_{-\ga-\gb}v_\gl.$$
Therefore $$e_{\gb+\gc }e_{-\gc} e_{-\ga-\gb } e_{-\gb-\gc
}v_\gl  =  (e_{-\gc}e_{\gb+\gc } +e_{\gb})e_{-\ga - \gb}e_{-\gb-\gc}v_\gl.$$
$$=[(a-1)e_{-\gc} e_{-\ga-\gb } + e_{-\ga} e_{-\gb-\gc}
+ e_{-\ga-\gb-\gc}]v_\gl.$$
and
$$e_{\gb+\gc }e_{-\ga-\gb-\gc} e_{-\gb-\gc}
v_\gl = [ae_{-\ga-\gb-\gc} +e_{-\ga} e_{-\gb-\gc}]v_\gl
\;\;\quad \mbox{using  (\ref{sstar})}.$$
Combining these equations gives
$$e_{\gb+\gc }x =
e_{\gb+\gc}e_{\ga + \gb }e_{-\ga-\gb-\gc}e_{-\ga- \gb}e_{-\gb-\gc }v_\gl$$
$$= [(a-1)e_{-\gc} e_{-\ga-\gb } + e_{-\ga} e_{-\gb-\gc}
+ e_{-\ga-\gb-\gc}]v_\gl $$%
$$-a[ae_{-\ga-\gb-\gc} +e_{-\ga} e_{-\gb-\gc}]v_\gl
$$%
$$=[(1-a^2)e_{-\ga-\gb-\gc} +(a-1)e_{-\gc} e_{-\ga-\gb }
+(1 - a)e_{-\ga} e_{-\gb-\gc} ]v_\gl.$$
$$ = (1-a)w.$$\epf
\bc \label{san} $e_{\ga + \gb }y =
(a-1)w.$
\ec
\bpf Immediate from \eqref{dog}, \eqref{cat} and Lemma \ref{quack}. \epf
\noi
\bc \label{cor33} If $a \neq -1,$ then $$U(\mathfrak{g}_0)w \; \oplus \;
U(\mathfrak{g}_0) f_1v_\gl \; \oplus \; U(\mathfrak{g}_0)f_2v_\gl  =
U(\mathfrak{g}_0)f_3v_\gl
 \; \oplus \; U(\mathfrak{g}_0) f_1v_\gl \; \oplus \;
U(\mathfrak{g}_0)f_2v_\gl.
$$  The LHS is a direct sum of $U(\fg_0)$ highest weight modules. \ec
\bpf By Lemma \ref{w}, $w=[(a+1)f_3- e_{-\gc}f_1 + e_{-\ga}f_2] v_\gl $ and this gives the first statement.  The second follows from Subsection  \ref{cob}.\epf
\bc \label{Ian} If $a\neq 1$ then $w, x, y$ and $v_\gk$ all generate the same submodule of $M(a).$
\ec
\bpf This follows from Lemmas \ref{wot}, \ref{quack}, \eqref{dog} and \eqref{cat}. \epf

\bl \label{l35} Suppose that $a \neq 0.$ Then we have \bi
\itema
\[ U(\mathfrak{g}_0) f_2f_1v_\gl \; \oplus
\; U(\mathfrak{g}_0) x
 = U(\mathfrak{g}_0) f_2f_1v_\gl \;
\oplus \; U(\mathfrak{g}_0) f_3f_2v_\gl,
\]
\itemb \[ U(\mathfrak{g}_0) f_2f_1v_\gl \; \oplus
\; U(\mathfrak{g}_0) y
 = U(\mathfrak{g}_0) f_2f_1v_\gl \;
\oplus \; U(\mathfrak{g}_0) f_3f_1v_\gl.
\] \itemc
The left sides  are direct sums of
$\fg_0$ highest weight modules.  \ei \el \bpf  (a) follows immediately from \eqref{quick}, and (b) is proved in a similar
 way. Note that $x, y$ are $\fg_0$ highest weight vectors since they appear in the second table in
Subsection \ref{cob}. \epf
\bl \label{r} If $a=1$, then $e_\gb v_\kappa$ is a non-zero multiple of $e_{-\ga}e_{-\gc}e_{-\gb}w.$
\el
\bpf
Note that
$$v_\kappa = e_{-(\ga + \gb + \gc)}e_{-(\gc + \gb) }e_{-(\ga + \gb) }v_\gl.$$
An easy computation shows  that
$e_\gb v_\kappa $ is a non-zero multiple of $e_{-\ga -\gb-\gc}w.$
If $a=1$ then $w$ has weight zero as a weight vector for $\fsl(2)\ti\fsl(2)$.  Hence $e_{-\ga} w = e_{-\gc}w =0$. Another computation shows that
$e_{-\ga}e_{-\gc}e_{-\gb}w$  is also a non-zero multiple of $e_{-\ga -\gb-\gc}w.$
\epf
\noi For the convenience of the reader, we record that the elements $v_\gk, w, x, y, z$ are defined in this Section in Equations \eqref{vk}, \eqref{ww},  \eqref{dog}, \eqref{cat}, \eqref{k11} respectively. All these elements except for $z$ are $\fg_0$ highest weight vectors that also appear in the table at the end of Section \ref{cob}. For the elements $f_1, f_2, f_3$ see \eqref{fff}. 
\subsection{$M^X(\lambda)$ as a factor module of $M^\gb(\lambda)$.} \label{3}
We assume throughout this section that $\gl \in\cH_X$, and $(\gl,\ga^\vee)=a.$
\noi Let $N=N^X(\gl)$ be the kernel of the natural map from $M^\gb(\gl)$ onto $M^X(\gl),$ and
let $\cN^X(\gl)$ be the submodule of $M^\gb(\gl)$ generated by $w$ and $v_\gk$. The goal of  this section is to show that
$N^X(\gl)=\cN^X(\gl)$.
\bt \label{gnn} The submodule $N^X(\gl)$ of $M^\gb(\gl)$ is equal to
 $\cN^X(\gl)$, and we have $x, y \in
\cN^X(\gl)$. Furthermore $N^X(\gl)$  is a highest weight module iff $a\neq 1$.
\et
\bpf
Recall the element ${\widetilde{w}} \in U(\fg)_B \gth_{\ga +\gb +\gc} v_{{\widetilde{\lambda}}}$ from \eqref{ing}.
 After setting $T=0$,  ${\widetilde{w}} $  specializes to $w$ as in \eqref{ww}.  Hence $w \in N^X(\gl)$. By \eqref{kwik} and \eqref{quick} $x, y \in N^X(\gl),$ and if $a\neq 1$ we have
$v_\gk \in N^X(\gl)$ by \eqref{kwak} or \eqref{quock}, so  by  Lemma \ref{w} $N^X(\gl)$ is a highest weight module generated by $w$.  Now suppose that $a= 1$. Then

\be \label{ing1} e_{-\ga -\gb}e_{-\gb -\gc}{\widetilde{w}} =\pm Tv_\gk \in U(\fg)_B \gth_{\ga +\gb +\gc} v_{{\widetilde{\lambda}}}.\ee
Hence $v_\gk \in U(\fg)_B \gth_{\ga +\gb +\gc} v_{{\widetilde{\lambda}}}$, so
$\cN^X(\gl) \subseteq N^X(\gl).$ 
 To complete the proof it suffices to show that the factor module $ M^\gb(\gl)/\cN^X(\gl)$ has character $\tte^{\gl}p_X.$  This is done in  Theorems \ref{3.7}, \ref{3.10} and \ref{3.11}. 
\epf
\noi For the remainder of this section we fix $\gl \in \cH_X$ and suppose $a$ is as in \eqref{sstar}.

\bt \label{3.7} Suppose $a \neq 0, -1$ then
\bi \itema $M^\gb(\gl)$ is a direct sum of 8 Vermas for
$\fg_0.$ The highest weight vectors are
$$v_\gl, \; f_1v_\gl, \; f_2v_\gl, \;  f_2f_1v_\gl, \; v_\gk  =f_3f_2f_1v_\gl, \; x, \; y, \; w.$$
\itemb The factor module $ M^\gb(\gl)/\cN^X(\gl) $ is the direct sum of
the 4 Vermas for
$\fg_0$ with highest weight vectors
\be \label{ext}v_\gl, \; f_1v_\gl, \; f_2v_\gl, \;  f_2f_1v_\gl.\ee
\itemc We have $\ch \cN^X(\gl)=\tte^{\gl - \gs_3}p_X$ and $\ch M^\gb(\gl)/\cN^X(\gl) =\tte^{\gl}p_X.$ Thus $N^X(\gl)=\cN^X(\gl)$.
\ei
\et \bpf (a) is immediate from
Corollary \ref{cor33} and Lemma
\ref{l35}. For (b) note that $\cN^X(\gl)= U(\fg) w$ contains the $U(\fg_0)$ highest weight vectors $x,y$ and $v_\gk$ by Lemma \ref{wot}.
The statement follows since $w$ is a highest weight vector for $\fg$. Finally (c) is a computation based on (b).
\epf

\subsubsection{The case $a=0$.} \label{A=0}
\bt \label{3.10} For the case $a = 0,$ we have as $\fg_0$-modules
\bi \itema $M^\gb(\gl) = U(\mathfrak{g}_0)
v_\gl \; \oplus \; U(\mathfrak{g}_0) f_1v_\gl \; \oplus \;
U(\mathfrak{g}_0)f_2v_\gl \; \oplus
\;U(\mathfrak{g}_0)w \; \oplus \; U(\mathfrak{g}_0) v_\gk \; \oplus \;
I
$ where $I$ is an indecomposable with Verma socle $U(\mathfrak{g}_0)f_1
f_2v_\gl$ and highest weight vectors $f_1f_3v_\gl$ and $f_2f_3v_\gl$ in the
top.
\itemb
We have 
$$\cN^X(\gl) = U(\mathfrak{g}_0)w
 \; \oplus \; U(\mathfrak{g}_0) v_\gk \; \oplus \;J$$ where $J$
 is an indecomposable $U(\mathfrak{g}_0)$-module containing 
 $J'=(U(\mathfrak{g}_0)x + U(\mathfrak{g}_0)y)$ and with $J/J'=U(\mathfrak{g}_0)z$, 
a Verma module with highest weight $\gl-2(\ga+\gb+\gc)$.
\itemc   The character of $\cN^X(\gl)$ is \[ \tte^{\gl - \gs_3}(1 +
\tte^{-\ga - \gb})(1 + \tte^{-\gb-\gc})/ \prod_{\gs \in
\Delta^+_{0}} (1 - \tte^{- \gs}) = \tte^{\gl -
\gs_3}p_X.\]
Thus $N^X(\gl)=\cN^X(\gl)$. \ei \et \bpf (a)  
Note that 
\be \label{k10}  e_{\gc}f_3f_2v_\gl =f_1f_2v_\gl \mbox{ and } e_{\ga}f_3f_1v_\gl = \pm f_1f_2v_\gl.\ee 
For example the first equation holds since $e_{- \gb}v_\gl = 0,$ and
$$[e_{\gc} , e_{-\ga -\gb - \gc}e_{-\gb-\gc}] = e_{-\ga - \gb}e_{-\gb-\gc} - e_{-\ga- \gb - \gc}e_{- \gb}.$$
In addition it is easy to see that 
\[e_{\ga}f_3f_2v_\gl =  e_{\gc}f_3f_1v_\gl = 0.\]
Now the first five summands in the
expression for $M^\gb(\gl)$ are generated by $\fg_0$-highest weight vectors, and their sum $K$ is direct.
Thus $M^\gb(\gl)/K$ has the same character as the direct sum of $\fg_0$ Vermas with highest weights $\gl- \gs_1-\gs_2$, $\gl- \gs_1-\gs_3$ and $\gl- \gs_2 -\gs_3$.  The claim follows. \\
\\
(b)  By Lemmas \ref{wot} and  Lemma \ref{yamo} (a)
all summands on the right side are contained in the left.  Also the sum is direct by comparing weights. (Note that the weights of $v_\gk$ involve $\gb$ three times, $x, y, z$ twice and $w$ only once). Lemma \ref{yamo} also implies that $z$ is a highest weight vector mod $J'$. The weight of $z$ is antidominant, so $z$ generates a $\fg_0$ simple Verma module. \\ \\
(c) If $a  = 0,$ then $$x = e_{-\gc} e_{-\ga-\gb } e_{-\gb-\gc }v_\gl,$$ and
$$y = \pm e_{-\ga} e_{-\ga-\gb }
e_{-\gb-\gc }v_\gl.$$
 Hence $U(\mathfrak{g}_0)x \cap
U(\mathfrak{g}_0)y = U(\mathfrak{g}_0)e_{-\ga} e_{-\gc}f_1f_2v_\gl.$
Therefore $U(\mathfrak{g}_0)x + U(\mathfrak{g}_0)y$ has character
\be \label{badc} \tte^{\gl - \gs_3}(
\tte^{-\gs_1} + \tte^{-\gs_2} -
\tte^{-\gs_3})/ \prod_{\gs \in
\Delta^+_{0}}(1 - \tte^{- \gs}). \ee
The result follows since $z$ has weight $\gl - 2\gs_3,$ $v_\gk $ has
weight $\gl - \gs_1 - \gs_2 - \gs_3,$ and $w$ has weight $\gl -
\gs_3.$\epf
\noi When
$a=0$ the elements $z, x, y$ from Subsection  \ref{cms} can be written as follows.
\be \label{fem} z= e_{-\gc}f_3 f_1v_\gl - e_{-\ga} f_3 f_2v_\gl, \quad
x =e_{-\gc}f_1 f_2v_\gl, \quad y= \pm e_{-\ga} f_1 f_2v_\gl.\ee
It follows that $J\subset I$ in Theorem \ref{3.10}.
Also \[e_{-\gc}y = \pm e_{-\ga} x\]
The table below gives a complete list of $\fg_0$-singular vectors in $M/N$. In order to show the Weyl group symmetry, for each $\fg_0$-singular vector we list the shifted weight (i.e. the $\gr$-shifted weight) of $s$ as the ordered pair \be \label{rot}((\wt \; s +\gr, \ga^\vee), (\wt \;s +\gr, \gc^\vee))\ee  This makes evident the singular vectors with the same central character. 
\[
 \begin{tabular} {|c|c|} \hline
$\fg_0$-singular vector $s$&
shifted weight of $s$\\ \hline
$v_\gl$ & $(1,1)$\\ \hline
$f_1f_2v_\gl$ & $(-1,-1)$\\ \hline
$f_3v_\gl$ &$(0,0)$\\   \hline
$f_1v_\gl$ & $(0, 2)$\\ \hline
$f_2v_\gl$ &$(2,0)$\\   \hline
$f_1f_2f_3v_\gl$ &$(0,0)$\\   \hline
$f_1f_3v_\gl$ & $(-1, 1)$\\ \hline
$f_2f_3v_\gl$ &$(1,-1)$\\   \hline
$x$ & $(1,-1)$\\ \hline
$y$ & $(-1,1)$\\ \hline
$z$ &$(-1,-1)$\\   \hline
$v_\gk$& $(0,0)$\\ \hline
 \end{tabular}
  \]
The highest weight vectors of the Verma composition factors of $I$ and their weights
for $[\fg_0,\fg_0] \cong \fsl(2)\ti\fsl(2)$ are displayed in the diagram  below, see \eqref{k10}.  We write $f_1f_2$ in place of $ f_1f_2 v_\gl $ etc.
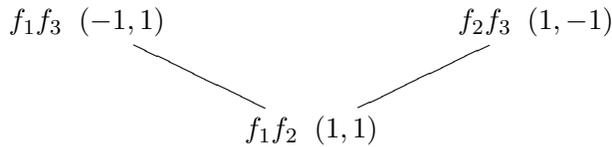
\begin{displaymath}
    \xymatrix{
& f_1f_3 \;\; (-1,1)&& f_2f_3 \;\;(1,-1)&\\
&& f_1 f_2 \;\;(1, 1)\ar@{-}[ul] \ar@{-}[ur]&&
   }
\end{displaymath}
The weights of the elements  $x, y, z \in J$ are displayed in the diagram below. We note that as  a $\fg_0$-module, $J$
does not have a Verma flag see \eqref{badc}, and also Corollary \ref{bdc}.
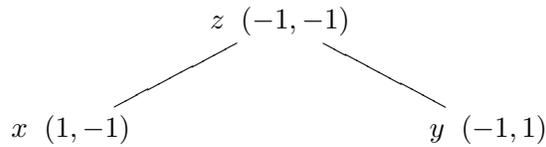
\begin{displaymath}
    \xymatrix{
&&z\;\; (-1,-1)&&\\
& x \;\; (1, -1)\ar@{-}[ur]&& y \;\; (-1,1)\ar@{-}[ul]
   }
\end{displaymath}
\subsubsection{ The case $a=-1$.}\label{3.111}
When $a=-1$, $(\gl+\gr,\ga^\vee)=0$, and as we shall see $M^X(\gl)$ is simple.
 \bt \label{3.11} For the case $a = -1,$ we have as $\fg_0$-modules
\bi \itema
$M^\gb(\gl) = U(\mathfrak{g}_0)
v_\gl \; \oplus \; U(\mathfrak{g}_0)x \; \oplus \;  U(\mathfrak{g}_0)y\; \oplus \;
U(\mathfrak{g}_0)f_2f_1v_\gl \;\oplus \; U(\mathfrak{g}_0) v_\gk
 \; \oplus \;
I$ where $I$ is an indecomposable having Verma socle
$\cV=  U(\mathfrak{g}_0)f_1
v_\gl \oplus U(\mathfrak{g}_0)f_2v_\gl$ and such that the image of $f_3v_\gl$ is a highest weight vector
which generates the $\fg_0$-Verma module $I/\cV$.
\itemb We have 
\be \label{wvk}\cN^X(\gl) = U(\mathfrak{g}_0)x \;\oplus \; U(\mathfrak{g}_0)y \; \oplus \;
 U(\mathfrak{g}_0)w
 \; \oplus \; U(\mathfrak{g}_0) v_\gk,
\ee  where $U(\mathfrak{g}_0)w \subset \cV$. \itemc
 The character of $\cN^X(\gl) $ is \[ \tte^{\gl - \gs_3}(1 +
\tte^{-\ga - \gb})(1 + \tte^{-\gb-\gc})/ \prod_{\gs \in
\Delta^+_{0}} (1 - \tte^{- \gs}) =\tte^{\gl -
\gs_3}p_X.\] Hence $\cN^X(\gl) = N^X(\gl)$.
\itemd 
Set $M=M^\gb(\gl)$, $N= N^X(\gl)$, $u_1=e_{-\gc}f_1v_\gl$ and $u_2 =e_{-\ga}f_2v_\gl$. 
Then as a $\fg_0$-module \[M/N = [U(\mathfrak{g}_0)
v_\gl \; \oplus
U(\mathfrak{g}_0)f_2f_1v_\gl \oplus \;
J],\] where $J$ is an indecomposable with  socle
$J_1= U(\mathfrak{g}_0)u_1 =
 U(\mathfrak{g}_0)u_2$, and unique maximal submodule $J_2 =U(\mathfrak{g}_0)f_1v_\gl + U(\mathfrak{g}_0)f_2v_\gl$ .
The quotient $J_2/J_1$ is the direct sum of two highest weight modules with weights  $(-2,0)$ and $(0,-2)$.  The second $($resp. first$)$
copy of $\fsl(2)$ acts trivially on the first $($resp. second$)$ module.
\ei
\et \bpf (a) For the structure of $I$, note that $e_\ga f_3v_\gl = - f_2v_\gl$
 and $e_\gc f_3v_\gl = f_1v_\gl$ by an easy calculation. For the other terms use Lemma \ref{l35}.\\
(b) Equation \eqref{wvk} follows easily from  Lemma \ref{wot}. By Lemma \ref{w},  
$$w = (e_{-\ga}e_{-\gb-\gc }-e_{-\gc}e_{-\ga-\gb } )v_\gl= (e_{-\ga}f_2-e_{-\gc}f_1)v_\gl= u_2-u_1\in \cV$$ 
is a highest weight vector  for $\fg$, and by Lemma \ref{wot},  $\cN^X(\gl)= U(\fg)w$. Finally (c) and (d) follow from (b).
\epf
\noi 
The table below gives a complete list of $\fg_0$-singular vectors in $M/N$. The weights are shifted as in \eqref{rot}.
\[
 \begin{tabular} {|c|c|} \hline
$\fg_0$-singular vector $s$&
shifted weight of $s$\\ \hline
$v_\gl$ & $(0,0)$\\ \hline
$f_1f_2v_\gl$ & $(0,0)$\\ \hline
$f_3v_\gl$ &$(-1,-1)$\\   \hline
$u=e_{-\ga}f_2v_\gl = e_{-\gc}f_1v_\gl$ & $(-1,-1)$\\ \hline
$f_1v_\gl$ & $(-1, 1)$\\ \hline
$f_2v_\gl$ &$(1,-1)$\\   \hline
 \end{tabular}
  \]
\bt \label{riz} The module $M^X(\gl)=M/N$ is simple. \et
\bpf It suffices to show that no linear combination of the singular vectors (other than $v_\gl$) in the table is a highest weight vector.  We have $e_\gb u= e_{-\ga}e_{-\gc}v_\gl$, as noted in the proof of Theorem \ref{3.11}. Since $e_\gb$ commutes with $f_3$, it follows that   $e_\gb f_3v_\gl = 0$.  However $f_3$ is not a $\fg_0$ highest weight vector.\\
\\
For the rest of the proof we show that no linear combination, other than $v_\gl$ is killed by $e_\gb$.
\\ \\
We have $[e_\gb,f_1f_2] = e_{-\ga}f_2 - f_1e_{-\gc},$ and
 it follows that $e_\gb f_1f_2v_\gl = (2u-f_3)v_\gl\neq 0.$ So there is no highest weight vector with weight $(-1, -1)$.\\
\\
Finally 
 we have $[e_\gb,f_1]v_\gl = e_{-\ga}v_\gl\neq 0$, and
$[e_\gb,f_2]v_\gl = e_{-\gc}v_\gl\neq 0.$
\epf

\subsection{The Maximal Finite Quotient.} \label{mfq}
Let $n=a+1>0$, and suppose $(\gl+\gr, \ga^\vee) = n$.  Then $M^X(\gl)/(U(\fg)e_{-\ga}^n v_\gl)$ is $\ga$-finite. If also $(\gl+\gr, \gc^\vee)=n$, then
$M^X(\gl)/(U(\fg)e_{-\ga}^n v_\gl +U(\fg)e_{-\gc}^n v_\gl)$ is finite (dimensional).  Now the maximal finite quotient of $M(\gl)$  is the Kac module $K(\gl)$. We use known results on the structure of $K(\gl)$ to determine the maximal finite quotient of $M^X(\gl).$
From \cite{SZ}, the Jantzen filtration on $K(\gl)$ is the unique Loewy filtration, and hence coincides with both the socle and radical series of $K(\gl).$ For $\fg =\fgl(2,2)$, this filtration takes the form
\[K(\gl) = K^0(\gl) \supset K^1(\gl) \supset K^2(\gl) \supset K^3(\gl) =0.\]
Setting $K_i(\gl) =K^i(\gl)/K^{i+1}(\gl)$ we have $K_0(\gl) = L(\gl)$. Concerning $K_1(\gl)$ and   $K_2(\gl)$ there are three cases.
\bi \itema If $n>2$, then \[K_1(\gl) = L(\gl-\gb) \op L(\gl-\ga-\gb-\gc),\quad
K_2(\gl) = L(\gl-\ga-2\gb-\gc).
\]
\itemb If $n=2$, then \[K_1(\gl) = L(\gl-\gb) \op L(\gl-\ga-\gb-\gc) \op L(\gl-2\ga-3\gb-2\gc),\]
and $K_2(\gl)$ is as in case (a). However $K(\gl)$ does not contain a highest weight vector of weight $\gl-2\ga-3\gb-2\gc$.
\itemc If $n=1$, then \[K_1(\gl) = L(\gl-\gb),\quad K_2(\gl) = L(\gl-2(\ga+\gb+\gc)).\]
\ei
These facts can also be deduced from \cite{MS}.  Until Theorem \ref{simp}, $v_\gl$ will denote the highest weight vector in $K(\gl)$. To describe the
highest weight and singular vectors in $K(\gl)$, recall the element $v_\gk$ from \eqref{vk}, now regarded as an element of $K(\gl)$ and consider the following products, (which are zero in some cases)
\[v_1 =e_{-\gb} v_\gl,\quad  v_2 = \gth_{\ga +\gb+\gc}v_\gl,\quad  v_3= \gth_{\ga +\gb+\gc}e_{-\gb} v_\gl,\quad
v_4=f_1 f_2 f_3 e_{-\gb} v_\gl.\]
Then we have
\bi
\itema In case (a), $v_1, v_2$ are highest weight vectors
 which together generate $K^1(\gl)$, and $v_3$ is a highest weight vector which generates $K^2(\gl)$.
\itemb In case (b), $v_1, v_2$ are highest weight vectors and $v_\gk$
is a highest weight vector mod $K^2(\gl)$, see Lemma \ref{r}.
Together these elements generate $K^1(\gl)$. Also $v_3$ is a highest weight vector which generates $K^2(\gl)$.
\itemc In case (c), $v_1$ is a highest weight vector
 which generates $K^1(\gl)$, and $v_4$ is a highest weight vector which generates $K^2(\gl)$. Note that $\gl-\ga-\gb-\gc$ is not  dominant in this case.
\ei
For convenience we record the dimensions of the finite dimensional irreducible $\fg$-modules. The same recursive  method can be used to give an explicit expression for their characters, see also  \eqref{chl}. First note that if $\gl, \mu\in \cH_X$ and  $(\gl+\gr, \ga^\vee) = (\mu+\gr, \ga^\vee) = n$, then $L(\gl)$ is isomorphic to a tensor product of $L(\mu)$ with a one dimensional module.  Hence $b_n= \dim  L(\gl)$ depends only on $n$. Also $\dim K(\gl) = 16n^2.$ Now if $n= 1,$ then $L(\gl)$ is trivial as a $[\fg_0,\fg_0]$-module, and so $b_1=1.$ Then from   (a)-(c) above we have
\by 16 &=& 2b_1 +b_2\nn\\
64 &=& 2b_1 + 2b_2 + b_3\nn\\
16n^2 &=& b_{n-1}+ 2b_{n}+b_{n+1} \mbox{ for } n\ge 3. \nn
\ey
Solving these equations we see that
\be \label{beq} b_{n} = 4n^2 -2\mbox{ for } n \ge2. \ee
The next result supports Conjecture \ref{conje}.
\bt \label{simp} If $\gl \in \cH_X$ is dominant, the maximal
finite quotient of $M^X(\gl)$ is simple.
Thus the maximal submodule of $M^X(\gl)$ is $U(\fg)e_{-\ga}^n v_\gl +U(\fg)e_{-\gc}^n v_\gl$.
\et

\bpf This is shown by an analysis of cases (a)-(c) above.  \epf
\noi Theorem \ref{simp} is perhaps surprising since
for any $\gl \in \cH_X$ we have in $K(\cO)$
\be \label{leg} [M(\gl)] = [M^X(\gl)] +[M^X(\gl- \ga-\gb-\gc)] +[M^X(\gl- \gb)] +[M^X(\gl- \ga-2\gb-\gc)].\ee
In fact $M^X(\gl)$ has a filtration with factors having the same characters as the modules on the right of this equation.
Thus since a Kac module $K(\gl)$  can have 3, 4 or 5 composition factors, it is worth observing what happens to these composition factors 
relative to this filtration. First we note that $M(\gl)$ has a submodule $M^\gb(\gl-\gb)$ with factor module $M^\gb(\gl).$  The first (resp. last) pair of summands on the right side of \eqref{leg} appear in a filtration of $M^\gb(\gl)$ (resp. $M^\gb(\gl-\gb)$). The most interesting case is (b) since if $n=2,$ the Kac module $K(\gl)$ has greatest length.  Here $v_1,$ and $v_3$ map to zero in $M^\gb(\gl)$ and $v_2, v_\gk \in N^X(\gl)$.   Thus the only composition factor of the Kac module that survives in $M^X(\gl)$ is $L(\gl),$ and $N^X(\gl)$ is  not a highest weight module. The analysis of cases (a) and (c) is straightforward once it is noted that in case (a), the weight $\gl- \ga-\gb-\gc$ is not dominant.
\subsection{On the structure of $M^X(\gl)$.}\label{smx}
\bt \label{ns} The module $M^X(a)$ is simple unless $a \in\N$. \et
\bpf We assume that $a\neq -1,$ since that case has just been covered. Then if $a$ is not a non-negative integer, then none of the vectors in \eqref{ext} are dominant.  Thus $M^X(\gl)$ is a direct sum of simple $\fg_0$-modules. The result now follows from easy calculations. For example, since $h_{\ga+\gb}v_\gl = av_\gl$, $f_1v_\gl$ does not generate a proper submodule.\epf

\bt \label{genst} If $(\gl+\gr,\ga^\vee)=n>1$, then
\bi \itema The socle $\cS$ of $M^X(\gl)$ is  $M^X(s_\ga s_\gc \cdot \gl) \cong L(s_\ga s_\gc \cdot \gl).$
\itemb We have
\[U(\fg)e_{-\ga}^n v_\gl \cong M^Y(s_\ga \cdot \gl), \quad U(\fg)e_{-\gc}^n v_\gl \cong M^Y(s_\gc \cdot \gl)\] and
\[L(s_\ga s_\gc \cdot \gl) =M^X(s_\ga s_\gc\cdot \gl).\]
\itemc \[M^X_1(\gl)=U(\fg)e_{-\ga}^n v_\gl+U(\fg)e_{-\gc}^n v_\gl 
 \quad  {\rm and } \;\; \cS \subseteq M^X_2(\gl).\]
\ei The lattice of submodules of $M^X(\gl)$ is as in Figure \ref{fig1} with $V_1= U(\fg)e_{-\ga}^n v_\gl$ and $V_2= U(\fg)e_{-\gc}^n v_\gl.$
\et
\bpf Omitted. It is similar to, but easier than the proof of Theorem \ref{stz1}.\epf
\noi Now we discuss the structure of $M^\gb(\lambda)$ and $M^X(\gl)$
in the exceptional cases.
\subsubsection{The case $a=0$.}\label{a=0}
We continue with the notation from Subsection  \ref{A=0}.
 If $n=(\gl+\gr,\ga^\vee) =1$, equivalently $a=0,$ then from Theorem \ref{3.10},
\be \label{1fad}  M^\gb(\gl)/N^X(\gl) = U(\mathfrak{g}_0)
v_\gl \; \oplus \; U(\mathfrak{g}_0) f_1v_\gl \; \oplus \;
U(\mathfrak{g}_0)f_2v_\gl \; \oplus
I/J.\ee
\bc \label{bdc} Let $\fg'=[\fg_0,\fg_0]$.  Then as a $\fg'$-module, $I/J$ has composition factors with highest weights $(0,0), (-2,0)$ and $(0,-2)$ each with
multiplicity one. \ec
\bl \label{shi} The element $p=f_1 f_2v_\gl$ is a $\fg$ highest weight vector in $M/N$ with weight $\gl -\ga-2\gb -\gc$ and spans a 1-dimensional trivial $\fg$-module. We have $p \in U(\fg)e_{-\ga}v_\gl  \cap U(\fg) e_{-\gc}  v_\gl$.\el
\bpf As noted in Subsection \ref{cob},  $f_1 f_2v_\gl=-f_2 f_1v_\gl$ is a highest weight vector for $\fb_5$.  In particular it is a $\fg_0$ highest weight vector.  Also we have in $M$ that
\[ e_{\gb}f_2f_1 v_\gl =\pm(e_{-\gb-\gc}e_{-\ga} - e_{-\gc} e_{-\ga-\gb }) v_\gl = \pm w \in N^X(\gl),\]
see Lemma \ref{w}. In addition $e_{-\ga}p = e_{-\gc}p= 0$ by \eqref{fem} and Theorem \ref{3.10}.  Also $f_3p =v_\gk \in N^X(\gl)$. Since $\fg_0$, $f_3$ and $e_\gb$ generate $[\fg,\fg]$ we obtain the first statement. The second follows since 
$p = \pm e_{-\gb-\gc}e_{-\gb} e_{-\ga}v_\gl =\pm e_{-\ga-\gb}e_{-\gb} e_{-\gc}v_\gl $. \epf
\br{\rm In contrast to the case $n>1$, when $n =1$ the submodule $U(\fg)e_{-\ga}^n v_\gl$  of $M^X(\gl)$ is not isomorphic to
$M^Y(s_\ga \cdot \gl)$. Indeed, if $u =e_{-\ga} v_\gl$ then $e_{-\gb}u$ is a highest weight vector for
 $\fb^{(2)}$
which generates the same  submodule of $M^X(\gl)$ as $u$. So from  $U(\fg)e_{-\ga} v_\gl \cong M^Y(s_\ga \cdot \gl)$ we would have $U(\fg)e_{-\ga} v_\gl \cong \Ind_\fp^\fg \;\ttk
e_{-\gb}e_{-\ga} v_\gl $, see Theorem \ref{2} (a) and hence   $e_{-\gb-\gc}u=0$.  However $e_{-\gb-\gc}u =\pm f_2f_1 v_\gl$.  What this shows is the following.
\bl 
\label{k9}
There are non-split exact sequences
\[0 \lra \ttk f_2f_1 v_\gl  \lra U(\fg)e_{-\ga} v_\gl \lra M^Y(s_\ga \cdot \gl) \lra 0\] and
\[0 \lra \ttk f_2f_1 v_\gl  \lra U(\fg)e_{-\gc} v_\gl \lra M^Y(s_\gc \cdot \gl) \lra 0.\]
\el }\er
\bl \label{yak}
\bi \itema $p_Y - p_X =
\tte^{-\gb}.$
\itemb When $(\gl, \ga) = (\gl, \gc)=(\gl, \gb) = 0$ we have
$$[L(\gl)] +
[L(\gl-\ga -2\gb-\gc)]+
[M^Y(s_\ga \cdot \gl)]
+[M^Y(s_\gc \cdot \gl)] -[M^X(w_0\cdot\gl)]=M^X(\gl).$$
\itemc When $(\gl, \ga^\vee) = (\gl, \gc^\vee) > 0=(\gl, \gb)$ we have
$$[L(\gl)] +
[M^Y(s_\ga \cdot \gl)]
+[M^Y(s_\gc \cdot \gl)] -[M^X(w_0\cdot\gl)]=M^X(\gl).$$
\ei
\el
\bpf (a) follows from Lemma \ref{ink} (a).
To prove (b) and (c), we identify the class of a module in the Grothendieck group $K(\cO)$ with its character. Then (b) follows since by (a),
\by && [M^Y(s_\ga \cdot \gl)]
+[M^Y(s_\gc \cdot \gl)] -[M^X(w_0\cdot\gl)]-M^X(\gl)\nn\\
&=& \tte^\gl[(\tte^{-\ga} +\tte^{-\gc})(p_X+ \tte^{-\gb}) -p_X -\tte^{-\ga-\gc}p_X]\nn\\
&=& \tte^\gl[(\tte^{-\ga}-\tte^{-\gc})\tte^{-\gb}-p_X
(1 - \tte^{-\ga}) (1 - \tte^{-\gc})]\nn\\
&=& -\tte^\gl -  \tte^{\gl-\ga-2\gb-\gc}\nn \\
&=& -[L(\gl)] -
[L(\gl-\ga -2\gb-\gc)].\ey
Similarly if $(\gl+\gr,\ga^\vee)=n>1$, then
setting $x=\tte^{-\ga}, y=\tte^{-\gb}, z =\tte^{-\gc}$ and $\Xi=(1 - x^{-\ga}) (1 - z^{-\gc}),$ we have
\by \label{chl}
  && [M^X(\gl)]+[M^X(w_0\cdot\gl)]
-[M^Y(s_\ga \cdot \gl)]
-[M^Y(s_\gc \cdot \gl)]\\
&=& \tte^\gl[(1+ x^nz^ n)(1+xy)(1+yz) - (x^n+z^n)(1+y)(1+xyz)]/\Xi\nn.\ey
Taking limits as $x,y,z\lra 1$ we obtain $4n^2-2$, which equals $\dim L(\gl)$ by \eqref{beq}.
\epf
\bt \label{stz1}
\bi \itema The socle $\cS$ of $M^X(\gl)$ satisfies
\[\cS = \ttk f_2f_1 v_\gl \oplus U(\fg)e_{-\ga}e_{-\gc}  v_\gl \; = \; U(\fg)e_{-\ga}v_\gl  \cap U(\fg) e_{-\gc}  v_\gl.  \]
\itemb We have
\[U(\fg)e_{-\ga}v_\gl/\ttk f_2f_1 v_\gl \cong M^Y(s_\ga\cdot \gl), \quad U(\fg)e_{-\gc}v_\gl/\ttk f_2f_1 v_\gl \cong M^Y(s_\gc\cdot \gl)
\] and
\[U(\fg)e_{-\ga}e_{-\gc}v_\gl\cong M^X(s_\ga s_\gc\cdot \gl).\]
\itemc $\cS \subseteq M_2^X(\gl).$\ei
\et
\bpf First observe that $U(\fg)e_{-\ga}e_{-\gc}  v_\gl \cong M^X(s_\ga s_\gc \cdot \gl)$ is simple  by Theorem \ref{ns}. The other statements in (b) hold by Lemma \ref{k9}.  
Also by Lemma \ref{shi}, $\ttk f_2f_1 v_\gl $ spans a one-dimensional submodule of $M^X(\gl)$. 
This can also be shown using Lemma \ref{yamo} (d).
Thus $\cS' =\ttk f_2f_1 v_\gl \oplus U(\fg)e_{-\ga}e_{-\gc}  v_\gl \subseteq \cS$.
Furthermore 
 $M^X(\gl)/(U(\fg)e_{-\ga}v_\gl  +U(\fg) e_{-\gc}  v_\gl)$ is the one dimensional module $L(\gl)$ by Theorem
\ref{simp}. Next we show that
$L= (U(\fg)e_{-\ga}v_\gl  \cap U(\fg) e_{-\gc}  v_\gl)/\cS'$ is zero.  
Set $M_\ga = U(\fg)e_{-\ga}v_\gl, M_\gc = U(\fg)e_{-\gc}v_\gl.$     
Then in the Grothendieck group $K(\cO)$ we have
\be \label{k7}
 [M_\ga /M_\ga\cap M_\gc] = [M^Y(s_\ga)/M^X(w_0\cdot\gl)] -[L], \ee 
\be \label{k8}
[M_\gc /M_\ga\cap M_\gc] = [M^Y(s_\gc)/M^X(w_0\cdot\gl)] -[L]
.\ee
For example using Lemma \ref{k9} 
\by \label{k6}
 [M_\ga /M_\ga\cap M_\gc] &=& [M_\ga]-  [L(\gl-\ga -2\gb-\gc)]-[L]-[M^X(w_0\cdot\gl)]\nn\\
&=& [M^Y(s_\ga)/M^X(w_0\cdot\gl)] -[L].\nn
\ey 
In addition  using \eqref{k7} and \eqref{k8},  
\be \label{k5}
 [M_\ga\cap M_\gc] = [L]+[L(\gl-\ga -2\gb-\gc)]+
[M^X(w_0\cdot\gl)].
\ee 
Finally 
\by \label{k4}
 [M_\ga+M_\gc] &-&[M_\ga\cap M_\gc]\\
&=& [M^Y(s_\ga \cdot \gl)/M^X(w_0\cdot\gl)]+
[M^Y(s_\gc \cdot \gl)/M^X(w_0\cdot\gl)].\nn\ey
Combining \eqref{k4} and \eqref{k5} with Theorem \ref{simp} and comparing the result with Lemma \ref{yak} (b), we obtain $L=0.$
Finally (c) holds since the module $M_1^X(\gl)/M_2^X(\gl)$ is self-dual. \epf
\subsection{The module $M^Y(\lambda)$.}\label{mmt}
The main result of this subsection is the following.
\bt \label{mnt}If $\gl\in \cH_Y$ and $(\gl+\gr,\ga^\vee) = a \in \Z$ , then
$M_1^Y(\gl)$ is simple and
\[M_1^Y(\lambda) \cong \left\{
\begin{array}{cl}
M^X(s_\ga\cdot \gl) & \mbox{ if } a \ge 0 \\
M^X(s_\gc\cdot \gl)& \mbox{ if } a<0.
\end{array}
\right.\]
\et
\noi \bpf This follows from Theorems \ref{gnat} and Corollary \ref{lcor} below.
\epf
\noi We show that if $a\neq 0,$ we can construct $M^Y(\lambda)$ using induction from the parabolic
subalgebra $\fp= \fb+\fb^{(2)}.$ As in the footnote in Subsection  \ref{cms}, there is no loss of generality in assuming that $\mu = a(\ga+\gb),$ with $a \in \mathbb{Z}.$ Let $\gl= \mu +\gb$.
 Then $(\mu, \ga +\gb) = (\mu, \gb+\gc) = 0,$ so $\mu$
defines a one dimensional
$\fp$-module $\ttk v_{\mu}.$ Also
\be \label{k3}(\gl +\gr, \ga^\vee) = -(\gl +\gr, \gc^\vee)= a.\ee

\noi 
Note also that \be \label{star}h_{\ga} v_{\mu} =  h_{\gc} v_{\mu} =  av_{\mu},
\quad \quad h_{\gb} v_{\mu} =  -av_{\mu}\ee and that
$$(\mu, \ga^\vee) = -(\mu, \gc^\vee).$$
Set   $R^Y(a)=\Ind_\fp^\fg \;\ttk v_\mu$   and $v_\gl=e_{\gb}v_\mu.$
Then $v_\gl, v_\mu$ are highest weight vectors for the Borels $\fb, \fb^{(2)} $ respectively.
If $a=0$, $v_\gl$ generates a proper  submodule.

\bl \label{Ma1} For all $a$ we have a direct sum of $U(\fm_0)$-modules
$$R^Y(a) =  U(\mathfrak{m}_0)v_{\mu} \; \oplus \;
U(\mathfrak{m}_0)e_{-\ga-\gb-\gc} v_{\mu} \; \oplus \;
U(\mathfrak{m}_0)e_{\gb} v_{\mu} \; \oplus \; U(\mathfrak{m}_0)e_{-(\ga +
\gb + \gc)} e_{\gb} v_{\mu} .$$ \el
\noi \bpf Since $R^Y(a)$ is induced from $\fp$ and $\fg = \fp \oplus \fm,$ this follows at once from the
PBW Theorem. \epf

\bl \label{sub} Suppose $(\mu,Y)=0,$ so that $\mu$ defines a character of $\fp.$ Then
\[\ch \Ind_\fp^\fg \;\ttk v_\mu = \tte^{\mu +\gb}p_Y.\]
Thus $R^Y(a)$  and $M^Y(a)$ have the same character.
\el
\bpf This follows from Lemma \ref{Ma1} because
\[\tte^\mu(1 + \tte^{\gb})(1 + \tte^{-\ga - \gb -\gc}) =
\tte^{\mu + \gb}
(1 + \tte^{
-\gb})(1 + \tte^{-\ga - \gb -\gc}).
\] \epf

\bl \label{udef} Set
$$u = (e_{-\ga}e_{-\gc} + a e_{-\ga-\gb-\gc}e_{\gb})v_\mu$$
\[u'=(e_{-\ga}e_{-\gc}
e_{-\gb}- a^2 e_{-\ga-\gb-\gc})v_{\gl} \]
Then $u, u'$ are $\fg_0$ highest weight  vectors.  \el \bpf A computation.
Note that the substitution $e_{-\gb}v_\gl = a v_\mu$ gives $au=u'.$ \epf

\bl \label{107}  The elements $e_{-\ga-\gb}e_{-\gb}v_{{{\lambda}}}$ and
$e_{-\gb-\gc}e_{-\gb}v_\gl=0$ are in the kernel of the natural map $M(\gl)\lra M^Y(a)$.
\el
\bpf By \eqref{a+b} we have in $M^Y(a)_A$ that
\[0 = \gth_{\ga +\gb}v_{{\widetilde{\lambda}}}  = (e_{-\gb}e_{-\ga}+ e_{-\ga -\gb}h_\ga)v_{{\widetilde{\lambda}}} =  (e_{-\gb}e_{-\ga}+ (T+a)e_{-\ga -\gb})v_{{\widetilde{\lambda}}}.\]
Multiplying both sides by $e_{-\gb}$ and using the fact $M^Y(a)_A$ is $\ttk[T]$ torsion-free, we obtain
$e_{-\gb} e_{-\ga-\gb}v_{{\widetilde{\lambda}}}=0$ in  $M^Y(a)_A$.
Hence
$e_{-\gb} e_{-\ga-\gb}v_{{{\lambda}}}=-e_{-\ga-\gb}e_{-\gb}v_{{{\lambda}}}=0$
in $M^Y(a)$.
Similarly using the expression for $\gth_{\gb+\gc}$ in \eqref{c+b}, we see that  $e_{-\gb-\gc}e_{-\gb}v_\gl=0$ in $M^Y(a)$.
\epf
\noi Now fix $a$ and set $M=R^Y(a)$. We will improve Lemma \ref{Ma1}  to obtain decompositions into $\fg_0$-modules, see Theorems
\ref{2} and \ref{2.5}.
\subsubsection{The case $a>0$.}
Because of \eqref{star} and
\eqref{k3}, we can restrict our attention to the cases  where $a\ge 0.$
If $a < 0$ we can carry out a similar
analysis using $\gc$ in place of $\alpha.$
Note that $v_\gl$ generates a proper submodule of $M$ iff $a = 0.$
 \bt \label{2}
\bi
\itema If $a \neq 0,$ then $M$ is a highest weight module for the distinguished Borel with highest weight
$\gl=\mu+\gb $.
\itemb Furthermore
\by M
&=& U(\mathfrak{m}_0)v_{\mu} \; \oplus \; U(\mathfrak{m}_0)u  \; \oplus \;
U(\mathfrak{m}_0)e_{-\ga-\gb-\gc} v_{\mu} \; \oplus \;
U(\mathfrak{m}_0)e_{\gb}v_{\mu},\nn\\
&=& U(\mathfrak{m}_0)e_{-\gb}v_{\gl}\; \oplus \; U(\mathfrak{m}_0)(e_{-\ga}e_{-\gc}
e_{-\gb}- a^2 e_{-\ga-\gb-\gc})v_{\gl} \nn\\      &\oplus&
U(\mathfrak{m}_0)e_{-\ga-\gb-\gc}e_{-\gb}v_{\gl}\; \oplus \;
U(\mathfrak{m}_0)v_{\gl},\ey
and these are direct sums of $\fg_0$-modules. In addition $M$ is generated as a
$\mathfrak{g}$-module by $v_{\gl}$ which is a highest weight vector for $\fb.$\ei
\et \bpf First (a) follows from Lemma \ref{Ma1}.
By Lemma \ref{udef} we have for $a \neq 0,$
$$U(\mathfrak{m}_0)v_{\mu} \; \oplus \;
U(\mathfrak{m}_0)e_{-\ga-\gb-\gc} e_{\gb} v_{\mu} =
U(\mathfrak{m}_0)v_{\mu} \; \oplus \; U(\mathfrak{m}_0)u.$$ The first decomposition in (b) follows from Lemma \ref{Ma1}, and the second from $e_{-\gb}v_\gl = a v_\mu$.\epf
\noi Since
$(\gl + \gr, \ga^\vee) = a $, $v_{\eta} =
e_{-\ga}^{a}v_{\gl}$ is a highest weight vector for
$\fb$ with weight  $\eta = s_{\ga}\cdot\gl.$
 Let $N$ be the submodule generated by $v_{\eta}$. Theorem \ref{2} gives a decomposition of $M$ into $\fg_0$-submodules, and we give a compatible  decomposition of
$N$.
\bt \label{gnat} Set $v_\gk= e_{-\ga}^{a+1}v_{\mu}$. Then
\bi \itema we have a direct sum of $\fg_0$-modules
$$ N = U(\fm_0)v_\gk
 \; \oplus \; U(\fm_0) e_{-\gb - \gc}v_\eta  \; \oplus \;
U(\fm_0)v_\eta  \; \oplus
\; U(\fm_0)e_{-\gb - \gc}v_\gk,$$
with
\by U(\fm_0)v_\gk
 &\subseteq& U(\mathfrak{m}_0)v_{\mu},\quad \; \;U(\fm_0) e_{-\gb - \gc}v_\eta  \subseteq\; U(\mathfrak{m}_0)u\label{one}\\
U(\fm_0)v_\eta  &\subseteq&
U(\mathfrak{m}_0)e_{\gb}v_{\mu},\;\;
 U(\fm_0)e_{-\gb - \gc}v_\gk\subseteq
U(\mathfrak{m}_0)e_{-\ga - \gb - \gc} v_{\mu}
.\label{two}\ey

\itemb The character of $N$ is \[ \tte^\gk(1 + \tte^{\ga +
\gb})(1 + \tte^{-\gb-\gc})/ \prod_{\gs \in \Delta^+_{0}} (1 -
\tte^{- \gs}) \]
\[ = \tte^{\eta}(1
+ \tte^{-\ga - \gb})(1 + \tte^{-\gb-\gc})/ \prod_{\gs \in
\Delta^+_{0}} (1 - \tte^{- \gs}) = \tte^{s_\ga\cdot\gl} p_X.\]
\itemc $N$ and $M/N$ are simple.
\ei \et \bpf
Note that  $v_\gk$
 is a highest weight vector for $\fb^{(2)}$ with weight $$\gk = s_\ga \cdot \mu = \mu - (a+1)\ga.$$
Also $e_{-\ga -\gb }v_\eta=-v_\gk$  and $e_{\ga + \gb }v_\gk = -(a+1)v_\eta,$ so $v_\eta$ and $v_\gk$ both generate the same submodule. Hence $N \cap U(\mathfrak{m}_0)v_{\gl}$  contains $v_\eta,$ which has weight \be \label{ad}\eta = \gk + \ga + \gb = s_\ga \cdot (\mu  + \gb) = s_\ga\cdot\gl.\ee
Next the inclusions in  \eqref{one} and \eqref{two} follow since
\by e_{-\gb - \gc}v_\eta &=& e_{-\gb-\gc}e_{-\ga}^{a}e_\gb v_\mu,\nn\\ 
&=&(e_{-\ga}^{a}e_{-\gb-\gc} + ae_{-\ga}^{a-1}e_{-\ga-\gb-\gc})e_\gb v_\mu\nn\\  
&=&e_{-\ga}^{a-1}(e_{-\ga}e_{-\gc} + ae_{-\ga-\gb-\gc}
e_{\gb})v_\mu\nn\\  &=& e_{-\ga}^{a-1}u.\nn\ey and
\by e_{-\gb-\gc}v_\kappa  &=&
e_{-\gb-\gc}e_{-\ga}^{a+1}v_\mu\nn\\
&=& (a+1)e_{-\ga}^{a}e_{-\ga-\gb-\gc}v_\mu.\nn\ey
From these computations we also see that $e_{-\gb-\gc}v_\eta$
is a $\fg_0$ highest weight  vector in $N \cap U(\mathfrak{m}_0)e_{-(\ga
+ \gb + \gc)} v_{\gl}$
 with weight $s_\ga \cdot (\mu  - \ga - \gc ) = \gk + \ga - \gc,$ 
and that
$e_{-\gb-\gc}v_\gk$
is a $\fg_0$ highest weight vector in $N \cap
U(\mathfrak{m}_0)e_{-\ga-\gb-\gc}v_\mu$ with weight $s_\ga \cdot
(\mu - \ga - \gb -\gc ) = \gk -\gb -\gc.$
Also $e_{\ga + \gb }e_{-\ga}^{a}e_{-\ga-\gb-\gc}v_\mu = e_{-\ga}^{a-1}u.$ Equality in \eqref{one} and \eqref{two} as well as the statements about simplicity follow easily by looking at the highest weights and using $\fsl(2)$ theory. Finally
(b) follows from (a) and (\ref{ad}). \epf
\bc If $a\neq 0$, then $R^Y(a)\cong M^Y(a)$.\ec
\bpf
\noi As $\fp = \fb\op \span\{e_{-\ga-\gb}, e_{-\gb-\gc}\}$, $M(\gl)=\Ind_\fb^\fg \;\ttk v_\gl$ and $R^Y(a)=\Ind_\fp^\fg \;\ttk e_{-\gb}v_\gl$,
it follows that $R^Y(a)$ is the (universal) module obtained from $M(\gl)$ by setting $e_{-\ga-\gb} e_{-\gb}v_\gl$ and
$e_{-\gb-\gc}e_{-\gb}v_\gl$ equal to zero. Since the same relations hold in $M^Y(a)$ by Lemma \ref{107}, there is an onto map from $R^Y(a)\lra M^Y(a)$.  Because  both modules have the same character, it must be an isomorphism.
\epf

\subsubsection{The case $a=0$.} \label{a=0v2}
From \eqref{a+b} we have
\[(e_{-\ga }e_{-\gb }+(T+a)e_{-\ga -\gb} v_{{\widetilde{\lambda}}}) =\gth_{\ga +\gb } v^Y_{{\widetilde{\lambda}}}=0.\]
\noi
Now set $M=R^Y(\gl)$, where $\gl=0$ and $v_\mu= e_{-\gb}v_\gl$.\
\bt \label{2.5} We have
\bi \itema  If $I=  U(\mathfrak{m}_0)v_{\mu} \; \oplus \;U(\mathfrak{m}_0)e_{-(\ga +\gb + \gc)} e_{\gb} v_{\mu},$ we have
\be \label{fyj} M = U(\mathfrak{m}_0) e_{\gb }v_\mu\; \oplus \; U(\mathfrak{m}_0) e_{-\ga-\gb-\gc}v_\mu \; \oplus \; I,\ee
a direct sum of $\fg_0$-modules. Also $I$ is an indecomposable containing the  Verma submodule $S=U(\mathfrak{g}_0)v_\mu$, such that $I/S$ is a Verma module with highest weight vector $e_{-\ga-\gb-\gc}e_{\gb}v_\mu$.
\itemb The unique maximal $\fg$-submodule of $M$ is
\[N =  U(\mathfrak{m}_0)e_\gb v_{\mu} \; \oplus \;
U(\mathfrak{m}_0)e_{-\ga-\gb-\gc} v_{\mu} \oplus J,\] where $J$ is an indecomposable, codimension one submodule of $I$ which fits into an exact sequence
\[0\lra U(\mathfrak{m}_0)\mathfrak{m}_0v_{\mu} \lra J\lra U(\mathfrak{m}_0)e_{-(\ga +
\gb + \gc)} e_{\gb} v_{\mu} \lra 0.\]
\itemc The character of $N$ is
$\tte^{\gb}p_X,$
and $M/N$ is the trivial module.
\itemd The module $N$ is a simple highest weight module with highest weight $\gb$. Thus $N\cong L(\gb)$.\ei \et
\bpf First \eqref{fyj} follows from Lemma \ref{Ma1}. Note that $s_\ga\cdot
\mu = \mu - \ga$, $s_\gc\cdot \mu = \mu - \gc$.
Also $e_{-\ga}v_{\gl} = e_{\gb}e_{-\ga}v_{\mu},\; e_{-\gc} v_{\gl} =
e_{\gb}e_{-\gc}v_{\mu} .$ It is easy to check that  $e_{\gb} v_{\mu} $ and $e_{-\ga-\gb-\gc} v_{\mu}$ are highest weight vectors for $\fg_0$ and that
\[e_\ga e_{-\ga-\gb-\gc} e_{\gb} v_{\mu} = -e_{-\gc} v_{\mu},\quad \quad
 e_\gc e_{-\ga-\gb-\gc} e_{\gb} v_{\mu} = -e_{-\ga} v_{\mu}\]
\noi
The character of $U(\mathfrak{m}_0)\mathfrak{m}_0v_{\mu}$ is \[
(\tte^{-\ga} + \tte^{-\gc} - \tte^{-\ga - \gc})/
\prod_{\gs \in \Delta^+_{0}}(1 - \tte^{- \gs}). \] In addition
$N$ contains the $\fg_0$ highest weight vectors $e_{-\ga-\gb-\gc}
v_{\mu},\; e_{\gb} v_{\mu}, e_{-\ga-\gb-\gc} e_{\gb} v_{\mu}$
with weights $ -\ga - \gb - \gc, \gb$ and $-\ga -\gc$ respectively.
Adding the characters of the $\fg_0$-modules generated by these elements
we obtain
\[\ch N =
(\tte^{-\ga} + \tte^{-\gc} + \tte^{\gb} +
\tte^{-\ga -\gb -\gc})/
\prod_{\gs \in \Delta^+_{0}}(1 - \tte^{- \gs}) = \tte^{\gb}p_X. \]
Now (c) follows since by (b), Lemma \ref{2} (a) and Lemma \ref{yak} we have
\[\ch M/N = \ch M -\ch N =
\tte^{\gb} ( p_Y - p_X) =1.\]

\noi Finally, using
$$e_{-\gb-\gc}e_\gb v_\mu  =\pm e_{-\gc}v_\mu, \;  e_{-\ga -\gb}e_{-\gc}v_\mu =\pm e_{-\ga}v_\mu, \;  e_{-\ga-\gb}e_{-\gb-\gc}e_{\gb}v_\mu = \pm e_{-(\ga +\gb + \gc)} v_{\mu},$$
 it is easy to check  that $e_\gb v_\mu$  is a highest weight vector which generates $N$.
Since  $(\gl+\gr, \ga^\vee) = \; (\gl +\gr, \gc^\vee) = 0.$
the results of Subsection   \ref{3.111}  apply, and in particular
(d) follows from Theorem \ref{riz}.
\epf
\bc In the Grothendieck group $K(\cO)$, we have $[R^Y(0)] = [M^Y(0)]$.  Thus $M^Y(0)$ has length two with unique maximal $\fg$-submodule isomorphic to the module
$N$ in Theorem \ref{2.5} $(b)$. However $R^Y(0)$ and $M^Y(0)$ are not isomorphic.
\ec
\bpf The first statement holds since $R^Y(0)$ and $M^Y(0)$ have the same character.
For the last statement note that $M^Y(0)$ is a highest weight module for the distinguished Borel but $R^Y(0)$ is not.
By Theorem \ref{2.5}, the submodule $N$ of $R^Y(0)$ has codimension one. On the other hand the elements
$e_{-\ga} v_\gl, e_{-\gb} v_\gl$ and $e_{-\gc} v_\gl$ generate a proper submodule of
$M^Y(0)$ with codimension one. The result follows since $R^Y(0)$ and $M^Y(0)$ have the same character. \epf
\bc \label{lcor} If $a=0,$ the unique maximal submodule of $M^Y(\gl)$ is simple and isomorphic to $M^X(\gl).$\ec
\bpf Combine the previous corollary and Theorem \ref{2.5}.\epf

\end{document}